\documentclass{amsart}
\usepackage{amsthm,amsmath,amssymb,amscd, mathrsfs}
\usepackage[latin1]{inputenc}
\usepackage[all]{xy}
\usepackage{url}
\usepackage{mathtools}
\usepackage{hyperref}

\usepackage{framed}
\usepackage{color}
\usepackage{pbox}
\usepackage{url}

\numberwithin{equation}{section}

\theoremstyle{plain}
\newtheorem{thm}{Theorem}[section]
\newtheorem*{thm*}{Theorem}
\newtheorem{prop}[thm]{Proposition}
\newtheorem{lem}[thm]{Lemma}
\newtheorem{cor}[thm]{Corollary}

\newtheorem{lem-defn}[thm]{Lemma-Definition}
\newtheorem{prop-defn}[thm]{Proposition-Definition}
\newtheorem{thm-defn}[thm]{Theorem-Definition}

\newtheorem*{claim*}{Claim}

\theoremstyle{definition}
\newtheorem{defn}[thm]{Definition}
\newtheorem{example}[thm]{Example}

\theoremstyle{remark}
\newtheorem{rem}[thm]{Remark}

\DeclareMathOperator{\Frac}{Frac}

\DeclareMathOperator{\Ker}{Ker}

\DeclareMathOperator{\Image}{Im}
\DeclareMathOperator{\Hom}{Hom}
\DeclareMathOperator{\End}{End}
\DeclareMathOperator{\Aut}{Aut}

\DeclareMathOperator{\GL}{GL}

\DeclareMathOperator{\rank}{rank}

\DeclareMathOperator{\Der}{Der}
\DeclareMathOperator{\Gal}{Gal}
\DeclareMathOperator{\Spec}{Spec}
\DeclareMathOperator{\Spa}{Spa}
\DeclareMathOperator{\Spf}{Spf}

\DeclareMathOperator{\Pic}{Pic}
\DeclareMathOperator{\Ob}{Ob}
\DeclareMathOperator{\Mor}{Mor}
\DeclareMathOperator{\Sh}{Sh}
\DeclareMathOperator{\PSh}{PSh}
\DeclareMathOperator{\Fix}{Fix}
\DeclareMathOperator*{\colim}{colim}

\newcommand{\lft}{\mathrm{lft}}

\newcommand{\mult}{\mathrm{mult}}

\newcommand{\Afd}{\mathbf{Afd}}
\newcommand{\Rig}{\mathbf{Rig}}
\newcommand{\AdSp}{\mathbf{AdSp}}

\newcommand{\EtSp}{\mathbf{EtSp}}
\newcommand{\Set}{\mathbf{Set}}

\newcommand{\Grpd}{\mathbf{Grpd}}
\newcommand{\Sch}{\mathbf{Sch}}
\newcommand{\FSch}{\mathbf{FSch}}
\newcommand{\Mod}{\mathrm{Mod}}
\newcommand{\Alg}{\mathrm{Alg}}
\newcommand{\Jac}{\mathrm{Jac}}

\newcommand{\id}{\mathrm{id}}
\newcommand{\pr}{\mathrm{pr}}
\newcommand{\sw}{\mathrm{sw}}
\newcommand{\Fr}{\mathrm{Fr}}

\newcommand{\Inf}{\mathrm{Inf}}

\newcommand{\can}{\mathrm{can}}
\newcommand{\cont}{\mathrm{cont}}

\newcommand{\an}{\mathrm{an}}

\newcommand{\red}{\mathrm{red}}
\newcommand{\rig}{\mathrm{rig}}
\newcommand{\cris}{\mathrm{cris}}

\newcommand{\Ncris}{\mathrm{Ncris}}
\newcommand{\st}{\mathrm{st}}
\newcommand{\dR}{\mathrm{dR}}

\newcommand{\Et}{\mathrm{\acute{E}t}}
\newcommand{\et}{\mathrm{\acute{e}t}}
\newcommand{\fppf}{\mathrm{fppf}}

\newcommand{\fin}{\mathrm{fin}}
\newcommand{\Isom}{\mathrm{Isom}}

\newcommand{\Z}{\mathbb{Z}}
\newcommand{\Q}{\mathbb{Q}}

\newcommand{\C}{\mathbb{C}}
\newcommand{\A}{\mathbb{A}}
\newcommand{\F}{\mathbb{F}}

\newcommand{\Isoc}{\mathrm{Isoc}}

\newcommand{\FIsoc}{F\text{-}\mathrm{Isoc}}

\newcommand{\Bun}{\mathrm{Bun}}
\newcommand{\Coh}{\mathrm{Coh}}
\newcommand{\calCoh}{\mathcal{C}oh}

\newcommand{\isoc}{\mathrm{isoc}}
\newcommand{\V}{\mathbf{V}}
\newcommand{\lf}{\mathrm{lf}}

\newcommand{\nilp}{\mathrm{nilp}}

\newcommand{\affine}{\mathrm{affine}}

\newcommand{\aff}{\mathrm{aff}}
\newcommand{\Gm}{\mathbb{G}_{\mathrm{m}}}
\newcommand{\irr}{\mathrm{irr}}

\newcommand{\rar}{\rightarrow}

\newcommand{\xrar}{\xrightarrow}

\newcommand{\wh}{\widehat}
\DeclareMathOperator{\Supp}{Supp}

\newcommand{\fT}{\mathfrak{T}}
\newcommand{\fp}{\mathfrak{p}}
\newcommand{\fq}{\mathfrak{q}}

\newcommand{\bA}{\mathbb{A}}
\newcommand{\bC}{\mathbb{C}}
\newcommand{\bF}{\mathbb{F}}
\newcommand{\bG}{\mathbb{G}}
\newcommand{\bN}{\mathbb{N}}
\newcommand{\bP}{\mathbb{P}}
\newcommand{\bQ}{\mathbb{Q}}

\newcommand{\bZ}{\mathbb{Z}}

\newcommand{\cA}{\mathcal{A}}
\newcommand{\cC}{\mathcal{C}}
\newcommand{\cE}{\mathcal{E}}
\newcommand{\cF}{\mathcal{F}}
\newcommand{\cG}{\mathcal{G}}
\newcommand{\cH}{\mathcal{H}}
\newcommand{\cI}{\mathcal{I}}
\newcommand{\cL}{\mathcal{L}}
\newcommand{\cM}{\mathcal{M}}
\newcommand{\cO}{\mathcal{O}}
\newcommand{\cS}{\mathcal{S}}
\newcommand{\cT}{\mathcal{T}}
\newcommand{\cX}{\mathcal{X}}
\newcommand{\cY}{\mathcal{Y}}
\newcommand{\cZ}{\mathcal{Z}}

\newcommand{\fkm}{\mathfrak{m}}
\newcommand{\fkn}{\mathfrak{n}}

\newcommand{\fkp}{\mathfrak{p}}

\newcommand{\fkX}{{\mathfrak X}}

\newcommand{\calA}{\mathcal{A}}

\newcommand{\calC}{\mathcal{C}}
\newcommand{\calD}{\mathcal{D}}
\newcommand{\calE}{\mathcal{E}}
\newcommand{\calF}{\mathcal{F}}
\newcommand{\calG}{\mathcal{G}}
\newcommand{\calH}{\mathcal{H}}
\newcommand{\calI}{\mathcal{I}}
\newcommand{\calJ}{\mathcal{J}}

\newcommand{\calL}{\mathcal{L}}
\newcommand{\calM}{\mathcal{M}}
\newcommand{\calN}{\mathcal{N}}
\newcommand{\calO}{\mathcal{O}}
\newcommand{\calP}{\mathcal{P}}

\newcommand{\calS}{\mathcal{S}}
\newcommand{\calT}{\mathcal{T}}
\newcommand{\calU}{\mathcal{U}}
\newcommand{\calV}{\mathcal{V}}
\newcommand{\calW}{\mathcal{W}}
\newcommand{\calX}{\mathcal{X}}
\newcommand{\calY}{\mathcal{Y}}
\newcommand{\calZ}{\mathcal{Z}}

\newcommand{\Vect}{\mathrm{Vect}}
\newcommand{\andR}{{\mathrm{an}\text{-}\mathrm{dR}}}

\newcommand{\univext}[1]{E(#1)}

\title{Moduli stacks of crystals and isocrystals}
\author{Gyujin Oh}
\address{Department of Mathematics, Columbia University, 2990 Broadway, New York, NY 10027, USA}
\email{gyujinoh@math.columbia.edu}
\author{Koji Shimizu}
\address{Yau Mathematical Sciences Center, Tsinghua University, Beijing 100084, China;
Beijing Institute of Mathematical Sciences and Applications, Beijing 101408, China}
\email{shimizu@tsinghua.edu.cn}

\date{\today}

\begin{document}

\begin{abstract}
Given a liftable smooth proper variety over $\mathbb{F}_p$, we construct the moduli stacks of crystals and isocrystals on it. We show that the former is a formal algebraic stack over $\mathbb{Z}_p$ and the latter is an adic stack---Artin stack in rigid geometry---over $\mathbb{Q}_p$. Both stacks come equipped with the Verschiebung endomorphism $V$ corresponding to the Frobenius pullback of (iso)crystals. {We study the geometry of the $V$-fixed points over the open substack of irreducible isocrystals, which we use to geometrically count the rank one $F$-isocrystals.
Along the way, we carefully develop the theory of adic stacks.}
\end{abstract}

\maketitle
\tableofcontents

\section{Introduction}

In \cite{Simpson1,Simpson2}, Simpson constructed and studied three coarse moduli spaces of objects on a smooth projective complex variety $Z$ over $\bC$: the Betti moduli space $M_B(Z)$ of $\C$-local systems on $Z(\C)$; the de Rham moduli space $M_\dR(Z)$ of vector bundles with integrable connections on $Z$; and the Dolbeault moduli space $M_{\mathrm{Dol}}(Z)$ of semistable Higgs bundles on $Z$ with vanishing rational Chern classes. Here, Betti, de Rham, and Dolbeault refer to the cohomology theories whose \emph{coefficient} objects are parametrized by the respective moduli spaces. 
Simpson also geometrically realized, for example, the Riemann--Hilbert correspondence $M_B(Z)^\an\cong M_\dR(Z)^\an$ and the $\bC^{\times}$-action on $M_B(Z)$ induced by the scaling of the Higgs field. 

In this paper, we study the analogous moduli spaces of coefficient objects for \emph{crystalline cohomology}. Let $Z$ be a liftable smooth proper variety over a finite field $k=\F_{q}$ of characteristic $p$, let $W\coloneqq W(k)$ be the ring of $p$-typical Witt vectors of $k$, and set $K\coloneqq \cO_{K}[p^{-1}]$. The crystalline cohomology $H^i_\cris(Z/W)$ of $Z$ is defined as the cohomology of the structure sheaf $\calO_{Z/W}$ on the crystalline site $(Z/W)_\cris$, and the rationalization $H^i_\cris(Z/K)\coloneqq H^i_\cris(Z/W)[p^{-1}]$ gives a good $p$-adic cohomology theory for $Z$.
The coefficient objects for $H^i_\cris(Z/W)$ and $H^i_\cris(Z/K)$ are called \emph{crystals} and \emph{isocrystals}, respectively. 
We construct the moduli spaces of crystals and isocrystals on $Z$ as reasonable geometric objects:

\begin{thm-defn}[Rough form of Definitions~\ref{defn:Mcris} and \ref{defn:Misoc}]
The moduli stack $\calM_\cris(Z)$ of crystals on $Z$ exists as a formal algebraic stack over $W$, and the moduli stack $\calM_\isoc(Z)$ of isocrystals on $Z$ exists as an adic stack over $K$.
\end{thm-defn}

It requires a certain amount of care to correctly define the geometric setup for the moduli stacks $\cM_{\cris}(Z)$ and $\cM_{\isoc}(Z)$. In particular, we develop a theory of adic stacks---an analogue of Artin stacks in the setting of rigid geometry over $K$. We describe the theory in more detail later in the introduction.

The moduli stack of crystals $\cM_{\cris}(Z)$ has the expected moduli description; namely, for a $p$-adic formal scheme $S$, $\Mor(S,\cM_{\cris}(Z))$ is a groupoid of crystals on the crystalline site $(Z_{S_{1}}/S)_{\cris}$, where $Z_{S_{1}}\coloneqq Z\times_{W}S$ (see Theorem~\ref{thm:moduli interpretation of Mcris}). Its generic fiber $\calM_\cris(Z)_\eta$ is well-defined as an adic stack and admits an epimorphism
\[
\calM_\cris(Z)_\eta \rightarrow \calM_\isoc(Z)
\]
of adic stacks; this may well be called the \emph{passage-to-isogeny} map, reflecting the fact that the category of isocrystals on $Z$ is defined as the isogeny category of crystals. Moreover, for any finite extension $L$ of $K$, $\Mor(\Spa(L,\cO_{L}),\cM_{\isoc}(Z))$ is the groupoid of isocrystals with $L$-structure on $Z$ (see Corollary~\ref{cor:moduli interpretation of classical points of Misoc}).

The construction of these moduli stacks is in line with the fact that crystalline cohomology of $Z$ can be computed by de Rham cohomology of a smooth proper lift $Z_W/W$ of $Z$. In particular, $\calM_\isoc(Z)$ is defined to be an open adic substack of the analytification of the moduli stack $\calM_\dR(Z_K/K)$ of integrable connections on the generic fiber $Z_K$. We will return to this later in the introduction.

\bigskip
\noindent
\textbf{Verschiebung endomorphism and $F$-isocrystals.}
More subtle and interesting is the \emph{Verschiebung endomorphism} on these moduli stacks. Recall that $H^i_\cris(Z/W)$ comes with the Frobenius endomorphism, whose characteristic polynomial computes the zeta function of $Z$. This endomorphism is  defined by the functoriality of the crystalline topos $\Sh((Z/W)_\cris)$ along the $q$-Frobenius $\Fr_{q,Z}$ on $Z$. In the same way, the categories of crystals and isocrystals on $Z$ admit the Frobenius pullback endofunctor $\Fr_{q,Z/W,\cris}^\ast$. Correspondingly, we realize $\Fr_{q,Z/W,\cris}^\ast$ as the endomorphisms
\[
V\colon \cM_{\cris}(Z)\rar\cM_{\cris}(Z)\quad\text{and}\quad  V\colon\cM_{\isoc}(Z)\rar\cM_{\isoc}(Z)
\]
of a formal algebraic stack over $W$ and an adic stack over $K$, respectively, and call them the Verschiebung endomorphisms (Theorems~\ref{thm:Verschiebung on Mcris} and \ref{thm:Verschiebung on Misoc}).

We believe that the study of the geometry of $\cM_{\isoc}(Z)$ and the Verschiebung endomorphism $V$ is both interesting and feasible in many cases. The importance comes from the fact that the fixed points of $V$ on $\cM_{\isoc}(Z)$ correspond to \emph{$F$-isocrystals} on $Z$.

An \emph{$F$-isocrystal} on $Z$ is an isocrystal $\calF$ on $Z$ together with an isomorphism $\varphi_\calF\colon \Fr_{q,Z/W,\cris}^\ast\calF\xrightarrow{\cong}\calF$. Here the extra structure $\varphi_\calF$ is crucial for arithmetic questions; it defines the Frobenius endomorphism on the cohomology $H^i_\cris(Z,\calF)$ and the Frobenius characteristic polynomial at every closed point of $Z$.

Deligne's \emph{companion conjecture} entails, among many other things, the statement that every irreducible $F$-isocrystal with $\overline{\Q}_p$-structure (also called $\overline{\Q}_p$-$F$-isocrystal) of finite order determinant is pure of weight zero, and, for every prime $\ell\neq p$, there exists a lisse $\overline{\Q}_\ell$-Weil sheaf on $Z$ that has the same Frobenius characteristic polynomial at every closed point.\footnote{Here $Z$ is assumed to be smooth and proper over $k$. In particular, the notion of $F$-isocrystals on $Z$ defined above agrees with that of \emph{overconvergent} $F$-isocrystals on $Z$. The latter is the coefficient for rigid cohomology and the correct one in the formulation of the companion conjecture; see also Remark~\ref{rem:F-isocrystals and overconvergent/convergent ones}.}
This $p$-to-$\ell$ companion statement is now a theorem: the curve case is proved by Abe \cite{Abe-Langlands}, following (and combined with) the works of Drinfeld and L.~Lafforgue on the Langlands correspondence for $\GL_n$ over the function field \cite{Drinfeld-Langlands-GL2, LLafforgue-Langlands}; the higher-dimensional case is obtained independently by Abe--Esnault \cite{Abe-Esnault} and Kedlaya \cite{Kedlaya-companionI}, deduced from the curve case using Drinfeld's theorem of constructing a lisse sheaf on $Z$ from those on curves over $Z$ \cite{Drinfeld-conjecture-of-Deligne}. For the detail, we refer the reader to these references and \cite{Deligne-WeilII, Crew-F-isocrystals-monodromy-groups, Cadoret-Bourbaki}.

Under this companion conjecture, our picture is analogous to the complex analytic picture of Simpson, where the fixed points of the scaling $\C^{\times}$-action on the Betti moduli space correspond to the local systems underlying complex variations of Hodge structures.

\bigskip

Let us now summarize the main theorems in this paper regarding the geometry of $\cM_{\isoc}(Z)$ and its Verschiebung endomorphism.

\bigskip
\noindent
\textbf{Irreducible locus and $V$-fixed locus.}
The adic stack $\cM_{\isoc}(Z)$ is a \emph{stack}, as each isocrystal has non-trivial endomorphisms. We show that the \emph{absolutely irreducible} isocrystals on $Z$ form an open substack $\cM_{\isoc,\irr}(Z)\subset\cM_{\isoc}(Z)$, which we call the \emph{irreducible locus}. The adic stack $\cM_{\isoc,\irr}(Z)$ has the advantage of being only mildly stacky, and we can discuss the geometry of the Verschiebung fixed points:

\begin{thm}[Rough form of Theorems~\ref{thm:Misocirr is Gm-gerbe} and \ref{thm:irreducible F-isocrystals discrete}]
The substack $\cM_{\isoc,\irr}(Z)$ of absolutely irreducible isocrystals on $Z$ admits a coarse moduli space $M_{\isoc,\irr}(Z)$. 
Furthermore, the space $M_{\isoc,\irr}^{V=\id}(Z)$ of the Verschiebung fixed points is discrete and of the form $\coprod_{i\in I}\Spa(L_{i},\cO_{L_{i}})$, where each $L_{i}$ is a finite extension of $K$. 
\end{thm}

We remark that the set of $\overline{\Q}_p$-valued points of $M_{\isoc,\irr}^{V=\id}(Z)$, namely, the set
\[
\varinjlim_{L/K\;\text{finite}}M_{\isoc,\irr}^{V=\id}(Z)(\Spa(L,\calO_L)) 
\]
is identified with the set of isomorphism classes of $\overline{\Q}_p$-isocrystals on $Z$ underlying irreducible $\overline{\Q}_p$-$F$-isocrystals (see Definition~\ref{def:MFisoc}).

In the above, $M_{\isoc,\irr}(Z)$ is an \emph{\'etale space}---a rigid geometry analogue of an algebraic space, and the formulation of the coarse moduli space involves the theory of \emph{$\Gm^{\an}$-gerbes} in the context of adic stacks. By constructing a coarse moduli space, we can argue geometrically, via deformation theory and the theory of weights for $F$-isocrystals, that the $V$-fixed points are discrete in the coarse moduli space.

\bigskip
\noindent
\textbf{Counting rank one $F$-isocrystals.}
There is a  stratification of $M_{\isoc,\irr}^{V=\id}(Z)$ by the rank of an isocrystal, 
\[
M_{\isoc,\irr}^{V=\id}(Z)=\coprod_{r\in\bN}M_{\isoc,\irr,\rank r}^{V=\id}(Z).
\]
For each $r$, we deduce from \cite{Abe-Esnault,Kedlaya-companionI} that $M_{\isoc,\irr,\rank r}^{V=\id}(Z)$ is \emph{finite} \'etale over $K$ since its degree $N(Z,r)$ is the number of absolutely irreducible $\overline{\Q}_{p}$-$F$-isocrystals of rank $r$ over $Z$ up to constant twists, which is also the number of rank $r$ absolutely irreducible $\overline{\Q}_{\ell}$-local systems on $Z$ up to constant twists by the companion conjecture. 

Starting from the pioneering work of Drinfeld   \cite{Drinfeld-counting}, there has been a considerable progress in establishing the formula computing $N(Z,r)$ in various contexts (e.g. \cite{DeligneFlicker,Yu-unramified}). One of the most intriguing features of these formulas is that $N(Z_{\bF_{q^{m}}},r)$ involves Weil numbers and thus resembles the number of $\bF_{q^{m}}$-points of an $\bF_{q}$-variety $V$. Deligne \cite[2.18]{Deligne-comptage} suggests that such a formula would follow from the existence of an analytic open subset of the de Rham moduli space of $Z_K$, equipped with a ``Verschiebung endomorphism'' $V$, such that $N(Z_{\bF_{q^{m}}},r)$ is given by the trace of $V^{m}$ acting on a suitable cohomology of that open subset. Our work realizes the first part at the level of stacks.

Elaborating further on our construction of $\cM_{\isoc}(Z)$ and its Verschiebung endomorphism, we geometrically derive a formula of the above form for the \emph{rank $1$ part} $M_{\isoc,1}^{V=\id}(Z)\coloneqq M_{\isoc,\irr,\rank 1}^{V=\id}(Z)\subset M_{\isoc,\irr}^{V=\id}(Z)$.

\begin{thm}[See Theorems~\ref{thm:counting rank 1 F-isocrystals},~\ref{thm:counting rank 1 F-isocrystals Fqm}]
Suppose that $Z$ is a curve or an abelian variety over $k=\bF_{q}$. Then, $M_{\isoc,1}^{V=\id}(Z)$  is finite \'etale over $K$ of degree equal to the number of $k$-points of $\Pic^{0}(Z/k)$. In particular, 
\[
N(Z_{\bF_{q^{m}}},1)=\lvert\Pic^{0}(Z/k)(\bF_{q^{m}})\rvert.
\]
\end{thm}

The virtue of the above theorem is that we arrive at the correct formula of $N(Z_{\bF_{q^{m}}},1)$ by only using  geometric arguments. Furthermore, the proof manifests another interesting feature, which we may call the \emph{dimension reduction} phenomenon: namely, the variety whose rational points count $N(Z_{\bF_{q^{m}}},1)$ is $\Pic^{0}(Z/k)$, whose dimension is half the dimension of the moduli space of integrable connections. We obtain these results by considering an analogous coarse moduli space $\calM_{\cris,1}(Z)$ as a formal model of the adic stack $\calM_{\isoc,1}(Z)$ of rank one isocrystals and studying the mod $p$ reduction of the latter via the \emph{Cartier descent}; see Theorem~\ref{thm:geometric description of M1cris,BunV=id}.

\bigskip
\noindent
\textbf{Adic stacks.}
Let us discuss the foundations of stacks. Our main theorem on the moduli spaces of crystals and isocrystals needs the theories of algebraic spaces and Artin stacks in formal geometry and non-archimedean geometry. For the former theory, we simply use the works of Fujiwara--Kato, Stacks Project, and Emerton \cite{Fujiwara-Kato, stacks-project, EmertonStack}. For the latter, several authors have already worked out such theories using different foundations of non-archimedean geometry; to name a few, see \cite{ConradTemkin1} (rigid spaces and Berkovich spaces),  \cite{Ulirsch} (Berkovich spaces), \cite{Hellmann, Warner} (adic spaces); see also \cite{Bambozzi-Ben-Bassat, Ben-Bassat-Kremnizer, Paugam, Porta-Yu-higheranalyticstacks, Porta-Yu-derived} for approaches with more emphasis on derived non-archimedean geometry as well as the recent works of Clausen--Scholze and Rodr\'iguez Camargo for the approach using condensed mathematics. 

We mainly follow Warner's work \cite{Warner} using Huber's adic spaces. 
{We will write $\cO_{K}\coloneqq W$, to emphasize that the theory of adic stacks applies over any complete discrete valuation field of mixed characteristic $(0,p)$.}
Let $\Afd_K$ denote the category of affinoid adic spaces $\Spa(A,A^+)$ such that the completion $(\widehat{A},\widehat{A}^+)$ is topologically of finite type over an analytic affinoid field $(L,L^+)$ over $(K,\calO_K)$. We equip it with \'etale topology and write $(\Afd_K)_\Et$ for the resulting site. There is a fully faithful embedding from the category $\Rig_K$ of \emph{rigid $K$-spaces} (by which we mean adic spaces locally of finite over $K$) to the category $\Sh((\Afd_K)_\Et)$ of sheaves. Then \'etale spaces and adic stacks are defined as follows (here we omit to define the notions and properties of representable morphisms).

\begin{defn}[Definitions~\ref{def:etale spaces}, \ref{def:adic stacks}]\hfill
\begin{enumerate}
 \item An \emph{\'etale space} over $K$ is a sheaf $F$ on $(\Afd_K)_\Et$ such that the diagonal $\Delta\colon F\rightarrow F\times F$ is representable by adic spaces and $F$ receives an \'etale surjection from a rigid $K$-space.
 \item An \emph{adic stack} over $K$ is a category $\calX\rightarrow (\Afd_K)_\Et$ fibered in groupoids such that $\calX$ is a stack in groupoids, the diagonal $\Delta\colon \calX\rightarrow \calX\times \calX$ is representable by \'etale spaces, and $\calX$ receives a $1$-morphism from a rigid $K$-space that is representable by \'etale spaces and is a smooth surjection.
\end{enumerate}
\end{defn}

There are several possible choices of the underlying category in place of $\Afd_K$. In fact, Warner \cite{Warner} defines the notion of an admissible category on which he develops the theory of \'etale spaces and adic stacks. Alternatively, one can work on the category $\Rig_K$.
There are mainly two reasons for our choice $\Afd_K$: firstly, it contains $\Spa(L,L^+)$ for any analytic affinoid field $(L,L^+)$ over $(K,\calO_K)$; secondly, the analytification functor and the generic fiber functor are well described on affinoids, which we now explain.

To a scheme $Y$ locally of finite type over $K$, one can functorially associate a rigid $K$-space $Y^\an$. Then $Y^\an$ as a sheaf on $\Afd_K$ is the sheafification of the presheaf sending $\Spa(A,A^+)$ to $Y(\Spec \widehat{A})$. We show that for an algebraic space $F$ locally of finite type over $K$, the sheafification $F^\an$ of the presheaf sending $\Spa(A,A^+)$ to $F(\Spec \widehat{A})$ is an \'etale space over $K$ (Proposition~\ref{prop:analytification of algebraic spaces}), which we call the \emph{analytification} of $F$. Similarly, we define the analytification $\calY^\an$ of an algebraic stack $\calY$ locally of finite type over $K$ as the stackification of the category fibered in groupoids over $(\Afd_K)_\Et$ whose fiber category over $\Spa(A,A^+)$ is $\calY_{\Spec \widehat{A}}$; then we show that $\calY^\an$ is an adic stack over $K$ (Proposition~\ref{prop:analytification of algebraic stacks}).

To define the generic fiber of formal algebraic spaces and stacks, we proceed as follows. Recall that to a formal scheme $Y$ locally formally of finite type over $\calO_K$, one can functorially associate a rigid $K$-space $Y_\eta$ called the generic fiber of $Y$. When regarded as a sheaf on $\Afd_K$, $Y_\eta$ is described as the sheafification of the presheaf 
\[
\Spa(A,A^+)\mapsto \varinjlim_{A_0\subset \widehat{A}^+}Y(\Spf A_0),
\]
where $A_0$ runs over the open and bounded $\calO_K$-subalgebras of  $\widehat{A}^+$. Similarly, we define the \emph{generic fiber} $\calY_\eta$ of a formal algebraic space or stack $\calY$ locally formally of finite type over $\calO_K$ and show that $\calY_\eta$ is an \'etale space or adic stack over $K$, respectively (Propositions~\ref{prop:generic fiber of formal algebraic spaces} and \ref{prop:generic fiber of formal algebraic stacks}).

We can compare these two constructions in the following situation. Let $\calY$ be an algebraic stack locally of finite type over $\calO_K$. To $\calY$, one can associate the $p$-adic completion $\widehat{\calY}$, which is a formal algebraic stack over $\calO_K$, and the generic fiber $\calY_K$, which is an algebraic stack over $K$. 
When $\calY$ is a scheme $Y$, we have a morphism $f\colon (\widehat{Y})_\eta\rightarrow (Y_K)^\an$ of rigid $K$-spaces, which is, locally on the source, an isomorphism; it is even an open immersion if $Y$ is separated over $\calO_K$. Similarly, we have a $1$-morphism $f\colon \widehat{\calY}_\eta \rightarrow \calY_K^\an$ of adic stacks over $K$, generalizing the case of schemes (see Definition~\ref{defn:comparison map from generic fiber to analytification}). However, $f$ can also exhibit a stacky phenomenon as in the following example.

\begin{example}
When $\calY=B\GL_r=[\ast/\GL_r]$, the classifying stack of the general linear group over $\calO_K$, we obtain
$[\ast/(\widehat{\GL_r})_\eta]\rightarrow [\ast/\GL_r^\an]$.
Intuitively, if $\calO_K=\Z_p$, this map is a passage-to-isogeny map sending a family of rank $r$ free $\Z_p$-modules to the associated family of $r$-dimensional $\Q_p$-vector spaces.
\end{example}

In general, $f\colon \widehat{\calY}_\eta \rightarrow \calY_K^\an$ is a mixture of these two examples; Proposition~\ref{prop:image of comparison map from generic fiber to analytification} shows that $f$ factors as an epimorphism followed by an open immersion. This is exactly the feature we want to use to construct the stack $\calM_\isoc(Z)$ of isocrystals on $Z$ mentioned earlier.

\bigskip
\noindent
\textbf{Construction of $\cM_{\cris}(Z)$ and $\cM_{\isoc}(Z)$.} 
Recall that $Z\rightarrow\Spec k$ is smooth proper and assumed to lift to a smooth proper scheme $Z_W\rightarrow \Spec W$ (now we write again $\calO_K=W\coloneqq W(k)$). In the body of this article, we also work on the logarithmic setting when $Z_W$ has a relative smooth normal crossings divisor $D_W$, but we assume $D_W=\emptyset$ here for simplicity. For any $p$-adic formal scheme $S$, set $S_1\coloneqq S\times_Wk$ and $Z_S\coloneqq Z_{W}\times_{W}S$.

Let $\calM_{\dR}(Z_W/W)\rightarrow (\Sch/W)_\fppf$ denote the moduli stack of integrable connections on $Z_W/W$ (Definition~\ref{def:moduli stack of integrable connections}). It is well known that the category of crystals of $\calO_{Z_{S_1}/S}$-modules on the crystalline site $(Z_{S_1}/S)_\cris$ is equivalent to the category of integrable $S$-connections on $Z_S$ that are topologically quasi-nilpotent (Proposition~\ref{prop:crystals and connections in p-adic case}). 
So the stack $\calM_\cris(Z)$ of crystals is a strictly full subcategory of $\calM_{\dR}(Z_W/W)$ cut out by an integrable connection $(\calE,\nabla)$ on $Z_S$ being topologically quasi-nilpotent. We use the fact that this is equivalent to the condition that the mod $p$ reduction of $(\calE,\nabla)$ is an integrable connection on $Z_{S_1}$ with \emph{nilpotent $p$-curvature} (or, in short, \emph{nilpotent} integrable connection; see Definition~\ref{def:p-curvature}).

Therefore, to define the moduli of crystals $\cM_{\cris}(Z)$, we need to analyze the strictly full subcategory $\cM_{\dR}^{\nilp}(Z/k)$ (called the \emph{nilpotent locus}) of $\calM_\dR(Z/k)$ consisting of nilpotent integrable $T$-connections on $Z_{T}$. While $\cM_{\dR}(Z/k)$ is an algebraic stack,  $\cM_{\dR}^{\nilp}(Z/k)$ is a \emph{formal algebraic stack} over $k$;
we will prove the latter result by first showing that its reduction $(\cM_{\dR}^{\nilp}(Z/k))_\red$ is a closed substack (Theorem~\ref{thm:reduced closed substack of nilpotent connections}). This requires a universal bound of the exponent of nilpotence for any nilpotent endomorphism on a fixed coherent sheaf on a reduced scheme; see Proposition~\ref{prop:bounding exponent of nilpotence} for the details.

Once these results are established, $\calM_\cris(Z)$ is defined to be the formal completion of $\calM_\dR(Z_W/W)$ along $(\cM_{\dR}^{\nilp}(Z/k))_\red$. The resulting formal algebraic stack $\calM_\cris(Z)$ has the desired crystalline moduli interpretation and is equipped with the Verschiebung endomorphism corresponding to the Frobenius pullback of crystals (Theorems~\ref{thm:moduli interpretation of Mcris} and \ref{thm:Verschiebung on Mcris}).

To construct the adic stack $\calM_\isoc(Z)$, recall the comparison map
\[
(\calM_{\dR}(Z_W/W)^\wedge)_\eta\rightarrow \calM_{\dR}(Z_K/K)^\an
\]
given by the general formalism of adic stacks discussed above. We see that $1$-morphism $\calM_\cris(Z)_\eta\rightarrow (\calM_{\dR}(Z_W/W)^\wedge)_\eta$ is an open immersion, and the essential image of the composite $\calM_\cris(Z)_\eta\rightarrow \calM_{\dR}(Z_K/K)^\an$ defines an open adic substack of $\calM_{\dR}(Z_K/K)^\an$, which is $\calM_\isoc(Z)$ (see Definition~\ref{defn:Misoc} for the precise definition). Then we are able to show that $\calM_\isoc(Z)$ has the desired moduli interpretation on $L$-valued points for $L/K$ finite (Corollary~\ref{cor:moduli interpretation of classical points of Misoc}) and that the endomorphism $V_\eta$ on $\calM_\cris(Z)_\eta$ descends and defines the Verschiebung $V$ on $\calM_\isoc$ 
(Theorem~\ref{thm:Verschiebung on Misoc}).

\bigskip
\noindent
\textbf{Related works.}
Our work on the rank $1$ moduli of isocrystals is largely inspired by the discussions of \cite{Deligne-IHES,Deligne-comptage}, where a proof of a result analogous to Theorem~\ref{thm:counting rank 1 F-isocrystals} is sketched in the case where $p>2$ and $Z_{W}/W$ is a relative curve; this is further elaborated in \cite{Katsigianni}. These works rely on the nilpotent crystalline site interpretation of the universal vector extensions.

Relevant but slightly different directions involving moduli stacks of crystals and isocrystals are pursed in other works. Kedlaya \cite{Kedlaya-companionII} works on the $\ell$-to-$p$ part of the companion conjecture by introducing a system of moduli stacks capturing $F$-crystals. Esnault and Groechenig \cite{EG-revisited} work on the cristallinity of rigid connections. One of their key ingredients is the set $\calM_{\dR}(W)^\mathrm{iso}$ of isomorphism classes on the $W$-points of $\calM_\dR(Z_W/W)$ for $p>2$, which is equipped with the natural $p$-adic topology and the Frobenius pullback map $F^\ast$ (the Verschiebung endomorphism).\footnote{By direct computation as in \emph{op.~cit.} or nilpotent crystalline interpretation.} Among other things, they prove that $F^\ast$ is an open mapping, in particular, injective; on the contrary, it is not difficult to see that the Verschiebung $V$ of our work is not a $1$-monomorphism of stacks (e.g., not a monomorphism on $\C_p$-points). These two results illustrate an interesting feature of $V$.

There are recent advances in $p$-adic Hodge theory where coefficient objects on a space $X$ (either a scheme over a finite field or a $p$-adic formal scheme) of a cohomology theory are realized as quasi-coherent sheaves on a \emph{stack} (in fpqc topology) associated to $X$. This idea, commonly referred to as the \emph{stacky approach} to coefficient objects, dates back to Simpson \cite{Simpson-dR}, where he introduces the \emph{de Rham stack} of a smooth variety $X/\bC$, whose category of quasi-coherent sheaves is naturally identified with the category of quasi-coherent sheaves with integrable connections on $X$. In particular, the stacky approach to crystalline (iso)crystals is realized in \cite{Drinfeld-stacky}; see also \cite{Drinfeld-Prismatization,BhattLurie-APC,BhattLurie-Prismatization,Bhatt-Fgauges} for the stacky approach for other types of crystals. These works primarily focus on realizing (iso)crystals as quasi-coherent sheaves on a certain stack, whereas our goal is to construct the moduli stack of (iso)crystals.

\bigskip
\noindent
\textbf{Outline.}
Section~\ref{sec:formal algebraic stacks} reviews the theory of formal algebraic spaces and formal algebraic stacks in the existing literature. In Section~\ref{sec:adic stacks}, following and complementing \cite{Warner}, we develop the theory of \emph{\'etale spaces} and \emph{adic stacks}, which are the analogues of algebraic spaces and algebraic stacks in non-archimedean geometry, respectively. Section~\ref{sec:analytification and generic fiber} defines the analytification of an algebraic stack and the generic fiber of a formal algebraic stack, which are both adic stacks. Furthermore, given an algebraic stack $\cY$ locally of finite type over $\bZ_{p}$, we study the $1$-morphism from the generic fiber $(\wh{\cY})_{\eta}$ of the $p$-adic completion to the analytification $(\cY_{\bQ_{p}})^{\an}$ of the algebraic generic fiber (Proposition~\ref{prop:image of comparison map from generic fiber to analytification}).

In Section~\ref{sec:Moduli stacks of integrable connections}, we first review the algebraic moduli stack of integrable connections $\cM_{\dR}$, and prove various results about it. For example, for a smooth proper variety $Z$ over a finite field $k$, we prove that the \emph{nilpotent locus} $\cM_{\dR}^{\nilp}(Z/k)$ is a formal algebraic stack (Theorem~\ref{thm:reduced closed substack of nilpotent connections}). Furthermore, we prove that, for a smooth proper variety $Z_K$ over a complete non-archimedean field $K$ of mixed characteristic, the analytification of $\cM_{\dR}(Z_K/K)$ is naturally identified with the analytic moduli stack of (logarithmic) integrable connections $\cM_{\andR}(Z_K^{\an}/K)$ (Theorem~\ref{thm:analytification of MdR}). In Section~\ref{sec:crystals and isocrystals on crystalline site}, we recall the relevant facts about crystals  and isocrystals on the (logarithmic) crystalline site. Then Section~\ref{section:moduli of crystals and isocrystals} defines the moduli stack of crystals $\cM_{\cris}(Z)$ and the moduli stack of isocrystals $\cM_{\isoc}(Z)$ for a liftable smooth proper variety $Z/k$. We show that they have the desired properties, notably the moduli descriptions and the existence of the Verschiebung endomorphisms.

As geometric applications of our constructions of $\cM_{\cris}(Z)$ and $\cM_{\isoc}(Z)$, we pursue two directions. Section~\ref{sec:geometry of irreducible locus} studies the open substack $\cM_{\isoc,\irr}(Z)\subset\cM_{\isoc}(Z)$ (called the \emph{irreducible locus}) corresponding to the \emph{absolutely irreducible} isocrystals, and show that it is a \emph{$\Gm^{\an}$-gerbe}. This defines the associated coarse moduli space $M_{\isoc,\irr}(Z)$, which is an \'etale space. We then show that the Verschiebung-fixed points $M_{\isoc,\irr}^{V=\id}(Z)$ is \'etale over $K$ (Theorem~\ref{thm:irreducible F-isocrystals discrete}). In Section~\ref{sec:rank 1}, we describe the geometry of $\cM_{\isoc}(Z)$ in the case of rank one isocrystals, when $Z/k$ is a curve or an abelian variety. This gives a geometric way to count $\overline{\Q}_{p}$-$F$-isocrystals of rank one up to constant twists (Theorems~\ref{thm:counting rank 1 F-isocrystals} and \ref{thm:counting rank 1 F-isocrystals Fqm}). For this, we review the universal vector extensions in Section~\ref{sec:universal vector extension} and discuss the Cartier descent in Section~\ref{sec:Moduli stacks of rank 1 crystals and isocrystals}. 

Finally, Appendix \ref{sec:appendix} features several general facts about sites and stacks.

\bigskip
\noindent
\textbf{Notation and conventions.}
By complete, we mean Hausdorff complete, namely, complete and separated.

For a ring $A$, let $\Mod_A$ denote the category of $A$-modules and write $\Mod_A^\mathrm{fg}$ for the full subcategory of finitely generated $A$-modules. Let $\Alg_A$ denote the category of $A$-algebras.

For a ringed space or ringed site $(X,\calO_X)$, let $\Mod(\calO_X)$ denote the category of sheaves of $\calO_X$-modules and write $\Coh(\calO_X)$ (resp.~$\Vect(\calO_X)$) for the full subcategory of coherent $\calO_X$-modules (resp.~ finite locally free $\calO_X$-modules); see \cite[01BV, 03DL]{stacks-project}. If there is no confusion, we also write $\Coh(X)$ (resp.~$\Vect(X)$).

We mainly follow \cite{stacks-project, EmertonStack} for the foundations of the theory of algebraic spaces, algebraic stacks (i.e., Artin stacks), and formal algebraic stacks. In particular, an algebraic space (resp.~algebraic stack) over a scheme $S$ is an fppf sheaf (resp.~stack) on $(\Sch/S)$; note that the theory using fppf topology coincides with the one using \'etale topology as in \cite[076A, 076V]{stacks-project}. However, we prefer to fix a universe $\mathscr{U}$; a site means a $\mathscr{U}$-site, and fibered categories are in $\mathscr{U}$. This convention is different from that of \cite{stacks-project, EmertonStack}, but we still cite these references.

\bigskip
\noindent
\textbf{Acknowledgments.}
We thank Andres Fernandez Herrero for several informative discussions, and especially for suggesting the proof of Lemma~\ref{lem:controlling torsion filtration} and how to deduce Theorem~\ref{prop:bounding exponent of nilpotence} from Corollary~\ref{cor:bounding exponent of nilpotence, projective}. 

A part of this material was completed while the  first author was in residence at the Mathematical Sciences Research Institute (MSRI/SLMath) in Berkeley, California, during the Spring 2023
semester, supported by the National Science
Foundation under Grant No.~DMS-1928930.
The second author was partially supported by AMS--Simons Travel Grant at the early stage of the project and by the NSFC Excellent Young Scientists Fund Program (Overseas) at the later stage.

\section{Formal algebraic spaces and stacks}\label{sec:formal algebraic stacks}

This section reviews the theory of formal algebraic spaces and formal algebraic stacks, following \cite{Fujiwara-Kato, stacks-project, EmertonStack}.

\subsection{Classical formal schemes}
The main references are \cite{EGAI, EGAI-Springer, Abbes-EGRI, Fujiwara-Kato, stacks-project}.
Recall that a \emph{linearly topologized ring} is a topological ring such that $0$ has a fundamental system of neighborhoods consisting of ideals. An ideal $I$ of a linearly topologized ring $A$ is said to be an \emph{ideal of definition} if it is open and every open neighborhood of $0$ contains $I^n$ for some $n$. We say that $A$ is \emph{admissible} if it is complete and admits an ideal of definition.

To an admissible topological ring $A$, one can associate a topologically locally ringed space $\Spf A$ (called the \emph{formal spectrum} of $A$) whose underlying topological space is the closed subset of $\Spec A$ consisting of open prime ideals and whose structure sheaf is given by $\calO_{\Spf A}\coloneqq \varprojlim_{I}\widetilde{A/I}|_{\Spf A}$ where $I$ runs through all the ideals of definition: see \cite[\S10.1]{EGAI}, \cite[\S10.1]{EGAI-Springer}, \cite[\S2.1]{Abbes-EGRI}, \cite[I.1.1.(b)]{Fujiwara-Kato}, or \cite[0AHY]{stacks-project} for details.

\begin{defn}
An \emph{affine formal scheme} is a topologically locally ringed space that is isomorphic to $\Spf A$ for an admissible topological ring $A$. A \emph{formal scheme} is a topologically locally ringed space that is locally isomorphic to an affine formal scheme. A morphism of formal schemes is defined to be a morphism of topologically locally ringed spaces. Let $\FSch$ denote the category of formal schemes. 
\end{defn}

\begin{rem}\label{rem:first remarks on formal schemes}
For admissible topological rings $A$ and $B$, we have a functorial identification $\Hom_{\cont}(A,B)=\Hom_{\FSch}(\Spf B,\Spf A)$ by \cite[Prop.~10.2.2]{EGAI} or \cite[0AHZ]{stacks-project}.
Every scheme is naturally regarded as a formal scheme, and the resulting functor $\Sch\rightarrow \FSch$ is fully faithful (see \cite[0AHY]{stacks-project}).
\end{rem}

\begin{lem}[{\cite[0AI1, 0AI2]{stacks-project}}]\label{lem:classical formal scheme is an fpqc sheaf}
For a formal scheme $X$, define the contravariant functor $h_\fkX\colon \Sch\rightarrow \Set$ by $h_X(S)\coloneqq\Hom_{\FSch}(S,X)$. Then, the assignment $X\rightarrow h_X$ yields a fully faithful functor $\FSch\rightarrow \PSh(\Sch)$. Moreover, $h_X$ is a sheaf with respect to fpqc topology.
\end{lem}

\begin{defn}[{\cite[Def.~I.1.1.3, I.1.1.6; I.1.1.14, I.1.1.16]{Fujiwara-Kato}}]\label{def:formal scheme is adic of finite ideal type}
\hfill
\begin{enumerate}
 \item An admissible ring $A$ is called \emph{adic} if the topology on $A$ is $I$-adic for some ideal $I$ (i.e., $\lbrace I^{n}\rbrace$ is a fundamental basis of $0\in A$). 
Moreover, it is called \emph{adic of finite ideal type} if $I$ can be taken to be finitely generated. 
 \item A formal scheme is called \emph{adic} (resp.~\emph{adic of finite ideal type}) if it admits an affine open covering consisting of $\Spf A$ with $A$ being adic (resp.~adic of finite ideal type).
\end{enumerate}
The affine formal scheme $\Spf A$ is adic of finite ideal type if and only if $A$ is adic of finite ideal type by \cite[Cor.~I.3.7.13]{Fujiwara-Kato} (or \cite[0GXQ]{stacks-project}); being adic of finite ideal type is also referred to as being \emph{locally adic*} in \cite[0AID]{stacks-project}.
\end{defn}

\begin{defn}[{\cite[Def.~I.1.1.18, Prop.~I.1.19]{Fujiwara-Kato}}]\label{def:ideal of definition for formal scheme}
Let $A$ be an admissible ring. For any open ideal $J\subset A$, we define the sheaf $J^\triangle$ on $\Spf A$ by
\[
J^\triangle=\varprojlim_{I\subset J}\widetilde{J/I},
\]
where $I$ runs over all ideals of definition in $J$. Then $J^\triangle$ is an sheaf of ideals of $\calO_{\Spf A}$ and satisfies $\Gamma(\Spf A,J^\triangle)=J$. Such a sheaf of ideals is called an \emph{ideal of definition} of $\Spf A$. Note that one can similarly define an $\calO_{\Spf A}$-module sheaf $M^\triangle$ for any $A$-module $M$.

An ideal sheaf $\calI$ on a formal scheme $X$ is said to be an \emph{ideal of definition} if for any $x\in X$, there exists an affine open neighborhood $\Spf A\subset X$ such that $\calI|_{\Spf A}=J^\triangle$ for some open ideal $J\subset A$. When $X$ is affine, this definition agrees with the one in the previous paragraph. An ideal of definition is called \emph{of finite type} if it is of finite type as an $\calO_X$-module. 
\end{defn}

\begin{defn}[{\cite[Def.~I.1.1.5.(2), I.1.3.1]{Fujiwara-Kato}, \cite[0GBR]{stacks-project}}]\label{def:adic morphisms between formal schemes}
\hfill
\begin{enumerate}
 \item A continuous map $A\rightarrow B$ between adic rings is called \emph{adic} if there exists an ideal definition $I$ of $A$ such that $IB$ is an ideal of definition of $B$.
 \item A morphism $f\colon X\rightarrow Y$ between adic formal schemes of finite ideal type is called \emph{adic} if $f$ is Zariski locally of the form $\Spf B\rightarrow \Spf A$ with the map $A\rightarrow B$ being adic.
\end{enumerate}
 By \cite[Prop.~I.1.3.2]{Fujiwara-Kato}, for adic rings $A,B$ of finite ideal type, $A\rightarrow B$ is adic if and only if the associated morphism $\Spf B\rightarrow \Spf A$ is adic.
\end{defn}

\begin{defn}[{\cite[Def.~I.1.5.(d)]{Fujiwara-Kato}}]\label{def:adically P}
Let $\calP$ be one of the following properties of morphisms between schemes: (locally) of finite type; (locally) of finite presentation; (faithfully) flat; finite; \'etale; smooth.

A morphism $f\colon X\rightarrow Y$ between adic formal schemes of finite ideal type is said to be \emph{adically $\calP$} if $f$ is adic and for any open covering $Y=\bigcup_\alpha Y_\alpha$ and any ideal $\calI_\alpha$ of definition of finite type of $Y_\alpha$, and for any $k\geq 0$, the induced morphism $(f^{-1}(Y_\alpha), \calO_{f^{-1}(Y_\alpha)}/\calI_\alpha^{k+1}\calO_{f^{-1}(Y_\alpha)})\rightarrow (Y_\alpha, \calO_{Y_\alpha}/\calI_\alpha^{k+1})$ of schemes satisfies $\calP$. This is equivalent to the condition with ``any'' replaced by ``some''. Due to the following Remark~\ref{rem:adically P vs P}, we often say $\calP$ instead of adically $\calP$ for an adic $f$.
\end{defn}

\subsection{Formal algebraic spaces}
Formal algebraic spaces are discussed in \cite{Knutson, Fujiwara-Kato, stacks-project}, and we mainly follow \cite{stacks-project}. Fix a base scheme $S$ and consider the big fppf site $(\Sch/S)_{\fppf}$. 

Recall the following in \cite[025V, 03XZ]{stacks-project}: let $\calP$ be a property of morphisms of algebraic spaces that is stable under base change and fppf local on the base. For a morphism $f\colon F\rightarrow G$ between presheaves of sets on $(\Sch/S)_{\fppf}$ that is representable by algebraic spaces, we say that $f$ \emph{has property $\calP$} if the morphism $F\times_{G,f}U\rightarrow U$ of algebraic spaces has $\calP$ for every $U\in (\Sch/S)_{\fppf}$ and $\xi\in G(U)$.

\begin{defn}[{\cite[0AI7]{stacks-project}}]
An \emph{affine formal algebraic space} over $S$ is a sheaf $U$ on $(\Sch/S)_{\fppf}$ that is expressed as a colimit of sheaves $U\cong\varinjlim_{\lambda\in \Lambda} U_{\lambda}$, where $\Lambda$ is a directed set, $U_{\lambda}$ is an $S$-affine scheme, and the transition morphisms are thickenings, i.e., closed immersions defined by nilpotent ideals.
\end{defn}

\begin{defn}[{\cite[0AIM]{stacks-project}}]\label{defn:formal algebraic spaces}
A \emph{formal algebraic space} over $S$ is a sheaf $F$ on $(\Sch/S)_{\fppf}$ such that there exist affine formal algebraic spaces $U_i$ over $S$ and morphisms $f_i\colon U_i\rightarrow F$ of sheaves satisfying the following properties:
\begin{enumerate}
\item each $f_i$ is representable by algebraic spaces and \'etale;
\item $\coprod_{i}U_{i}\rightarrow F$ is an epimorphism of sheaves.
\end{enumerate}
A morphism of formal algebraic spaces is defined to be a morphism of sheaves. 
\end{defn}

\begin{rem}
Strictly speaking, \cite{stacks-project} defines the big fppf site $(\Sch/S)_{\fppf}$ using a limit ordinal to avoid set-theoretic issues. Correspondingly, it also imposes a set-theoretic condition in the definition of affine formal schemes or formal schemes; namely, $X$ is a coequalizer of two morphisms between representable sheaves (see \cite[021R, 021S, 0AIS]{stacks-project}). We fix a universe $\mathcal{U}$ instead and ignore these aspects.
\end{rem}

\begin{rem}\label{rem:morphism from affine formal algebraic space}
If a morphism from an affine formal algebraic space to a presheaf on $(\Sch/S)_\fppf$ is representable by algebraic spaces and \'etale, then it is necessarily representable (i.e., representable by schemes) by Lemma~\ref{lem:representability of diagonal of formal algebraic space} below.
Hence, in the definition of a formal algebraic space, one may replace (i) with the condition that $U_i\rightarrow F$ is representable and \'etale. This result also implies that an affine formal algebraic space that is an algebraic space is an affine scheme.
\end{rem}

\begin{lem}[{\cite[0AIP]{stacks-project}}]\label{lem:representability of diagonal of formal algebraic space}
If $F$ is a formal algebraic space over $S$, then the diagonal morphism $\Delta\colon F\rightarrow F\times_SF$ is representable. Moreover, it is a monomorphism, locally quasi-finite, locally of finite type, and separated.
\end{lem}

\begin{lem}
Let $X$ be a formal scheme and consider the corresponding fppf sheaf $h_X\in\Sh((\Sch)_\fppf)$ in Lemma~\ref{lem:classical formal scheme is an fpqc sheaf}. Then $h_X$ is a formal algebraic space over $\Spec \Z$. Moreover, if $X$ is affine, then $h_X$ is affine.
\end{lem}

By abuse of language, we will simply refer to $h_X$ regarded as a sheaf on $(\Sch/S)_\fppf$ as a formal scheme $X$ over $S$ and write $F(X)$ for $\Mor_{\Sh((\Sch/S)_\fppf)}(h_X,F)$ for any sheaf $F$. An affine formal scheme over $S$ is also referred to as a \emph{classical} affine formal space (see \cite[0AID, 0AIE]{stacks-project}). 

\begin{proof}
It is easy to see that if $X\rightarrow X'$ is an open immersion between formal schemes, then $h_X\rightarrow h_{X'}$ is representable and an open immersion. Hence it suffices to treat the affine case, which follows from Remark~\ref{rem:first remarks on formal schemes} (see \cite[0AIE]{stacks-project}).
\end{proof}

\begin{defn}[{\cite[0AKX, 0AKY; 0AQ3]{stacks-project}}]
\hfill
\begin{enumerate}
 \item A formal algebraic space $F$ on $(\Sch/S)_\fppf$ is said to be \emph{locally adic*} (resp.~\emph{locally Noetherian}) if for every (equivalently, some) covering $\{U_i\rightarrow F\}$ from affine formal schemes $U_i$ such that $U_i\rightarrow F$ is representable by algebraic spaces and \'etale, $U_i$ is classical and of the form $\Spf A_i$ with $A_i$ being adic of finite ideal type (resp.~Noetherian). 
 \item  A morphism of locally adic* formal algebraic spaces on $(\Sch/S)_\fppf$ is said to be \emph{adic} if it is representable by algebraic spaces. 
 \item Let $A$ be an adic topological ring of finite ideal type with $I$ a finitely generated ideal of definition.
An \emph{$I$-adic formal algebraic space} is a locally adic* formal algebraic space equipped with an adic morphism $\calX\rightarrow \Spf A$. 
\end{enumerate}
\end{defn}

\begin{rem}\label{rem:adically P vs P}
The above definition of being locally adic* agrees with Definition~\ref{def:formal scheme is adic of finite ideal type} when $F$ is a formal scheme; a locally adic* affine formal scheme is classical.

For a morphism $f\colon X\rightarrow Y$ between locally adic* formal schemes, the above definition of $f$ being adic agrees with Definition~\ref{def:adic morphisms between formal schemes} (cf.~Remark~\ref{rem:morphism from affine formal algebraic space}).
In particular, for a property $\calP$ in Definition~\ref{def:adically P}, $f$ is adically $\calP$ if and only if $f$ is adic and $\calP$ (as a morphism representable by algebraic spaces).
\end{rem}

\begin{lem}\label{lem:representable over formal scheme}
Let $f\colon F\rightarrow G$ be an adic morphism of locally adic* formal algebraic spaces on $(\Sch/S)_\fppf$. If $G$ is a locally adic* formal scheme and $f$ is representable, then $F$ is also a locally adic* formal scheme. 
\end{lem}

\begin{proof}
This follows from \cite[Prop.~I.1.4.3, Exer.~I.1.4]{Fujiwara-Kato}.
\end{proof}

\begin{lem-defn}[{\cite[0AIN]{stacks-project}}]
Let $F$ be a formal algebraic space over $S$. Then there exists a unique reduced algebraic space $F_\red$ together with a representable thickening $F_\red\rightarrow F$ such that every morphism $U\rightarrow F$ from a reduced algebraic space $U$ factors uniquely through $F_\red$. We call $F_\red$ the \emph{reduction} of $F$.
\end{lem-defn}

\begin{lem}\label{lem:affine formal space iff reduction is affine}
Let $F$ be a formal algebraic space over $S$. Then $F$ is an affine formal scheme if and only if $F_\red$ is affine. If $F$ is locally adic*, then $F$ is a formal scheme if and only if $F_\red$ is a scheme.
\end{lem}

\begin{proof}
The first assertion is \cite[0DE8]{stacks-project}. The second follows from the first: here we impose the locally adic* condition to use the fact that a locally adic* affine formal scheme is classical.
\end{proof}

\begin{defn}
Let $U$ be a locally adic* formal scheme over $S$. An \emph{\'etale equivalence relation} on $U$ is an equivalence relation $(R,s,t\colon R\rightrightarrows  U)$ such that $R$ is a locally adic* formal scheme over $S$ and such that $s$ and $t$ are adic and \'etale.
\end{defn}

\begin{prop}\label{prop:presentation of formal algebraic space}
Let $F$ be a locally adic* formal algebraic space over $S$ and let $f\colon U\rightarrow F$ be an epimorphism of fppf sheaves from a locally adic* formal scheme $U$. Assume that $f$ is adic and \'etale, and set $R\coloneqq U\times_FU$. Then $R$ is a locally adic* formal scheme $U$ and $s\coloneqq \pr_2, t\coloneqq \pr_1\colon R\rightrightarrows  U$ define an \'etale equivalence relation on $U$. Moreover, the diagram $R\rightrightarrows  U\rightarrow F$ is a coequalizer diagram of sheaves on $(\Sch/S)_\fppf$. Hence $f$ induces an isomorphism $U/R\xrightarrow{\cong}F$.
\end{prop}

Such an \'etale equivalence relation is called a \emph{presentation} of $F$. Since $f$ is an epimorphism of sheaves, it is automatically surjective (as a representable morphism). We also remark that for every locally adic* formal algebraic space $F$, there always exists such $f\colon U\rightarrow F$.

\begin{proof}
Since $j$ is the base change of the diagonal $\Delta\colon F\rightarrow F\times_SF$ along $U\times_SU\rightarrow F\times_SF$, it is representable by Lemma~\ref{lem:representability of diagonal of formal algebraic space}. Hence $R$ is also a locally adic* formal scheme by Lemma~\ref{lem:representable over formal scheme}. Since $U\rightarrow F$ is \'etale, so are the base changes $s$ and $t$. Now it is straightforward to check that $s,t\colon R\rightrightarrows  U$ is an \'etale equivalence relation. The remaining assertions follow from Example~\ref{eg:coequalizer diagram of sheaves}.
\end{proof}

\begin{rem}\label{rem:presentation of formal algebraic space as etale sheaf}
With the notation as in Proposition~\ref{prop:presentation of formal algebraic space}, let $U/_\et R$ denote the \'etale sheaf on $\Sch/S$ associated to the presheaf $T\mapsto U(T)/R(T)$. We claim that the induced morphism $\alpha \colon U/_\et R\rightarrow U/R$ in $\PSh(\Sch/S)$ is an isomorphism. In fact, since $f$ is representable, \'etale, and surjective, the morphism $U\rightarrow U/R$ is an epimorphism of \'etale sheaves, hence so is $\alpha$. On the other hand, $U\times_{(U/_\et R)}U=R=U\times_{(U/R)}U$ (note that the fiber product of sheaves agrees with the fiber product taken in the category of presheaves), which implies that $\alpha$ is a monomorphism of presheaves.
\end{rem}

\subsection{Small \'etale site and a.q.c. sheaves}\label{sec:etale site and aqc sheaves}

Let us define the small \'etale site of a formal algebraic space. Fix a base scheme $S$ and let $X$ be a locally adic* formal algebraic space over $S$.

\begin{defn}
The \emph{small \'etale site} $X_\et$ and its variants $X_{\affine,\et}, X_{\mathrm{space},\et}$ are defined as follows: an object of $X_{\mathrm{space},\et}$ is a formal algebraic space $U$ over $S$ together with a morphism $\varphi\colon U\rightarrow X$ that is adic and \'etale; a morphism $U\rightarrow U'$ in $X_{\mathrm{space},\et}$ is a morphism of formal algebraic spaces that is compatible with structure morphisms to $X$. A family of morphisms $\{f_i\colon (U_i\rightarrow X)\rightarrow (U\rightarrow X)\}$ is a covering if $\coprod U_i\rightarrow U$ is surjective (see the following Remark~\ref{rem:small etale site of formal algebraic space}(i)). 

Let $X_\et$ (resp.~$X_{\affine,\et}$) denote the full subcategory of $X_{\mathrm{space},\et}$ consisting of formal schemes (resp.~affine formal schemes) equipped with the induced topology. 
Then the fully faithful embeddings $X_{\affine,\et}\rightarrow X_\et\rightarrow X_{\mathrm{space},\et}$ are special cocontinuous functors of sites and induce equivalences of topoi 
\[
\Sh(X_{\affine,\et})\xrightarrow{\cong}\Sh(X_\et)\xrightarrow{\cong}\Sh(X_{\mathrm{space},\et}).
\]
\end{defn}

\begin{rem}\label{rem:small etale site of formal algebraic space}
\hfill
\begin{enumerate}
\item A morphism $U\rightarrow U'$ in $X_{\mathrm{space},\et}$ is necessarily adic and \'etale, which follows from the graph decomposition $U\rightarrow U'\times_XU\rightarrow U'$. Hence for a family of morphisms $\{f_i\colon (U_i\rightarrow X)\rightarrow (U\rightarrow X)\}$, the resulting morphism $\coprod U_i\rightarrow U$ is surjective as a morphism representable by algebraic spaces if and only if it is an epimorphism of sheaves on $(\Sch/S)_\fppf$.
\item By Lemma~\ref{lem:affine formal space iff reduction is affine} and \cite[0DEG]{stacks-project}, the functor
\[
X_\et\rightarrow (X_\red)_\et: \quad (U\rightarrow X)\mapsto (U\times_X X_\red\rightarrow X_\red)
\]
is an equivalence of sites. In particular, $X_\et$ defined as above coincides with the small \'etale site of $X$ in the sense of \cite[0DEA(1)]{stacks-project}. A similar remark applies to $X_{\affine,\et}$ and $X_{\mathrm{space},\et}$ (see \cite[0DEA, 0DEG, 0DEH]{stacks-project}).
\item If $X$ is locally adic* \emph{formal scheme}, then it is an adic formal scheme of finite ideal type in the sense of \cite{Fujiwara-Kato} (cf.~Definition~\ref{def:formal scheme is adic of finite ideal type}). Then $X_\et$ is canonically identified with the small \'etale site defined in \cite[Def.~I.6.2.8]{Fujiwara-Kato}.
\end{enumerate}
\end{rem}

\begin{defn}
Let $\calO_X$ denote the structure sheaf of topological rings on $X_\et$ defined in \cite[0DEI, 0DEJ]{stacks-project} via Remark~\ref{rem:small etale site of formal algebraic space}(ii): for $U\in X_{\affine, \et}$, we have $\calO_X(U)=\Gamma(U, \calO_U)$ (recall that $U$ is classical).
\end{defn}

\begin{rem}
One can define the notion of \emph{log formal algebraic spaces} based on the locally ringed topos $(X_\et,\calO_X)$ as in \cite[\S12.1]{Gabber-Ramero-foundations}, which coincides with the existing one (e.g.~\cite[\S1.1]{Shiho-II}) in the case of formal schemes.
Moreover, we can define several types of morphisms (e.g.~\'etale and smooth) between fine log formal algebraic spaces similar to the ones in \cite{Kato-log}. Since we only need the minimum formalism in \S\ref{sec:crystals and isocrystals on crystalline site}, we omit the detail.
\end{rem}

\smallskip
\noindent
\textbf{A.q.c.\ sheaves and admissible ideals}.
We briefly review the notions of adically quasi-coherent sheaves and admissible ideals, following \cite{Fujiwara-Kato}.

Let $U$ be a locally adic* formal scheme. We first work on the Zariski site $U$. Let $\calF$ be an $\calO_U$-module. If $U$ admits an ideal $\calI$ of definition (which always exists locally on $U$), set $\widehat{\calF}\coloneqq \varprojlim_n \calF/\calI^n\calF$. This is independent of the choice of $\calI$. In particular, $\widehat{\calF}$ is well-defined even without assuming that $U$ admits an ideal of definition globally. We call $\widehat{\calF}$ the \emph{completion} of $\calF$. Note that it comes with the natural map $\calF\rightarrow\widehat{\calF}$. We say that $\calF$ is \emph{complete} if the map $\calF\rightarrow\widehat{\calF}$ is an isomorphism (see \cite[I.3.1.(a)]{Fujiwara-Kato}).

An $\calO_U$-module $\calF$ (on the Zariski site $U$) is said to be \emph{adically quasi-coherent} (\emph{a.q.c.} for short) if it is complete and if, for every open $V\subset U$ and every ideal $\calI$ of definition of finite type on $V$, the sheaf $\calF|_V/\calI(\calF|_V)$ is a quasi-coherent sheaf on the scheme $(V,\calO_V/\calI)$ (see \cite[Def.~I.3.1.3]{Fujiwara-Kato}). If, moreover, $\calF$ is of finite type as an $\calO_U$-module, we say that $\calF$ is a.q.c.\ sheaf \emph{of finite type}. When $U=\Spf A$ is affine, the global section functor gives an equivalence between the category of a.q.c.\ sheaves (of finite type) and the category of (finitely generated) complete $A$-modules with the quasi-inverse given by $M\mapsto M^\triangle$ in Definition~\ref{def:ideal of definition for formal scheme} (see \cite[Thm.~I.3.2.8]{Fujiwara-Kato}). 
An \emph{admissible ideal} is an ideal sheaf of $\calO_U$ that is a.q.c.\ of finite type and contains, Zariski locally, an ideal of definition (see \cite[Def.~I.3.7.4]{Fujiwara-Kato}). 

Next we introduce the similar notions on the \'etale site $U_\et$ (see Remark~\ref{rem:small etale site of formal algebraic space}(iii)): if $U$ admits an ideal $\calI$ of definition, set $(U_n,\calO_{U_n})\coloneqq (U,\calO_U/\calI^{n+1})$ and consider the morphism of ringed topoi $i_n=(i_{n,\ast},i_n^\ast)\colon (U_{n,\et},\calO_{U_n})\rightarrow (U_\et,\calO_U)$ associated to the closed immersion $i_n\colon U_n\hookrightarrow U$. For an $\calO_U$-module (on $U_\et$), define the completion $\widehat{\calF}$ by $\widehat{\calF}\coloneqq \varprojlim i_{n,\ast}i_n^\ast\calF$. This construction is independent of the choice of $\calI$ and globalizes (note that $U_\et$ admits an ideal of definition of finite type \'etale locally). Now one can define the notion of a.q.c. sheaves (of finite type) on $U_\et$ in a similar manner as in the Zariski site case: see \cite[I.6.2.(b)]{Fujiwara-Kato} for details.

\begin{prop}[{\cite[Prop.~I.6.2.12]{Fujiwara-Kato}}]\label{prop:aqc sheaves on Zariski site and etale site}
Let $U$ be a locally adic* formal scheme and let $\varepsilon=(\varepsilon_\ast,\varepsilon^\ast)\colon (U_\et,\calO_U)\rightarrow (U,\calO_U)$ denote the natural morphism of ringed topoi. Then $\varepsilon_\ast$ induces an equivalence from the category of a.q.c. sheaves of finite type on $U_\et$ to the category of a.q.c. sheaves of finite type on $U$. The quasi-inverse functor is given by $\calF\mapsto \widehat{\varepsilon^\ast\calF}$.
\end{prop}

We extend these notions to the case of locally adic* formal algebraic spaces using the small \'etale site. We mainly follow \cite[Chap.~I.6.3]{Fujiwara-Kato} but need to make necessary adjustments since the definition of formal algebraic spaces in \cite{stacks-project} and this paper is different from that of \cite{Fujiwara-Kato}.

\begin{defn}[{cf.~\cite[Def.~I.6.3.16]{Fujiwara-Kato}}]
Let $X$ be a locally adic* formal algebraic space over some scheme $S$. An \emph{ideal of definition of finite type} is an ideal sheaf $\calI\subset \calO_X$ of finite type on $X_\et$ such that if for any $(f\colon U\rightarrow X)\in X_\et$, we have $\widehat{f^\ast\calI}=\calI\calO_U$ as $\calO_U$-modules on $U_\et$, and $\calI\calO_U$ is an ideal of definition of finite type of $U$ via Proposition~\ref{prop:aqc sheaves on Zariski site and etale site} (note $f^{-1}\calO_X=\calO_U$ and thus $f^\ast\calI=f^{-1}\calI$). When $X$ is a formal scheme, this definition agrees with the one in Definition~\ref{def:ideal of definition for formal scheme} by adically flat descent (see \cite[Prop.~I.6.2.12, Exer.~I.6.3]{Fujiwara-Kato}).
\end{defn}

\begin{lem}\label{lem:formal space as inductive limi of algebraic spaces}
Let $X$ be a locally adic* formal algebraic space over a scheme $S$ and assume that $X_\et$ admits an ideal of definition $\calI$. For each $n\geq 0$, let $X_n$ denote the subfunctor of $X\colon (\Sch/S)_\fppf\rightarrow \Set$ defined by for $T\in \Ob(\Sch/S)$,
\[
X_n(T)=\{(t\colon T\rightarrow X)\in X(T)\mid t^\ast\calI^{n+1}=0 \;\text{on}\;T_\et\},
\]
where $(t_\ast,t^\ast)\colon (T_\et,\calO_{T})\rightarrow (X_\et,\calO_X)$ denotes the associated morphism of small \'etale ringed topoi. Then $X_n$ is an algebraic space over $S$, the induced map $X_n\rightarrow X_{n+1}$ is a thickening, and $X=\varinjlim_n X_n$ as sheaves on $(\Sch/S)_\fppf$.
\end{lem}

We often write $X_n=(X,\calO_X/\calI^{n+1})$ for simplicity.

\begin{proof}
Take a presentation $s,t\colon R\rightrightarrows  U$ given by an epimorphism $f\colon U\rightarrow X$ from a locally adic* formal scheme $U$ that is representable by algebraic spaces and \'etale (see Proposition~\ref{prop:presentation of formal algebraic space}). Let $U_n$ (resp.~$R_n)$ denote the scheme $(U,\calO_U/\calI^{n+1}\calO_U)$ (resp.~$(R,\calO_R/\calI^{n+1}\calO_R)$) over $S$. Then the induced diagram $s_n,t_n\colon R_n\rightrightarrows  U_n$ is an \'etale equivalence relation on $U_n$, and it is straightforward to see that $X_n$ represents the quotient sheaf, $X_n\rightarrow X_{n+1}$ is a thickening, and $X=\varinjlim_n X_n$.
\end{proof}

Let $X$ be a locally adic* formal algebraic space over $S$. If $X_\et$ admits an ideal $\calI$ of definition, consider $\calF\rightarrow\widehat{\calF}\coloneqq \varprojlim_n \calF/\calI^n\calF$ for an $\calO_X$-module $\calF\in \Sh(X_{\mathrm{space},\et})=\Sh(X_\et)$. This is independent of the choice of $\calI$ and for $(f\colon X'\rightarrow X)\in X_{\mathrm{space},\et}$, we have $f^\ast\widehat{\calF}=\varinjlim_n f^\ast\calF/\calI^nf^\ast\calF$. Hence $\calF\rightarrow\widehat{\calF}$ is well-defined even without assuming that $X_\et$ admits an ideal of definition.

\begin{defn}\label{def:aqc sheaves in etale topology}
We say that an $\calO_X$-module $\calF$ on $X_\et$ is \emph{a.q.c.~(of finite type or of finite presentation)} if $\calF$ is complete, namely, $\calF\rightarrow\widehat{\calF}$ is an isomorphism, and for any $U\in X_{\mathrm{space},\et}$ and any ideal $\calI$ of definition of finite type of $U$, the sheaf $\calF|_U/\calI(\calF|_U)$ is a quasi-coherent sheaf (of finite type or of finite presentation) on the algebraic space $(U,\calO_U/\calI)$ (cf.~\cite[Def.~I.6.3.23]{Fujiwara-Kato}). 

An \emph{admissible ideal} is an ideal sheaf of $\calO_X$ that is a.q.c.\ of finite type and contains an ideal of definition \'etale locally (cf.~\cite[Def.~I.3.7.4]{Fujiwara-Kato}).

One can easily deduce from adically flat descent \cite[Chap.~I.6.1]{Fujiwara-Kato} that the notion of a.q.c. sheaves (of finite type) or admissible ideals agrees with the previously defined one when $X$ is a locally adic* formal scheme.
\end{defn}

\begin{rem}\label{rem:aqc sheaves in etale topology}
In the work of Fujiwara--Kato, the notion of a.q.c.~sheaves of finite presentation is only discussed for locally universally rigid-Noetherian formal schemes and algebraic spaces (see \cite[\S I.3.5, I.6.5]{Fujiwara-Kato}). We need to use this notion in \S\ref{sec:crystals and isocrystals on crystalline site} for more general (but $p$-adic) formal algebraic spaces, for which the above definition suffices. In fact, when $X$ is adic over $\Spf \Z_p$, it follows from fpqc descent that an a.q.c.~sheaf $\calF$ on $X_\et$ is of finite presentation in the above sense if and only if $\calF/p^n\calF$ is of finite presentation as an $\calO_X/p^n$-module for every $n\geq 1$.
\end{rem}

\begin{example}
With the notation and assumption as in Lemma~\ref{lem:formal space as inductive limi of algebraic spaces}, the defining ideal of a closed subspace of $X_1$ of finite type gives an admissible ideal on $X_\et$ (cf.~\cite[03MB]{stacks-project} and \cite[Cor.~I.3.7.8]{Fujiwara-Kato}). If $Z$ is another formal space over $S$ and $Z\rightarrow X$ is an adic morphism of finite type, one can define the \emph{sheaf of differentials} $\Omega_{Z/X}^1$, which is an a.q.c.~sheaf of finite type on $Z_\et$ (cf.~\cite[Thm.~I.5.2.1, Prop.~I.5.2.4]{Fujiwara-Kato}). A similar construction yields the a.q.c.~sheaf of finite type $\omega_{Z/X}^1$, the \emph{sheaf of log differentials}, for an adic morphism $(Z,M_Z)\rightarrow (X,M_X)$ of finite type between coherent log formal algebraic spaces locally adic* over $S$.
\end{example}

\subsection{Formal algebraic stacks}
We continue to fix a base scheme $S$ and work on $(\Sch/S)_{\fppf}$. Every formal algebraic space $U$ over $S$ gives a stack $\calS_U\rightarrow (\Sch/S)_\fppf$ in groupoids as in Remark~\ref{rem:category fibered in setoids associated to presheaf}. We say that a stack in groupoids over $(\Sch/S)_\fppf$ is \emph{representable by a formal algebraic space} if it is equivalent to $\calS_U$ for a formal algebraic space $U$ over $S$.

\begin{defn}[{\cite[Def.~5.3, 5.16, 8.11]{EmertonStack}}]\label{defn:formal algebraic stacks}
A \emph{formal algebraic stack} over $S$ is a stack in groupoids $\calX\rightarrow (\Sch/S)_{\fppf}$ such that there exist a formal algebraic space $U$ over $S$ and a $1$-morphism $\calS_U\rightarrow \calX$ that is representable by algebraic spaces, smooth, and surjective. If $U$ can be taken to be locally adic* (resp.~locally Noetherian), we call $\calX$ \emph{locally adic*} (resp.~\emph{locally Noetherian}). 

Note that the diagonal of a formal algebraic stack over $S$ is representable by algebraic spaces and locally of finite type \cite[Lem.~5.12]{EmertonStack}. For a formal algebraic space $U$ over $S$, the formal algebraic stack $\calS_U$ is locally adic* (resp.~locally Noetherian) if and only if $U$ is locally adic* (resp.~locally Noetherian): see \cite[Ex.~5.7, Rem.~5.17. Rem.~8.12]{EmertonStack}. We will write $U$ for $\calS_U$ when there is no confusion.
\end{defn}

See Definition~\ref{def:groupoid for site} for the definition of a groupoids in sheaves on a site.

\begin{defn}
A \emph{smooth groupoid in formal algebraic spaces} over $S$ is a groupoid $(U,R,s,t,c)$ in sheaves on $(\Sch/S)_\fppf$ such that $U$ and $R$ are formal algebraic spaces and such that $s,t\colon R\rightrightarrows  U$ are representable by algebraic spaces and smooth.
\end{defn}

\begin{prop}\label{prop:presentation of formal algebraic stack}
Let $\calX$ be a formal algebraic stack over $S$ and let $\pi\colon U\rightarrow\calX$ be a $1$-morphism from a formal algebraic space $U$ that is representable by algebraic spaces, smooth, and surjective. Then $R=U\times_\calX U$ is representable by a formal algebraic space. Moreover, $R\rightrightarrows  U$ defines a smooth groupoid in formal algebraic spaces and induces an equivalence $[U/R]\xrightarrow{\cong}\calX$.
\end{prop}

\begin{proof}
This is \cite[Lem.~5.25]{EmertonStack}.
\end{proof}

\begin{defn}\label{def:presentation of formal algebraic stack}
A \emph{presentation} of a formal algebraic stack $\calX$ over $S$ consists of a smooth groupoid $(U,R,s,t,c)$ in formal algebraic spaces and an equivalence $[U/R]\xrightarrow{\cong}\calX$. 
\end{defn}

\begin{prop}\label{prop:formal algebraic space given by smooth groupoid}
Let $(U,R,s,t,c)$ be a smooth groupoid in formal algebraic spaces over $S$. Then the quotient stack $[U/R]$ is a formal algebraic stack over $S$. Moreover, $U\rightarrow[U/R]$ is representable by algebraic spaces, smooth, and surjective, and $U\times_{[U/R]}U$ is representable by $R$.
\end{prop}

\begin{proof}
Recall from Lemma~\ref{lem:2-coequalizer property of quotient stack}(i) that there is a natural $1$-morphism $\pi\colon \calS_U\rightarrow [U/ R]$ and that the induced $1$-morphism $\calS_R\rightarrow \calS_U\times_{[U/R]}\calS_U$ is an equivalence. 

We also remark that for any $x,y\in \Ob([U/R]_T)$, there exist an fppf covering $T'\rightarrow T$ and $\tilde{x}',\tilde{y}'\in U(T')$ such that $\pi(\tilde{x}')=x|_{T'}$ and $\pi(\tilde{y}')=y|_{T'}$ by Lemma~\ref{lem:criteria for equivalence of stacks}(ii), which implies that $\Isom_{[U/R]}(x,y)$ is an algebraic space by \cite[04S6]{stacks-project}. Hence the diagonal map $[U/R]\rightarrow [U/R]\times [U/R]$ is representable by algebraic spaces by \cite[045G]{stacks-project}.

It remains to show that $\pi$ is representable by algebraic spaces, smooth, and surjective. Take any $T\in \Ob(\Sch/S)$ and any $1$-morphism $t\colon \Sch/T\rightarrow [U/R]$. By the previous paragraph, $t$ is representable by algebraic spaces. Hence $\calS_U\times_{\pi,[U/R],t}(\Sch/T)$ is representable by a formal algebraic space $V$ by \cite[04SC, 02YR, 0AJI]{stacks-project}. We need to show that $V$ is indeed an algebraic space and the morphism $V\rightarrow T$ is smooth and surjective. For this, take an fppf covering $T'\rightarrow T$ and a morphism $\tilde{t}'\colon T'\rightarrow U$ with $t|_{T'}=\pi(\tilde{t}')$ as before. Using $\calS_U\times_{[U/R]}\calS_U\cong \calS_R$, we can identify the morphism $V\times_TT'\rightarrow T'$ with $R\times_UT'\rightarrow T'$, which is representable by algebraic spaces, smooth, and surjective. In particular, $V\times_TT'=R\times_UT'$ is an algebraic space, and the assertions follow from \cite[04S6, 0429, 041Q]{stacks-project}.
\end{proof}

\begin{lem}\label{lem:presentation of formal algebraic stack in etale topology}
Let $\calX$ be a formal algebraic stack over $S$ together with a presentation given by $(U,R,s,t,c)$. Let $[U/_\et R]\rightarrow \Sch/S$ denote the \'etale stackification of presheaf in groupoids $T\mapsto (U(T),R(T),s(T),t(T),c(T))$. Then the induced $1$-morphism $[U/_\et R]\rightarrow [U/R]\xrightarrow{\cong}\calX$ over $\Sch/S$ is an equivalence.
\end{lem}

\begin{proof}
Before staring the proof, we remark the following: the notion of $2$-fiber products or $1$-morphisms representable by algebraic spaces does not depend on the choice of topology; for a morphism $\alpha\colon F\rightarrow G$ of presheaves on $\Sch/S$ that is representable by algebraic spaces, the induced $1$-morphism $\calS_F\rightarrow\calS_G$ over $\Sch/S$ is representable by algebraic spaces. Moreover, it is smooth (resp.~surjective) if and only if $\alpha$ is smooth  (resp.~surjective) by \cite[045A]{stacks-project}.

Let $\pi\colon U\rightarrow [U/_\et R]$ be the natural $1$-morphism;  $\calS_U\times_{[U/_\et R]}\calS_U$ is naturally equivalent to $\calS_R$ by Lemma~\ref{lem:2-coequalizer property of quotient stack}(i). We also claim that (i) the diagonal map $[U/_\et R]\rightarrow [U/_\et R]\times [U/_\et R]$ is representable by algebraic spaces; (ii) $\pi$ is representable by algebraic spaces, smooth, and surjective. These are proved exactly as in the proof of Proposition~\ref{prop:formal algebraic space given by smooth groupoid}: for (i), apply an \'etale topology analogue of \cite[045G]{stacks-project}; for (ii), consider an \'etale covering $T'\rightarrow T$ instead.

Now let us show that $f\colon [U/_\et R]\rightarrow [U/R]$ is an equivalence, following the arguments in \cite[076U]{stacks-project}. Since $\calS_U\times_{[U/_\et R]}\calS_U\cong \calS_R$ and $\calS_U\times_{[U/R]}\calS_U\cong \calS_R$, we deduce from Lemma~\ref{lem:explicit description of quotient stack} that $f$ is fully faithful. Then we conclude that $f$ is an equivalence by the argument in the last paragraph of \cite[076U]{stacks-project}.
\end{proof}

\begin{defn}[{\cite[Def.~3.27, Rem.~3.30, Lem.~5.26]{EmertonStack}}]\label{defn:reduced substack}
Let $(\Sch/S)_\red$ denote the full subcategory of $\Sch/S$ consisting of reduced $S$-schemes. For a stack $\calX\rightarrow (\Sch/S)_\fppf$, let $\calX_\red$ denote the substack of $\calX$ generated by the full subcategory $\calX|_{(\Sch/S)_\red}$ of $\calX$ (cf.~\cite[Rem.~3.24]{EmertonStack}). We call $\calX_\red$ the \emph{underlying reduced substack} of $\calX$. If $\calX=\calS_U$ for a formal algebraic space $U$, then $(\calS_U)_\red=\calS_{(U_\red)}$ as substacks of $\calS_U$.
If $\calX$ is a formal algebraic stack over $S$, then $\calX_\red$ is a closed and reduced substack of $\calX$, and the inclusion $\calX_\red\hookrightarrow\calX$ is a thickening; any $1$-morphism $\calY\rightarrow \calX$ from a reduced algebraic stack factors through $\calX_\red$.
\end{defn}

\begin{rem}\label{rem:locally Noetherian formal algebraic stack as colimit}
Let $\calX$ be a locally Noetherian formal algebraic stack over $S$. For each $n\geq 1$, define a strictly full subcategory $\calX_n$ of $\calX$ as follows: set $\calX_1=\calX_\red$. For every $S$-scheme $T$ and $x\in \Ob(\calX_T)$, the fiber product $\calX_{1}\times_{\calX,x}T$ is represented by a closed immersion $T_1\hookrightarrow T$. Define $(\calX_n)_T$ as the full subcategory of $\calX_T$ consisting of the objects $x$ such that the $n$th power of the ideal sheaf of $T_1\hookrightarrow T$ is zero.
We claim that $\calX_n$ is a locally Noetherian algebraic stack over $S$, the transition map $\calX_{n}\rightarrow \calX_{n+1}$ is a thickening, and the induced $1$-morphism $\lim_n\calX_n\rightarrow \calX$ is an equivalence. To see this, first observe that the condition $x\in \Ob((\calX_n)_T)$ can be checked smooth locally on $T$. Hence by taking a presentation of $\calX$, the assertions are reduced to analogues of locally Noetherian formal algebraic spaces, which are exactly \cite[0GHM]{stacks-project}. 
\end{rem}

\begin{defn}[{\cite[Ex.~5.9]{EmertonStack}}]\label{def:formal completion}
Let $\calX$ be an algebraic stack over $S$, and let $Z\subset|\calX|$ be a locally closed subset, where $\lvert \calX\rvert$ denotes the underlying topological space of $\calX$ (cf.~\cite[04XG, 04XL]{stacks-project}). The \emph{completion} $\widehat{\calX}_{|Z}$ of $\calX$ along $Z$ is a substack of $\calX$ defined as
\[
\widehat{\calX}_{|Z}(U)=\lbrace f\colon U\rightarrow \calX\mid f(|U|)\subset Z\rbrace.
\]
This is a formal algebraic stack.
The proof therein shows that this construction is consistent with the (formal) completion of schemes (\cite[I.1.4.(b)]{Fujiwara-Kato}) and formal algebraic spaces (\cite[0AIX]{stacks-project}). More generally, one can still obtain a formal algebraic stack when $\calX$ is a formal algebraic stack over $S$ and $Z$ is a locally closed subset of $\lvert \calX_\red\rvert$; use \cite[0GVR]{stacks-project}.
\end{defn}

\begin{defn}[{\cite[Def.~7.3, Rem.~7.4, 7.5]{EmertonStack}}]
A morphism $f\colon \calX\rightarrow \calY$ of stacks over $S$ with $\calY$ a locally adic* formal algebraic stack is called \emph{adic} if it is representable by algebraic stacks. In this case, $\calX$ is also a locally adic* formal algebraic stack. When $f$ arises as $\calS_X\rightarrow\calS_Y$ for a morphism $g\colon X\rightarrow Y$ of formal algebraic spaces, $f$ is adic if and only if $g$ is adic.
\end{defn}

\begin{defn}[{\cite[Def.~7.6, Ex.~7.7]{EmertonStack}}]\label{def:I-adic completion of algebraic stack}
Let $A$ be an adic topological ring of finite ideal type with $I$ a finitely generated ideal of definition.

An \emph{$I$-adic formal algebraic stack} is a formal algebraic stack $\calX$ over $\Spec A$ that admits an adic morphism $\calX\rightarrow \Spf A$. By definition, for each $n\geq 1$, the fiber product $\calX_n\coloneqq \calX\times_{\Spf A}\Spec A/I^n$ is an algebraic stack over $\Spec A/I^n$. Moreover, we have $\calX=\varinjlim_n \calX_n$ by \cite[Lem.~4.8, Rem.~4.9.]{EmertonStack}.

For an algebraic stack $\calX$ over $\Spec A$, the \emph{$I$-adic completion} $\widehat{\calX}_{I}$ of $\calX$ is defined to be the fiber product $\calX \times_{\Spec A} \Spf A$. This is an $I$-adic formal algebraic stack and is identified with the completion of $\calX$ along the closed substack $\calX\times_{\Spec A}\Spec(A/I)$.
\end{defn}

\begin{rem}
The following basic properties hold: if $\calY$ is an $I$-adic formal algebraic stack and $\calX\rightarrow\calY$ is a morphism of stacks that is representable by algebraic stacks, then $\calX$ is also an $I$-adic formal algebraic stack \cite[Lem.~7.9]{EmertonStack}; a morphism between $I$-adic formal algebraic stacks is adic \cite[Lem.~7.10]{EmertonStack}.
\end{rem}

We end this subsection with a few remarks about formal algebraic stacks over a complete discrete valuation ring. Let $\calO$ be a complete discrete valuation ring with residue field $k$ and maximal ideal $\fkm$.

\begin{rem}\label{rem:formal scheme valued point of formal algebraic stack}
Let $\calX\rightarrow (\Sch/\calO)_\fppf$ be a formal algebraic stack. For every $\fkm$-adic formal scheme $S$ over $\calO$, set $S_n\coloneqq S\times_{\Spec\calO}\Spec (\calO/\fkm^n)$. Let $\calX_{S}$ denote the groupoid $\varprojlim_n\calX_{S_n}$: more precisely, an object is $(x_n,i_n\colon x_n\rightarrow x_{n-1})_n$ where $x_n\in \calX_{S_n}$ and $i_n$ is a morphism in $\calX$, and a morphism from $(x_n,i_n)$ to $(x_n',i_n')$ consists of a morphism $f_n\colon x_n\rightarrow x_n'$ for each $n$ such that $i_n'\circ f_n=f_{n-1}\circ i_n$. One can check that $\calX_{S}$ defines a category fibered in groupoids over the category of $\fkm$-adic formal schemes over $\calO$.

More generally, let $S$ be an $\fkm$-adic formal algebraic space over $\calO$ and define $S_n$ similarly. Since $S=\varinjlim_n S_n$ as sheaves, we also have $\calS_S=\varinjlim_n\calS_{S_n}$ as stacks by Lemma~\ref{lem:sheafification and stackification} and the construction in \cite[Lem.~2.48]{EmertonStack}. 
It follows that the natural functor $\Mor_{\Sch/\calO}(\calS_S,\calX)\rightarrow \varprojlim_n\Mor_{\Sch/\calO}(\calS_{S_n},\calX)$ between the categories of $1$-morphisms over $(\Sch/\calO)$ is an equivalence, and it is further identified with the above $\calX_S$ when $S$ is an $\fkm$-adic formal scheme. We often write $\calX_S$ or $\calX(S)$ for $\Mor_{\Sch/\calO}(\calS_S,\calX)$ by abuse of notation.
\end{rem}

\begin{rem}\label{rem:enhanced moduli interpretation of completion of algebraic stack}
Let $(\Sch/\calO)_{\fkm\text{-}\nilp}$ denote the full subcategory of $\Sch/\calO$ consisting of $\calO$-schemes in which $\fkm$ is locally nilpotent. 
Let $\calX\rightarrow (\Sch/\calO)_\fppf$ be an algebraic stack. It is straightforward to see that the $\fkm$-adic completion $\widehat{\calX}_\fkm$ coincides with $\calX|_{(\Sch/\calO)_{\fkm\text{-}\nilp}}$ as substacks of $\calX$.
Set $\calX_k\coloneqq \calX\times_{\Spec \calO}\Spec k$, which we identify with $\calX|_{(\Sch/k)}$ and regard as a substack of $\calX$. We use the same convention for any substack of $\calX$.

Let $Z\subset \lvert\calX_k\rvert$ be a closed subset and write $\calZ$ for the corresponding reduced closed substack of $\calX_k$ as in \cite[050C]{stacks-project}. Observe $(\widehat{\calX}_{|Z})_k=(\widehat{\calX_k})_{|Z}$ as substacks of $\calX$.
Let $\calY\rightarrow(\Sch/\calO)_\fppf$ be an algebraic stack that factors through $(\Sch/\calO)_{\fkm\text{-}\nilp}$. We claim that the category $\Mor_{\Sch/\calO}(\calY,\widehat{\calX}_{|Z})$ is precisely the full subcategory of $\Mor_{\Sch/\calO}(\calY,\calX)$ consisting of $1$-morphisms $\calY\rightarrow\calX$ such that the composite with $\calY_\red\rightarrow \calY$ factors through the substack $\calZ$; any $1$-morphism $\calS_U\rightarrow\widehat{\calX}_{|Z}$ satisfies this property after post-composing $\widehat{\calX}_{|Z}\hookrightarrow \calX$ by the fact $(\widehat{\calX}_{|Z})_\red=\calZ$ (see \cite[Ex.~5.30]{EmertonStack}). For the converse, observe that the assertion holds if $\calY$ comes from a scheme in $(\Sch/\calO)_{\fkm\text{-}\nilp}$ by definition. By taking a presentation, we see that the same holds for algebraic spaces in which $\fkm$ is locally nilpotent and for the general $\calY$. Finally, note that the same claim also holds for
any $\fkm$-adic formal algebraic stack as well as any locally Noetherian formal algebraic stack that admits a map to $\Spf\calO$; write such a formal algebraic stack as a $2$-colimit and apply the original claim (see Definition~\ref{def:I-adic completion of algebraic stack} and Remark~\ref{rem:locally Noetherian formal algebraic stack as colimit}).
\end{rem}

\section{\'Etale spaces and adic stacks}\label{sec:adic stacks}

This section introduces a theory of algebraic spaces and Artin stacks in non-archimedean geometry (\emph{\'etale spaces} and \emph{adic stacks}).
We use adic spaces for our foundation and mainly follow \cite{Warner}.

\subsection{Generalities on adic spaces}
Fix a complete non-archimedean field $K$ (i.e., complete with respect to a rank-$1$ valuation) of mixed characteristic $(0,p)$ with ring of integers $\calO_K$. Let $\pi_K$ be a uniformizer of $\calO_K$.

We use Huber's theory of adic spaces for the foundation of non-archimedean geometry and refer the reader to \cite{Huber-cont, Huber-german, Huber-gen, Huber-etale} for the details. For an adic space $S$, let $\AdSp_S$ denote the category of adic spaces over $S$ and write $\AdSp_S^\aff$ (resp.~$\AdSp_S^{\aff,\lft}$) for the full subcategory of affinoid adic spaces (locally of finite type over $S$).
Write $\AdSp_K^{(\aff,\lft)}$ for $\AdSp_{\Spa(K,\calO_K)}^{(\aff,\lft)}$.

\begin{rem}[{\cite{Huber-cont}}]\label{rem:open bounded subring}
For an f-adic ring $A$, a subring $A_0\subset A$ is a ring of definition if and only if $A_0$ is open and bounded in $A$ (Prop.~1(ii), \emph{op.~cit.}). If $(A,A^+)$ is an affinoid ring over $(K,\calO_K)$, then it is Tate. Moreover, every open and bounded subring $A_0\subset A$ admits $\pi_KA_0$ as an ideal of definition (Prop.~1.5(ii), \emph{op.~cit.}), and $A^+$ is the filtered union of open and bounded $\calO_K$-subalgebras of $A^+$ (Cor.~1.3, \emph{op.~cit.}). We usually write $(\widehat{A},\widehat{A}^+)$ for the completion of $(A,A^+)$.
A morphism $f\colon (A,A^+)\rightarrow (B,B^+)$ of affinoid rings over $(K,\calO_K)$ is adic (Prop.~1.10, \emph{op.~cit.}), and for every open and bounded $\calO_K$-subalgebra $A_0\subset A^+$, there exists an open and bounded $\calO_K$-subalgebra $B_0\subset B^+$ with $f(A_0)\subset B_0$ (Lem.~1.8(iii), \emph{op.~cit.}).

In particular, every morphism in $\AdSp_K$ is adic. Hence for morphisms $X\rightarrow Y$ and $Y'\rightarrow Y$ in $\AdSp_K$, the fiber product $X\times_YY'$ exists if $X\rightarrow Y$ is locally of weakly finite type by \cite[Prop.~1.2.2(b)]{Huber-etale}.
\end{rem}

\smallskip
\noindent
\textbf{Rigid $K$-spaces}.
Write $\Spa(K)\coloneqq \Spa(K,\calO_K)$ when there is no confusion.
In this paper, a \emph{rigid $K$-space} refers to an adic space locally of finite type over $\Spa K$. We let $\Rig_K$ denote the category of rigid $K$-spaces. By \cite[Thm.~II.A.5.2]{Fujiwara-Kato}, $\Rig_K$ is naturally equivalent to the category of locally of finite type rigid spaces over $(\Spf \calO_K)^\rig$ developed in \cite{Fujiwara-Kato}, and we often refer its results via this equivalence.

\smallskip
\noindent
\textbf{Coherent sheaves}.
See Notation and conventions after Introduction for our notation.

\begin{prop}\label{prop:basics of coherent sheaves}
 Let $X=\Spa(A,A^+)$ be a strongly Noetherian affinoid adic space over $(K,\calO_K)$ with $A$ complete. 
\begin{enumerate}
\item  Let $M\in\Mod_A^{\mathrm{fg}}$. For any rational subset $W\subset X$, consider the $\calO_X(W)$-module $M\otimes_A\calO_X(W)$ with the natural $\calO_X(W)$-topology, and for any open $U\subset X$, set
\[
(M\otimes\calO_X)(U)=\varprojlim_{W\subset U}M\otimes_A\calO_X(W),
\]
with the projective limit topology, where the limit is taken over all the rational subsets of $W$ of $X$ contained in $U$.
Then $M\otimes\calO_X$ is a sheaf of complete topological groups on $X$, and $H^i(W,M\otimes\calO_X)=0$ for every $i\geq 1$ and every rational open $W$ of $X$, where the cohomology is taken in the category of abelian sheaves on $X$ (which agrees with the one in $\Mod(\calO_X)$). Moreover, the resulting functor $\Mod_A^{\mathrm{fg}}\rightarrow\Mod(\calO_X)$ is exact.
\item The global sections functor $\Gamma(X,-)\colon \Mod(\calO_X)\rightarrow\Mod_A$ induces an equivalence of categories
\[
\Coh(\calO_X)\xrightarrow{\cong}\Mod_A^{\mathrm{fg}}
\]
with quasi-inverse given by $M\mapsto M\otimes\calO_X$.
Moreover, finite locally free $\calO_X$-modules correspond to finite projective $A$-modules.
\item For the morphism $f\colon X'\coloneqq \Spa(A,A^\circ)\rightarrow X$, the pullback $f^\ast\colon \Mod(\calO_X)\rightarrow\Mod(\calO_{X'})$ restricts to an equivalence of categories $\Coh(\calO_X)\xrightarrow{\cong}\Coh(\calO_{X'})$.
\end{enumerate}
\end{prop}

\begin{proof}
Part (i) is given in \cite[Satz 3.3.17]{Huber-german}, \cite[Thm~2.5]{Huber-gen}; (ii) follows from \cite[Satz 3.3.18]{Huber-german} and \cite[Cor.~2.5.19(a)]{Kedlaya-Liu-I}; (iii) is a consequence of (ii).
\end{proof}

\begin{cor}\label{cor:coherent sheaf does not depend on ring of integral elements}
 Let $(L,L^+)$ be an analytic affinoid field over $(K,\calO_K)$. Let $X$ be an adic space locally of finite type over $(L,L^+)$ and set 
 \[
 f\colon X'\coloneqq X\times_{\Spa(L,L^+)}\Spa(L,\calO_L)\rightarrow X.
 \]
 Then the pullback $f^\ast\colon \Mod(\calO_X)\rightarrow\Mod(\calO_{X'})$ restricts to an equivalence of categories $\Coh(\calO_X)\xrightarrow{\cong}\Coh(\calO_{X'})$.
\end{cor}

\begin{proof}
This follows from Proposition~\ref{prop:basics of coherent sheaves}(iii).
\end{proof}

\begin{prop}\label{prop:coherent sheaves on etale site of adic space}
Let $X$ be an adic space locally of finite type over $\Spa(K,\calO_K)$ and let $\lambda\colon X_\et\rightarrow X$ be the morphism from the \'etale site of $X$.
\begin{enumerate}
\item Assume that $X$ is affinoid. Let $\calF\in\Coh(\calO_X)$ and set $\calF_\et\coloneqq \lambda^\ast\calF$. Then $\calF_\et(U)=\calF(X)\otimes_{\calO_X(X)}\calO_U(U)$ for every affinoid $U\in X_\et$, and $H^i(X_\et,\calF_\et)=0$ for every $i\geq 1$.
\item The pullback $\lambda^\ast\colon \Mod(\calO_X)\rightarrow \Mod(\calO_{X_\et})$ gives an equivalence of categories $\Coh(\calO_X)\xrightarrow{\cong}\Coh(\calO_{X_\et})$.
The same holds for the full subcategories of finite locally free modules.
\end{enumerate}
\end{prop}

\begin{proof}
Part (i) is \cite[Prop.~9.2(i)]{scholze-p-adic-hodge}. The first part of (ii) follows from (i) and Proposition~\ref{prop:basics of coherent sheaves}(ii), and the second is \cite[Lem.~7.3]{scholze-p-adic-hodge} (note that the proof therein does not use the smoothness assumption on $X$).
\end{proof}

\begin{rem}\label{rem:flatnesss of coherent sheaf and morphism for adic spaces}
Let $f\colon X\rightarrow Y$ be a morphism locally of finite type between locally strongly Noetherian adic spaces over $K$ and let $\calF$ be a coherent module on $X$. Consider the following conditions:
\begin{enumerate}
 \item[(a)] for every $x\in X$, the stalk $\calF_x$ is a flat $\calO_{Y,f(x)}$-module;
 \item[(b)] for every pair of affinoid opens $U=\Spa(B,B^+)\subset X$ and $V=\Spa(A,A^+)\subset Y$ with $f(U)\subset V$ such that $A,B$ are complete and the $B$-module $\calF(U)$ is a flat $A$-module;
 \item[(b')] there exist coverings $X=\bigcup_i\Spa(B_i,B_i^+)$ and $Y=\bigcup_i\Spa(A_i,A_i^+)$ by affinoid opens with $f(\Spa(B_i,B_i^+))\subset \Spa(A_i,A_i^+)$ such that $A_i,B_i$ are complete and the $B_i$-module $\calF(\Spa(B_i,B_i^+))$ is flat over $A_i$ for every $i$.
\end{enumerate}

We say that $\calF$ is \emph{$Y$-flat} if (a) holds; $f$ is \emph{flat} if (a) holds for $\calF=\calO_X$ (\cite[p.~78]{Huber-etale}; see \cite[\S2]{Lazda-descent}, \cite[\S2.5]{Warner}, and \cite[B.4]{Zavyalov-quotients} for relevant discussions).

Note that (b) are (b') are equivalent: obviously, (b) implies (b'); for the converse, recall that if $U=\Spa(B,B^+)\in X$ is an affinoid open and $U=\Spa(B_j,B_j^+)$ is an open covering by rational subsets with $B,B_j$ complete, then $B\rightarrow \prod_j B_j$ is faithfully flat by \cite[Lem.~1.4, (II.1)(iv)]{Huber-gen} and \cite[00HQ]{stacks-project}. One can deduce (b) from (b') using this observation and Proposition~\ref{prop:basics of coherent sheaves}.

Moreover, (a) implies (b'). This is \cite[Lem.~B.4.2]{Zavyalov-quotients} if $\calF=\calO_X$, and the same argument works for any coherent sheaf $\calF$: we may shrink $X$ and $Y$. With the notation of \emph{loc.~cit.}, let $M$ be a finitely generated $B$-module and write $\calF$ for the corresponding coherent sheaf on $X=\Spa(B,B^+)$. For any maximal ideal $\fkm$ of $B$, take a point $v\in X$ with $\operatorname{supp}(v)=\fkm$. Then $\calF_v=M_\fkm\otimes_{B_\fkm}\calO_{X,v}$. Set $w=f(v)\in Y=\Spa(A,A^+)$ and write $\fkp$ for the prime ideal of $A$ below $\fkm$. If $\calF_v$ is flat over $\calO_{Y,w}$, then $M_\fkm$ is flat over $A_\fkp$ since $A_\fkp\rightarrow \calO_{Y,w}$ and $B_\fkm\rightarrow \calO_{X,v}$ are faithfully flat. We conclude that $M$ is flat over $A$ by \cite[00HT]{stacks-project}.
\end{rem}

\subsection{The category \texorpdfstring{$\V_K$}{VK}}

In the rest of this section, fix a complete non-archimedean field $K$ of mixed characteristic $(0,p)$ and let $\calO_K$ denote its ring of integers. 
In \cite{Warner}, Warner introduces the category $\V$ of adic spaces locally of finite type over some complete affinoid field and defines the notions of \'etale spaces and adic stacks relative to an admissible subcategory of $\V$ (see Remark~\ref{rem:Warner's admissible category} below). In this paper, we consider the following variant.

\begin{defn}
Define the category $\V_K$ as follows: 
\begin{itemize}
 \item an object is an adic space locally of finite type over $(L,L^+)$ for some analytic affinoid field $(L,L^+)$ over $(K,\calO_K)$;
 \item  a morphism is a morphism of adic spaces over $(K,\calO_K)$. 
\end{itemize}
We equip $\V_K$ with the \'etale topology: a family $(f_i\colon S_i\rightarrow S)$ of morphisms is a \emph{covering} if each $f_i$ is \'etale and $S=\bigcup_i f_i(S_i)$. This satisfies the axioms of a site and let $\V_{K,\Et}$ denote the resulting big \'etale site. 

Let $\Afd_{K}$ denote the full subcategory of $\V_{K}$, consisting of the affinoids.\footnote{Recall that $\AdSp_K^{\aff(,\lft)}$ denotes the category of affinoid adic spaces (locally of finite type) over $(K,\calO_K)$.}
More generally, for any $S\in \V_K$ (not necessarily an affinoid), we write $\Afd_K/S$ for the full subcategory of the localized category $\V_K/S$ consisting of affinoids. 
With respect to the induced \'etale topology, the fully faithful embedding $(\Afd_{K}/S)_\Et\rightarrow (\V_{K}/S)_\Et$ is a special cocontinuous functor and gives an equivalence of categories $\Sh((\Afd_{K}/S)_\Et)\xrightarrow{\cong}\Sh((\V_{K}/S)_\Et)$. For a morphism $T\rightarrow S$ in $\V_K$, we have a morphism of topoi $\Sh((\V_K/T)_\Et)\rightarrow \Sh((\V_K/S)_\Et)$ by \cite[00Y0]{stacks-project}. 
\end{defn}

\begin{rem}\label{rem:Warner's admissible category}
In the definition of an object of $\V_K$, the affinoid field  $(L,L^+)$ is \emph{not} part of the datum as opposed to Warner's $\V$.
The category $\V_K$ satisfies the following two properties:
\begin{enumerate}
 \item If $S\in\Ob(\V_K)$ and $S'\rightarrow S$ is a morphism of adic spaces that is locally of finite type, then $S'\in\Ob(\V_K)$.
 \item If $S'\rightarrow S$ and $S''\rightarrow S$ are morphisms in $\V_K$ with $S'\rightarrow S$ locally of finite type, the fiber product $S'\times_SS''$ as adic space exists and belongs to $\V_K$.
\end{enumerate}
A subcategory of Warner's $\V$ satisfying these two conditions is called \emph{admissible} in \cite[Def.~2.2.6]{Warner}.
Essentially, general results on admissible categories in \cite{Warner} continue to hold for $\V_K$ with the same proof. In what follows, we cite Warner's results without mentioning this subtle difference.
\end{rem}

\begin{prop}\label{prop:subcanonical}
For every $S\in\V_K$, the site $(\V_{K}/S)_{\Et}$ is subcanonical. 
\end{prop}

\begin{proof}
This is \cite[Prop.~2.3.9]{Warner}.
\end{proof}

\begin{prop}\label{prop:descent of coherent sheaves}
Let $f\colon Y'\rightarrow Y$ be a morphism in $\V_K$ that is locally of finite type,  flat, and surjective. Then, the category of coherent sheaves on $Y$ is equivalent to the category of coherent sheaves on $Y'$ together with descent data. 
\end{prop}

See Remark~\ref{rem:flatnesss of coherent sheaf and morphism for adic spaces} for the definition of flatness.

\begin{proof}
This is \cite[Lem.~3.1.1]{Warner} based on \cite[Thm~3.1]{Bosch-Gortz}, \cite[Thm.~4.2.8]{ConradAmpleness}, and \cite[Thm.~2.8]{Lazda-descent}.
\end{proof}

\smallskip
\noindent
\textbf{Properties of morphisms in $\V_K$}.

\begin{defn}\label{def:etale-local-on-target-in-V}
 Let $\calP$ be a property of morphisms locally of weakly finite type in $\V_K$. 
 We say that $\calP$ is \emph{\'etale local on the target}  if $\calP$ satisfies the following conditions:
\begin{enumerate}
 \item[(a)] $\calP$ is stable under base change;
 \item[(b)] for every morphism $f\colon X\rightarrow Y$ in $\V_K$ and every \'etale surjection $Z\rightarrow Y$ in $\V_K$, the morphism $f'\colon X\times_YZ\rightarrow Z$ has $\calP$ if and only if $f$ has $\calP$.
\end{enumerate}
We similarly define the notion of \emph{smooth local on the target in $\V_K$} by replacing ``\'etale surjection'' in (b) with ``smooth surjection''. 
\end{defn}

\begin{lem}\label{lem:reduction of smooth locality to finite etale case}
A property $\calP$ of morphisms locally of weakly finite type in $\V_K$ is smooth local on the target if it satisfies the following conditions:
\begin{enumerate}
 \item[(a)] $\calP$ is stable under base change;
 \item[(b)] a morphism $f\colon X\rightarrow Y$ in $\V_K$ has $\calP$ if there exists an open covering $Y=\bigcup_i Y_i$ such that the restriction $X\times_YY_i\rightarrow Y_i$ has $\calP$ for every $i$;
 \item[(c)] a morphism $f\colon X\rightarrow Y$ in $\V_K$ satisfies $\calP$ if there exists a finite \'etale surjection $Z\rightarrow Y$ in $\V_K$ such that $X\times_YZ\rightarrow Z$ has $\calP$.
\end{enumerate}
\end{lem}

\begin{proof}
Take any morphism $f\colon X\rightarrow Y$ in $\V_K$ and assume that there exists a smooth surjection $g\colon Z\rightarrow Y$ in $\V_K$ such that $f_Z\coloneqq f\times \id_Z\colon X\times_YZ\rightarrow Z$ has $\calP$. We need to show that $f$ has $\calP$. Let us start with reduction steps to impose additional conditions on $g$:  by \cite[Prop.~3.1.11]{Warner}, there exists a morphism $Z'\rightarrow Z$ in $\V_K$ such that the composite $Z'\rightarrow Z\rightarrow Y$ is \'etale and surjective. Hence by (a), we may replace $Z$ with $Z'$ and assume that $g$ is an \'etale surjection. We know from \cite[Lem.~2.2.8]{Huber-etale} that there exists an open covering $Z=\bigcup_i Z_i$ such that for every $i$, the morphism $g|_{Z_i}\colon Z_i\rightarrow Y$ factors as $Z_i\hookrightarrow  V_i\rightarrow Y_i\hookrightarrow Y$ where $Z_i\hookrightarrow  V_i$ and $Y_i\hookrightarrow Y$ are open immersions, and $V_i\rightarrow Y_i$ is finite \'etale. By shrinking $Y_i$ to $g(Z_i)$, we may further assume that $Z_i\rightarrow Y_i$ is surjective. Using (b), we are finally reduced to the following case: $g\colon Z\rightarrow Y$ is an \'etale surjection between quasi-compact adic spaces such that $g$ factors as $Z\hookrightarrow V\rightarrow Y$ where $Z\hookrightarrow V$ is an open immersion and $V\rightarrow Y$ is finite \'etale.

We now show that $f_Z$ having $\calP$ implies $f$ having $\calP$ by induction on the degree $d$ of $V\rightarrow Y$. We follow \cite[Prop.~3.2.2, Pf.]{deJong-vanderPut}. When $d=1$, then $g$ is an isomorphism. Assume $d>1$ and observe that the diagonal $\Delta(V)\subset V\times_YV$ is a union of connected components of $V\times_YV$ (which follows from the scheme case by \cite[Ex.~1.6.6(ii)]{Huber-etale}). Hence $\pr_2\colon W\coloneqq V\times_YV\smallsetminus \Delta(V)\rightarrow V$ is finite \'etale of degree less than $d$. Set $Z'\coloneqq (Z\times_YV)\cap W$ and $Y'\coloneqq \pr_2(Z')$. Then $g'\coloneqq \pr_2|_{Z'}\coloneqq Z'\rightarrow Y'$ is an \'etale surjection between quasi-compact adic spaces and it factors as $Z'\hookrightarrow W\cap \pr_2^{-1}(Y')=W\times_VY'\rightarrow Y'$ where the first map is an open immersion and the second is finite \'etale of degree less than $d$. Assume that $f_Z$ has $\calP$. Then so does $f_{Z'}$ by (a). We know by induction that $f_{Y'}$ has $\calP$. Now observe that $V=Y'\cup Z$ is an open covering: for $v\in V\smallsetminus Z$, choose $z\in Z$ with the same image in $Y$; then $(z,v)\in Z'$ and thus $v\in Y'$. We deduce from (b) that $f_V$ has $\calP$. Finally, we conclude by (c) that $f$ has $\calP$.
\end{proof}

\begin{prop}\label{prop:smooth-local-on-target}
The following properties are smooth local on the target in $\V_K$: locally of finite type; universally closed; universally open; separated; proper; finite; closed immersion (defined by a square-zero coherent ideal sheaf); open immersion; isomorphism; surjective; monomorphism; quasi-finite; unramified; \'etale; smooth; flat and locally of finite type.
\end{prop}

\begin{proof}
Each of the listed properties is known to be stable under base change and (analytic-)local on the target. Hence it remains to check (c) in Lemma~\ref{lem:reduction of smooth locality to finite etale case}, which is done in \cite[Prop.~3.1.12, Pf.]{Warner}.  
\end{proof}

\begin{defn}\label{def:etale local on source-and-target}
A property $\calP$ of morphisms in $\V_K$ is said to be \emph{\'etale (resp.~smooth) local on the source-and-target} if $\calP$ is \'etale (resp.~smooth) local on the target and satisfies the following property: for every map $f\colon X\rightarrow Y$ in $\V_K$ and every \'etale (resp.~smooth) surjection $g\colon X'\rightarrow X$, $f$ has $\calP$ if and only if $f\circ g$ has $\calP$.
\end{defn}

\begin{prop}
\hfill
\begin{enumerate}
\item The following properties are smooth local on the source-and-target in $\V_K$: surjective, locally of finite type, flat and locally of finite type, smooth.
\item The following properties are \'etale local on the source-and-target in $\V_K$: quasi-finite, unramified, \'etale.
\end{enumerate}    
\end{prop}

We present a modified version of \cite[Prop.~3.1.14, Pf.]{Warner} since we are unable to check a reduction argument therein in the smooth case.

\begin{proof}
Keep the notation as in Definition~\ref{def:etale local on source-and-target}. Each property $\calP$ above is stable under composition, and $g$ has $\calP$. By Proposition~\ref{prop:smooth-local-on-target}, it suffices to show that $f$ has $\calP$ if $f\circ g$ has $\calP$. Moreover, an argument similar to the second part of the proof of Lemma~\ref{lem:reduction of smooth locality to finite etale case} shows that in checking the \'etale local on the source-and-target property, we may further assume that $g$ is finite \'etale.

(i) The claim is obvious for surjectivity. For local finite typeness, we know from \cite[Prop.~3.1.11]{Warner} that there exists a morphism $g'\colon X''\rightarrow X'$ such that $g\circ g'$ is an \'etale surjection. Since $f\circ g\circ g'$ is still locally of finite type, we may reduced to the case where $g$ is an \'etale surjection. By the above remark, we may further assume that $g$ is finite \'etale, and this case is proved in the third paragraph of \cite[Prop.~3.1.14, Pf.]{Warner}. Moreover, it implies the case where $\calP$ is flat and locally of finite type since $g$ is flat and surjective (see Remark~\ref{rem:flatnesss of coherent sheaf and morphism for adic spaces}).

If $\calP$ is smoothness, it is \'etale local on the source-and-target in $\V_K$: by the above remark, we may assume that $g$ is finite \'etale, and this case is proved in the fifth paragraph of \cite[Prop.~3.1.14, Pf.]{Warner}. Now suppose that $g$ is a smooth surjection and $f\circ g$ is smooth. By \cite[Prop.~3.1.11]{Warner}, there exists an \'etale covering $\{j_i\colon U_i\rightarrow X\}$ together with a closed immersion $s_i\colon U_i\hookrightarrow X_i'$ into an open subspace $X_i'\subset X'$ such that $j_i$ is separated and $j_i=g\circ s_i$ for each $i$. Then $g_i\colon U_i\times_XX'_i\rightarrow U_i$ is a smooth surjection such that it admits a section induced by $s_i$ and the composite $f\circ j_i\circ g_i$ is smooth. Note that $s_i\colon U_i\hookrightarrow U_i\times_XU_i\hookrightarrow U_i\times_XX'_i$ is a closed immersion.
If we can show that $f\circ j_i\colon U_i\rightarrow X$ is smooth, we can conclude that $f$ is smooth by the \'etale locale property of smoothness. This means that we can reduce to the case where $g$ admits a section $s\colon X\hookrightarrow X'$ by a closed immersion.

Let $\calI\subset \calO_{X'}$ be the coherent ideal sheaf defining $s$. Since $g$ and $g\circ f$ are smooth, the sheaves of differentials $\Omega_{X'/Y}$ and $\Omega_{X'/X}$ are locally free $\calO_{X'}$-modules and sit in the exact sequence
$0\rightarrow g^\ast\Omega_{X/Y}\rightarrow\Omega_{X'/Y}\rightarrow\Omega_{X'/X}\rightarrow 0$
by \cite[Prop.~1.6.3(i), 1.6.9(i)(iii)]{Huber-etale}. In particular, $g^\ast\Omega_{X/Y}$ is also locally free.
Consider the following commutative diagram
\[
\xymatrix{
\calI \ar[r]^-{d_{X'/Y}}\ar[rd]_-{d_{X'/X}} &\Omega_{X'/Y}\ar[d]\\
& \Omega_{X'/X}.
}
\]
Take any $x\in X$. Since $g$ and $\id_X$ are smooth, there exist generators $a_1,\ldots,a_n\in \calI_x$ such that $d_{X'/X}(a_1),\ldots, d_{X'/X}(a_n)\in (\Omega_{X'/X})_{s(x)}$ can be extended to a basis of the $\calO_{X',s(x)}$-module $(\Omega_{X'/X})_{s(x)}$ by \cite[Prop.~1.6.9(ii)]{Huber-etale}. Then the previous discussion shows that $d_{X'/Y}(a_1),\ldots, d_{X'/Y}(a_n)\in (\Omega_{X'/Y})_{s(x)}$ can be extended to a basis of the $\calO_{X',s(x)}$-module $(\Omega_{X'/Y})_{s(x)}$. Since $g\circ f$ is smooth, we conclude from \cite[Prop.~1.6.9(ii)]{Huber-etale} that $f=s\circ (g\circ f)$ is smooth at $x$, finishing the proof of the smooth case.

(ii) By the beginning remark, we may assume that $g$ is finite \'etale, and this case is proved in the sixth paragraph of \cite[Prop.~3.1.14, Pf.]{Warner}.
\end{proof}

\smallskip
\noindent
\textbf{\'Etale equivalence relation}.

\begin{defn}\label{def:etale equivalence relation}
Let $U\rightarrow S$ be a morphism in $\V_K$ that is locally of finite type.
An \emph{\'etale equivalence relation} on $U$ over $S$ is an equivalence relation $(R,s,t\colon R\rightrightarrows  U)$ such that $R\in \Ob(\V_{K}/S)$ and such that $s$ and $t$ are \'etale.

A \emph{quotient} of an \'etale equivalence relation $s,t\colon R\rightrightarrows  U$ on $U$ over $S$ consists of an adic space $Q$ locally of type over $S$ together with an 
\'etale surjection $q\colon U\rightarrow Q$ such that $s\circ q=t\circ q$ and the resulting map $R\rightarrow U\times_QU$ is an isomorphism.
\end{defn}

\begin{thm}[Warner]\label{thm:effective-etale-equivalence-relation}
Let $s,t\colon R\rightrightarrows  U$ be an \'etale equivalence relation over $S$ in $\V_K$. If $(t,s)\colon R\rightarrow U\times_SU$ is a closed immersion, then a quotient exists and the resulting adic space is locally of finite type and separated over $S$.
\end{thm}

Conrad and Temkin proved a similar result on the existence of quotients for \'etale equivalence relations with closed diagonals in the framework of Berkovich spaces \cite[Thm.~1.2.2]{ConradTemkin1}. Warner adopted their proof in the adic space setup.

\begin{proof}
This is \cite[Thm.~3.1.5]{Warner}. The assertion that the resulting quotient is locally of finite type over $S$ is not stated in \textit{loc.~cit.}, but it follows from the proof therein and \cite[Cor.~2.10.6]{Warner}.
\end{proof}

\begin{cor}\label{cor:Warner-3.1.9}
Let $T\rightarrow S$ be an \'etale surjection in $\V_K$. Given a morphism $f\colon X\rightarrow T$ locally of finite type and separated and an isomorphism $\phi\colon X\times_ST\rightarrow T\times_SX$ over $T\times_ST$ satisfying the cocycle condition $\pr_{13}^\ast\phi=\pr_{23}^\ast\phi\circ \pr_{12}^\ast\phi$ over $T\times_ST\times_ST$ (where $\pr_{ij}\colon T\times_ST\times_ST\rightarrow T\times_ST$ are the obvious projections), there exists an adic space $Y$ locally of finite type and separated over $S$ such that $(X,\phi)$ is isomorphic to the canonical descent datum $(Y\times_ST, Y\times_S(T\times_ST)\xrightarrow{
\cong} (T\times_ST)\times_SY)$.
\end{cor}

\begin{proof}
We follow \cite[\S5]{Kiehl-aquivalenzrelationen} and \cite[Cor.~3.1.9]{Warner}. 
Set 
\[
R\coloneqq T\times_SX,\quad s\coloneqq\pr_2\colon R\rightarrow X,\quad \text{and}\quad t\coloneqq \pr_{1}\circ\phi^{-1}\colon R\rightarrow X\times_ST\rightarrow X. 
\]
Thanks to Theorem~\ref{thm:effective-etale-equivalence-relation}, it suffices to show that $s,t\colon R\rightrightarrows  X$ is an \'etale equivalence relation and $(t,s)\colon R\rightarrow X\times_SX$ is a closed immersion. The former is easy. For the latter, it is enough to show that $(t,s)\times \id_T\colon R\times_ST\rightarrow (X\times_SX)\times_ST$ is a closed immersion: since $T\rightarrow S$ is \'etale and surjective, we can apply the descent of closed immersions in Proposition~\ref{prop:smooth-local-on-target}.

By assumption, $\Delta\colon X\rightarrow X\times_TX$ is a closed immersion. By base-changing $\Delta$ along $R\times_TR\coloneqq R\times_{f\circ s,T,f\circ t}R\rightarrow X\times_TX$, we see that so is $a\colon R\times_{s,X,t}R\rightarrow R\times_TR$. It now suffices to construct a commutative diagram of the form\footnote{A similar commutative diagram appears in \cite[\S5]{Kiehl-aquivalenzrelationen} and \cite[(3.1.1)]{Warner}, but we are unable to verify that their right vertical map is an isomorphism, so we provide an alternative diagram.}
\[
\xymatrix{
R\times_{s,X,t}R\ar[r]^-{a}\ar[d]^-{\cong}_-{b} & R\times_TR \ar[d]_-{\cong}^-{c} \\
R\times_ST\ar[r]^-{(s,t)\times \id_T} & (X\times_SX)\times_ST.
}
\]
We define the isomorphism $b$ as the composite
\[
b\colon R\times_{s,X,t}R\underset{\cong}{\xrightarrow{\id_R\times \phi^{-1}}}(T\times_SX)\times_{\pr_2,X,\pr_1}(X\times_ST)
\underset{\cong}{\xrightarrow{b_1}}(T\times_SX)\times_ST=R\times_ST,
\]
where $b_1$ is given in terms of a functor of points by 
\[
b_1((\tau_1,\xi),(\xi,\tau_2))\coloneqq \bigl(\bigl(\tau_2,t(\tau_1,\xi)\bigr),f(\xi)\bigr).
\]
Using $\phi^{-1}=(t,f\circ s)$ and $\id_R=(\pr_1,s)$, observe that $b_1$ is an isomorphism with the inverse given by 
\[
((\tau,\xi),\tau')\mapsto \bigl(\bigl(\pr_1(\phi(\xi,\tau')),s(\phi(\xi,\tau'))\bigr),\bigl(s(\phi(\xi,\tau')),\tau\bigr)\bigr)=\bigl(\phi(\xi,\tau'),\bigl(s(\phi(\xi,\tau')),\tau\bigr)\bigr).
\]
Let $\sw\colon X\times_ST\xrightarrow{\cong} T\times_SX$ and $\sw'\colon T\times_SX\xrightarrow{\cong} X\times_ST$ be the maps swapping the components (e.g., $\sw(\xi,\tau)=(\tau,\xi)$).
The isomorphism $c$ is defined as the composite
\begin{align*}
c\colon R\times_TR&\underset{\cong}{\xrightarrow{\id_R\times \phi^{-1}}}(T\times_SX)\times_{f\circ\pr_2,T,f\circ\pr_1}(X\times_ST)\\
&\underset{\cong}{\xrightarrow{c_1}}
(X\times_ST)\times_{\pr_2,T,\pr_1}(T\times_SX)=X\times_ST\times_SX\\
&\underset{\cong}{\xrightarrow{c_2}}(X\times_SX)\times_ST,
\end{align*}
where $c_1\coloneqq \phi^{-1}\times(\sw'\circ \phi^{-1}\circ\sw)$ and $c_2\coloneqq \id_X\times \sw'$; equivalently, they are given in terms of functors of points by 
\begin{align*}
c_1((\tau_1,\xi_1),(\xi_2,\tau_2))&\coloneqq\bigl(\phi^{-1}(\tau_1,\xi_1),\bigl(f(\xi_2),t(\tau_2,\xi_2)\bigr)\bigr);\\
c_2(\xi_1,\tau,\xi_2)&\coloneqq \bigl(\bigl(\xi_1,\xi_2\bigr),\tau\bigr).
\end{align*}
By construction, $c$ is an isomorphism.%\]
It remains to show the commutativity of the diagram, or equivalently, $c\circ a\circ (\id_R\times\phi)=\bigl((s,t)\times\id_T\bigr)\circ b_1$. Again using $\phi^{-1}=(t,f\circ s)=(t,f\circ \pr_2)$, we compute
\begin{align*}
c(a((\id_R\times\phi)((\tau_1,\xi),(\xi,\tau_2))))&=c_2(c_1((\tau_1,\xi),(\xi,\tau_2)))=\bigl(\bigl(t(\tau_1,\xi),t(\tau_2,\xi)\bigr),f(\xi)\bigr),\\
((s,t)\times\id_T)(b_1((\tau_1,\xi),(\xi,\tau_2)))&=\bigl(\bigl(t(\tau_1,\xi),t(\tau_2,t(\tau_1,\xi))\bigr),f(\xi)\bigr).
\end{align*}
Hence we need to verify $t(\tau_2,t(\tau_1,\xi))=t(\tau_2,\xi)$.
For this, we apply the inverse of the cocycle condition $\pr_{13}^\ast\phi^{-1}=\pr_{12}^\ast\phi^{-1}\circ \pr_{23}^\ast\phi^{-1}$ to $(\tau_2,\tau_1,\xi)$:
\begin{align*}
\pr_{13}^\ast\phi^{-1}(\tau_2,\tau_1,\xi)&=\bigl(t(\tau_2,\xi),\tau_1,f(\xi)\bigr),\\
\pr_{12}^\ast\phi^{-1}\bigl(\pr_{23}^\ast\phi^{-1}(\tau_2,\tau_1,\xi)\bigr)&=\pr_{12}^\ast\phi^{-1}\bigl(\tau_2,t(\tau_1,\xi),f(\xi)\bigr)\\
&=\bigl(t(\tau_2,t(\tau_1,\xi)),f(t(\tau_1,\xi)),f(\xi)\bigr).
\end{align*}
Comparing the first component gives the desired equality.
\end{proof}

\smallskip
\noindent
\textbf{Representable sheaves and representable morphisms}. Let $S\in \Ob(\V_K)$.

\begin{defn}
\hfill
\begin{enumerate}
\item For $Z\in \Ob(\V_K/S)$, define $h_{Z/S}\in \Sh((\Afd_{K}/S)_\Et)$ to be the sheaf $T\mapsto \Hom_{S}(T,Z)$ (recall Proposition~\ref{prop:subcanonical}). 
If $Z$ is locally of finite type over $S$, then for any morphism $S'\rightarrow S$ in $\V_K$, the restriction of $h_{Z/S}$ to $(\Afd_K/S')_\Et$ is given by $h_{Z\times_SS'/S'}$. We usually write $h_Z$ for $h_{Z/S}$.
\item  A sheaf $F$ on $(\Afd_{K}/S)_\Et$ is said to be \emph{representable} if it is isomorphic to $h_{Z/S}$ for some $Z\in \Ob(\V_{K}/S)$ that is \emph{locally of finite type} over $S$. Note that we do not require $Z$ be an affinoid.
\end{enumerate}
\end{defn}

\begin{defn}\label{def:representable by adic spaces}
Let $f\colon F\rightarrow G$ be a morphism of presheaves on $\Afd_{K}/S$.
\begin{enumerate}
\item  We say that $f$ is \emph{representable by adic spaces} if, for every $T\in\Ob(\V_{K}/S)$ and every morphism $h_T\rightarrow G$, the presheaf pullback $h_{T}\times_{G}F\rightarrow h_T$ is isomorphic to $h_Z\rightarrow h_T$ for some $Z\in \Ob(\V_K/S)$ and a morphism $Z\rightarrow T$ that is \emph{locally of finite type}. 
 \item Let $\calP$ be a property of morphisms in $\V_{K}$, which is \'etale local on the target in $\V_K$ (see Definition~\ref{def:etale-local-on-target-in-V}). We say that $f$ \emph{has property $\mathcal{P}$} if it is representable by adic spaces, and for every $T\in\Ob(\V_{K}/S)$ and every morphism $h_T\rightarrow G$, the morphism $h_T\times_{G}F\rightarrow h_T$ is represented by a morphism $Z\rightarrow T$ of adic spaces having property $\calP$.
\end{enumerate}
\end{defn}

Note that if $F$ and $G$ are sheaves on $(\Afd_K/S)_\Et$, then the sheaf pullback $h_{T}\times_{G}F\rightarrow h_T$ agrees with the presheaf pullback. For a morphism $U\rightarrow V$ of adic spaces locally of finite type over $S$, the corresponding morphism $h_U\rightarrow h_V$ is representable by adic spaces. Moreover, $U\rightarrow V$ has $\calP$ if and only if $h_U\rightarrow h_V$ has $\calP$.

\begin{lem}\label{lem:DiagonalRepresentable}
Let $F$ be a sheaf on $(\Afd_{K}/S)_{\Et}$. If the diagonal morphism $\Delta\colon F\rightarrow F\times_{h_S}F$ is representable by adic spaces, then any morphism $h_{U}\rightarrow F$ for $U\in\Ob(\V_{K}/S)$, locally of finite type over $S$, is representable by adic spaces.
\end{lem}

\begin{proof}
This is \cite[Lem.~3.2.3]{Warner}.
\end{proof}

\subsection{\'Etale spaces}
Let $S\in \V_K$.

\begin{defn}\label{def:etale spaces}
An \emph{\'etale space} over $S$ is a sheaf $F$ on $(\Afd_K/S)_\Et$ such that the following properties hold. 
\begin{enumerate}
\item The diagonal morphism $\Delta\colon F\rightarrow F\times_{h_{S}} F$ is representable by adic spaces. 
\item There exist $U\in\Ob(\V_{K}/S)$, locally of finite type over $S$, and a morphism $h_{U}\rightarrow F$ that is surjective and \'etale. Such a morphism $h_{U}\rightarrow F$ is called an \emph{\'etale cover} of $F$ and is surjective as a morphism of sheaves.
\end{enumerate}
Note that $h_U\rightarrow F$ in (ii) is necessarily representable by (i) and Lemma~\ref{lem:DiagonalRepresentable}.
A morphism of \'etale spaces over $S$ is a morphism of sheaves. 
Let $\EtSp_S\subset \Sh((\Afd_{K}/S)_\Et)$ denote the strictly full subcategory of \'etale spaces over $S$.
\end{defn}

\begin{lem}
 The fiber product exists in $\EtSp_{S}$ and agrees with the one in $\Sh((\Afd_{K}/S)_\Et)$.
 \end{lem}

\begin{proof}
The proof of \cite[Prop.~5.4.6]{Olsson-book} works in our situation.
\end{proof}

\begin{lem} 
For $F\in\Ob(\EtSp_{S})$ and $S'\in\Ob(\V_{K}/S)$, the restriction of $F$ to $\Sh((\Afd_K/S')_\Et)$ is an \'etale space over $S'$.
\end{lem}

\begin{proof}
Set $F_{S'}\coloneqq F\times_{h_{S}}h_{S'}\in\Sh((\Afd_{K}/S)_{\Et})$. By definition, $F_{S'}$ and $F$ have the same restrictions on $(\Afd_{K}/S')_{\Et}$. The diagonal map $\Delta_{S'}\colon F_{S'}\rightarrow F_{S'}\times_{h_{S'}}F_{S'}$ is the base change of the diagonal $\Delta\colon F\rightarrow F\times_{h_{S}}F$ over $h_S$ via $h_{S'}\rightarrow h_S$, so it is still representable by adic spaces. Furthermore, if $h_{U}\rightarrow F$ is an \'etale cover, then $h_{U}\times_{h_{S}}h_{S'}\rightarrow F_{S'}$ is also \'etale and surjective. Note that the fiber product $U\times_{S}S'$ exists in $\V_K$ and is locally of finite type over $S'$ since $U$ is locally of finite type over $S$. Hence $h_{U\times_{S}S'}=h_{U}\times_{h_{S}}h_{S'}\rightarrow F_{S'}$ is an \'etale cover over $(\Afd_{K}/S')_{\Et}$.
\end{proof}

\smallskip
\noindent
\textbf{Presentation of \'etale spaces}

\begin{defn}\label{def:presentation of etale space}
A \emph{presentation} of an \'etale space $F$ over $S$ is an \'etale equivalence relation $R\rightrightarrows  U$ on $U\in \Ob(\V_K/S)$, locally of finite type over $S$, together with an isomorphism $U/R\xrightarrow{\cong}F$.
\end{defn}

The following justifies the notion of presentation.

\begin{prop}\label{prop:presentation of etale adic spaces}
\hfill
\begin{enumerate}
\item For $F\in\Ob(\EtSp_{S})$ and an \'etale cover $f\colon h_U\rightarrow F$, the fiber product $h_U\times_{F}h_U$ is representable by an adic space $R$ locally of finite type over $S$. Moreover, the natural morphisms $R\rightrightarrows  U$ define an \'etale equivalence relation on $U$, and the resulting diagram $h_R\rightrightarrows  h_U\rightarrow F$ is a coequalizer diagram in $\Sh((\Afd_K/S)_\Et)$.
\item Given an \'etale equivalence relation $R\rightrightarrows  U$ on $U\in \Ob(\V_K/S)$, locally of finite type over $S$, the quotient sheaf $U/R$ is an \'etale space $F$ over $S$. Moreover, $h_U\rightarrow F$ is an \'etale cover, and $h_R\xrightarrow{\cong}h_U\times_Fh_U$.
\end{enumerate}
\end{prop}

\begin{proof}
(i) For every $T\in \Ob(\V_K/S)$ and every morphism $h_T\rightarrow F$, the fiber product $h_U\times_Fh_T$ is representable by an adic space \'etale surjective over $T$ by Lemma~\ref{lem:DiagonalRepresentable} and our assumption. Now the assertions follow from this and Example~\ref{eg:coequalizer diagram of sheaves}.

(ii) This is \cite[Prop.~3.2.11]{Warner}.
\end{proof}

\begin{prop}\label{prop:separated etale space is adic space}
Let $F\in \Ob(\EtSp_S)$. If the diagonal $F\rightarrow F\times_S F$ is a closed immersion, then $F$ is an adic space locally of finite type and separated over $S$.
\end{prop}

\begin{proof}
This is \cite[Cor.~3.2.12]{Warner} together with Theorem~\ref{thm:effective-etale-equivalence-relation}.
\end{proof}

\smallskip
\noindent
\textbf{Properties of morphisms of \'etale spaces}.
We introduce several properties of morphisms of \'etale spaces that are not necessarily representable by adic spaces.

\begin{defn}\label{def:morphism of etale spaces has property etale local on source-and-target} 
Let $\calP$ be a property of morphisms in $\V_{K}$ which is \'etale local on the source-and-target (Definition~\ref{def:etale local on source-and-target}): e.g.~surjective, \'etale, or smooth. A morphism $f\colon X\rightarrow Y$ in $\EtSp_S$ is said to \emph{have $\mathcal{P}$} if there exist \'etale covers $h_U\rightarrow X$ and $h_V\rightarrow Y$ such that the induced map $h_U\times_{Y}h_V\rightarrow h_V$ corresponds to a morphism in $\V_K$ having $\calP$ (note that $h_U\times_{Y}h_V$ is representable by Lemma~\ref{lem:DiagonalRepresentable}).
\end{defn}

\begin{lem}\label{lem:sanity check of etale local on source-and-target}
Let $\calP$ be a property of morphisms in $\V_{K}$ which is \'etale local on the source-and-target, and let $f\colon X\rightarrow Y$ be a morphism of \'etale spaces over $S$. 
\begin{enumerate}
\item If $f$ has $\calP$, then $h_U\times_{Y}h_V\rightarrow h_V$ corresponds to a morphism  
 in $\V_K$ having $\calP$ for any \'etale covers $h_U\rightarrow X$ and $h_V\rightarrow Y$.
\item Assume that $f$ is representable by adic spaces. Then $f$ has $\calP$ in the sense of Definition~\ref{def:morphism of etale spaces has property etale local on source-and-target} 
if and only if $f$ has $\calP$ in the sense of Definition~\ref{def:representable by adic spaces}(ii).
\end{enumerate}
\end{lem}

\begin{proof}
Part (ii) follows from (i).
For (i), consider a commutative diagram
\[
\xymatrix{
h_{U'}\times_Y h_{V'}\ar[r]\ar[d]&h_{U}\times_Y h_{V'}\ar[rr]\ar[d] && h_{V'}\ar[d]^-d\\
h_{U'}\times_Y h_{V}\ar[r]\ar[d]& h_{U}\times_Y h_{V} \ar[rr]\ar[d] && h_V\ar[d]^-b\\
h_{U'}\ar[r]^-c & h_U\ar[r]^-a & X\ar[r]^-f & Y,
}
\]
where $a$ and $b$ are \'etale covers, and $c\colon U'\rightarrow U$ and $d\colon V'\rightarrow V$ are \'etale surjective morphisms in $\V_K$. Note that the top two squares are Cartesian of representable sheaves. So take $Z,Z'\in \Ob(\V_K/S)$ with $h_{U}\times_Y h_{V}\cong h_Z$ and $h_{U'}\times_Y h_{V'}\cong h_{Z'}$. Since $c$ and $d$ are \'etale surjective, $Z\rightarrow V$ has $\calP$ if and only if $Z'\rightarrow V'$ has $\calP$.

Now given two pairs of \'etale covers $(h_{U_1}\rightarrow X,h_{V_1}\rightarrow Y)$ ($i=1,2$), consider $U=U_i$. $V=V_i$,  $h_{U'}=h_{U_1}\times_X h_{U_2}$, and $h_{V'}=h_{V_1}\times_X h_{V_2}$. The preceding argument shows that $h_{U_1}\times_{Y}h_{V_1}\rightarrow h_{V_1}$ corresponds to a morphism in $\V_K$ having $\calP$ if and only if $h_{U_2}\times_{Y}h_{V_2}\rightarrow h_{V_2}$ corresponds to a morphism in $\V_K$ having $\calP$, which yields (i). 
\end{proof}

\begin{rem}\label{rem:smooth surjective morphism of etale spaces is epimorphism}
Let $f\colon X\rightarrow Y$ be a morphism in $\EtSp_S$. If $f$ is an epimorphism of sheaves, then it is surjective in the sense of Definition~\ref{def:morphism of etale spaces has property etale local on source-and-target}. Conversely, if $f$ is smooth and surjective in this sense, then it is an epimorphism of sheaves by the existence of \'etale quasi-sections to smooth surjections \cite[Prop.~3.1.11]{Warner}.
\end{rem}

\begin{lem}\label{lem:image of smooth morphism of etale spaces}
Let $f\colon X\rightarrow Y$ be a smooth morphism in $\EtSp_S$ and write $f(X)$ for the image of $f$ as a sheaf. Then $f(X)$ is an \'etale space over $S$, and the inclusion $f(X)\hookrightarrow Y$ is representable by adic spaces and is an open immersion.
\end{lem}

\begin{proof}
Take \'etale covers $h_U\rightarrow X$ and $h_V\rightarrow Y$ such that the induced map $h_U\times_{Y}h_V\rightarrow h_V$ is smooth. Since $h_U\times_{Y}h_V\rightarrow X$ is an epimorphism of sheaves, we may replace $X$ with $h_U\times_{Y}h_V$ and assume that $X$ is an adic space and $f$ factors as $h_X\rightarrow h_V\rightarrow Y$ where $V$ is also an adic space, $X\rightarrow V$ is smooth, and $h_V\rightarrow Y$ is an \'etale cover. For any $T\in \Ob(\V_K/S)$ and any morphism $h_T\rightarrow Y$, consider the fiber product $f(X)\times_Yh_T$ of sheaves. It is straightforward to see that $f(X)\times_Yh_T \rightarrow h_T$ is a monomorphism of sheaves and agrees with the sheaf image of $f_T\colon h_X\times_Yh_T\rightarrow h_T$. By replacing $f$ by $f_T$, it suffices to show that the sheaf image $f(X)$ is an open subspace of $Y$ under the additional assumption that $Y$ is also an adic space. Now $f\colon X\rightarrow Y$ is a smooth morphism of adic spaces and thus an open map. Let $W\subset Y$ be the image of $f$ as a topological space. It remains to show  $f(X)=h_W$. The inclusion $f(X)\subset h_W$ is obvious; since $X\rightarrow W$ is smooth and surjective, we deduce from \cite[Prop.~3.1.11]{Warner} that $f(X)=h_W$. 
\end{proof}

\smallskip
\noindent
\textbf{Bootstrap}.
The following definitions and results will be used in the discussion of adic stacks.

\begin{defn}
A morphism $F\rightarrow G$ of presheaves on $\Afd_{K}/S$ is said to be \emph{representable by \'etale spaces} if, for every $T\in\Ob(\V_{K}/S)$ and every morphism $h_T\rightarrow G$, the presheaf fiber product $h_{T}\times_{G}F$ is an \'etale space over $T$. 
\end{defn}

\begin{lem}\label{lem:bootstrap for morphisms represetanble by etale spaces}
Let $F\rightarrow G$ be a morphism of presheaves on $\Afd_{K}/S$ that is representable by \'etale spaces. 
\begin{enumerate}
\item If $G$ is an \'etale space over $S$, then so is $F$.
\item  Let $T\in \Ob(\V/S)$ and regard $F\rightarrow G$ as a morphism of presheaves on $\Afd_K/T$. For every \'etale space $X$ over $T$ and every morphism $X\rightarrow G$ of presheaves on $\Afd_K/T$, the presheaf fiber product $F\times_GX$ is an \'etale space over $T$.
\end{enumerate}
\end{lem}

\begin{proof}
Part (i) is \cite[Lem.~3.2.20]{Warner}.
For (ii), it is formal to check that $F\times_GX\rightarrow X$ is representable by \'etale spaces (as a morphism of presheaves on $\Afd_K/T$). Hence $F\times_GX$ is an \'etale sheaf over $T$ by (i). 
\end{proof}

\begin{cor}\label{cor:bootstrap (composite of morphisms representable by etale spaces)}
If $a\colon F\rightarrow G$ and $b\colon G\rightarrow H$ are morphisms of presheaves on $\Afd_{K}/S$ that are representable by \'etale spaces, then so is $b\circ a\colon F\rightarrow H$.
\end{cor}

\begin{proof}
This follows from Lemma~\ref{lem:bootstrap for morphisms represetanble by etale spaces}(ii).
\end{proof}

\begin{defn}
Let $F\rightarrow G$ be a morphism of presheaves on $\Afd_{K}/S$ that is representable by \'etale spaces. For a property $\calP$ of morphisms in $\V_K$ that is \'etale local on the source-and-target, we say that $F\rightarrow G$ \emph{has $\calP$} if for every $T\in\Ob(\V_{K}/S)$ and every morphism $h_T\rightarrow G$, the morphism $h_{T}\times_{G}F\rightarrow h_T$ comes from a morphism of \'etale spaces over $T$ having $\calP$ in the sense of Definition~\ref{def:morphism of etale spaces has property etale local on source-and-target}.
\end{defn}

This definition applies to, for example, surjective, \'etale, or smooth morphisms.

\begin{prop}\label{prop:bootstrap for etale spaces (sheaf admitting etale surjective map from etale sapce)}
Let $a\colon F\rightarrow G$ be a morphism of sheaves on $(\Afd_{K}/S)_\Et$. Assume that $F$ is an \'etale space over $S$ and that $a$ is representable by \'etale spaces and \'etale surjective.
Then $G$ is an \'etale space over $S$.
\end{prop}

This is an analogue of \cite[03Y3]{stacks-project} and implies \cite[Lem.~3.2.23]{Warner}. When $a$ is representable by adic spaces and \'etale surjective, it follows from Proposition~\ref{prop:presentation of etale adic spaces}(ii).

\begin{proof}
Since $F$ is an \'etale space over $S$, take an \'etale cover $U\rightarrow F$. Then the composite $U\rightarrow F\rightarrow G$ is representable by \'etale spaces and \'etale surjective by Corollary~\ref{cor:bootstrap (composite of morphisms representable by etale spaces)}. Now it suffices to show that the diagonal $G\rightarrow G\times_S G$ is representable by adic spaces: then $U\rightarrow G$ is representable by adic spaces Lemma~\ref{lem:DiagonalRepresentable} and is \'etale surjective by Lemma~\ref{lem:sanity check of etale local on source-and-target}.
Note that $R\coloneqq U\times_GU$ is an \'etale space over $S$ by Lemma~\ref{lem:bootstrap for morphisms represetanble by etale spaces}(ii). 

Take any $W\in\Ob(\V_K/S)$ and any morphism $W\rightarrow G\times_S G$. We need to show that the sheaf $H\coloneqq G\times_{G\times_S G}W$ on $(\Afd_K/S)_\Et$ is an adic space locally of finite type over $W$.
Since $U\rightarrow G$ is representable by \'etale spaces and \'etale surjective, it is an epimorphism of sheaves. Hence there is an \'etale surjection $W'\rightarrow W$ of adic spaces such that $W'\rightarrow W\rightarrow G\times_S G$ factors as $W'\rightarrow U\times_SU\rightarrow G\times_S G$. Set $H'\coloneqq 
G\times_{G\times_S G}W'\in \Sh((\Afd_K/S)_\Et)$. Then $H'$ is identified with $R\times_{U\times_SU}W'$ and thus is an \'etale space over $W'$. Since $R\rightarrow U\times_SU$ is a monomorphism, so is $H'\rightarrow W'$. Hence $H'\rightarrow H'\times_{W'}H'$ is an isomorphism. We know from Proposition~\ref{prop:separated etale space is adic space} that $H'$ is an adic space locally of finite type over $W'$.
Since $W'\rightarrow W$ is an \'etale surjection of adic spaces, $H'\rightarrow H$ is an \'etale surjective representable morphism of sheaves. In particular, $H'\times_HH'$ is an adic space locally of finite type over $W$, and $H'\times_HH'\rightrightarrows  H'$ is an \'etale equivalence relation with $H$ being the quotient sheaf. So $H$ is an \'etale space over $W$ by Proposition~\ref{prop:presentation of etale adic spaces}(ii). Since $H'\rightarrow W'$ is a monomorphism and $W'\rightarrow W$ is \'etale surjective, $H\rightarrow W$ is a monomorphism, which implies that $H\rightarrow H\times_WH$ is an isomorphism. Hence $H$ is an adic space locally of finite type over $W$ by Proposition~\ref{prop:separated etale space is adic space}.
\end{proof}

\subsection{Adic stacks}
We continue to fix $S\in \V_K$ unless otherwise stated.

\smallskip
\noindent
\textbf{Representable stacks and morphisms by \'etale spaces}.
For a presheaf $F$ on $\Afd_K/S$, let $\calS_{F}\rightarrow \Afd_{K}/S$ denote the category fibered in groupoids given in Definition~\ref{def:category fibered in setoids associated to presheaf}.

\begin{defn}
Let $\calX\rightarrow\Afd_{K}/S$ be a category fibered in groupoids.
\begin{enumerate}
\item We say that $\calX$ is \emph{representable by an adic space} if there exists $T\in\Ob(\V_{K}/S)$, locally of finite type over $S$, and an equivalence of categories $\calX\xrightarrow{\cong}\Afd_{K}/T$  over $\Afd_{K}/S$. 
\item We say that $\calX$ is \emph{representable by an \'etale space} if there exists $F\in\Ob(\EtSp_{S})$ and an equivalence of categories $\calX\xrightarrow{\cong} \calS_{F}$ over $\Afd_{K}/S$.
\end{enumerate}
\end{defn}

\begin{defn}
A $1$-morphism $f\colon \calX\rightarrow\calY$ of categories fibered in groupoids over $\Afd_{K}/S$ is \emph{representable by \'etale spaces} (resp.~\emph{by adic spaces}) if, for every $T\in\Ob(\V_{K}/S)$ and every $1$-morphism $y\colon (\Afd_{K}/T)\rightarrow\calY$ over $\Afd_{K}/S$, the category fibered in groupoids $(\Afd_{K}/T)\times_{y,\calY}\calX\rightarrow \Afd_{K}/T$
is representable by an \'etale space (resp.~by an adic space).
\end{defn}

\begin{lem}\label{lem:sanity check for representable by adic spaces for etale spaces and associated categories fibered in groupoids}
A morphism $f\colon F\rightarrow G$ of presheaves on $\Afd_{K}/S$ is representable by adic spaces (resp.~\'etale spaces) if and only if the associated $1$-morphism $f'\colon \calS_{F}\rightarrow\calS_{G}$ over $\Afd_{K}/S$ is representable by adic spaces (resp.~\'etale spaces).
\end{lem}

\begin{proof}
We only treat the representability by adic spaces as the representability by \'etale spaces is similar.
Suppose that $f\colon F\rightarrow G$ is representable by adic spaces. Let $T\in\Ob(\V_{K}/S)$ be an adic space together with a $1$-morphism $y'\colon (\Afd_{K}/T)\rightarrow\calS_G$ of categories fibered in groupoids over $\Afd_{K}/S$. By Remark~\ref{rem:category fibered in setoids associated to presheaf}, $y'$ gives a morphism of sheaves $y\colon h_{T}\rightarrow G$.
By assumption, the fiber product $F\times_{f,G,y}h_{T}$ is identified with $h_{V}$ for some $V\in\Ob(\V_{K}/T)$ locally of finite type over $T$. By using Remark~\ref{rem:category fibered in setoids associated to presheaf} again, we see that $y$ yields an equivalence $p\colon \calS_{F}\times_{f',\calS_{G},y'}(\Afd_{K}/T)\xrightarrow{\cong}\Afd_{K}/V$ over $\Afd_{K}/T$. Thus $f'\colon \calS_{F}\rightarrow\calS_{G}$ is representable by adic spaces. The converse direction follows similarly from Remark~\ref{rem:category fibered in setoids associated to presheaf}.
\end{proof}

Below, we also use $f$ to denote the associated $1$-morphism $\calS_{F}\rightarrow\calS_{G}$.

\begin{lem}\label{lem:DiagonalRepresentable-Stacks}
Let $f\colon\calX\rightarrow(\Afd_{K}/S)_{\Et}$ be a category fibered in groupoids. Consider the following three properties:
\begin{enumerate}
\item the diagonal $1$-morphism $\Delta\colon\calX\rightarrow\calX\times\calX$ is representable by \'etale spaces;
\item for every $T\in\Ob(\V_{K}/S)$ and every pair of $1$-morphisms $x,y\colon (\Afd_K/T)\rightarrow \calX$ over $\Afd_K/S$, the presheaf $\Isom(x,y)$ on $\Afd_{K}/T$ defined by 
\[
\Isom(x,y)(p\colon V\rightarrow T)\coloneqq\Isom_{\calX_{V}}(x\circ p,y\circ p)
\]
is an \'etale space over $T$, where $1$-morphisms $x\circ p,y\circ p\colon (\Afd_K/V)\rightarrow \calX$ over $\Afd_K/S$ are regarded as objects in $\calX_V$ by the $2$-Yoneda lemma;
\item for every $T\in\Ob(\V_{K}/S)$ locally of finite type over $S$ and every $1$-morphism $x\colon (\Afd_K/T)\rightarrow\calX$ over $\Afd_{K}/S$, $x$ is representable by \'etale spaces. 
\end{enumerate}
Then \emph{(i)} and \emph{(ii)} are equivalent, and these imply \emph{(iii)}.
\end{lem}

\begin{proof}
(i)$\iff$(ii). 
Take any $T\in\Ob(\V_{K}/S)$ and $1$-morphisms $x,y\colon (\Afd_K/T)\rightarrow \calX$ over $\Afd_K/S$. Then the fibered category $\calX\times_{\Delta,\calX\times\calX,(x,y)}(\Afd_K/T)\rightarrow \Afd_K/T$ is fibered in setoids and naturally identified with $\calS_{\Isom(x,y)}\rightarrow \Afd_K/T$ by \cite[04SI]{stacks-project}. Hence the former is representable by an \'etale space if and only if $\Isom(x,y)$ is an \'etale space over $T$.

(i)$\implies$(iii). Take any $Z\in \Ob(\V_K/S)$. Since $T$ is locally of finite type over $S$, the fiber product $Z\times_ST$ exists in $\V_K$, and $\Afd_K/(Z\times_ST)$ is identified with $(\Afd_K/Z)\times(\Afd_K/T)$ as fibered categories over $\Afd_K/S$.
 Given a $1$-morphism $y\colon (\Afd_{K}/Z)\rightarrow\calX$ over $\Afd_{K}/S$, the $2$-fiber product $(\Afd_K/Z)\times_{y,\calX,x}(\Afd_K/T)$ is identified with $((\Afd_K/Z)\times(\Afd_K/T))\times_{(y,x),\calX\times\calX,\Delta}\calX$ by \cite[02XE, 0H2E]{stacks-project}. We know from (i) that 
$(\Afd_K/Z)\times_{y,\calX,x}(\Afd_K/T)\rightarrow \Afd_K/(Z\times_ST)$ is representable by an \'etale space. Since $Z\times_ST\rightarrow Z$ is locally of finite type,  $(\Afd_K/Z)\times_{y,\calX,x}(\Afd_K/T)\rightarrow \Afd_K/Z$ is also representable by an \'etale space.
\end{proof}

\smallskip
\noindent
\textbf{Properties of morphisms represented by \'etale spaces}.

\begin{defn}\label{def:1-morphism representable by etale spaces has property P}
Let $\calP$ be a property of morphisms in $\V_{K}$ which is smooth local on the source-and-target. A $1$-morphism $f\colon\calX\rightarrow\calY$ over $\Afd_{K}/S$ representable by \'etale spaces is said to \emph{have $\calP$} if, for every $T\in\Ob(\V_{K}/S)$ and every $1$-morphism $y\colon(\Afd_{K}/T)\rightarrow\calY$, the $1$-morphism $(\Afd_{K}/T)\times_{y,\calY}\calX\rightarrow(\Afd_{K}/T)$ corresponds to a morphism of \'etale spaces having $\calP$ in the sense of Definition~\ref{def:morphism of etale spaces has property etale local on source-and-target}.
\end{defn}

\begin{lem}\label{lem:sanity check of 1-morphism representable by etale spaces}
Let $f\colon \calX\rightarrow\calY$ be a $1$-morphism over $(\Afd_{K}/S)_\Et$ such that $\calX$ is a stack and $f$ is representable by \'etale spaces. Then for every \'etale space $F$ over $S$ and every $1$-morphism $\calS_F\rightarrow \calY$ over $\Afd_{K}/S$, the $2$-fiber product $\calS_F\times_{\calY}\calX\rightarrow \Afd_K/S$ is representable by an \'etale space.
\end{lem}

\begin{proof}
Take an \'etale cover $T\rightarrow F$. Then $(\Afd_K/T)\times_{\calY}\calX\rightarrow \Afd_K/T$ is representable by an \'etale space. In particular, it is a category fibered in setoids. Since $T\rightarrow F$ is representable by adic spaces and \'etale surjective, we deduce that $\calS_F\times_{\calY}\calX\rightarrow \Afd_K/S$ is also a category fibered in setoids. By Remark~\ref{rem:category fibered in setoids associated to presheaf}, the $1$-morphism $\calS_F\times_{\calY}\calX\rightarrow \calS_F$ over $\Afd_K/S$ comes from a morphism of presheaves $F'\rightarrow F$ on $\Afd_K/S$  for some $F'$. Moreover, $F'\rightarrow F$ is representable by \'etale spaces by Lemma~\ref{lem:sanity check for representable by adic spaces for etale spaces and associated categories fibered in groupoids}. Hence $F'$ is an \'etale space over $S$ by Lemma~\ref{lem:bootstrap for morphisms represetanble by etale spaces}(i).
\end{proof}

\begin{lem}\label{lem:1-morphism from etale space to stack is representable by etale spaces}
Let $\calX\rightarrow (\Afd_{K}/S)_{\Et}$ be a stack in groupoids such that the diagonal $1$-morphism $\Delta\colon \calX\rightarrow \calX\times\calX$ over $\Afd_K/S$ is representable by \'etale spaces.
Then every $1$-morphism $\calS_F\rightarrow \calX$ over $\Afd_K/S$ for $F\in \Ob(\EtSp_S)$ is representable by \'etale spaces.
\end{lem}

\begin{proof}
Take any $T\in \Ob(\V_K/S)$ and $1$-morphism $x\colon \Afd_K/T\rightarrow\calX$ over $\Afd_K/S$. We need to show that $f\colon (\Afd_K/T)\times_{x,\calX}\calS_F\rightarrow \Afd_K/T$ is representable by an \'etale space.
Take an \'etale cover $Z\rightarrow F$ from an adic space $Z$. By assumption and Lemma~\ref{lem:DiagonalRepresentable-Stacks}, $(\Afd_K/T)\times_{x,\calX}\calS_Z\rightarrow \Afd_K/T$ is representable by an \'etale space, say $G\in\Ob(\EtSp_T)$. 
We deduce from Lemma~\ref{lem:sanity check for representable by adic spaces for etale spaces and associated categories fibered in groupoids} that the $1$-morphism $a\colon (\Afd_K/T)\times_{x,\calX}\calS_Z\rightarrow(\Afd_K/T)\times_{x,\calX}\calS_F$ over $\Afd_K/T$ is representable by adic spaces and \'etale surjective. In particular, $(\Afd_K/T)\times_{x,\calX}\calS_F\rightarrow \Afd_K/T$ is a category fibered in setoids, and by Remark~\ref{rem:category fibered in setoids associated to presheaf}, $a$ comes from a morphism $G\rightarrow G'$ of presheaves on $\Afd_K/T$ that is representable by adic spaces and \'etale surjective. Since $f$ is a stack, $G'$ is a sheaf on $(\Afd_K/T)_\Et$ by \cite[0432]{stacks-project}. Hence $G'$ is an \'etale space over $T$ by Proposition~\ref{prop:bootstrap for etale spaces (sheaf admitting etale surjective map from etale sapce)}. 
\end{proof}

\smallskip
\noindent
\textbf{Adic stacks}.

\begin{defn}\label{def:adic stacks}
An \emph{adic stack} over $S$ is a category fibered in groupoids $\calX\rightarrow (\Afd_{K}/S)_{\Et}$ satisfying the following properties:
\begin{enumerate}
 \item it is a stack in groupoids; 
 \item the diagonal $\Delta\colon \calX\rightarrow \calX\times\calX$ is representable by \'etale spaces;
 \item there exist $X\in\Ob(\V_{K}/S)$, locally of finite type over $S$, and a $1$-morphism $f\colon (\Afd_{K}/X)\rightarrow \calX$ over $\Afd_K/S$ that is smooth and surjective (note that $f$ is representable by \'etale spaces by Lemma~\ref{lem:DiagonalRepresentable-Stacks}). Such an $f$ is called a \emph{smooth cover} of $\calX$.
\end{enumerate}
A \emph{$1$-morphism} of adic stacks over $S$ is a $1$-morphism over $(\Afd_{K}/S)_{\Et}$.
\end{defn}

\begin{lem}\label{Fiber-Product-Stacks}
The $2$-fiber product of adic stacks as fibered categories is an adic stack.
\end{lem}

\begin{proof}
The proof of \cite[04TF]{stacks-project} also works in our case, thanks to Lemma~\ref{lem:DiagonalRepresentable-Stacks}.
\end{proof}

\smallskip
\noindent
\textbf{Presentations of adic stacks}.
All adic stacks can be given as quotient stacks.

\begin{defn}\label{def:smoot groupoid in etale spaces}
A \emph{smooth groupoid in \'etale spaces} over $S$ is a groupoid $(U,R,s,t,c)$ in sheaves on $(\Afd_{K}/S)_\Et$ such that $U$ and $R$ are \'etale spaces over $S$ and such that $s,t\colon R\rightarrow U$ are smooth morphisms of \'etale spaces (cf.~Definition~\ref{def:groupoid for site}).
\end{defn}

\begin{defn}\label{def:presentation of adic stack}
A \emph{presentation} of an adic stack $\calX$ over $S$ is a smooth groupoid $(U,R,s,t,c)$ in \'etale spaces over $S$ together with an equivalence $[U/R]\xrightarrow{\cong}\calX$.
\end{defn}

A smooth surjective\footnote{As every $1$-morphism from an \'etale space to an adic stack is representable by \'etale spaces by Lemma~\ref{lem:1-morphism from etale space to stack is representable by etale spaces}, we can discuss smoothness and surjectivity as in Definition~\ref{def:1-morphism representable by etale spaces has property P}.} $1$-morphism from an \'etale space to an adic stack yields a presentation.

\begin{prop}\label{prop:smooth cover to smooth groupoid}
Let $f\colon\calX\rightarrow(\Afd_{K}/S)_{\Et}$ be an adic stack over $S$. For every $U\in\Ob(\EtSp_{S})$ and every smooth surjective  $1$-morphism $\calS_{U}\rightarrow\calX$ over $(\Afd_K/S)_\Et$, the $2$-fiber product $\calS_{U}\times_{\calX}\calS_{U}$ is represented by an \'etale space $R$ over $S$. Furthermore, the projections $t\coloneqq\pr_1, s\coloneqq\pr_2\colon R\rightarrow U$ and the morphism $c\colon R\times_{s,U,t}R\rightarrow R$ given by
\[
\calS_{R\times_{s,U,t}R}=(\calS_{U}\times_{\calX}\calS_{U})\times_{\pr_{2},\calS_{U},\pr_{1}}(\calS_{U}\times_{\calX}\calS_{U})=\calS_{U}\times_{\calX}\calS_{U}\times_{\calX}\calS_{U}\xrar{\pr_{13}}\calS_{R},
\]
define a smooth groupoid $(U,R,s,t,c)$ in \'etale spaces and induce an equivalence $[U/R]\xrightarrow{\cong}\calX$ over $(\Afd_K/S)_{\Et}$.
\end{prop}

\begin{proof}
The $2$-fiber product $\calS_{U}\times_{\calX}\calS_{U}$ is representable by an \'etale space by Lemma~\ref{lem:sanity check of 1-morphism representable by etale spaces}. 
The remaining assertions are proved just as in the algebraic case \cite[04T4]{stacks-project}.
\end{proof}

\begin{thm}\label{thm:smooth groupoid gives adic stack}
Let $(U,R,s,t,c)$ be a smooth groupoid in \'etale spaces over $S$. Then the category fibered in groupoids $[U/R]\rightarrow(\Afd_{K}/S)_{\Et}$ is an adic stack, and $1$-morphism $\calS_U\rightarrow [U/R]$ is smooth and surjective. \end{thm}

\begin{proof}
We first prove that the diagonal $\Delta\colon[U/R]\rightarrow[U/R]\times[U/R]$ is representable by \'etale spaces. By Lemma~\ref{lem:DiagonalRepresentable-Stacks}, it suffices to show that for every $T\in\Ob(\V_K/S)$ and every pair of $1$-morphisms $x,y\colon(\Afd_{K}/T)_{\Et}\rightarrow[U/R]$ over $(\Afd_{K}/S)_{\Et}$, the presheaf $\Isom(x,y)$ defined in \textit{loc.~cit.} is an \'etale space over $T$. For this, we may assume that $T$ is an affinoid. Since $[U/R]$ is a stack by definition, $\Isom(x,y)$ is a sheaf on $(\Afd_K/T)_\Et$, and the following holds by Lemma~\ref{lem:stackification}: there exists an \'etale covering $\{ T_{i}\rightarrow T\}$ by finitely many affinoids $T_i$ such that $x|_{T_{i}},y|_{T_{i}}\colon (\Afd_{K}/T_{i})_{\Et}\rightarrow [U/R]$ come from morphisms  $x_{i},y_{i}\colon T_{i}\rightarrow U$ of sheaves on $(\Afd_K/S)_\Et$. Then $\Isom(x,y)|_{T_i}$ is isomorphic to $T_{i}\times_{(y_i,x_i),U\times_SU}R$, which is an \'etale space over $T_i$ (cf.~\cite[044V]{stacks-project}). Set $T'\coloneqq \coprod_i T_i\in \Ob(\Afd_K/T)$. Then $T'\rightarrow T$ is \'etale surjective, and thus $\Isom(x,y)|_{T'}=\Isom(x,y)\times_TT'\rightarrow \Isom(x,y)$ is a morphism of sheaves on $(\Afd_K/T)_\Et$ which is representable by adic spaces and \'etale surjective. Hence $\Isom(x,y)$ is an \'etale space over $T$ by Proposition~\ref{prop:bootstrap for etale spaces (sheaf admitting etale surjective map from etale sapce)}. 

Next we show that $\calS_{U}\rightarrow[U/R]$ is smooth and surjective. Take any $T\in\Ob(\V_{K}/S)$ and any $1$-morphism $f\colon(\Afd_{K}/T)_{\Et}\rightarrow[U/R]$ over $(\Afd_K/S)_\Et$. Then $(\Afd_{K}/T)_{\Et}\times_{f,[U/R]}\calS_{U}\rightarrow(\Afd_{K}/T)_{\Et}$ comes from a morphism in $\EtSp_{T}$, and we need to show that the resulting morphism is smooth and surjective.  Again, as above, there is an \'etale covering $\lbrace T_{i}\rightarrow T\rbrace$ by affinoids such that the composite $(\Afd_{K}/T_{i})_{\Et}\rightarrow(\Afd_K/T)_{\Et}\rightarrow[U/R]$
comes from a morphism $T_{i}\rightarrow U$ of sheaves on $(\Afd_K/S)_\Et$. As both smoothness and surjectivity are \'etale local on the source-and-target, we may assume that $f$ itself comes from a morphism $\tilde{f}\colon T\rightarrow U$ of sheaves. Then we see
\begin{align*}
(\Afd_{K}/T)_{\Et}\times_{f,[U/R]}\calS_{U}
&\cong(\Afd_{K}/T)_{\Et}\times_{\tilde{f},\calS_U}(\calS_{U}\times_{[U/R]}\calS_{U})\\
&\cong(\Afd_{K}/T)_{\Et}\times_{\tilde{f},\calS_U,t}\calS_{R}
\cong \calS_{T\times_{\tilde{f},U,t}R},
\end{align*}
where the second equivalence follows from Lemma~\ref{lem:2-coequalizer property of quotient stack}(i).
Since $t$ is smooth and surjective, so is the base change $T\times_{\tilde{f},U,t}R\rightarrow T$.
\end{proof}

\smallskip
\noindent
\textbf{Properties of morphisms of adic stacks}.
Let us define the properties of morphisms of adic stacks that is not necessarily smooth local on the source-and-target: monomorphisms, epimorphisms, and open immersions.

\begin{defn}\label{def:1-morphism being mono or epi}
Let $f\colon\calX\rightarrow\calY$ be a $1$-morphism of adic stacks over $S$. 
We say that $f$ is a \emph{monomorphism} if the underlying functor $\calX\rightarrow \calY$ is fully faithful.
We say that $f$ is an \emph{epimorphism} if, for every $T\in\Ob(\Afd_{K}/S)$ and $y\in\Ob(\calY_{T})$, there exist an \'etale covering $\lbrace T_{i}\rightarrow T\rbrace$ and $x_{i}\in\Ob(\calX_{T_{i}})$ for each $i$ such that $f(x_{i})$ and $y\rvert_{T_{i}}$ are isomorphic in $\calY_{T_{i}}$.
\end{defn}

\begin{defn}
A $1$-morphism $\calU\rightarrow \calX$ of adic stacks over $S$ is said to an \emph{open immersion} if it is representable by adic spaces and is an open immersion.\footnote{Recall that a morphism of \'etale spaces over $S$ is said to be an open immersion if it is representable by adic spaces and is an open immersion.}

An \emph{open substack} of an adic stack $\calX$ over $S$ is a strictly full subcategory $\calW\subset\calX$ such that $\calW$ is also an adic stack and the $1$-morphism $\calW\rightarrow\calX$ is an open immersion.
\end{defn}

\begin{rem}\label{rem:open immersion and open substack}
Every open immersion $\calU\rightarrow \calX$ of adic stacks over $S$ uniquely factors as $\calU\rightarrow \calW\subset \calX$ where $\calW$ is an open substack of $\calX$ and $\calU\rightarrow \calW$ is an epimorphism: in fact, $\calW$ is the full subcategory of $\calX$ consisting of objects that are isomorphic to $f(u)$ for some $u\in \Ob(\calU)$. More generally, if a $1$-morphism $f\colon \calX'\rightarrow \calX$ of adic stacks over $S$ admits a factorization $\calX'\rightarrow\calW\subset \calX$ such that $\calW$ is an open substack of $\calX$ and $\calX'\rightarrow\calW$ is an epimorphism, then such an open substack $\calW$ is uniquely determined by $f$.
\end{rem}

\begin{lem}\label{lem:epi-open immersion decompositon}
Let $(U',R',s',t',c')\rightarrow (U,R,s,t,c)$ be a morphism of smooth groupoids in \'etale spaces over $S$ such that  $U'\rightarrow U$ is an open immersion, and set $W\coloneqq s(t^{-1}(U'))\subset U$. Then $W\hookrightarrow U$ is an open immersion and $s^{-1}(W)=t^{-1}(W)$ as subsheaves of $R$, which we denote by $R_W$. Moreover, $R_W\rightrightarrows  W$ gives a smooth groupoid in \'etale spaces over $S$, the induced $1$-morphism $[W/R_W]\rightarrow [U/R]$ is an open immersion, and the $1$-morphism $[U'/R']\rightarrow [U/R]$ factors as $[U'/R']\rightarrow [W/R_W]\rightarrow [U/R]$ with the first $1$-morphism being an epimorphism.
\end{lem}

\begin{proof}
We know from Lemma~\ref{lem:image of smooth morphism of etale spaces} and the inverse structure of a groupoid (cf.~\cite[0230]{stacks-project}) that $W\hookrightarrow U$ is an open immersion and $s^{-1}(W)=t^{-1}(W)$. It follows that $R_W\rightrightarrows  W$ becomes a smooth groupoid in \'etale spaces over $S$ by restricting $(U,R,s,t,c)$ along $W\hookrightarrow U$. 
Moreover, the proof of \cite[04ZN]{stacks-project} also works in our setting, showing that the following $2$-commutative diagram is a $2$-fiber product:
\[
\xymatrix{
 \calS_W \ar[r]\ar[d] & [W/R_W]\ar[d]\\
 \calS_U \ar[r] & [U/R].
}
\]
Since $U\rightarrow [U/R]$ is smooth and surjective, we deduce from Proposition~\ref{prop:smooth-local-on-target} (together with Theorem~\ref{thm:effective-etale-equivalence-relation} and \cite[Prop.~3.1.14]{Warner}) that $[W/R_W]\rightarrow [U/R]$ is an open immersion.

By construction, we have a factorization $[U'/R']\rightarrow [W/R_W]\rightarrow [U/R]$. It remains to show that the first $1$-morphism is an epimorphism. For this, take any $T\in\Ob(\Afd_K/S)$. 
By Lemma~\ref{lem:explicit description of quotient stack}, an object of $[W/R_W]_T$ is given by a $[W/R_W]$-descent datum $(w_i,r_{ij})$ (where $w_i\in W(T_i)$ and $r_{ij}\in R_W(T_i\times_TT_j)$) relative to some covering $\calT=\{T_i\rightarrow T\}$. Since $s\colon t^{-1}(U')\rightarrow W$ is an epimorphism of sheaves by Remark~\ref{rem:smooth surjective morphism of etale spaces is epimorphism}, there exists an \'etale surjection $f_i\colon T_i'\rightarrow T_i$ of affinoids for each $i$ such that the composite $T_i'\xrightarrow{f_{i}} T_i\xrightarrow{w_{i}} W$ factors through $t^{-1}(U')$ as in the left square in the following commutative diagram:
\[
\xymatrix{
 && U' \ar@{^{(}->}[r] & W \\
T_i'\ar[d]_-{f_i}\ar@{..>}[r] \ar@/^.6pc/@{..>}[urr]^-{u_i'} \ar@/^-1.0pc/@{..>}[rr]_(.3){r_i'} & t^{-1}(U')\ar@{^{(}->}[r]\ar[d]\ar[ru] & R_W \ar[d]^-s\ar[ru]^-t & \\
T_i\ar[r]^-{w_i} & W\ar@{=}[r] & W. &
}
\]
Define $u_i'\in U'(T_i')$ and $r_i'\in R_W(T_i')$ so that the remaining part of the diagram becomes commutative. We claim that the compositions of $r_i'$s and $r_{ij}'\coloneqq r_{ij}\circ (f_i\times f_j)$s yield elements $\tilde{r}_{ij}'\in R_W(T_i'\times_TT_j')$, making the system $(u_i',\tilde{r}_{ij}')$ a $[W/R_W]$-descent datum relative to $\calT'=\{T_i'\rightarrow T_i\rightarrow T\}$:
more precisely, set $w_i'\coloneqq w_i\circ f_i\colon T_i'\rightarrow W$ and let $\iota\colon R_W\rightarrow R_W$ (resp.~$c\colon R_W\times_{s,W,t}R_W\rightarrow R_W$) denote the inverse (resp.~the composition). We define
\[
\tilde{r}_{ij}\coloneqq c\circ (r_{ij}',\iota\circ r_i'\circ \pr_1), \quad
\tilde{r}_{ij}'\coloneqq c\circ (r_j'\circ \pr_2,\tilde{r}_{ij})\colon T_i'\times_TT_j'\rightarrow R_W.
\]
For example, $\tilde{r}_{ij}$ sits in the following commutative diagram:
\[
\xymatrix{
&&T_i'\times_TT_j'\ar[d]^-{(r_{ij}',\iota\circ r_i'\circ \pr_1)}\ar@/^-2.0pc/[ddll]_-{\pr_1}\ar@/^2.0pc/[ddrr]^-{\pr_2} \ar@/^-3.0pc/@{..>}[dd]_-{\tilde{r}_{ij}} &&\\
&&R_W\times_{s,W,t}R_W\ar[dl]_-{\pr_2}\ar[d]^-c\ar[dr]^-{\pr_1} &&\\
T_i'\ar[d]_-{u_i'}& R_W\ar[d]_-s& R_W\ar[dl]_-s\ar[dr]^-t& R_W\ar[d]^-t& T_j'\ar[dl]^-{w_j'} \\
U'\ar@{^{(}->}[r]&W & & W. &
}
\]
Then we can check
\[
\tilde{r}_{ij}'\circ s=u_i'\circ \pr_1,\quad \tilde{r}_{ij}'\circ t=u_j'\circ \pr_2,\quad\text{and}\quad
c\circ (\tilde{r}_{jk}'\circ \pr_{23},  \tilde{r}_{ij}'\circ \pr_{12})=\tilde{r}_{ij}'\circ \pr_{13},
\]
and obtain the desired $[W/R_W]$-descent datum $(u_i',\tilde{r}_{ij}')$.
Moreover, $(r_i')$ defines a morphism from $(w_i,r_{ij})|_{\calT'}=(w_i',r_{ij}')$ to $(u_i',\tilde{r}_{ij}')$ by construction. Since  $u_i'\in U'(T_i')$, we conclude that $[U'/R']\rightarrow [W/R_W]$ is an epimorphism.
\end{proof}

\begin{cor}\label{cor:epi-open immersion factorization}
Let $f\colon \calX'\rightarrow \calX$ be a $1$-morphism of adic stacks over $S$ and suppose that there exists a morphism $(U',R',s',t',c')\rightarrow (U,R,s,t,c)$ of smooth groupoids  that gives a presentation of $f$ as in the following $2$-commutative diagram:
\[
\xymatrix{
[U'/R']\ar[r]\ar[d]_-\cong & [U/R] \ar[d]^-\cong\\
\calX'\ar[r]^-f & \calX.\\ 
}
\]
Assume further that $U'\rightarrow U$ an open immersion.
Then there exists a unique open substack $\calW$ of $\calX$ such that $f$ factors as $\calX'\rightarrow \calW\subset \calX$ with $\calX'\rightarrow \calW$ being an epimorphism. 
In fact, $\calW$ is the full subcategory of $\calX$ consisting of objects that are isomorphic to $f(x')$ for some $x'\in \Ob(\calX')$. In particular, the factorization is independent of choice of such a presentation of $f$.
\end{cor}

\begin{proof}
With the notation as in Lemma~\ref{lem:epi-open immersion decompositon}, define $\calW\subset \calX$ to be the open substack given by the open immersion $[W/R_W]\rightarrow [U/R]\rightarrow \calX$ as in Remark~\ref{rem:open immersion and open substack}.
It is straightforward to see that $f$ factors as $\calX'\rightarrow \calW\subset \calX$ with $\calX'\rightarrow \calW$ being an epimorphism. The remain assertions also follow easily.
\end{proof}

\subsection{Formal deformation theory}

In this subsection, we define the notions of tangent spaces and infinitesimal automorphisms of adic stacks. 
Fix $S\in\Ob(\V_{K})$.

\begin{defn}
 A morphism $X\rightarrow Y$ in $\EtSp_{S}$ is called a \emph{first-order thickening} if it is representable by adic spaces and is a closed immersion defined by a square-zero coherent ideal sheaf. A first-order thickening is said to be \emph{split} if it admits a retraction (called a \emph{splitting}) that is \emph{representable by adic spaces}. 
\end{defn}

It is clear that (split) first-order thickenings are stable under base change.

\begin{lem}\label{lem:pushout etale space}
Let $X\rightarrow Y$ and $X\rightarrow Z$ be first-order thickenings in $\EtSp_{S}$ with $X\rightarrow Z$ split. Then, there exists a pushout $Y\amalg_{X}Z$ in $\EtSp_{S}$.
\end{lem}

\begin{proof}
We follow \cite[07VX]{stacks-project}. Take an \'etale cover $h_{V}\rightarrow Y$ with $V\in\Ob(\V_{K}/S)$. Then, $h_{V}\times_{Y}X$ is representable by an adic space, say, $U$. Fix a splitting $Z\rightarrow X$. The fiber product $h_{U}\times_{X}Z$ along the splitting is representable by an adic space, say, $V'$. By construction, the projection $h_{V'}\rightarrow Z$ is an \'etale cover, and $h_{U}=h_{V'}\times_{Z}X$. These constructions give two first-order thickenings $U\rightarrow V$, $U\rightarrow V'$ in $\V_{K}/S$. By \cite[Prop~2.10.4]{Warner}, there exists a pushout $W\coloneqq V\amalg_{U}V'$ in $\V_{K}/S$.

Let $R_{U}\coloneqq U\times_{X}U$, $R_{V}\coloneqq V\times_{Y}V$, and $R_{V'}\coloneqq V'\times_{Z}V'$. We have $R_{V}\times_{Y}X=R_{U}=R_{V'}\times_{Z}X$ and thus obtain first-order thickenings $R_{U}\rightarrow R_{V}$ and $R_{U}\rightarrow R_{V'}$. By \emph{loc.~cit.}, a pushout $R\coloneqq R_{V}\amalg_{R_{U}}R_{V'}$ exists in $\V_{K}/S$. Now the proof of \cite[07VX]{stacks-project} goes through with \cite[07RT]{stacks-project} replaced by \cite[Prop.~2.10.4]{Warner} (and the explicit description of the pushout in the first paragraph of the proof therein): one can obtain an etale equivalence relation $R\rightrightarrows  W$ and show that the quotient $R/W$ represents the desired pushout.
\end{proof}

Let $(L,L^{+})$ be a complete analytic affinoid field over $S$. For a finite-dimensional $L$-vector space $V$, let $L[V]$ be the $L$-algebra whose underlying $L$-module is $L\oplus V$ and whose multiplication is given by $(a,v)(a',v')=(aa',av'+a'v)$ for $a,a'\in L$ and $v,v'\in V$. The submodule $L[V]^{+}\coloneqq L^{+}\oplus V$ agrees with the integral closure of $L^{+}$ in $L[V]$, and the pair $(L[V],L[V]^{+})$ is a complete affinoid ring of topologically finite type over $(L,L^+)$ whose formation is functorial in $V$. Set $q_V\colon D_{V}\coloneqq \Spa(L[V],L[V]^{+})\rightarrow \Spa(L,L^+)$ and simply write $D$ and $q$ if $V=L$. The $L$-algebra map $L[V]\rightarrow L$ sending $(a,v)$ to $a$ yields a split first-order thickening $\Spa(L,L^+)\hookrightarrow D_V$.

\begin{lem}\label{lem:pushout square of D}
Let $V_1$ and $V_2$ be finite-dimensional $L$-vector spaces. The commutative diagram 
\begin{equation}\label{eq:pushout-square}
\xymatrix{
\Spa(L,L^{+})\ar[r]\ar[d]&D_{V_{1}}\ar[d]\\
D_{V_{2}}\ar[r]&D_{V_{1}\times V_{2}},
}    
\end{equation}
is a pushout square in $\EtSp_{S}$ and induces an equivalence of categories 
\begin{equation}\label{eq:equivalence-categories}
\cX_{D_{V_{1}\times V_{2}}}\xrightarrow{\cong}\cX_{D_{V_{1}}}\times_{\cX_{\Spa(L,L^{+})}}\cX_{D_{V_{2}}}.
\end{equation}
\end{lem}

\begin{proof}
Note $L[V_{1}\times V_{2}]=L[V_{1}]\times_{L}L[V_{2}]$ and $L[V_{1}\times V_{2}]^{+}=L[V_{1}]^{+}\times_{L^{+}}L[V_{2}]^{+}$. Hence \eqref{eq:pushout-square} is a pushout square in $\Afd_{K}/S$ by \cite[Prop~2.10.4, Pf.]{Warner}. This implies that $h_{D_{V_{1}\times V_{2}}}$ agrees with the pushout $h_{D_{V_{1}}}\amalg_{h_{\Spa(L,L^{+})}}h_{D_{V_{2}}}$ in $\Sh((\Afd_{K}/S)_{\Et})$, so \eqref{eq:pushout-square} is necessarily a pushout square in $\EtSp_{S}$.

For the second assertion, we follow \cite[07WN]{stacks-project}. 
Since $\calX$ is an adic stack, the diagonal $\calX\rightarrow \calX\times\calX$ is representable by \'etale spaces. Hence the full faithfulness of \eqref{eq:equivalence-categories} follows from the first assertion and Lemma~\ref{lem:DiagonalRepresentable-Stacks}.
To show the essential surjectivity, take any $y_{1}\in\Ob(\cX_{D_{V_{1}}})$, $y_{2}\in\Ob(\cX_{D_{V_{2}}})$ and an isomorphism $f\colon y_{1}\rvert_{\Spa(L,L^{+})}\xrar{\sim} y_{2}\rvert_{\Spa(L,L^{+})}$. We need to construct an object $y\in\Ob(\cX_{D_{V_{1}\times V_{2}}})$ such that $y\rvert_{D_{V_{1}}}\cong y_{1}$ and $y\rvert_{D_{V_{2}}}\cong y_{2}$ and such that $f$ is given by the composition $y_{1}\rvert_{\Spa(L,L^{+})}\cong y\rvert_{\Spa(L,L^{+})}\cong y_{2}\rvert_{\Spa(L,L^{+})}$. In fact, we may construct such an $y$ smooth locally on $D_{V_1\times V_2}$ thanks to Lemma~\ref{lem:criteria for equivalence of stacks} and \cite[Prop.~3.1.11]{Warner}.
Take a smooth cover $\cS_{h_{U}}\rightarrow\cX$ with $U\in\Ob(\V_{K}/S)$. Set $y_{0}\coloneqq y_{1}\rvert_{\Spa(L,L^{+})}$.
The $2$-fiber products $\cS_{h_{U}}\times_{\cX,y_{0}}\cS_{h_{\Spa(L,L^{+})}}$, $\cS_{h_{U}}\times_{\cX,y_{1}}\cS_{h_{D_{V_{1}}}}$, $\cS_{h_{U}}\times_{\cX,y_{2}}\cS_{h_{D_{V_{2}}}}$
are representable by \'etale spaces, which we denote by $U_{0},U_{1},U_{2}$, respectively. As $\Spa(L,L^{+})\rightarrow D_{V_{i}}$ is a split first-order thickening, so is $U_{0}\rightarrow U_{i}$. Therefore, a pushout $U_{1}\amalg_{U_{0}}U_{2}$ exists in $\EtSp_{S}$ by Lemma~\ref{lem:pushout etale space}, from which the desired $y$ can be constructed as in the proof of \cite[07WN]{stacks-project}.
\end{proof}

\begin{defn}\label{Tangent-Spaces-Defn}
Let $\calX\rightarrow(\Afd_{K}/S)_{\Et}$ be an adic stack over $S$ and let $x_{0}\in\Ob(\calX_{\Spa(L,L^{+})})$. The \emph{tangent groupoid} $\fT_{\calX,x_{0}}$ of $\calX$ at $x_{0}$ is defined to be the subgroupoid of $\calX_{D}$ consisting of objects $x$ satisfying $x\rvert_{\Spa(L,L^{+})}\cong x_{0}$. The set $T_{\calX,x_{0}}$ of isomorphism classes of $\fT_{\calX,x_{0}}$ is called the \emph{tangent space}. The space of \emph{infinitesimal automorphisms} is defined as 
$\Inf_{\calX,x_{0}}\coloneqq\Ker(\Aut_{\calX_{D}}(x_{0}')\rightarrow\Aut_{\calX_{\Spa(L,L^{+})}}(x_{0}))$, where $x_0'$ denotes the \emph{trivial deformation} $q^\ast x_0\in\Ob(\calX_D)$.
\end{defn}

\begin{lem}
The set $T_{\calX,x_{0}}$ and the group $\Inf_{\calX,x_{0}}$ have the natural structures of $L$-vector spaces. Moreover, every $1$-morphism $f\colon\cX\rightarrow\cY$ of adic stacks induces $L$-linear maps $df_{x_0}\colon T_{\cX,x_{0}}\rightarrow T_{\cY,f(x_{0})}$ and $df_{x_0}\colon\Inf_{\cX,x_{0}}\rightarrow \Inf_{\cY,f(x_{0})}$.
\end{lem}

\begin{proof}
Let $\mathbf{Vect}_{L}^{\fin}$ be the category of finite-dimensional $L$-vector spaces. Define the tangent space functor $\calT_{\cX,x_{0}}\colon\mathbf{Vect}_{L}^{\fin}\rightarrow\Set$ as
\[
\calT_{\cX,x_{0}}(V)\coloneqq \lbrace\text{isomorphism classes of $x_{V}\in\Ob(\cX_{D_{V}})$ such that $x_{V}\rvert_{\Spa(L,L^{+})}\cong x_{0}$}\rbrace.
\]
It follows from Lemma~\ref{lem:pushout square of D} that for $V_{1}, V_{2}\in \mathbf{Vect}_L^\fin$, the natural map
\[
\cT_{\cX,x_{0}}(V_{1}\times V_{2})\rightarrow \cT_{\cX,x_{0}}(V_{1})\times \cT_{\cX,x_{0}}(V_{2})
\]
is an isomorphism. Now \cite[06IA]{stacks-project} gives the $L$-vector space structure on $T_{\cX,x_{0}}=\cT_{\cX,x_{0}}(L)$. Similarly, consider a functor $\cI_{\cX,x_{0}}\colon\mathbf{Vect}_{L}^{\fin}\rightarrow\Set$ defined by
\[
\cI_{\cX,x_{0}}(V)\coloneqq\Ker(\Aut_{\cX_{D_{V}}}(q_V^\ast x_{0})\rightarrow\Aut_{\cX_{\Spa(L,L^{+})}}(x_{0})).
\]
Again Lemma~\ref{lem:pushout square of D} and \cite[06IA]{stacks-project} give the $L$-vector space structure on $\Inf_{\cX,x_{0}}=\cI_{\cX,x_{0}}(L)$.
Finally, the second assertion follows from \cite[06ID, Pf.~of 07W6]{stacks-project}.
\end{proof} 

\begin{lem}\label{lem:T exact sequence}
Let $f_1\colon\cX_1\rightarrow\cY$ and $f_2\colon\cX_2\rightarrow\cY$ be $1$-morphisms of adic stacks over $S$. Set $\cZ\coloneqq\cX_1\times_{\cY}\cX_2$ and let $\pr_{i}$ denote the projection from $\cZ$ to $\cX_i$. For every $z_{0}\in\Ob(\cZ_{\Spa(L,L^{+})})$, we have an exact sequence of $L$-vector spaces
\[
0\rightarrow\Inf_{\cZ,z_{0}}\rightarrow\Inf_{\cX_1,x_{1,0}}\oplus\Inf_{\cX_2,x_{2,0}}\rightarrow\Inf_{\cY,y_{0}}\rightarrow T_{\cZ,z_{0}}\rightarrow T_{\cX_1,x_{1,0}}\oplus T_{\cX_2,x_{2,0}}\rightarrow T_{\cY,y_{0}},
\]
where $x_{i,0}\coloneqq \pr_{i}(z_{0})$ and $y_{0}\coloneqq f_1(x_{1,0})=f_2(x_{2,0})$.
\end{lem}

\begin{proof}
The construction and exactness of the above sequence formally follow from the fact that $\cY$ is the $2$-fiber product of categories: see \cite[06L5, Pf.]{stacks-project}.
\end{proof}

\begin{lem}\label{lem:tangent space of etale spaces}
Let $L$ be a finite extension over $K$ and write $\Spa L=\Spa(L,\calO_L)$.
\begin{enumerate}
 \item Let $U\in\Ob(\Rig_K)$. Take a morphism $x_0\colon \Spa L\rightarrow U$ and write $u\in U$ for the image of $x_0$. Then there are canonical isomorphisms
\[
T_{\calS_U,x_0}\cong \Hom_{\calO_{U,u}}((\Omega_{U/K})_u, L)\cong \Hom_{L}((\Omega_{U/K})_u\otimes_{\calO_{U,u}}L,L)
\]
of $L$-vector spaces. Moreover, if $f\colon U\rightarrow U'$ is a morphism in $\Rig_K$, the induced map $T_{\calS_U,x_0}\rightarrow T_{\calS_{U'},f\circ x_0}$ corresponds via the above identification to the $L$-linear dual of $f^\ast\colon (\Omega_{U'/K})_{f(u)}\otimes_{\calO_{U',f(u)}}L\rightarrow (\Omega_{U/K})_u\otimes_{\calO_{U,u}}L$.
 \item Let $X\in\Ob(\EtSp_K)$ and take a presentation $U/R\cong X$. Assume that a morphism $x_0\colon \Spa L\rightarrow X$ lifts to $\tilde{x}_0\colon \Spa L\rightarrow U$. Then the induced $L$-linear map $T_{\calS_U,\tilde{x}_0}\rightarrow T_{\calS_X,x_0}$ is an isomorphism. 
\end{enumerate}
\end{lem}

\begin{proof}
(i) We follow \cite[0B2D]{stacks-project}. 
Write $D=\Spa(L\oplus L\epsilon,(L\oplus L\epsilon)^\circ)$ and identify $\Gamma(D,\Omega_{D/L})=L\,d\epsilon\cong L$.
By Remark~\ref{rem:category fibered in setoids associated to presheaf}, every element of $T_{\calS_U,x_0}$ corresponds to a morphism $x\colon D\rightarrow U$ such that $x|_{\Spa(L,\calO_L)}=x_0$. Then the composite $x^\ast\Omega_{U/K}\rightarrow \Omega_{D/K}\rightarrow \Omega_{D/L}$ yields an $\calO_{U,u}$-linear map $\xi_x\colon (\Omega_{U/K})_u
\rightarrow \Gamma(D,\Omega_{D/L})\cong L$. Conversely, take an $\calO_{U,u}$-linear map $\xi\colon (\Omega_{U/K})_u
\rightarrow L$. There exists an affinoid open neighborhood $\Spa(A,A^\circ)$ of $u\in U$ such that $\xi$ is induced from an $R$-linear map $\xi_A\colon \Omega_{A/K}\rightarrow L$. Note that $\xi_A$ corresponds to a continuous $K$-derivation $\delta\colon A\rightarrow L$. Let $x_0^\ast\colon A\rightarrow k(u)\rightarrow L$ denote the map corresponding to $x_0$ and define a $K$-linear map $\alpha\colon A\rightarrow L\oplus L\epsilon$ by $\alpha(a)=x_0^\ast (a)+\delta(a)\epsilon$. It is immediate to see that $\alpha$ is a $K$-algebra map and defines a morphism $x_\xi\colon D\rightarrow U$ such that $x_\xi|_{\Spa(L,\calO_L)}=x_0$. One can also check that $x_\xi$ is independent of the choice of $A$ and that the maps $x\mapsto \xi_x$ and $\xi\mapsto x_\xi$ define $L$-linear isomorphism $T_{\calS_U,x_0}\cong \Hom_{\calO_{U,u}}((\Omega_{U/K})_u, L)$. The remaining assertions are easy to verify.

(ii) This is deduced from the unique existence of a lift of $x\colon D\rightarrow X$ in the following commutative diagram
\[
\xymatrix{
U\ar[r] & X \\
\Spa L \ar[u]^-{\tilde{x}_0}\ar[r]& D \ar[u]_-x\ar@{.>}[ul]_{\exists  ! \tilde{x}},
}
\]
which follows from the definition of \'etale morphisms of adic spaces and Lemma~\ref{lem:DiagonalRepresentable}.
\end{proof}

\section{Analytification and generic fiber}
\label{sec:analytification and generic fiber}

In this section, we provide a detailed exposition of the analytification of algebraic stacks (Proposition~\ref{prop:analytification of algebraic stacks}) and the generic fiber of formal algebraic stacks (Proposition~\ref{prop:generic fiber of formal algebraic stacks}). In particular, given an algebraic stack $\calY$ locally of finite type over $\Z_p$, we study the $1$-morphism from the generic fiber $\widehat{\calY}_\eta$ of $p$-adic completion of $\calY$ to the analytification $\calY_{\Q_p}^\an$ of generic fiber of $\calY$ (Proposition~\ref{prop:image of comparison map from generic fiber to analytification}).

Let $\calO_K$ be a complete discrete valuation ring of mixed characteristic $(0,p)$ with fraction field $K$ and a uniformizer $\pi$.
There are two major constructions to obtain rigid $K$-spaces: the \emph{analytification} of a scheme of finite type over $K$ and the \emph{generic fiber} of a formal scheme of finite type over $\calO_K$. Moreover, given a scheme $Y$ of finite type $\calO_K$, there is a morphism of rigid $K$-spaces $\widehat{Y}_\eta\rightarrow Y_K^\an$ from the generic fiber of the $p$-adic completion $\widehat{Y}$ to the analytification of the generic fiber $Y_K$. The goal of this section is to extend these constructions to (formal) algebraic stacks (Definitions~\ref{def:analytification of algebraic stack}, \ref{def:generic fiber of formal algebraic stack}, and \ref{defn:comparison map from generic fiber to analytification}).

\subsection{Analytification}\label{sec:analytification}
For an adic space $T$, we also use $T$ to denote its underlying locally ringed space.

\begin{prop-defn}\label{prop-defn:analytificaiton}
Let $Y\rightarrow S$ be a morphism of schemes that is locally of finite type and let $T\rightarrow S$ be a morphism  of locally ringed spaces from a locally strongly Noetherian analytic adic space $T$. Then there exist an adic space $Z$ over $T$ and a morphism $Z\rightarrow Y$ of locally ringed spaces such that the diagram
\[
\xymatrix{
Z \ar[r]\ar[d] & T\ar[d]\\
Y\ar[r] & S
}
\]
of locally ringed spaces
commutes and it is universal among such diagrams.
Moreover, $Z\rightarrow T$ is locally of finite type, and the construction is functorial in $Y$ and $T$.
We call it \emph{the fiber product of $Y$ and $T$ over $S$} and denote it by $Y\times_ST$.
\end{prop-defn}

\begin{proof}
 This is a special case of \cite[Prop.~3.8]{Huber-gen}.
\end{proof}

\begin{lem}\label{lem:first properties of analytification of schemes}
Keep the notation as in Proposition-Definition~\ref{prop-defn:analytificaiton}.
\begin{enumerate}
\item The functor $Y\mapsto Y\times_ST$ commutes with fiber products. Moreover, it preserves the following properties of morphisms: open immersions; surjections; \'etale morphisms; smooth morphisms.
\item If $T'\rightarrow T$ is a morphism of locally strongly Noetherian analytic adic spaces, then $(Y\times_ST)\times_TT'\cong Y\times_ST'$. 
\item  Assume $T=\Spa(A,A^+)$ with $A$ complete and $S=\Spec A$. Regard $Y\times_ST$ as a set-valued presheaf on $\AdSp_T^\aff$ (which is a sheaf with respect to analytic topology).
If $Y=\Spec R$ is affine, then $Y\times_ST$ is given by
\[
 \Spa(B,B^+)\mapsto Y(\Spec \widehat{B})=\Hom_{\mathrm{Alg}_{A}}(R,\widehat{B}).
\]
For a general $Y$, the adic space $Y\times_ST$ as a presheaf on $\AdSp_T^\aff$ agrees with the sheafification with respect to analytic topology of the presheaf
\[
 \Spa(B,B^+)\mapsto Y(\Spec \widehat{B}).
\]
\end{enumerate}
\end{lem}

\begin{proof}
For (i), the assertions for fiber products and open immersions follow easily from the universality of the fiber product. The statement about surjections and \'etale morphisms follow from \cite[Lem.~3.9(ii)]{Huber-gen} and \cite[Cor.~1.7.3(i)]{Huber-etale}, respectively. Finally, the case of smooth morphisms follows from the cases of open immersions and \'etale morphisms since $\mathbb{A}^n_S\times_ST$ is smooth over $T$ (cf.~\cite[Prop.~3.8, Pf.]{Huber-gen}).
Part (ii) also follows from the universality, and part (iii) is deduced from the fact that $\Spec$ is right adjoint to the global sections functor on the category of locally ringed spaces \cite[01I1]{stacks-project}.
\end{proof}

\begin{defn}
Let $\Sch_K^\lft$ denote the category of schemes locally of finite type over $K$. 
The \emph{analytification functor} 
\[
(\cdot)^\an\colon \Sch_K^\lft\rightarrow \Rig_K
\]
is defined as $Y\mapsto Y^\an\coloneqq Y\times_{\Spec K}\Spa(K,\calO_K)$.

More generally, for a locally strongly Noetherian adic space $T=\Spa(A,A^+)$ over $K$ and an $A$-scheme $Y$ locally of finite type, we often write $Y^\an$ for the adic space $Y\times_{\Spec A}T$ locally of finite type over $T$.
\end{defn}

\begin{rem}
One also associates to a scheme locally of finite type over $K$ a rigid analytic variety as in \cite[9.3.4]{BGR}. See \cite[Rem.~4.6(i)]{Huber-gen} for the comparison of these two analytification functors.
\end{rem}

\begin{lem}\label{lem:analytification and base change}
For $Y\in \Sch_K^\lft$ and $T=\Spa(A,A^+)$ a locally strongly Noetherian adic space over $K$, there is a natural isomorphism of adic spaces over $T$:
\[
Y^\an\times_{\Spa(K,\calO_K)}T\cong (Y\times_{\Spec K}\Spec A)\times_{\Spec A}T.
\]
\end{lem}

\begin{proof}
This follows easily from the universality of the fiber product.
\end{proof}

\begin{prop}\label{prop:comparison of stalks for analytification}
Let $T=\Spa(A,A^+)\in \Afd_K$ with $A$ complete and set $S=\Spec A$. For every scheme $Y$ locally of finite type over $S$, the structure morphism $q\colon Y^\an\rightarrow Y$ of locally ringed spaces is surjective. Moreover, for $y\in Y^\an$, the induced map $\calO_{Y,q(y)}\rightarrow \calO_{Y^\an, y}$ of local rings is flat.
\end{prop}

\begin{proof}
By definition of $\Afd_K$, there exists a complete analytic affinoid field $(L,L^+)$ over which $T$ is of finite type. We use results from Berkovich spaces over $L$.

By \cite[Lem.~3.9(ii)]{Huber-gen}, the surjectivity of $q$ is reduced to that of $T\rightarrow \Spec A$. The latter map factors as $T\rightarrow \calM(A)\rightarrow \Spec A$ where $\calM(A)$ denotes the Berkovich spectrum of $A$ (cf.~\cite[Def.~2.4.6]{Kedlaya-Liu-I}). Since the second map is surjective by \cite[Prop.~2.1.1]{Berkovich-etale}, the desired surjectivity follows.

The assertion on flatness is local on $Y$, so we may assume that $Y=\Spec B$ is affine of the form $B=A[x_1,\ldots,x_n]/I$ and write $Y^\an=\bigcup_{l\geq 0}\Spa(C_l,C_l^+)$ where 
$C_l\coloneqq A\langle p^lx_1,\ldots,p^lx_n\rangle/IA\langle p^lx_1,\ldots,p^lx_n\rangle$ and $C_l^+$ is the integral closure of the image of $A^+$ in $C_l$ (see \cite[Prop.~3.8, Pf.]{Huber-gen}). Fix $l$ and write $(C,C^+)\coloneqq (C_l,C_l^+)$. Take any $y\in \Spa(C,C^+)$ and set $\fkn\coloneqq \operatorname{supp}(y)\in \Spec C$ and $\fkm\coloneqq q(y)\in \Spec B$. Then the map in question is written as
\[
\calO_{Y,q(y)}=B_\fkm\rightarrow C_\fkn\rightarrow \calO_{Y^\an,y}.
\]
The second map is faithfully flat by \cite[Lem.~B.4.2]{Zavyalov-quotients}. Recall that $\Spa(C,C^+)\rightarrow \calM(C)$ has a canonical section (cf.~\cite[Def.~2.4.6]{Kedlaya-Liu-I}).
If $y$ is in $\calM(C)\subset \Spa(C,C^+)$, then the flatness of $\calO_{Y,q(y)}\rightarrow \calO_{Y^\an, y}$ is proved in \cite[Prop.~2.6.2]{Berkovich-etale}. For a general $y$, choose $y'\in \calM(C)$ with $\operatorname{supp}(y')=\fkn$. Since $B_\fkm\rightarrow \calO_{Y^\an,y'}$ is flat and $C_\fkn\rightarrow \calO_{Y^\an,y'}$ is faithfully flat, we deduce that $B_\fkm\rightarrow C_\fkn$ is flat, whence $B_\fkm\rightarrow C_\fkn\rightarrow \calO_{Y^\an,y}$ is also flat.
\end{proof}

\smallskip
\noindent
\textbf{Notation}.
Let $Y$ be a scheme locally of finite type over $K$.
For a locally strongly Noetherian adic space $T=\Spa(A,A^+)$  over $K$ with $A$ complete, set 
\[
Y_A\coloneqq Y\times_{\Spec K}\Spec A\quad\text{and}\quad Y_T^\an\coloneqq Y^\an\times_{\Spa(K,\calO_K)}T=Y_A\times_{\Spec A}T,
\]
and write $q_T\colon Y_T^\an\rightarrow Y_A$ for the structure morphism of locally ringed spaces. Since $Y_T^\an$ is locally strongly Noetherian, $\calO_{Y_T^\an}$ is coherent. Hence the pullback
\[
q_T^\ast\colon \Sh(Y_A,\calO_{Y_A})\rightarrow \Sh(Y_T^\an,\calO_{Y_T^\an}):\quad\calF\mapsto \calF^\an\coloneqq q_T^{-1}\calF\otimes_{q_T^{-1}\calO_{Y_A}}\calO_{Y_T^\an},
\]
restricts to $\Coh(\calO_{Y_A})\rightarrow\Coh(\calO_{Y_T^{\an}})$ by \cite[(0.5.3.11)]{EGAI}, 
which we call the \emph{analytification functor}.

\begin{prop}\label{prop:analytification and flatness}
With the above notation, $\calF\in\Coh(\calO_{Y_A})$ is $A$-flat (resp.~locally free) if and only if $\calF^\an$ is $T$-flat in the sense of Remark~\ref{rem:flatnesss of coherent sheaf and morphism for adic spaces} (resp.~locally free).
\end{prop}

\begin{proof}
The assertions are local on $Y_A$, so we assume that $Y_A=\Spec B$ is affine. Set $M=\Gamma(Y_A,\calF)$; it is a finitely generated $B$-module, and $\calF=M^\sim$.
Take any $y\in Y_T^\an$. Let $t\in T$ denote the image of $y$, and set $\fkm\coloneqq q_T(y)\in \Spec B$ and $\fkp\coloneqq \operatorname{supp}(t)\in\Spec A$ (which is below $\fkm$). Then $\calF_{q_T(y)}=M_\fkm$ and $\calF^\an_y=M\otimes_B\calO_{Y_T^\an,y}=M_\fkm\otimes_{B_\fkm}\calO_{Y_T^\an,y}$.
We know from Proposition~\ref{prop:comparison of stalks for analytification} that $B_\fkm\rightarrow \calO_{Y_T^\an,y}$ is faithfully flat. Hence $\calF_{y}^\an$ is flat over $A$ if and only if $\calF_{q_T(y)}$ is flat over $A$.

Assume that $\calF$ is flat over $A$. Then $\calF_{y}^\an$ is flat over $A$. By running the same argument for any rational localization $\Spa(A',A'^+)$ of $T$ containing $t$, we see that $\calF_{y}^\an$ is flat over any such $A'$. Since $\calO_{T,t}=\varinjlim A'$, we conclude that $\calF_{y}^\an$ is flat over $\calO_{T,t}$ by \cite[05UU]{stacks-project}. As this is true for any $y$, $\calF^\an$ is $T$-flat.

Conversely, assume that $\calF^\an$ is $T$-flat. Then for any $\fkm\in \Spec B$, one can find $y\in Y_T^\an$ with $\fkm=q_T(y)$ by Proposition~\ref{prop:comparison of stalks for analytification}. With the notation as above, since $A\rightarrow \calO_{T,t}$ is flat, we conclude that $\calF_{q_T(y)}$ is flat over $A$. Hence $\calF$ is $A$-flat.

The assertion for local freeness follows similarly, noting that $\calF_y^\an$ is finite projective over $\calO_{Y_T^\an,y}$ if and only if $\calF_{q_T(y)}$ is finite projective over $B_\fkm$.
\end{proof}

\begin{thm}\label{thm:rigid relative GAGA}
With the notation as above, assume that $Y$ is proper over $K$ and $T=\Spa(A,A^+)\in \Afd_K$.  
Then the analytification along $Y_T^\an\rightarrow Y_A$ defines an equivalence of categories 
\[
\Coh(Y_A)\xrightarrow{\cong}\Coh(Y_T^\an).
\]
Moreover, for every $\calF\in \Coh(Y_A)$ with the analytification $\calF^\an$, there is a canonical isomorphism of $A$-modules
\[
H^i(Z_A,\calF)\xrightarrow{\cong}H^i(Z_T^\an,\calF^\an).
\]
\end{thm}

\begin{proof}
Since $T\in\Afd_K$, there exists an analytic affinoid field $(L,L^+)$ over which $T$ is of finite type.
By Corollary~\ref{cor:coherent sheaf does not depend on ring of integral elements}, we may assume $L^+=\calO_L$ and $A^+=A^\circ$. Then the assertions follow from \cite{Kopf} (cf.~\cite[Thm.~9.1]{scholze-p-adic-hodge}).
\end{proof}

\smallskip
\noindent
\textbf{Analytification of algebraic spaces}.
Fix $S=\Spa(A,A^+)\in \Afd_K$ with $A$ complete.

\begin{defn}\label{def:analytification of algebraic space}
Let $F$ be an algebraic space locally of finite type over $\Spec A$. Define the sheaf $F^\an$ on $(\Afd_K)_\Et$ to be the sheafification of the presheaf
\[
\Ob(\Afd_K/S)\ni \Spa(B,B^+)\mapsto F(\Spec \widehat{B}).
\]
Proposition~\ref{prop:analytification of algebraic spaces} below shows that $F^\an$ is an \'etale space over $S$. We call $F^\an$ the \emph{analytification} of $F$.
\end{defn}

\begin{prop}\label{prop:analytification of algebraic spaces}
Let $F$ be an algebraic space locally of finite type over $\Spec A$ and take a presentation of $F$ given by an \'etale equivalence relation $s,t\colon R\rightrightarrows  U$.
Then $s^\an, t^\an\colon R^\an\rightrightarrows  U^\an$ is an \'etale equivalence relation on $U^\an$, and there is a natural isomorphism
\[
F^\an\cong U^\an/R^\an
\]
of sheaves on $(\Afd_K/S)_\Et$. In particular, $F^\an$ is an \'etale space over $S$.
\end{prop}

\begin{proof}
By Lemma~\ref{lem:first properties of analytification of schemes}(i), we see that $s^\an, t^\an\colon R^\an\rightrightarrows  U^\an$ is an \'etale equivalence relation on $U^\an$.
Define the presheaves $F_0, F_1, F_2$ on $\Afd_K/S$ by associating to $T=\Spa(B,B^+)$
\[
F_0(T)\coloneqq U(\Spec \widehat{B})/R(\Spec \widehat{B}),\quad
F_1(T)\coloneqq F(\Spec \widehat{B}), \quad
F_2(T)\coloneqq U^\an(T)/R^\an(T).
\]
Note that $F^\an$ (resp.~$U^\an/R^\an$) is the \'etale sheafification of $F_1$ (resp.~$F_2$). We will show that natural morphisms of presheaves $F_0\rightarrow F_1$ and $F_0\rightarrow F_2$ are isomorphisms after sheafification. By the universal property of sheafification, it is enough to show that given any sheaf $G\in \Sh((\Afd_K/S)_\Et)$, every morphism $F_0\rightarrow G$ uniquely factors through $F_0\rightarrow F_1$ (resp.~$F_0\rightarrow F_2$). For $F_0\rightarrow F_2$, this follows easily from Lemma~\ref{lem:first properties of analytification of schemes}(iii) and \cite[00WK]{stacks-project}. 
For $F_0\rightarrow F_1$, recall that the \'etale sheafification of the presheaf quotient of $U$ by $R$ on $\Sch/A$ is already an fppf sheaf and agrees with $F$ (cf.~\cite[076M, Pf.]{stacks-project}).
We also remark that if $\{\Spec B_i\rightarrow \Spec \widehat{B}\}$ is an \'etale covering, then $\{\Spec B_i\times_{\Spec \widehat{B}}T\rightarrow T\}$ is an \'etale covering in $\V_K/S$. By these observations and \cite[00WK]{stacks-project}, we can check that every morphism $F_0\rightarrow G$ uniquely factors through $F_0\rightarrow F_1$.
The last assertion follows from Proposition~\ref{prop:presentation of etale adic spaces}(ii).
\end{proof}

\smallskip
\noindent
\textbf{Analytification of algebraic stacks}.
Fix $S=\Spa(A,A^+)\in \Afd_K$ with $A$ complete.

\begin{defn}\label{def:analytification of algebraic stack}
Let $\calX$ be an algebraic stack locally of finite type over $\Spec A$. Define the stack $\calX^\an\rightarrow (\Afd_K)_\Et$ to be the stackification of the category fibered in groupoids over $(\Afd_K)_\Et$ whose fiber category over $\Spa(B,B^+)\in \Ob(\Afd_K/S)$ is given by
\[
\calX_{\Spec \widehat{B}}.
\] 
Proposition~\ref{prop:analytification of algebraic stacks} below shows that $\calX^\an$ is an adic stack over $S$. We call $\calX^\an$ the \emph{analytification} of $\calX$.
\end{defn}

\begin{prop}\label{prop:analytification of algebraic stacks}
Let $\calX$ be an algebraic stack locally of finite type over $\Spec A$ and take a presentation of $\calX$ given by a smooth groupoid $(U,R,s,t,c)$ in algebraic spaces.
Then $(U^\an,R^\an,s^\an,t^\an,c^\an)$ is a smooth groupoid in \'etale spaces, and there is a natural equivalence
\[
\calX^\an\cong [U^\an/R^\an]
\]
over $(\Afd_K/S)_\Et$. In particular, $\calX^\an$ is an adic stack over $S$.
\end{prop}

\begin{proof}
It follows from Definition~\ref{def:morphism of etale spaces has property etale local on source-and-target}, Lemma~\ref{lem:first properties of analytification of schemes}(i), and Proposition~\ref{prop:analytification of algebraic spaces} that $(U^\an,R^\an,s^\an,t^\an,c^\an)$ is a smooth groupoid in \'etale spaces.
Since $\calX\cong [U/R]$, we may assume that $\calX$ is split and regard it as a contravariant functor $\Sch/A\rightarrow \Grpd$ (see Remark~\ref{rem:equivalent to a split category fibered in groupoids}).
Define the split categories $\calX_0, \calX_1, \calX_2$ fibered in groupoids over $\Afd_K/S$ by associating to $T=\Spa(B,B^+)\in \Ob(\Afd_K/S)$
\[
\calX_0(T)\coloneqq [U(\Spec \widehat{B})/R(\Spec \widehat{B})],\;
\calX_1(T)\coloneqq \calX(\Spec \widehat{B}), \;
\calX_2(T)\coloneqq [U^\an(T)/R^\an(T)].
\]
Write $\calX_j'$ for the stackification of $\calX_j$. 
Note $\calX^\an=\calX_1'$ and $[U^\an/R^\an]=\calX_2'$. We will show that the natural $1$-morphisms $f_1\colon \calX_0\rightarrow\calX_1$ and $f_2\colon\calX_0\rightarrow\calX_2$ over $\Afd_K/S$ induce equivalences after stackification by using Lemma~\ref{lem:criteria for equivalence of stacks}.

First we show that $f_j'\colon \calX_0'\rightarrow\calX_j'$ is fully faithful by verifying Lemma~\ref{lem:criteria for equivalence of stacks}(i). Take any $T=\Spa(B,B^+)\in \Ob(\Afd_K/S)$ and any $x',y'\in \Ob((\calX_0')_T)$. By Lemma~\ref{lem:stackification}, there exists an \'etale covering $\{T_i=\Spa(B_i,B_i^+)\rightarrow T\}$ such that $x'|_{T_i},y'|_{T_i}$ arise from $x_i,y_i\in U(\Spec \widehat{B}_i)$. Then the morphism 
\begin{equation}\label{eq:analytification of stack}
f_1'\colon \Isom_{\calX_0'}(x'|_{T_i},y'|_{T_i})\rightarrow \Isom_{\calX_1'}(f_1'(x')|_{T_i},f_1'(y')|_{T_i})
\end{equation}
of sheaves on $(\Afd_K/T_i)_\Et$ is the sheafification of the assignment sending $T'=\Spa(B',B'^+)\in\Ob(\Afd_K/T_i)$ to the map
\[
(R\times_{U\times_AU,(y_i,x_i)}\Spec \widehat{B}_i)(\Spec \widehat{B}')\rightarrow \Isom_{\calX}(x_i,y_i)(\Spec \widehat{B}')
\]
obtained as the evaluation of sheaves on $(\Sch/\widehat{B}_i)_\fppf$ at $\Spec \widehat{B}'$, and the latter map is a bijection by \cite[044V]{stacks-project}.
Hence \eqref{eq:analytification of stack} is an isomorphism, and therefore so is $f_1'\colon \Isom_{\calX_0'}(x',y')\rightarrow \Isom_{\calX_1'}(f_1'(x'),f_1'(y'))
$.
Similarly, consider the sheafification of the presheaf morphism over $T_i$ sending $T'=\Spa(B',B'^+)$ to the map
\[
(R\times_{U\times_AU,(y_i,x_i)}\Spec \widehat{B}_i)(\Spec \widehat{B}')\rightarrow (R^\an\times_{U^\an\times_SU^\an,(y_i^\an,x_i^\an)}T_i)(T').
\]
By recalling Definition~\ref{def:analytification of algebraic space}, we see that $f_2'$ is fully faithful.

It remains to verify Lemma~\ref{lem:criteria for equivalence of stacks}(ii) for $f_j'\colon \calX_0'\rightarrow\calX_j'$. For $f_1'$, take any $T=\Spa(B,B^+)\in \Ob(\Afd_K/S)$ and any $x'\in \Ob((\calX_1')_T)$. Then there exists an \'etale covering $\{T_i=\Spa(B_i,B_i^+)\rightarrow T\}$ such that $x'|_{T_i}$ arises from $x_i\in \calX(\Spec \widehat{B}_i)$. Recall from \cite[076V]{stacks-project} that the (fppf) algebraic stack $\calX$ is identified with the stackification of $(U,R,s,t,c)$ with respect to \'etale topology. Hence there exists an \'etale surjective morphism $\Spec B_i'\rightarrow \Spec \widehat{B}_i$ such that $x_i|_{\Spec B_i'}$ arises from an object of $U(\Spec B_i')$. Now $\{\Spec B_i'\times_{\Spec \widehat{B}_i}T_i\rightarrow T\}$ is an \'etale covering by Lemma~\ref{lem:first properties of analytification of schemes}(i). So refine it to an \'etale covering $\{T_{il}\rightarrow T_i\rightarrow T\}$ by affinoids. By our choice, each $x'|_{T_{il}}$ is in the essential image of $f_1'\colon (\calX_0')_{T_{il}}\rightarrow (\calX_1')_{T_{il}}$. Similarly, we can verify Lemma~\ref{lem:criteria for equivalence of stacks}(ii) for $f_2'$ using Definition~\ref{def:analytification of algebraic space}.
\end{proof}

\subsection{Generic fiber}\label{sec:generic fiber}
\begin{prop}[{\cite[Prop.~4.1]{Huber-gen}}]
For a locally Noetherian formal scheme $Y$, there exists an adic space $t(Y)$ and a morphism of locally and topologically ringed spaces $(t(Y),\calO^+_{t(Y)})\rightarrow (Y,\calO_Y)$ satisfying the universal property (see \textit{loc.~cit.} for the precise formulation).
The association $Y\mapsto t(Y)$ is functorial, sending open immersions to open immersions, and if $Y=\Spf R$, then $t(Y)$ is given by $\Spa(R,R)$.
\end{prop}

\begin{rem}\label{rem:t functor commutes with fiber product}
Let $f\colon X\rightarrow Z$ and $g\colon Y\rightarrow Z$ be morphisms of locally Noetherian formal schemes. Then $f$ is adic (resp.~locally of finite type) if and only if $t(f)$ is adic (resp.~locally of finite type) by \cite[Prop.~4.2(i)]{Huber-gen}. If $f$ is adic and $g$ is locally of finite type, then $t(X\times_ZY)=t(X)\times_{t(Z)}t(Y)$. For this, one can reduce to the affine case, which follows by comparing the constructions of the fiber product in  \cite[Prop.~10.7.2]{EGAI} and \cite[Prop.~3.7]{Huber-gen}.
\end{rem}

If $Y$ is a formal scheme locally of finite type over $\Spf \calO_K$, the open subspace $t(Y)_a$ of analytic points of the adic space $t(Y)$ is identified with $t(Y)\times_{\Spa(\calO_K,\calO_K)}\Spa(K)$ and called the \emph{(Raynaud) generic fiber} of $Y$ as in \cite[Rem.~4.6(ii)(iii)]{Huber-gen}.

\begin{rem}\label{rem:explicit formula of generic fiber}
The constructions show that, if $Y=\Spf A_0$ is affine, then $t(Y)_a$ is given by $\Spa(A,A^+)$, where $A=A_0[p^{-1}]$ and $A^+$ is the integral closure of $A_0$ in $A$.
\end{rem}

Let us recall formal schemes locally formally of finite type over $\calO_K$.

\begin{defn}[{\cite[Def.~II.9.6.5]{Fujiwara-Kato}, \cite[Prop.~2.1, Def.~2.2]{Achinger-Lara-Youcis-specialization}}]
\label{def:locally formally of finite type}
\hfill
\begin{enumerate}
 \item An adic $\calO_K$-algebra of finite ideal type $R$ is said to be \emph{topologically formally of finite type} over $\calO_K$ if for every ideal of definition $I\subset R$ containing $\pi$, $R/I$ is an $\calO_K/\pi$-algebra of finite type, or equivalently, if it admits a continuous adic $\calO_K$-algebra surjection of the form $\calO_K\langle x_1,\ldots,x_m\rangle[\![y_1,\ldots,y_n]\!]\twoheadrightarrow R$ for some $m,n\geq 0$, where $\calO_K\langle x_1,\ldots,x_m\rangle[\![y_1,\ldots,y_n]\!]$ denotes the completion of $\calO_K[x_1,\ldots,x_m,y_1,\ldots,y_n]$ with respect to $(\pi,y_1,\ldots,y_n)$.
 \item A formal scheme $Y$ over $\Spf \calO_K$ is said to be \emph{locally formally of finite type} if $Y$ admits an affine open covering $\{U_i=\Spf R_i\}$ such that each $R_i$ is topologically formally of finite type over $\calO_K$. Let $\mathbf{FSch}_{\calO_K}^\mathrm{lfft}$ denote the category of formal schemes locally formally of finite type over $\Spf\calO_K$.
\end{enumerate}
\end{defn}

\begin{rem}
One can verify that the fiber product of formal schemes locally formally of finite type over $\calO_K$ is also locally formally of finite type over $\calO_K$.
\end{rem}

We define the generic fiber $Y_\eta$ of such a formal scheme $Y$ as follows:

\begin{prop}\label{prop:generic fiber of lfft formal scheme as fiber product}
For $Y\in \mathbf{FSch}_{\calO_K}^\mathrm{lfft}$, the fiber product $t(Y)\times_{\Spa(\calO_K,\calO_K)}\Spa(K,\calO_K)$ exists in the category of adic spaces and is representable by a rigid $K$-space, denoted $Y_\eta$.
Moreover, if $Y=\Spf R$ is affine, then for any $\Spa(A,A^+)\in\AdSp_K^\aff$ with $A$ complete, there are functorial identifications 
\[
 \Hom_{\AdSp_K}(\Spa(A,A^+),Y_\eta)\xrightarrow{\cong}\Hom_{\calO_K}(R,A^+)=\varinjlim_{A_0\subset A^+}Y(\Spf A_0),
\]
where $\Hom_{\calO_K}(-,-)$ is the set of continuous $\calO_K$-algebra maps, and $A_0$ runs over the open and bounded $\calO_K$-subalgebras of $A^+$. We call $Y_\eta$ the \emph{generic fiber} of $Y$.
\end{prop}

\begin{proof}
Fujiwara and Kato \cite[II.9.6]{Fujiwara-Kato} construct a functor $\mathbf{FSch}_{\calO_K}^\mathrm{lfft}\rightarrow \Rig_K$, which we denote by $Y\mapsto Y_{\eta'}$ for the time being, satisfying the following properties:
\begin{enumerate}
    \item the functor sends open immersions to open immersions;
    \item if $Y=\Spf R$ is affine, then the functorial identifications in the statement hold for $Y_{\eta'}$ in place of $Y_\eta$.
\end{enumerate}
Here (i) follows from the construction given in \cite[II.9.6(a)]{Fujiwara-Kato} and (ii) is a special case of \cite[Lem.~2.9]{Achinger-Lara-Youcis-specialization}. We claim that $Y_{\eta'}$ represents the fiber product $t(Y)\times_{\Spa(\calO_K,\calO_K)}\Spa(K,\calO_K)$. When $Y$ is affine, this follows from (ii) and the fact that $t(Y)=\Spa(R,R)$. The general case is deduced from the affine case by (i).
\end{proof}

\begin{cor}\label{cor:functor of points of generic fiber}
For $Y\in \mathbf{FSch}_{\calO_K}^\mathrm{lfft}$, the generic fiber $Y_\eta$ as a presheaf on $\AdSp_K^\aff$ is identified with the sheafification with respect to analytic topology of the presheaf
\[
\Spa(A,A^+)\mapsto \varinjlim_{A_0\subset \widehat{A}^+}Y(\Spf A_0),
\]
where $A_0$ runs over the open and bounded $\calO_K$-subalgebras of the completion $\widehat{A}^+$.
\end{cor}

\begin{proof}
This follows easily from Proposition~\ref{prop:generic fiber of lfft formal scheme as fiber product}.
\end{proof}

\begin{rem}\label{rem:description of generic fiber of formal scheme over L}
Let $(L,L^+)$ be a complete analytic affinoid field over $(K,\calO_K)$.
\begin{enumerate}
 \item Let $(A,A^+)$ be a complete affinoid ring of topologically finite type over $(L,L^+)$. Then $A^+$ is the filtered union of open $L^+$-subalgebras of $A^+$ that are $\pi$-adically of finite type over $L^+$ (hence bounded). To see this, take a ring $A_0$ of definition such that $A_0$ is $\pi$-adically of finite type over $L^+$\footnote{Note that the following properties mean the same: being $\pi$-adically of finite type over $L^+$, being adically of finite type over $L^+$, and being of topologically finite type over $L^+$.} and $A^+$ is the integral closure of $A_0$ in $A$. For any finitely many elements $a_1,\ldots,a_n\in A^+$, the subalgebra $A_0[a_1,\ldots,a_n]$ is finitely generated as an $A_0$-module and thus $\pi$-adically complete by \cite[Cor.~0.8.2.20]{Fujiwara-Kato}.
 \item With the notation as in Corollary~\ref{cor:functor of points of generic fiber}, assume further that $Y$ is locally adically of finite type over $\calO_K$. Then the generic fiber $Y_\eta$ as a presheaf on $\AdSp_{\Spa(L,L^+)}^{\aff,\lft}$ is identified with the sheafification with respect to analytic topology of the presheaf
\[
\Spa(A,A^+)\mapsto \varinjlim_{A_0\subset \widehat{A}^+}Y(\Spf A_0),
\]
where $A_0$ runs over the open $L^+$-subalgebras of $\widehat{A}^+$ that are $\pi$-adically of finite type over $L^+$; we may reduce it to the affine case; then use (i) to check that the equality $\Hom_{\calO_K}(R,\widehat{A}^+)=\varinjlim_{A_0\subset \widehat{A}^+}Y(\Spf A_0)$ in Proposition~\ref{prop:generic fiber of lfft formal scheme as fiber product} continues to hold when $A_0$ runs only over such $L^+$-subalgebras. We also remark that the same holds for \'etale sheafification since analytic topology is coarser than \'etale topology.
\end{enumerate}
\end{rem}

\begin{prop}\label{prop: generic fiber of etale/smooth map}
If a morphism $Y\rightarrow Y'$ between formal schemes locally formally of finite type over $\Spf\calO_K$ is \'etale (resp.~smooth), then the induced map $Y_\eta\rightarrow Y'_\eta$ is \'etale (resp.~smooth).
\end{prop}

\begin{proof}
The \'etale case is proved in \cite[Prop.~2.10]{Achinger-Lara-Youcis-specialization}, and the smooth case follows from the \'etale case since one can easily see that $(Y\times_{\Spf\calO_K}\Spf R)_\eta=Y_\eta\times_{\Spa(K,\calO_K)}\Spa(R[p^{-1}],R)$  for $R=\calO_K\langle x_1,\ldots,x_m\rangle$.
\end{proof}

\begin{prop}\label{prop:generic fiber of formal completion along subset in the case of formal schemes}
Let $Y$ be a formal scheme locally of finite type over $\Spf\calO_K$, let $Z\subset \lvert Y\rvert$ be a closed subset, and write $\widehat{Y}_{|Z}$ for the formal completion along $Z$. Then the induced map $(\widehat{Y}_{|Z})_\eta\rightarrow Y_\eta$ is an open immersion.
\end{prop}

\begin{proof}
This is part of \cite[Prop.~2.14]{Achinger-Lara-Youcis-specialization}, and the quasi-compactness hypothesis therein is not used for this assertion.
\end{proof}

\begin{rem}
In the above construction, we assume that $\calO_K$ is a \emph{discrete} valuation ring so that every formal scheme $Y$ locally formally of finite type over $\calO_K$ is locally Noetherian and thus $t(Y)$ is defined. One can still construct the generic fiber without the discreteness assumption: see \cite[\S1.9]{Huber-etale} (where the generic fiber is denoted by $d(Y)$) and \cite[\S2.2]{Achinger-Lara-Youcis-specialization}.
\end{rem}

\begin{lem}\label{lem:properties of A0 to (A,A+)}
For $T=\Spa(A,A^+)\in \AdSp_K$ with $A$ complete, let $A_0\subset A^+$ be an open and bounded $\calO_K$-subalgebra. Given a map $f\colon \Spf A_0'\rightarrow \Spf A_0$ of affine formal schemes adically of finite type, the fiber product $Z\coloneqq \Spa(A_0',A_0')\times_{\Spa(A_0,A_0)}T$ exists in the category of adic spaces and is representable by an affinoid adic space of finite type over $T$. If $f$ is \'etale (resp. an \'etale covering), then $Z\rightarrow T$ is \'etale (resp. an \'etale covering).
\end{lem}

\begin{proof}
Note that $\pi A_0$ is an ideal of definition of $A_0$ by Remark~\ref{rem:open bounded subring}. 
The first assertion follows from \cite[Prop.~1.2.2(b)]{Huber-etale}. In fact, the proof therein shows that it is represented by $\Spa(A',A'^+)$ (or its completion) where $A'\coloneqq A_0'[\pi^{-1}]$ is topologized by making $A_0'$ a ring of definition, and $A'^+$ denotes the integral closure of the image of $A_0'\otimes_{A_0}A^+$ in $A'$.
With the notation in \cite[Prop.~1.9.1]{Huber-etale}, we also see $\Spa(A',A'^+)=d(\Spf A_0')\times_{d(\Spf A_0)}T$. Hence the second assertion follows from \cite[Lem.~3.5.1, Prop.~1.6.7(iv)]{Huber-etale}.
\end{proof}

\smallskip
\noindent
\textbf{Generic fiber of formal algebraic spaces}.

\begin{defn}\label{def:locally formally of finite type for formal algebraic space}
A formal algebraic space $F$ over $\calO_K$ is said to be \emph{locally formally of finite type over $\Spf \calO_K$} if there exist a formal scheme locally formally of finite type $U$ over $\Spf \calO_K$ and  a morphism $U\rightarrow F$ that is representable, \'etale, and surjective. Note that if $F$ is a formal scheme, then this definition agrees with Definition~\ref{def:locally formally of finite type}(ii) by \cite[02KX]{stacks-project}.
\end{defn}

\begin{defn}\label{def:generic fiber of formal algebraic space}
Let $F$ be a formal algebraic space over $\calO_K$ that is locally formally of finite type over $\Spf\calO_K$. Define the sheaf $F_\eta$ on $(\Afd_K/S)_\Et$ to be the sheafification of the presheaf
\[
(\Afd_K)/S\ni \Spa(B,B^+)\mapsto \varinjlim_{B_0\subset \widehat{B}^+}F(\Spf B_0),
\]
where $B_0$ runs over the open and bounded $\calO_K$-subalgebras of the completion $\widehat{B}^+$ of $B^+$.
When $F$ is a formal scheme, $F_\eta$ agrees with the sheaf associated to the generic fiber by Corollary~\ref{cor:functor of points of generic fiber}. Proposition~\ref{prop:analytification of algebraic spaces} below shows that $F^\an$ is an \'etale space over $K$. We call $F_\eta$ the \emph{generic fiber} of $F$.
\end{defn}

We use a presentation of $F$ by an \'etale equivalence relation (Proposition~\ref{prop:presentation of formal algebraic space}) to show that $F_\eta$ is an \'etale space over $K$:

\begin{prop}\label{prop:generic fiber of formal algebraic spaces}
Let $F$ be a formal algebraic space locally formally of finite type over $\Spf\calO_K$ and take a presentation of $F$ given by an \'etale equivalence relation $s,t\colon R\rightrightarrows  U$ .
Then the diagram $s_\eta, t_\eta,\colon R_\eta\rightrightarrows  U_\eta$ is an \'etale equivalence relation on $U_\eta$ and induces a natural isomorphism
\[
F_\eta\cong U_\eta/R_\eta
\]
of sheaves on $(\Afd_K)_\Et$. In particular, $F_\eta$ is an \'etale space over $K$, and the functor $F\mapsto F_\eta$ sends open immersions to open immersions. 
\end{prop}

\begin{proof}
Note that $U$ and $R$ are necessarily formal schemes locally formally of finite type over $\Spf\calO_K$. 
By construction, $R_\eta\rightarrow U_\eta\times U_\eta$ is a monomorphism of sheaves. It follows from Proposition~\ref{prop: generic fiber of etale/smooth map} that the diagram $s_\eta, t_\eta,\colon R_\eta\rightrightarrows  U_\eta$ is an \'etale equivalence relation on $U_\eta$.
Define the presheaves $F_0, F_1, F_2$ on $\Afd_K$ by associating to $T=\Spa(B,B^+)$
\[
F_0(T)\coloneqq \varinjlim_{B_0\subset \widehat{B}^+}U(\Spf B_0)/R(\Spf B_0),\quad
F_1(T)\coloneqq \varinjlim_{B_0\subset \widehat{B}^+}F(\Spf B_0), 
\]
where $B_0$ runs over the open and bounded $\calO_K$-subalgebras of $\widehat{B}^+$, and $F_2(T)\coloneqq U_\eta(T)/R_\eta(T)$.
Note that $F_\eta$ (resp.~$U_\eta/R_\eta$) is the \'etale sheafification of $F_1$ (resp.~$F_2$). We will show that natural morphisms of presheaves $F_0\rightarrow F_1$ and $F_0\rightarrow F_2$ are isomorphisms after sheafification. By the universal property of sheafification, it is enough to show that given any sheaf $G\in \Sh((\Afd_K)_\Et)$, every morphism $F_0\rightarrow G$ uniquely factors through $F_0\rightarrow F_1$ (resp.~$F_0\rightarrow F_2$). For $F_0\rightarrow F_2$, this follows easily from Corollary~\ref{cor:functor of points of generic fiber} since \'etale topology is finer than analytic topology. For $F_0\rightarrow F_1$, the assertion follows from  Remark~\ref{rem:presentation of formal algebraic space as etale sheaf}, Lemma~\ref{lem:properties of A0 to (A,A+)}, and \cite[00WK]{stacks-project}.
The last assertions follow from the presentation $F_\eta\cong U_\eta/R_\eta$.
\end{proof}

\begin{prop}\label{prop:generic fiber of formal algebraic space restricted to adic spaces over L}
Let $F$ be a formal algebraic space locally of finite type over $\Spf\calO_K$. Fix a complete analytic affinoid field $(L,L^+)$ over $(K,\calO_K)$ and an affinoid $(A,A^+)$ topologically finite type over $(L,L^+)$. Then the restriction of $F_\eta$ to the small (affinoid) \'etale site of $T\coloneqq \Spa(A,A^+)$ agrees with the sheafification of
\[
\Spa(B,B^+)\mapsto \varinjlim_{B_0\subset \widehat{B}^+} F(\Spf R_0)
\]
where $B_0$ runs over the open $L^+$-subalgebras of $\widehat{B}^+$ of topologically finite type.
\end{prop}

\begin{proof}
Taking a presentation of $F$ and using Remark~\ref{rem:description of generic fiber of formal scheme over L}(ii), one can prove the assertion as in Proposition~\ref{prop:generic fiber of formal algebraic spaces}.
\end{proof}

\begin{prop}\label{prop: properties of generic fiber of formal algebraic spaces}
\hfill
\begin{enumerate}
\item The generic fiber functor commutes with fiber products and sends monomorphisms (resp.~epimorphisms) of sheaves to monomorphisms (resp.~epimorphisms) of sheaves.
\item If a morphism $F_1\rightarrow F_2$ between formal algebraic spaces locally formally of finite type over $\Spf\calO_K$ is representable by algebraic spaces and \'etale (resp.~smooth), then the induced morphism $F_{1,\eta}\rightarrow F_{2,\eta}$ of \'etale spaces is \'etale (resp.~smooth). 
\end{enumerate}
\end{prop}

\begin{proof}
(i) Observe that the association of the presheaf $(\Afd_K)/S\ni \Spa(B,B^+)\mapsto \varinjlim_{B_0\subset \widehat{B}^+}F(\Spf B_0)$ commutes with fiber products.
Now the assertion on fiber products follows since sheafification commutes with finite limits. The assertion on monomorphisms is obvious and the assertion on epimorphisms follows from Lemma~\ref{lem:properties of A0 to (A,A+)}.

(ii)
For $i=1,2$, choose a formal scheme $U_i$ locally formally of finite type over $\calO_K$ and a morphism $U_i\rightarrow F_i$ that is representable by algebraic spaces, \'etale and surjective. 
Consider the commutative diagrams
\[
\xymatrix{
U_1\times_{F_2}U_2\ar[rr]\ar[d]& &U_2\ar[d]\\
U_1\ar[r]&F_1\ar[r] &F_2.
}
\quad\text{and}\quad
\xymatrix{
(U_1\times_{F_2}U_2)_\eta\ar[rr]\ar[d]& &(U_2)_\eta\ar[d]\\
(U_1)_\eta\ar[r]&(F_1)_\eta\ar[r] &(F_2)_\eta.
}
\]
By Lemma~\ref{lem:representability of diagonal of formal algebraic space}, $U_1\times_{F_2}U_2$ is also a formal scheme locally formally of finite type over $\calO_K$, and its map to $U_2$ is \'etale (resp.~smooth).
By Proposition~\ref{prop: generic fiber of etale/smooth map}, we see that $(U_1\times_{F_2}U_2)_\eta=(U_1)_\eta\times_{(F_2)_\eta}(U_2)_\eta$, the map $(U_i)_\eta\rightarrow (F_i)_\eta$ is \'etale surjective, and $(U_1\times_{F_2}U_2)_\eta\rightarrow (U_2)_\eta$ is \'etale (resp.~smooth). Hence $F_{1,\eta}\rightarrow F_{2,\eta}$ is \'etale (resp.~smooth) according to Definition~\ref{def:morphism of etale spaces has property etale local on source-and-target}.
\end{proof}

\smallskip
\noindent
\textbf{Generic fiber of formal algebraic stacks}.

\begin{defn}
Let $\calX$ be a formal algebraic stack over $\calO_K$ with a $1$-morphism $\calX\rightarrow \Spf \calO_K$ of stacks. We say that $\calX$ is \emph{locally formally of finite type over $\Spf \calO_K$} if there exist a formal algebraic space $U$ that is locally formally of finite type over $\Spf\calO_K$ and a $1$-morphism $U\rightarrow\calX$ of stacks over $\Spf\calO_K$ that is representable by algebraic spaces, smooth, and surjective. If $\calX=\calS_X$ for a formal algebraic space $X$, this is equivalent to $X$ being locally formally of finite type over $\Spf\calO_K$.
\end{defn}

\begin{defn}\label{def:generic fiber of formal algebraic stack}
Let $\calX$ be a formal algebraic stack that is locally formally of finite type over $\Spf \calO_K$. Consider a category over $\Afd_K$  whose fiber category over $\Spa(B,B^+)\in \Ob(\Afd_K)$ is the filtered colimit of groupoids
\[
\colim_{B_0\subset \widehat{B}^+}\calX_{\Spf B_0},
\]
where $B_0$ runs over the open and bounded $\calO_K$-subalgebras of the completion $\widehat{B}^+$ of $B^+$ (see Remark~\ref{rem:formal scheme valued point of formal algebraic stack} for the definition of $\calX_{\Spf B_0}$).
It is a category fibered in groupoids by Remark~\ref{rem:open bounded subring}. Define the stack 
\[
\calX_\eta\rightarrow (\Afd_K)_\Et
\]
to be its stackification.
Proposition~\ref{prop:generic fiber of formal algebraic stacks} below shows that $\calX_\eta$ is an adic stack over $S$. We call $\calX_\eta$ the \emph{generic fiber} of $\calX$. If $\calX=\calS_F$ for a formal algebraic space $F$, one can see that $\calX_\eta$ is identified with $\calS_{F_\eta}$.
\end{defn}

\begin{prop}\label{prop:generic fiber of formal algebraic stacks}
Let $\calX$ be a formal algebraic stack locally formally of finite type over $\Spf\calO_K$ and take a presentation of $\calX$ given by a smooth groupoid $(U,R,s,t,c)$ in formal algebraic spaces.
Then $(U_\eta,R_\eta,s_\eta,t_\eta,c_\eta)$ is a smooth groupoid in \'etale spaces, and there is a natural equivalence
\[
\calX_\eta\cong [U_\eta/R_\eta]
\]
over $(\Afd_K/S)_\Et$. In particular, $\calX_\eta$ is an adic stack over $\Spa(K,\calO_K)$.
\end{prop}

See Definition~\ref{def:presentation of formal algebraic stack} (resp.~\ref{def:presentation of adic stack}) for the generalities of presentations of formal algebraic stacks (resp.~adic stacks).

\begin{proof}
Proposition~\ref{prop: properties of generic fiber of formal algebraic spaces}(ii) shows that $(U_\eta,R_\eta,s_\eta,t_\eta,c_\eta)$ is a smooth groupoid in \'etale spaces.
Since $\calX\cong [U/R]$, we may assume that $\calX$ is split and regard it as a contravariant functor $\Sch/\calO_K\rightarrow \Grpd$ (cf.~Remark~\ref{rem:equivalent to a split category fibered in groupoids}).
Define the split categories $\calX_0, \calX_1, \calX_2$ fibered in groupoids over $\Afd_K$ by associating to $T=\Spa(B,B^+)\in \Ob(\Afd_K)$
\[
\calX_0(T)\coloneqq \colim_{B_0\subset \widehat{B}^+} [U(\Spf B_0)/R(\Spf B_0)],\quad
\calX_1(T)\coloneqq \colim_{B_0\subset \widehat{B}^+}\calX(\Spf B_0), 
\]
where $B_0$ runs over the open and bounded $\calO_K$-subalgebras of $\widehat{B}^+$, and $\calX_2(T)\coloneqq [U_\eta(T)/R_\eta(T)]$.
Write $\calX_j'$ for the stackification of $\calX_j$. 
Note $\calX_\eta=\calX_1'$ and $[U_\eta/R_\eta]=\calX_2'$. We will show that each natural $1$-morphisms $f_j\colon \calX_0\rightarrow\calX_j$ over $\Afd_K$ induces an equivalence after stackification by using Lemma~\ref{lem:criteria for equivalence of stacks}.

First we show that $f_j'\colon \calX_0'\rightarrow\calX_j'$ is fully faithful by verifying Lemma~\ref{lem:criteria for equivalence of stacks}(i). Take any $T=\Spa(B,B^+)\in \Ob(\Afd_K)$ and any $x',y'\in \Ob((\calX_0')_T)$. By Lemma~\ref{lem:stackification}, there exists an \'etale covering $\{T_i=\Spa(B_i,B_i^+)\rightarrow T\}$ such that $x'|_{T_i},y'|_{T_i}$ arise from $x_i,y_i\in U(\Spf B_{i0})$ for some open and bounded $\calO_K$-subalgebra $B_{i0}$ of $\widehat{B}_i$. Then the morphism 
\begin{equation}\label{eq:generic fiber of formal algebraic stack}
f_1'\colon \Isom_{\calX_0'}(x'|_{T_i},y'|_{T_i})\rightarrow \Isom_{\calX_1'}(f_1'(x'|_{T_i}),f_1'(y'|_{T_i}))
\end{equation}
of sheaves on $(\Afd_K/T_i)_\Et$ is the sheafification of the assignment sending $T'=\Spa(B',B'^+)\in\Ob(\Afd_K/T_i)$ to the map
\[
\colim_{B'_0\subset \widehat{B}'^+}(R\times_{U\times U,(y_i,x_i)}\Spf B_{i0})(\Spf B_0')\rightarrow \colim_{B'_0\subset \widehat{B}'^+}\Isom_{\calX}(x_i,y_i)(\Spf B_0').
\]
where $U\times U\coloneqq U\times_{\Spf \calO_K}U$. By Remark~\ref{rem:open bounded subring}, $\Spf B_0'=\varinjlim_n\Spec B_0'/(\pi^n)$. So $(R\times_{U\times U,(y_i,x_i)}\Spf B_{i0})|_{\Spec B_0'/(\pi^n)}\rightarrow \Isom_{\calX}(x_i,y_i)|_{\Spec B_0'/(\pi^n)}$ is an isomorphism of algebraic spaces (cf.~\cite[044V, Pf.]{stacks-project}). Hence \eqref{eq:generic fiber of formal algebraic stack} is an isomorphism, which implies that $\Isom_{\calX_0'}(x',y')\rightarrow \Isom_{\calX_1'}(f_1'(x'),f_1'(y'))$ is also an isomorphism. Similarly, consider the presheaf morphism over $T_i$ sending $T'=\Spa(B',B'^+)$ to 
\begin{equation}\label{eq:generic fiber of formal algebraic stack: presheaf map for f2}
\colim_{B'_0\subset \widehat{B}'^+}(R\times_{U\times U,(y_i,x_i)}\Spf B_{i0})(\Spf B_0')\rightarrow (R_\eta\times_{U_\eta\times U_\eta,((y_i)_\eta,(x_i)_\eta)}T_i)(T'). \end{equation}
By recalling Definition~\ref{def:generic fiber of formal algebraic space}, we see that the morphism \eqref{eq:generic fiber of formal algebraic stack: presheaf map for f2} is an isomorphism after sheafification, from which we conclude that $f_2'$ is fully faithful.

It remains to verify Lemma~\ref{lem:criteria for equivalence of stacks}(ii) for $f_j'\colon \calX_0'\rightarrow\calX_j'$. For $f_1'$, choose any $T=\Spa(B,B^+)\in \Ob(\Afd_K)$ and any $x'\in \Ob((\calX_1')_T)$. Take an \'etale covering $\{T_i=\Spa(B_i,B_i^+)\rightarrow T\}$ and an open and bounded $\calO_K$-subalgebra $B_{i0}$ of $\widehat{B}_i$ such that $x'|_{T_i}$ arises from $x_i\in \calX(\Spf B_{i0})$. Lemma~\ref{lem:presentation of formal algebraic stack in etale topology} identifies the (fppf) formal algebraic stack $\calX$ with the \'etale stackification of $(U,R,s,t,c)$. 
Write $x_i=(x_i^{(n)})$ via $\calX(\Spf B_{i0})=\varprojlim_n \calX(\Spec B_{i0}/(\pi^n))$. Take an \'etale covering $\{\Spec \overline{B}_{i0,l}\rightarrow \Spec B_{i0}/(\pi)\}$ and $x_{il}^{(1)}\in U(\Spec \overline{B}_{i0,l})$ such that $x_i^{(1)}|_{\Spec B_{i0,l}/(\pi)}$ is the image of $x_{il}^{(1)}$. This \'etale covering lifts uniquely to an \'etale covering $\{\Spf B_{i0,l}\rightarrow \Spf B_{i0}\}$. Since $U\rightarrow \calX$ is representable by algebraic spaces and smooth, $x_{il}^{(1)}$ lifts to an element $x_{il}=(x_{il}^{(n)})\in U(\Spf B_{i0,l})=\varprojlim_nU(\Spec B_{i0,l}/(\pi^n))$ such that $x_i|_{\Spf B_{i0,l}}$ is the image of $x_{il}$ by \cite[04AM]{stacks-project}. By Lemma~\ref{lem:properties of A0 to (A,A+)}, this \'etale covering of $\Spf B_{i0}$ yields an \'etale covering $\{T_{il}\rightarrow T_i\}$ by affinoids, and $x'|_{T_{il}}$ is in the essential image of $f_1'\colon \calX_0'|_{T_{il}}\rightarrow \calX_1'|_{T_{il}}$, which verifies Lemma~\ref{lem:criteria for equivalence of stacks}(ii) for $f_1'$. The verification for $f_2'$ directly follows from Definition~\ref{def:locally formally of finite type for formal algebraic space}.
\end{proof}

\begin{rem}\label{rem:generic fiber of formal algebraic stacks}
Assume that $\calX$ is a formal algebraic stack locally of finite type over $\Spf\calO_K$. Let $T=\Spa(B,B^+)\in \Afd_K$ with $B$ complete. By Lemma~\ref{lem:stackification}, we know that for any $x'\in \Ob(\calX_{\eta,T})$, there exists an \'etale covering $\{T_i=\Spa(B_i,B_i^+)\rightarrow T\}$ with $B_i$ complete such that, for every $i$, the pullback $x'|_{T_i}$ is in the essential image of $\colim_{B_{i,0}\subset B_i^+}\calX_{\Spf B_{i,0}}\rightarrow\calX_{\eta,T_i}$ where $B_0$ runs over the open and bounded $\calO_K$-subalgebras of $B_i^+$.
In fact, we can deduce a slightly stronger result from the proof of Proposition~\ref{prop:generic fiber of formal algebraic stacks}: take a complete affinoid field $(L,L^+)$ over which $(B,B^+)$ is topologically of finite type. Then one can let $B_{i,0}$ run only over the open $L^+$-subalgebras of $B_i^+$ of topologically finite type. To see this, observe that $\calX_\eta$ is also a stackification of $\calX_2$ introduced in the proof and apply Proposition~\ref{prop:generic fiber of formal algebraic space restricted to adic spaces over L} to $U$.
\end{rem}

\begin{prop}\label{prop:generic fiber of formal completion along subset in the case of formal algebraic stacks}
Let $\calX\rightarrow (\Sch/\calO_K)_\fppf$ be a formal algebraic stack locally of finite type over $\Spf\calO_K$ and let $Z\subset\lvert \calX_\red\rvert$ be a closed subset. Then the induced $1$-morphism $(\widehat{\calX}_{|Z})_\eta\rightarrow \calX_\eta$ is an open immersion.
\end{prop}

\begin{proof}
The case when $\calX$ comes from a formal scheme is Proposition~\ref{prop:generic fiber of formal completion along subset in the case of formal schemes}. By taking a presentation, it is straightforward to see that the assertion continues to hold first for formal algebraic spaces and then for formal algebraic stacks.
\end{proof}

\subsection{Comparison of generic fiber and analytification}
For a scheme/algebraic space/algebraic stack $S$ locally of finite type over $\calO_K$, set $S_K\coloneqq S\times_{\Spec \calO_K}\Spec K$ and 
write $\widehat{S}$ for its $p$-adic completion (cf.~\cite[0AMC]{stacks-project} and Definition~\ref{def:I-adic completion of algebraic stack}).

\begin{prop}[{\cite[Rem.~4.6(iv)]{Huber-gen}}]
For every scheme $Y$ locally of finite type over $\calO_K$, there is a morphism of adic spaces over $\Spa(K,\calO_K)$
\[
\widehat{Y}_\eta \rightarrow Y_K^\an,
\]
which is a local isomorphism and functorial in $Y$. Moreover, if $Y$ is separated (resp.~proper) over $\calO_K$, then it is an open immersion (resp.~an isomorphism).
 \end{prop}

In general, $\widehat{Y}_\eta \rightarrow Y_K^\an$ is not an open immersion: consider two copies of $\A^1_{\calO_K}$ glued along the open $\A^1_K$, which gives the disjoint union of two closed unit disks mapping to the adic affine line.

\begin{rem}\label{rem:description of map from generic fiber to analytification} 
If $f\colon X\rightarrow Y$ is a morphism of schemes locally of finite type over $\calO_K$, the two composites of morphisms $\widehat{X}_\eta\xrightarrow{\widehat{f}_\eta}\widehat{Y}_\eta \rightarrow Y_K^\an $ and $\widehat{X}_\eta\rightarrow X_K^\an \xrightarrow{f_K^\an}Y_K^\an$ coincide.
Similarly, we have the following variant: let $(L,L^+)$ be a complete analytic affinoid field over $(K,\calO_K)$ and let $A_0$ be a $p$-adically complete ring of topologically finite type over $L^+$. Set $A\coloneqq A_0[p^{-1}]$ and let $A^+$ be the integral closure of $A_0$ in $A$. Then $\Spa(A,A^+)$ is an affinoid of finite type over $\Spa(L,L^+)$.
Assume that a morphism $f\colon \Spec A_0\rightarrow Y$ of schemes is given. 
On one hand, the base change $f_K\colon \Spec A\rightarrow Y_K$ yields $f_K^\an\colon \Spa(A,A^+)\rightarrow Y_K^\an$ by Lemma~\ref{lem:analytification and base change}. On the other hand, $f$ gives rise to a morphism $\Spf A_0\rightarrow \widehat{Y}$ of formal schemes, which in turn yields a morphism $\widehat{f}_\eta\colon \Spa(A,A^+)\rightarrow \widehat{Y}_\eta$ (cf.~\cite[Prop.~1.9.1]{Huber-etale}). Then the diagram
\[
\xymatrix{
\Spa(A,A^+)\ar[d]_-{\widehat{f}_\eta}\ar@{=}[r] & \Spa(A,A^+) \ar[d]^-{f_K^\an} \\
\widehat{Y}_\eta \ar[r] & Y_K^\an
}
\]
commutes. To see this, we can reduce it to the case where $Y$ is affine of the form $\Spec \calO_K[x_1,\ldots,x_m]/(f_1,\ldots,f_r)$, in which the commutativity follows from the constructions.
\end{rem}

\begin{defn}\label{defn:comparison map from generic fiber to analytification}
Define the natural transformation $\widehat{Y}_\eta \rightarrow Y_K^\an$ between the functors $Y\mapsto \widehat{Y}_\eta$ and $Y\mapsto Y_K^\an$ from the category of algebraic spaces locally of finite type over $\calO_K$ to $\EtSp_K$ as follows: for each $Y$, take a presentation $Y=U/R$ by an \'etale equivalence relation $R\rightrightarrows  T$ by $\calO_K$-schemes. Then we have a commutative diagram in $\V_K$
\[
\xymatrix{
\widehat{R}_\eta\ar[r]\ar@<-.5ex>[d]\ar@<.5ex>[d] &  R_K^\an \ar@<-.5ex>[d]\ar@<.5ex>[d]\\
\widehat{U}_\eta\ar[r] & U_K^\an,
}
\]
which induces a morphism $\widehat{Y}_\eta \rightarrow Y_K^\an$ in $\V_K$. It is straightforward to see that the induced map is independent of the choice of presentation and functorial in $Y$. Note that $\widehat{Y}_\eta \rightarrow Y_K^\an$ is, locally on the source, an open immersion.

Similarly, we obtain a $1$-morphism $\widehat{\calY}_\eta \rightarrow \calY_K^\an$ of adic stacks for an algebraic stack $\calY$ locally of finite type over $\calO_K$, which is unique up to unique $2$-isomorphisms.
\end{defn}

If $\calY$ has a moduli interpretation for which a GFGA type theorem holds, one can give an explicit description of $\widehat{\calY}_\eta \rightarrow \calY_K^\an$: see \S\ref{sec:generic fiber of p-adic completion of de Rham moduli} in the de Rham moduli case.

\begin{prop}\label{prop:image of comparison map from generic fiber to analytification}
Let $\calY\rightarrow (\Sch/\calO_K)_\fppf$ be an algebraic stack locally of finite type over $\calO_K$. Then there exists a unique open substack $\calW\subset \calY_K^\an$ such that 
the $1$-morphism $\widehat{\calY}_\eta \rightarrow \calY_K^\an$ factors as $\widehat{\calY}_\eta \rightarrow \calW\rightarrow \calY_K^\an$ with $\widehat{\calY}_\eta \rightarrow \calW$ being an epimorphism.
More generally, such an open substack exists when $\widehat{\calY}_\eta$ is replaced with an open immersion $\calV\hookrightarrow \widehat{\calY}_\eta$.
\end{prop}

\begin{proof}
Take a presentation $\calY=[U/R]$ such that $U$ is an algebraic space separated over $\calO_K$. Then we have a $2$-commutative diagram
\[
\xymatrix{
[\widehat{U}_\eta/\widehat{R}_\eta]\ar[d]_-\cong \ar[r]& [U_K^\an/R_K^\an] \ar[d]^-\cong\\
\widehat{\calY}_\eta\ar[r] & \calY_K^\an\\ 
}
\]
with $\widehat{U}_\eta\rightarrow U_K^\an$ being an open immersion. Hence the assertions follow from Corollary~\ref{cor:epi-open immersion factorization} when $\calV=\widehat{\calY}_\eta$. In the general case, we may assume that $\calV$ is an open substack of $\widehat{\calY}_\eta$ by Remark~\ref{rem:open immersion and open substack}. Then $\calV\times_{\widehat{\calY}_\eta}\widehat{U}_\eta$ is representable by an open \'etale subspace $V\subset \widehat{U}_\eta$ and it gives a presentation $[V/(\widehat{R}_\eta|_V)]\xrightarrow{\cong}\calV$. Now we can repeat the above argument to conclude.
\end{proof}

\section{Moduli stacks of integrable connections}\label{sec:Moduli stacks of integrable connections}

We study the algebraic stack $\cM_{\dR}$ of integrable (logarithmic) connections. In this section, fix a Noetherian scheme $S$  and a smooth proper $S$-scheme $Z$. We also fix a (possibly empty) normal crossings divisor $D\subset Z$ relative to $S$, and let $M_{Z}$ be the associated log structure on $Z$ (i.e., $M_{Z}=\calO_{Z}\cap j_\ast \calO_{Z\smallsetminus D}^\times$ where $j\colon Z\smallsetminus D\hookrightarrow Z$).
For each $S$-scheme $T$, write $(Z_T,D_{T})$ (resp. $M_{Z_{T}}$) for the base change of $(Z,D)$ along $T\rightarrow S$ (resp. the pullback log structure on $Z_T$ from $M_Z$). Let $\omega^{i}_{Z_{T}/T}$ denote the sheaf of log differential $i$-forms on $(Z_{T},M_{Z_{T}})$ over $T$ (see \cite[\S4]{Katz_Nilpotent} or \cite[(1.7)]{Kato-log}). If $D$ is empty, this coincides with the sheaf of differential $i$-forms $\Omega^{i}_{Z_{T}/T}$.

\subsection{Definition of \texorpdfstring{$\cM_{\dR}$}{MdR} and its variants}\label{section:MdR}
\begin{defn}\label{def:moduli stack of integrable connections}
We define the \emph{moduli stack of integrable connections},
\[
\calM_\dR((Z,D)/S)\rightarrow (\Sch/S)_\fppf,
\]
to be the category fibered in groupoids whose fiber over $T\rightarrow S$ is the groupoid of pairs $(\calE,\nabla)$ where $\calE$ is a $T$-flat $\calO_{Z_T}$-modules of finite presentation, and $\nabla\colon \calE\rightarrow\calE\otimes_{\cO_{Z_{T}}} \omega^{1}_{Z_{T}/T}$ is an integrable $T$-connection (see \cite[(4.1)]{Katz_Nilpotent}).  
For simplicity, we call such $(\calE,\nabla)$ an \emph{integrable $T$-connection} and, when $D=\emptyset$, write $\calM_\dR(Z/S)$ for $\calM_\dR((Z,D)/S)$. 
\end{defn}

By definition, $\cM_{\dR}((Z_{T},D_{T})/T)\cong\cM_{\dR}((Z,D)/S)\times_{S}T$ for any $S$-scheme $T$. 

\begin{defn}
Let $\calCoh(Z/S)\rightarrow (\Sch/S)_\fppf$ denote the stack of $S$-flat sheaves of finite presentation defined in \cite[08KB]{stacks-project}\footnote{Every $\calO_{Z_T}$-module $\calF$ of finite presentation is quasi-coherent; since $Z\rightarrow S$ is proper, the support of $\calF$ is proper over $T$ by \cite[0CYN]{stacks-project}.}.
By \cite[08WC, 08KD, 0CMY, 0DLY]{stacks-project}, $\calCoh(Z/S)$ is an algebraic stack with affine diagonal, locally of finite presentation over $S$. 
There is a natural forgetful $1$-morphism over $(\Sch/S)_\fppf$
\[
\calM_{\dR}((Z,D)/S)\rightarrow\calCoh(Z/S).
\]
\end{defn}

\begin{prop}\label{prop:MdR algebraic}
The category fibered in groupoids $\calM_\dR((Z,D)/S)$ is an algebraic stack with affine diagonal and is locally of finite presentation over $S$.
\end{prop}

\begin{proof}
Let $\Lambda=\bigcup_{n\geq 0}\Lambda_n$ denote the sheaf of rings of log PD-differential operators on $Z$ relative to $S$, where $\Lambda_n\subset \Lambda$ is the left $\calO_Z$-submodule of log PD-differential operators of order at most $n$ (cf.~\cite[Def.~4.4]{BO} and \cite[Def.~1.1.3]{Ogus-F-crystals}\footnote{Strictly speaking, the sheaves used therein are Zariski quasi-coherent sheaves, whereas we work with quasi-coherent sheaves on the \'etale site. This does not affect the argument.}). Then, $\Lambda$ satisfies \cite[2.1.1-2.1.6]{Simpson1} (i.e., $\Lambda$ is a sheaf of rings of differential operators in the sense of Simpson): when $D=\emptyset$, this follows from \cite[Cor.~I.4.5.3(i) ($m=1$), Prop.~II.4.2.5, Cor.~II.4.2.6]{Berthelot-book}; when $D\neq\emptyset$, the same argument works by using \cite[Prop.~6.5]{Kato-log} in place of  of \cite[Cor.~I.4.5.3(i)]{Berthelot-book} (see also \cite[\S1.1]{Ogus-F-crystals}).  
Moreover, giving a left $\Lambda$-module structure on a quasi-coherent $\calO_{Z}$-module is the same as equipping it with an integrable $T$-connection: when $D=\emptyset$, this follows from \cite[Cor.~II.4.2.12(i)]{Berthelot-book}; when $D\neq\emptyset$, this follows from \cite[Thm.~II.4.2.11]{Berthelot-book} by using a more general notion of PD-adic groupoids.\footnote{In the special case when $D=\emptyset$ and $m\calO_S=0$ for some $m>0$, these assertions can be also deduced from \cite[Prop.~3.32, Cor.~4.7, Thm.~4.8]{BO}.}  

Hence, $\cM_{\dR}((Z,D)/S)$ is equivalent to the category fibered in groupoids 
\[
\Lambda\calCoh(Z/S)\rightarrow(\Sch/S)_{\fppf}
\]
whose fiber over $T\rightarrow S$ is the groupoid of $T$-flat $\calO_{Z_T}$-modules $\cE$ of finite presentation equipped with the structure of left $\Lambda\rvert_{Z_{T}}$-modules, and the forgetful functor $\cM_{\dR}((Z,D)/S)\rightarrow\calCoh(Z/S)$ under this equivalence corresponds to the forgetful functor $\Lambda\calCoh(Z/S)\rightarrow\calCoh(Z/S)$. Recall that $\calCoh(Z/S)$ is an algebraic stack with affine diagonal, locally of finite presentation over $S$. 
Since $S$ is Noetherian, one can apply the discussions of \cite[\S2]{HLHJ-affineGrass} (see Remark 2.18 of \emph{op.~cit.}). In particular, the proof of Proposition 2.23 of \emph{op.~cit.}, which does not use the purity of the underlying sheaf, verbatim carries over and shows that the forgetful functor $\Lambda\calCoh(Z/S)\rightarrow\calCoh(Z/S)$ is representable by schemes, affine and of finite presentation. This implies that $\cM_{\dR}((Z,D)/S)\cong\Lambda\calCoh(Z/S)$ is an algebraic stack with affine diagonal, locally of finite presentation over $S$.
\end{proof}

We will also work on various open substacks of $\cM_{\dR}$.

\begin{defn}\label{def:moduli stack of vector bundles with integrable connections}
The \emph{moduli stack of vector bundles with integrable connections}
\[
\cM_{\dR,\Bun}((Z,D)/S)\rightarrow(\Sch/S)_{\fppf}
\]
is defined to be the strictly full subcategory of $\cM_{\dR}((Z,D)/S)$ whose objects over $T\rightarrow S$ are the pairs $(\cE,\nabla)\in\Ob(\cM_{\dR}((Z,D)/S)_{T})$ such that $\cE$ is a locally free $\cO_{Z_{T}}$-module. By \cite[0CZR]{stacks-project}, this condition is the same as requiring $\cE\rvert_{Z_{t}}$ to be a locally free $\cO_{Z_{t}}$-module for every $t\in T$.

For $r\in\bN$, let $\cM_{\dR,\Bun_{r}}((Z,D)/S)\rar(\Sch/S)_{\fppf}$ be the strictly full subcategory of objects of $\cM_{\dR,\Bun}((Z,D)/S)$ whose underlying vector bundles are of rank $r$.
\end{defn}

The formations of the above moduli categories  commute with base change. 

\begin{prop}\label{prop:MdRBun}
The natural functors
\[
\cM_{\dR,\Bun_{r}}((Z,D)/S)\rightarrow\cM_{\dR,\Bun}((Z,D)/S)\rightarrow\cM_{\dR}((Z,D)/S)
\]
are open immersions. In particular, $\cM_{\dR,\Bun}((Z,D)/S)$ and $\cM_{\dR,\Bun_{r}}((Z,D)/S)$ are algebraic stacks with affine diagonal and is locally of finite presentation over $S$. 
We have a decomposition 
\[
\cM_{\dR,\Bun}((Z,D)/S)=\coprod_{r\in\bN}\cM_{\dR,\Bun_{r}}((Z,D)/S).
\]
\end{prop}

\begin{proof}
It suffices to prove the analogous statements for $\calCoh(Z/S)$. More precisely, we define the \emph{moduli of vector bundles}
\[
\Bun(Z/S)\rightarrow(\Sch/S)_{\fppf}
\]
to be the strictly full subcategory of $\calCoh(Z/S)$ whose objects over $T\rightarrow S$ are the locally free $\cO_{Z_{T}}$-modules $\calE$. We first show that for any $S$-scheme $T$ and $1$-morphism $f\colon(\Sch/T)\rar\calCoh(Z/S)$ over $\Sch/S$, the $2$-fiber product $(\Sch/T)\times_{\calCoh(Z/S)}\Bun(Z/S)$ is represented by an open subscheme of $T$. Let $\pi_{T}\colon Z_{T}\rightarrow T$ be the projection and $\cE$ the sheaf on $Z_T$ corresponding to $f$. By \cite[0CZR]{stacks-project}, the locus $U$ consisting of points $z\in Z_{T}$ such that $\cE_{z}$ is a free $\calO_{Z_T,z}$-module is open in $Z_{T}$, and the formation commutes with arbitrary base change $T'\rightarrow T$. As $\pi_{T}$ is closed, $V\coloneqq T\smallsetminus \pi_{T}(Z_{T}\smallsetminus U)$ is open in $T$. Since $\pi_{T}$ is universally closed, the same argument applied to arbitrary base change of $T$ shows that $(\Sch/T)\times_{\calCoh(Z/S)}\Bun(Z/S)$ is represented by the open subscheme $V\subset T$. Hence $\Bun(Z/S)\rightarrow\calCoh(Z/S)$ is an open immersion.
As the rank of a vector bundle is locally constant, we have a decomposition 
 $\Bun(Z/S)=\coprod_{r\in\bN}\Bun_{r}(Z/S)$ where $\Bun_{r}(Z/S)$ classifies the rank $r$ locally free sheaves, and $\Bun_{r}(Z/S)\rightarrow\Bun(Z/S)$ is an open immersion.
\end{proof}

\begin{rem}\label{rem:connection is locally free over char 0}
Assume either $D=\emptyset$ or  $\dim(Z/S)=1$. Then the underlying module of an integrable $T$-connection is locally free when $T$ is the spectrum of a field of characteristic $0$ (cf. \cite[Prop.~8.8]{Katz_Nilpotent}). It follows from the fiberwise criterion of flatness that the same holds for any $\Q$-scheme $T$. Thus, if $S$ is a $\bQ$-scheme and either $D=\emptyset$ or $\dim(Z/S)=1$, we have $\cM_{\dR,\Bun}=\cM_{\dR}$. In general, $\cM_{\dR,\Bun}\hookrightarrow\cM_{\dR}$ is not an equivalence. 
\end{rem}

\begin{rem}\label{rem:morphism from algebraic space to MdR}
Let $T$ be an algebraic space over $S$. Then 
\[
\Mor_{(\Sch/S)_\fppf}(\calS_T,\calM_{\dR}((Z,D)/S))
\]
is still functorially identified with the groupoid of integrable $T$-connections, namely, pairs $(\calE,\nabla)$ where $\calE$ is a $T$-flat $\calO_{Z_T}$-modules of finite presentation, and $\nabla\colon \calE\rightarrow\calE\otimes_{\cO_{Z_{T}}} \omega^{1}_{Z_{T}/T}$ is an integrable $T$-connection. This follows from the descent result of integrable connections: when $D=\emptyset$, it is  \cite[Cor.~2.2.15]{Olsson-crystalline}, and the same proof works in the general case, using the first log infinitesimal neighborhood of the diagonal (cf.~Lemma~\ref{lem:connection and P1-linear isom}). The same holds for $\cM_{\dR,\Bun}$ and $\cM_{\dR,\Bun_{r}}$.
\end{rem}

\subsection{Nilpotent locus in characteristic \texorpdfstring{$p$}{p}}\label{sec:nilpotent locus}

In this subsection, fix a perfect field $k$ of characteristic $p>0$ and assume $S=\Spec k$. Write $\cM_{\dR,k}$ for $\cM_{\dR}((Z,D)/S)$.

For a $k$-scheme $T$, write $\Fr_{T}\colon T\rightarrow T$ for the absolute (i.e., $p$th power) Frobenius, and let $\Der((Z_{T},D_{T})/T)=\underline{\Hom}_{\calO_{Z_T}}(\omega_{Z_T/T}^1,\calO_{Z_T})$  denote the sheaf of log derivations with values in $\cO_{Z_{T}}$ relative to $T$. 

Let $(\cE,\nabla)\in\Ob\left(\left(\cM_{\dR,k}\right)_{T}\right)$ be an integrable $T$-connection.

\begin{defn}\label{def:p-curvature}
We define the \emph{$p$-curvature} of $(\cE,\nabla)$
\[
\psi(\cE,\nabla)\in\Hom_{\cO_{Z_{T}}}\left(\Der((Z_{T},D_{T})/T),\Fr_{Z_{T},*}\underline{\End}_{Z_{T}}(\cE)\right),
\]
by the formula
\[
\psi(\cE,\nabla)(\partial)\coloneqq (\nabla(\partial))^{p}-\nabla(\partial^{(p)})
\]
for a local log derivation $\partial$ of $\cO_{Z_{T}}$ into $\cO_{Z_{T}}$ relative to $T$ (see \cite[(5.0), (5.2), (6.0)]{Katz_Nilpotent}\footnote{Here we only consider the log structure coming from a normal crossings divisor. See also \cite[\S1.1, 1.2]{Ogus-F-crystals} for the case of a more general log structure.}).
We say that $(\cE,\nabla)$ is \emph{nilpotent} if, for every local log derivation $\partial$ and every local section $e$ of $\cE$ on a quasi-compact open, $\left(\psi(\cE,\nabla)(\partial)\right)^{n}(e)=0$ for some $n$ (see \cite[Rem.~1.2.2]{Ogus-F-crystals} and \cite[Thm.~6.2]{Kato-log}). We say that $(\cE,\nabla)$ has  \emph{exponent of nilpotence} $\le N$ if  $(\psi(\cE,\nabla)(\partial))^{N}=0$ for every local log derivation $\partial$.
Since $\underline{\End}_{\calO_{Z_T}}(\calE)$ is a quasi-coherent sheaf on $Z_T$, one immediately sees that $(\calE,\nabla)$ being nilpotent is fppf local on $T$.
\end{defn}

\begin{lem}\label{lem:any composition of p-curvature is zero}
The integrable $T$-connection $(\cE,\nabla)$ is nilpotent and has exponent of nilpotence $\le N$ if and only if, for $m\ge (N-1)\dim(Z)+1$, any local log derivations $\partial_{1},\hdots, \partial_{m}$ and any local section $e$ of $\cE$, we have 
\[
(\psi(\cE,\nabla)(\partial_{1}))(\psi(\cE,\nabla)(\partial_{2}))\cdots(\psi(\cE,\nabla)(\partial_{m}))(e)=0.
\]
\end{lem}

\begin{proof}
The if part is obvious. For the only-if part, we know from the smoothness assumption that $\Der((Z_{T},D_{T})/T)$ is a locally free $\cO_{Z_T}$-module of rank $d$: in fact, it is locally generated by $x_{1}\frac{\partial}{\partial x_{1}},\ldots,x_{c}\frac{\partial}{\partial x_{c}},\frac{\partial}{\partial x_{c+1}},\ldots,\frac{\partial}{\partial x_{d}}$ if $Z$ has local coordinates $x_{1},\ldots,x_{d}$ and if $D$ is defined by $x_{1}\cdots x_{c}=0$.
So the proof of \cite[Prop.~5.2, Cor.~5.5]{Katz_Nilpotent} works verbatim.
\end{proof}

\begin{defn}[Nilpotent locus]
Let $\cM_{\dR,k}^{\nilp}\rightarrow(\Sch/\Spec k)_{\fppf}$ be the strictly full subcategory of $\cM_{\dR,k}$ such that the objects over $T\rightarrow\Spec k$ consist of the \emph{nilpotent} integrable $T$-connections (see Definition~\ref{def:p-curvature}). It is immediate to see that $\cM_{\dR,k}^{\nilp}$ is a substack of $\cM_{\dR,k}$. 
We let $(\cM_{\dR,k}^{\nilp})_\red\rightarrow (\Sch/\Spec k)$ denote the underlying reduced substack of $\cM_{\dR,k}^{\nilp}$ in the sense of \cite[Def.~3.27]{EmertonStack}.

Similarly, let $\cM_{\dR,k}^{\psi=0}\rightarrow(\Sch/\Spec k)_{\fppf}$ denote the substack of $\cM_{\dR,k}$ such that the objects over $T\rightarrow\Spec k$ consist of the integrable $T$-connections with vanishing $p$-curvature.
\end{defn}

Let $\Fr_{Z_T/T}\colon Z_T\rar Z_T^{(p)}\coloneqq Z_T\times_{T,\Fr_T}T$ denote the relative Frobenius map and set $Z^{(p)}\coloneqq Z^{(p)}_k$. The Cartier descent describes $\cM_{\dR,k}^{\psi=0}$ in terms of $\calCoh(Z^{(p)}/k)$.

\begin{prop}[Cartier descent]\label{prop:Cartier descent for general MdR}
Assume $D=\emptyset$. Then the $1$-morphism $C_{Z/k}\colon\calCoh(Z^{(p)}/k)\rightarrow \calM_{\dR,k}$ sending $\calF\in \Ob(\calCoh(Z^{(p)}/k)_T)$ over a $k$-scheme $T$ to $(\Fr_{Z_T/T}^{\ast}\calF,\nabla_{\can})\in \Ob((\calM_{\dR,k})_T)$, where $\nabla_{\can}$ is the canonical connection in \cite[Thm.~5.1]{Katz_Nilpotent}, induces an equivalence of stacks over $k$
\[
C_{Z/k}\colon\calCoh(Z^{(p)}/k)\xrightarrow{\cong} \calM_{\dR,k}^{\psi=0}.
\]
In particular, $\cM_{\dR,k}^{\psi=0}$ is an algebraic stack.
\end{prop}

\begin{proof}
This is exactly the Cartier descent theorem \cite[Thm.~5.1]{Katz_Nilpotent}.
\end{proof}

\begin{rem}
When $D\neq\emptyset$, the Cartier descent theorem also involves conditions on the residue; see \cite[Thm.~1.3.4]{Ogus-F-crystals}.
\end{rem}

In contrast, the stack $\calM_{\dR,k}^\nilp$ is not algebraic, and we will mainly study $(\cM_{\dR,k}^{\nilp})_{\red}$ due to the following lemma.

\begin{lem}\label{lem:nilpotence can be checked after restricting to reduction}
An integrable $T$-connection $(\calE,\nabla)\in\Ob\left((\cM_{\dR,k})_{T}\right)$ is nilpotent if and only if its pullback to the reduction $T_{\red}$ of $T$ is nilpotent. 
\end{lem}

\begin{proof}
This is easily reduced to the following lemma in commutative algebra.
\end{proof}

\begin{lem}\label{lem:commutative algebra about nilpotent endomorphism}
Let $A$ be a ring and write $I$ for the nilradical.  An $A$-linear endomorphism $\psi\in \End(M)$ on a finitely generated $A$-module $M$ is nilpotent if and only if so is $\psi\otimes_A A/I\in \End(M/IM)$.
\end{lem}

\begin{proof}
The only-if part is obvious. Conversely, assume that $\psi\otimes_A A/I$ is nilpotent.
Take $A$-generators $m_1,\ldots,m_n$ of $M$ and $r\geq 1$ such that $(\psi\otimes_A A/I)^r(m_i\bmod{I})=0$ for every $1\leq i\leq n$. Write $\psi^r(m_i)=\sum_{1\leq j\leq n}a_{ij}m_j$ with $a_{ij}\in I$ and let $I'$ be the ideal of $A$ generated by $a_{ij}$ ($1\leq i,j\leq n$). Since $I'$ is a finitely generated subideal of the nilradical $I$, there exists $s\geq 1$ such that $(I')^s=0$. Then $\psi^{rs}=0$.
\end{proof}

\begin{rem}\label{rem:commutative algebra about nilpotent endomorphism}
A simpler variant of the proof of Lemma~\ref{lem:commutative algebra about nilpotent endomorphism} shows the following: Let $A$ be a reduced ring, $M$ a finitely generated $A$-module, and $\psi\in \End_A(M)$.
If $r\geq 1$ satisfies $(\psi\otimes_A \Frac(A/\fkp))^r=0$ for every $\fkp\in\Spec A$, then $\psi^r=0$.
\end{rem}

For a scheme $T$, let $(\Sch/T)_\red$ denote the full subcategory of $(\Sch/T)$ consisting of reduced $T$-schemes. 
Recall from Definition~\ref{defn:reduced substack} that the underlying reduced substack $(\cM_{\dR,k}^{\nilp})_{\red}\rightarrow (\Sch/\Spec k)_\fppf$ is defined to be the stack generated by $\cM_{\dR,k}^{\nilp}|_{(\Sch/\Spec k)_\red}\rightarrow(\Sch/\Spec k)_{\red}$ . 
We now prove  that $(\cM_{\dR,k}^{\nilp})_{\red}$ is a (reduced) closed algebraic substack of $\cM_{\dR,k}$. Combined with this, Lemma~\ref{lem:nilpotence can be checked after restricting to reduction} implies that  $\cM_{\dR,k}^{\nilp}$ is the  formal completion of $\cM_{\dR,k}$ along $(\cM_{\dR,k}^{\nilp})_{\red}$. This crucially relies on Proposition~\ref{prop:bounding exponent of nilpotence}, which we will prove later in this subsection.

\begin{thm}\label{thm:reduced closed substack of nilpotent connections} 
The $1$-morphism $(\cM_{\dR,k}^{\nilp})_\red\rightarrow\cM_{\dR,k}$ is a closed immersion of finite presentation. In particular, $(\cM_{\dR,k}^{\nilp})_{\red}$ is a reduced closed algebraic substack of $\cM_{\dR,k}$.
\end{thm}

\begin{proof}
Take any $k$-scheme $\pi\colon T\rightarrow\Spec k$ and $(\calE,\nabla)\in \Ob(((\cM_{\dR,k})_T)$, and form the $2$-fiber product $\calT^\nilp\coloneqq (\Sch/T)\times_{\cM_{\dR,k}}\cM_{\dR,k}^{\nilp}$. We will show that the induced  category fibered in groupoids $\calT^\nilp|_{(\Sch/T)_\red}\rightarrow (\Sch/T)_\red$ is represented by a unique reduced closed subscheme of finite presentation $T^\nilp_\red\hookrightarrow T$. Once this is proved, one obtains a similar result for any $1$-morphism $\calS_F\rightarrow \cM_{\dR,k}$ from an algebraic space $F$ over $k$. By considering a presentation of $\calM_{\dR,k}$, the latter in turn yields a reduced closed algebraic substack $(\cM_{\dR,k}^{\nilp})_\red'$ of $\cM_{\dR,k}$ whose restriction to $(\Sch/\Spec k)_\red$ consists of nilpotent integrable connections. Arguing as in \cite[Rem.~3.29]{EmertonStack}, we can conclude $(\cM_{\dR,k}^{\nilp})_\red'=(\cM_{\dR,k}^{\nilp})_\red$.

Since $\cM_{\dR,k}$ is locally of finite presentation over $k$, we may assume that $T=\Spec B$ is connected, affine and Noetherian. Moreover, we may also assume that $T$ is reduced. 
In the following, we simply write $\Der(Z_{T}/T)$ for $\Der((Z_{T},D_{T})/T)$ and so on; for a $k$-morphism $g\colon U'\rightarrow U$, let $g_Z$ denote the base change $Z_{U'}\rightarrow Z_{U}$ (except for $\Fr$). Observe the following identifications
\begin{align*}
\Hom_{\cO_{Z_{T}}}(\Der(Z_{T}/T),\Fr_{Z_{T},*}\underline{\mathrm{\End}}_{Z_{T}}(\cE)) 
& =
\Hom_{\cO_{Z_{T}}}(\Fr_{Z_{T}}^{\ast}\Der(Z_{T}/T),\underline{\mathrm{\End}}_{Z_{T}}(\cE))
\\
& =
\Hom_{\cO_{Z_{T}}}(\Fr_{Z_{T}}^{\ast}\pi_Z^{\ast}\Der(Z/k),\underline{\mathrm{\End}}_{Z_{T}}(\cE))
\\
& =
\Hom_{\cO_{Z_{T}}}(\pi_Z^{\ast}\Fr_{Z}^{\ast}\Der(Z/k),\underline{\mathrm{End}}_{Z_{T}}(\cE))
\\
& =
\Hom_{\cO_{Z_{T}}}(\pi_Z^{\ast}\Fr_{Z}^{\ast}\Der(Z/k)\otimes_{\cO_{Z_{T}}}\cE,\cE).
\end{align*}
Hence we regard $\psi(\cE,\nabla)\in\Hom_{\cO_{Z_{T}}}(\pi_Z^{\ast}\Fr_{Z}^{\ast}\Der(Z/k)\otimes_{\cO_{Z_{T}}}\cE,\cE)$. Consider the functor $\cH\colon(\Sch/T)^{\mathrm{op}}\rightarrow\Set$ sending a $T$-scheme $f\colon T'\rightarrow T$ to
\[
\cH(T')\coloneqq\Hom_{\cO_{Z_{T'}}}(f_Z^{\ast}\pi_Z^{\ast}\Fr_{Z}^{\ast}\Der(Z/k)\otimes_{\cO_{Z_{T'}}}f_{Z}^{\ast}\cE,f_Z^{\ast}\cE).
\]
Note that $\cE$ is a coherent $\cO_{Z_{T}}$-module, so its support is a closed subset of $Z_{T}$. As $Z_{T}\rar T$ is proper, by \cite[0CYN]{stacks-project}, the scheme-theoretic support of $\cE$ is proper over $T$. As $\cE$ is $T$-flat, we can apply \cite[08K6]{stacks-project} and conclude that $\cH$ is representable by an affine scheme of finite presentation over $T$, which we still denote by $\calH$. Since the formation of $p$-curvature commutes with base change, we may regard it as a $T$-scheme morphism $\psi(\cE,\nabla)\colon T\rightarrow \cH$. 

Let $d=\dim(Z/k)$ and $N$ be the universal bound on the exponent of nilpotence for $\cE$ as in Proposition~\ref{prop:bounding exponent of nilpotence} below. Consider the functor $\cH^{d(N-1)+1}\colon (\Sch/T)^{\mathrm{op}}\rightarrow\Set$ defined as 
\[
\cH^{d(N-1)+1}(T')\coloneqq\Hom_{\cO_{Z_{T'}}}((f_Z^{\ast}\pi_Z^{\ast}\Fr_{Z}^{\ast}\Der(Z/k))^{\otimes_{\cO_{Z_{T'}}} d(N-1)+1}\otimes_{\cO_{Z_{T'}}}f_{Z}^{\ast}\cE,f_Z^{\ast}\cE)
\]
for a $T$-scheme $f\colon T'\rightarrow T$. By the same reason as above, we can apply \cite[08K6]{stacks-project} and conclude that $\cH^{d(N-1)+1}$ is representable by an affine scheme of finite presentation over $T$. For each $T$-scheme $f\colon T'\rightarrow T$ and $\phi\in \calH(T')$, we define the $\cO_{Z_{T'}}$-linear map 
\[
c(T')(\phi)\colon (f_Z^{\ast}\pi_Z^{\ast}\Fr_{Z}^{\ast}\Der(Z/k))^{\otimes_{\cO_{Z_{T'}}} d(N-1)+1}\rar f_Z^{\ast}\underline{\mathrm{End}}_{Z_{T}}(\cE)
\]
by the formula 
\[
c(T')(\phi)(e_{1}\otimes \cdots \otimes e_{d(N-1)+1})\coloneqq \phi(e_{1})\circ\cdots\circ\phi(e_{d(N-1)+1})
\]
for local sections $e_{1},\ldots,e_{d(N-1)+1}$ of $f_Z^{\ast}\pi_Z^{\ast}\Fr_{Z}^{\ast}\Der(Z/k)$. Here $\phi(e_i)$'s are local sections of $f_Z^{\ast}\underline{\End}_{Z_{T}}(\cE)$ via the identification
$
\Hom_{\cO_{Z_{T'}}}(f_Z^{\ast}\pi_Z^{\ast}\Fr_{Z}^{\ast}\Der(Z/k)\otimes_{\cO_{Z_{T'}}}f_{Z}^{\ast}\cE,f_Z^{\ast}\cE)=\Hom_{\cO_{Z_{T'}}}(f_Z^{\ast}\pi_Z^{\ast}\Fr_{Z}^{\ast}\Der(Z/k),f_Z^{\ast}\underline{\mathrm{End}}_{Z_{T}}(\cE))
$, 
so we can compose them. This defines a $T$-morphism $c\colon\cH\rightarrow\cH^{d(N-1)+1}$ via the similar identification.

Set $\psi'\coloneqq c\circ\psi(\cE,\nabla)$, and let $0\colon T\rightarrow \cH^{d(N-1)+1}$ be the zero section; both $\psi'$ and $0$ are sections of the structure morphism $\cH^{d(N-1)+1}\rightarrow T$. Define $\widetilde{T}$ as the fiber product sitting in the following Cartesian square
\[
\xymatrix{
\widetilde{T}\ar[r]\ar[d]
&
\cH^{d(N-1)+1}\ar[d]^-{\Delta}
\\
T\ar[r]^-{(0,\psi')}
&
\cH^{d(N-1)+1}\times_{T}\cH^{d(N-1)+1}.
}
\]
Recall that $\cH^{d(N-1)+1}$ is affine of finite presentation over $T$. Hence the diagonal morphism $\Delta$ is a closed immersion of finite presentation and so is the base change $\widetilde{T}\rightarrow T$. By construction, for every $T$-scheme $f\colon T'\rightarrow T$,  we have
\[
\widetilde{T}(T')=\begin{cases}
 & 
\text{ if $(\psi(\cE_{Z_{T'}},\nabla_{Z_{T'}})(\partial_{1}))\circ\cdots\circ(\psi(\cE_{Z_{T'}},\nabla_{Z_{T'}})(\partial_{d(N-1)+1}))(e)=0$  }
\\ 
\{f\} & \text{ for any local log derivations $\partial_{1},\hdots,\partial_{d(N-1)+1}$ on $Z_{T'}$}\\&\text{ and for any local section $e$ of $\cE_{Z_{T'}}$,}\\
\;\emptyset & \text{ otherwise.}\end{cases}
\]
By Lemma~\ref{lem:any composition of p-curvature is zero}, $\widetilde{T}(T')=\lbrace f\rbrace$ if and only if $(\cE_{Z_{T'}},\nabla_{Z_{T'}})$ is nilpotent and has exponent of nilpotence $\le N$. If $T'$ is \emph{reduced}, then by Proposition~\ref{prop:bounding exponent of nilpotence}, any nilpotent endomorphism on $\cE_{Z_{T'}}$ has exponent of nilpotence $\le N$. Therefore, if $T'$ is reduced and $(\cE_{Z_{T'}},\nabla_{Z_{T'}})$ is nilpotent, then $(\cE_{Z_{T'}},\nabla_{Z_{T'}})$ has exponent of nilpotence $\le N$. This implies that $\widetilde{T}$ represents $\cT^{\nilp}\rvert_{(\Sch/T)_{\red}}$. Note that $\widetilde{T}$ is Noetherian, as $T$ is Noetherian. Thus, $\widetilde{T}_{\red}$ is a reduced closed subscheme of finite presentation of $T$ that represents $\cT^{\nilp}$. The uniqueness of a reduced closed subscheme of $T$ representing $\cT^{\nilp}$ is clear.
\end{proof}

Let $\widehat{\cM_{\dR,k}}_{|\lvert(\cM_{\dR,k}^{\nilp})_{\red}\rvert}\rar(\Sch/\Spec k)_{\fppf}$ denote the completion of $\cM_{\dR,k}$ along $\lvert(\cM_{\dR,k}^{\nilp})_{\red}\rvert$ (see Definition~\ref{def:formal completion}). 

\begin{cor}\label{cor:nilpotent locus is formal completion}
The nilpotent locus $\cM_{\dR,k}^{\nilp}$ is equal to $\widehat{\cM_{\dR,k}}_{|\lvert(\cM_{\dR,k}^{\nilp})_{\red}\rvert}$ as a substack of $\cM_{\dR,k}\rar(\Sch/\Spec k)_{\fppf}$. In particular, $\cM_{\dR,k}^{\nilp}$ is a locally Noetherian formal algebraic stack  over $\Spec k$ whose diagonal is representable (by schemes) and affine.
\end{cor}

\begin{proof}
Choose a $k$-scheme $T$ and write $i\colon T_{\red}\rar T$ for the natural closed immersion.
By Lemma~\ref{lem:nilpotence can be checked after restricting to reduction}, if $f\colon T\rar \cM_{\dR,k}$ is an object of $(\cM_{\dR,k}^{\nilp})_{T}$, the composition $T_{\red}\xrightarrow{i} T\xrightarrow{f}\cM_{\dR,k}$ factors through $(\cM_{\dR,k}^{\nilp})_{\red}$ (this makes sense as $(\cM_{\dR,k}^{\nilp})_{\red}$ is a strictly full subcategory of $\cM_{\dR,k}$). Since $(\cM_{\dR,k}^{\nilp})_{\red}$ is a reduced closed algebraic substack of $\cM_{\dR,k}$ by Theorem~\ref{thm:reduced closed substack of nilpotent connections}), it follows that $f(|T_{\red}|)=f(|T|)\subset \lvert(\cM_{\dR,k}^{\nilp})_{\red}\rvert$. 

Conversely, if $f\colon T\rar\cM_{\dR,k}$ is an object of $(\cM_{\dR,k})_{T}$ such that $f(|T|)\subset\lvert(\cM_{\dR,k}^{\nilp})_{\red}\rvert$, then $f(|T_{\red}|)\subset\lvert(\cM_{\dR,k}^{\nilp})_{\red}\rvert$. By \cite[050B]{stacks-project}, this implies that the composition $T_{\red}\xrightarrow{i}T\xrightarrow{f}\cM_{\dR,k}$ factors through the closed algebraic substack $(\cM_{\dR,k}^{\nilp})_{\red}$. Hence the corresponding integrable $T$-connection $(\cE,\nabla)$ satisfies the property that $(\cE_{Z_{T_{\red}}},\nabla_{Z_{T_{\red}}})$ is nilpotent, which by Lemma~\ref{lem:nilpotence can be checked after restricting to reduction} means that $f$ is an object of $(\cM_{\dR,k}^{\nilp})_{T}$.
\end{proof}

The goal of the rest of the subsection is to prove the following, which was a crucial ingredient in the proof of Theorem~\ref{thm:reduced closed substack of nilpotent connections}.

\begin{prop}[Universal bound on exponent of nilpotence]\label{prop:bounding exponent of nilpotence}
Let $S$ be a reduced Noetherian $k$-scheme, and let $\cE$ be a coherent $\cO_{Z_{S}}$-module which is $S$-flat. Then, there exists a positive integer $N>0$ such that, for any \emph{reduced} $S$-scheme $S'\rar S$, if an $\cO_{Z_{S'}}$-module endomorphism $\alpha\colon\cE_{Z_{S'}}\rar\cE_{Z_{S'}}$ is nilpotent (i.e., $\alpha^{n}=0$ for some $n>0$), then $\alpha^{N}=0$. Such an $N$ will be called a \emph{universal bound on the exponent of nilpotence} for $\cE$.
\end{prop}

Note that if $\calE$ is locally free of rank $r$, one may simply take $N=r$.

\begin{rem}
We cannot expect a universal bound as in Proposition~\ref{prop:bounding exponent of nilpotence} if we allow \emph{non-reduced} $S$-schemes $S'\rar S$. For example, if $Z=S=\Spec k$, $\calE=\calO_Z$, and $S'=\Spec k[\epsilon]/(\epsilon^{M+1})$, then the multiplication-by-$\epsilon$ endomorphism  $\alpha$ on any nonzero coherent sheaf clearly satisfies $\alpha^{M}\ne0$ and $\alpha^{M+1}=0$. 
\end{rem}

We prove this via a series of lemmas. 
We thank Andres Fernandez Herrero for suggesting the proof of Lemma~\ref{lem:controlling torsion filtration} and how to deduce the general case from the projective case (Corollary~\ref{cor:bounding exponent of nilpotence, projective}), and for allowing us to include them in the paper. 

\begin{rem}
If $\dim_{k}Z=1$, one may alternatively prove Proposition~\ref{prop:bounding exponent of nilpotence} using the Fitting ideals and the structure theory of modules over a Dedekind domain.
\end{rem}

We first prove Proposition~\ref{prop:bounding exponent of nilpotence} under an extra assumption that $Z$ is \emph{projective} over $k$. To achieve this, we first recall the standard concepts in the moduli theory: see \cite[\S1.1, 1.2]{HL} for example.

\begin{defn}\label{defn:pure}
Let $X$ be a Noetherian scheme and $\calF$ a coherent sheaf on $X$. We say that $\calF$ is \emph{pure} of dimension $d$ if $\cF$ is nonzero and, for every nonzero coherent subsheaf $\cF'\subset\cF$, $\dim\Supp(\cF')=d$. Note that this implies $\dim\Supp(\cF)=d$.

In general, $\calF$ admits the unique filtration (\emph{torsion filtration})
\[
0\subset T_{0}(\cF)\subset T_{1}(\cF)\subset\cdots\subset  T_{\dim\Supp(\cF)}(\cF)=\cF,
\]
where $T_{i}(\cF)\subset\cF$ is the maximal subsheaf whose support is of dimension $\le i$ (with convention $T_{-1}(\cF)=0$). By definition, $T_{i}(\cF)/T_{i-1}(\cF)$ is either $0$ or pure of dimension $i$.
\end{defn}

\begin{defn}\label{defn:Hilbert polynomial}
Let $X$ be a projective scheme over a field, and fix an ample line bundle $\cO(1)$ on $X$. For a nonzero coherent sheaf $\cF$ over $X$, the \emph{Hilbert polynomial} $P_{\cF}(t)\in\bQ[t]$ is defined as $P_{\cF}(n)\coloneqq\chi(\cF(n))$ for $n\in\bZ$. The Hilbert polynomial is of the form $P_{\cF}(t)=\sum_{k=0}^{\dim\Supp(\cF)}\frac{a_{k}(\cF)}{k!}t^{k}$ for $a_{k}(\cF)\in\bZ$. We define the \emph{multiplicity} of $\cF$ as $\mult(\cF)\coloneqq a_{\dim\Supp(\cF)}(\cF)$, which is in fact always a positive integer. The \emph{normalized Hilbert polynomial} $p_\calF(t)$ is defined as $p_\calF(t)\coloneqq P_\calF(t)/\mult(\calF)$. We also define both the (normalized) Hilbert polynomial and the multiplicity of the zero sheaf to be $0$.
\end{defn}

\begin{rem}\label{rem:rank}
Keep the setup as above. The \emph{rank} of a coherent sheaf $\cF$ on $X$ is defined as 
\[
\rank(\cF)\coloneqq
\begin{cases}
\frac{\mult(\cF)}{\mult(\cO_{X})} 
&
\text{if $\dim\Supp(\cF)=\dim X$,}
\\
0
&
\text{otherwise.}
\end{cases}
\]
This definition of rank coincides with the usual notion of rank for locally free sheaves. When $X$ is integral, this coincides with the generic rank. However, in general, $\rank(\cF)$ may not even be an integer.
\end{rem}

The following lemma is the starting point of the proof of Proposition~\ref{prop:bounding exponent of nilpotence}.

\begin{lem}\label{lem:bounding exponent of nilpotence, pure}
Let $X$ be a projective scheme over a field and fix an ample line bundle $\cO(1)$ on $X$. Let $\cF$ be a pure sheaf on $X$. If $\alpha\colon\cF\rar\cF$ is a nilpotent $\cO_{X}$-linear endomorphism, then $\alpha^{\mult(\cF)}=0$.
\end{lem}

\begin{proof}
Consider the decreasing filtration $\cF=\cF_{0}\supset\cF_{1}\supset\cF_{2}\supset\cdots$ of coherent subsheaves given by $\cF_{i}\coloneqq\Image(\alpha^{i}\colon\cF\rar\cF)$. 
Observe that if $\cF_{i}=\cF_{i+1}$, then $\cF_{n}=\cF_{i}$ for any $n\ge i$ as $\cF_{i+1}=\Image(\alpha\rvert_{\cF_{i}}\colon\cF_{i}\rar\cF)$.
Since $\alpha$ is nilpotent, there exists $N\geq 1$ such that $\alpha^N=0$ and $\cF=\cF_{0}\supsetneq\cF_{1}\supsetneq \cF_{2}\supsetneq\cdots\supsetneq\cF_{N-1}\supsetneq\cF_{N}=0$. Therefore, for $0\le i\le N-1$, $\cF_{i}$ is pure of dimension $d\coloneqq\dim\Supp(\cF)$. In particular, $\mult(\cF_{i})=a_{d}(\cF_{i})>0$ for $0\le i\le N-1$. 

We claim $\mult(\cF_{i})>\mult(\cF_{i+1})$ for $0\le i\le N-1$, which would imply $\mult(\cF)\ge N$ and thus $\alpha^{\mult(\cF)}=0$. Set $\cG_{i}=\Ker(\alpha\rvert_{\cF_{i}}\colon\cF_{i}\rar\cF)$. Then the short exact sequence $0\rar \cG_{i}\rar\cF_{i}\rar\cF_{i+1}\rar0$ gives $P_{\cF_{i}}(t)=P_{\cG_{i}}(t)+P_{\cF_{i+1}}(t)$. In particular $a_{d}(\cF_{i})=a_{d}(\cG_{i})+a_{d}(\cF_{i+1})$. If $\cG_{i}=0$, then $\cF_{i}\cong\cF_{i+1}$ via $\alpha$, which implies  $P_{\cF_{i}}(t)=P_{\cF_{i+1}}(t)$, or $P_{\cF_{i}/\cF_{i+1}}(t)=P_{\cF_{i}}(t)-P_{\cF_{i+1}}(t)=0$. However, as $\cF_{i}/\cF_{i+1}$ is a nontrivial coherent sheaf, its Hilbert polynomial cannot be zero (for example, its leading coefficient is a positive number). Therefore, $\cG_{i}$ must be nonzero, which implies by the purity of $\cF_{i}$ that $\cG_{i}$ itself is also pure of dimension $d$. Hence $a_{d}(\cG_{i})>0$, and we conclude $a_{d}(\cF_{i})>a_{d}(\cF_{i+1})$, or $\mult(\cF_{i})>\mult(\cF_{i+1})$.
\end{proof}

\begin{cor}\label{cor:exponent of nilpotence bound by multiplicities}
Let $X$ be a projective scheme over a field and fix an ample line bundle $\cO(1)$ on $X$. Let $\cF$ be a nonzero coherent sheaf on $X$. If $\alpha\colon\cF\rar\cF$ is a nilpotent $\cO_{X}$-linear endomorphism, then $\alpha^{\sum_{i=0}^{\dim\Supp(\cF)}\mult(T_{i}(\cF)/T_{i-1}(\cF))}=0$. 
\end{cor}

\begin{proof}
Note that the image of $\alpha\rvert_{T_{i}(\cF)}\colon T_{i}(\cF)\rar \cF$ is a quotient of $T_{i}(\cF)$, so is of support dimension $\le i$. Therefore, $\Image(\alpha\rvert_{T_{i}(\cF)}\colon T_{i}(\cF)\rar\cF)\subset T_{i}(\cF)$ by the definition of torsion filtration. We may now induct on $\dim\Supp(\cF)$. When $\dim\Supp(\cF)=0$, the assertion is the pure dimension $0$ case of Lemma~\ref{lem:bounding exponent of nilpotence, pure}. Now assume the assertion for $\dim\Supp(\cF)\le d-1$, and suppose that we have a coherent sheaf $\cG$ with $\dim\Supp(\cG)=d$. Then nilpotent $\alpha\colon\cG\rar\cG$ induces nilpotent $\overline{\alpha}\colon\cG/T_{d-1}(\cG)\rar\cG/T_{d-1}(\cG)$. As $\cG/T_{d-1}(\cG)$ is pure of dimension $d$, by Lemma~\ref{lem:bounding exponent of nilpotence, pure}, $\overline{\alpha}^{\mult(\cG/T_{d-1}(\cG))}=0$. This implies that the image of $\alpha^{\mult(\cG/T_{d-1}(\cG))}\colon\cG\rar\cG$ is contained in $T_{d-1}(\cG)$. By the induction hypothesis, $\left(\alpha\rvert_{T_{d-1}(\cG)}\right)^{\sum_{i=0}^{d-1}\mult(T_{i}(\cG)/T_{i-1}(\cG))}=0$. Therefore, $\alpha^{\sum_{i=0}^{d}\mult(T_{i}(\cG)/T_{i-1}(\cG))}=0$, proving the induction step.
\end{proof}

\begin{rem}
Using Remark~\ref{rem:rank}, if $X$ is integral, we can improve the bound in Corollary~\ref{cor:exponent of nilpotence bound by multiplicities} to $\rank(\cF)+\sum_{i=0}^{\dim\Supp(\cF)-1}\mult(T_{i}(\cF)/T_{i-1}(\cF))$, which coincides with the natural bound on the exponent of nilpotence when $\cF$ is locally free.
\end{rem}

Thus, we may obtain a controlled bound on the exponent of nilpotence by controlling how the torsion filtration of a fiber varies over a flat family of coherent sheaves. Unfortunately, the torsion filtration of a fiber can jump within a flat family of coherent sheaves.

\begin{example}\label{example:torsion jump}
Let $Z=\bP_{k}^{1}$ and $S=\bA_{k}^{1}$. Let $\cE$ be the ideal sheaf defining the closed point $(0,0)\in Z_{S}$. As $\cE$ is torsion-free, $\cE$ is $S$-flat. On one hand, the inclusion $\cE\rar\cO_{Z_{S}}$ is an isomorphism away from $(0,0)$, so for $s\in S$ with $s\ne0$, $\cE_{Z_{s}}\cong\cO_{Z_{s}}$ is torsion-free, i.e., $T_{0}(\cE_{Z_{s}})=0$. On the other hand, $T_{0}(\cE_{Z_{0}})$ is the skyscraper sheaf at $0\in Z_{0}$.
\end{example}

On the other hand, there is a finite stratification of the base such that the torsion filtration does not jump inside each stratum.

\begin{lem}\label{lem:controlling torsion filtration}
Let $X$ be a proper scheme over a field $k$, $S$ a reduced Noetherian scheme over $k$, and $\cE$ an $S$-flat coherent sheaf over $X_S$. Then, there is  a finite stratification $S=\bigcup_{i\in I}S_{i}$ by disjoint locally closed reduced subschemes $S_{i}\subset S$ such that, for each $i\in I$, there exists a filtration of $\cE_{X_{S_{i}}}$ by coherent sheaves over $X_{S_{i}}$,
\[
0=\cT_{-1,i}\subset \cT_{0,i}\subset\cT_{1,i}\subset\cdots\subset\cT_{\dim X,i}=\cE_{X_{S_{i}}},
\]
satisfying the following properties: for every $0\le j\le \dim X$, 
\begin{itemize}
\item $\cT_{j,i}/\cT_{j-1,i}$ is $S_{i}$-flat;
\item and for every $s\in S_{i}$, $(\cT_{j,i}/\cT_{j-1,i})_{X_{s}}$ is either $0$ or pure of dimension $j$.
\end{itemize}
In particular, for $s\in S_{i}$, taking the fiber at $s$ of the above filtration gives rise to the torsion filtration of $\cE_{X_{s}}$; namely, for $0\le j\le \dim X$, $(\cT_{j,i})_{X_{s}}\rar\cE_{X_{s}}$ is an injection whose image coincides with $T_{j}(\cE_{X_{s}})$. 
\end{lem}

\begin{proof}
It is enough to show the existence of a stratification and filtrations satisfying the two bullet points as the rest follows easily from them. When $S$ is the spectrum of a field, the torsion filtration of $\calE$ gives the desired filtration. For a general $S$, we prove the existence by Noetherian induction. We may assume that $S$ is irreducible by stratifying $S$. Let $\eta$ be the generic point of $S$. Spread out the torsion filtration $0=T_{-1}\subset T_{0}(\cE_{X_{\eta}})\subset\cdots\subset T_{\dim X}(\cE_{X_{\eta}})=\cE_{X_{\eta}}$ over $X_{\eta}$ to a filtration by coherent subsheaves $0=\cF_{-1}\subset \cF_{0}\subset\cdots\subset\cF_{\dim X}=\cE_{X_{U}}$ of $\cE_{X_{U}}$ for a nonempty open subscheme $U$ of $X$ such that if $T_{j-1}(\cE_{X_{\eta}})=T_{j}(\cE_{X_{\eta}})$, then $\cF_{j-1}=\cF_{j}$. By the generic flatness \cite[052B]{stacks-project}, we may shrink $U$ and assume that $\cF_{j}/\cF_{j-1}$ is $U$-flat for each $0\le j\le \dim X$. Then we know from \cite[Thm.~12.2.1(iii)]{EGAIV3} that each $\cF_{j}/\cF_{j-1}$ being pure of a certain dimension is an open condition. Hence after further shrinking $U$, we may also assume that $(\cF_{j}/\cF_{j-1})_{X_{s}}$ is either $0$ or pure of dimension $j$ for each $s\in U$. Now the first assertion follows by writing $X=U\cup (X\smallsetminus U)$ and applying the Noetherian induction hypothesis to $X\smallsetminus U$.
\end{proof}

\begin{cor}\label{cor:bounding exponent of nilpotence, projective}
Proposition~\ref{prop:bounding exponent of nilpotence} is true if $Z$ is projective over $k$.
\end{cor}

\begin{proof}
Take a finite stratification $S=\bigcup_{i\in I}S_{i}$ and a filtration $0=\cT_{-1,i}\subset\cT_{0,i}\subset\cdots\subset\cT_{\dim Z,i}=\cE_{Z_{S_{i}}}$ for each $i\in I$ as in Lemma~\ref{lem:controlling torsion filtration}. By further stratifying $S_i$, we may also assume that each $S_i$ is connected.

We fix an ample line bundle $\cO(1)$ on $Z$. For each $i\in I$ and $0\le j\le \dim Z$, $\cT_{j,i}/\cT_{j-1,i}$ is $S_{i}$-flat. Hence the Hilbert polynomial of $\cT_{j,i}/\cT_{j-1,i}$ is constant over $S_i$ by \cite[Thm. 7.9.4]{EGAIII2}; there exists a polynomial $F_{j,i}(t)\in\bQ[X]$ such that the Hilbert polynomial of $(\cT_{j,i}/\cT_{j-1,i})_{Z_{s}}$ is equal to $F_{j,i}(t)$ for any $s\in S_{i}$. In particular, the multiplicity of $(\cT_{j,i}/\cT_{j-1,i})_{Z_{s}}$ is constant on $s\in S_i$, say $m_{j,i}\in \Z_{\geq 0}$ as $(\cT_{j,i}/\cT_{j-1,i})_{Z_{s}}$ is either $0$ or pure of dimension $j$ (see Definition~\ref{defn:Hilbert polynomial}). Let $N_{i}\coloneqq \sum_{j=0}^{\dim Z}m_{j,i}$. By the $S_{i}$-flatness of $\cT_{j-1,i}$, we have 
\[
(\cT_{j,i}/\cT_{j-1,i})_{Z_{s}}=(\cT_{j,i})_{Z_{s}}/(\cT_{j-1,i})_{Z_{s}}=T_{j}(\cE_{Z_{s}})/T_{j-1}(\cE_{Z_{s}}).
\]
Therefore, we know from Corollary~\ref{cor:exponent of nilpotence bound by multiplicities} that for any $s\in S_{i}$, if $\alpha\colon\cE_{Z_{s}}\rar\cE_{Z_{s}}$ is a nilpotent $\cO_{Z_{s}}$-linear endomorphism, then $\alpha^{N_{i}}=0$. 

Set $N=\max_{i\in I}N_{i}$. We will show that this $N$ works; namely, for any reduced $S$-scheme $S'$, if an $\cO_{Z_{S'}}$-module endomorphism $\alpha\colon\cE_{Z_{S'}}\rar\cE_{Z_{S'}}$ is nilpotent, then $\alpha^{N}=0$. By construction, this is true if $S'=\Spec k(s)$ for $s\in S$. Hence it also holds for the spectrum of any field over $S$. For a general reduced $S$-scheme $S'$, note that $Z_{S'}$ is reduced as $k$ is perfect. It suffices to show that for every affine open $U=\Spec R$ of $Z_{S'}$, the restriction $\alpha_U\in \End(\calE_U)$ satisfies $\alpha_U^N=0$. The latter statement is deduced from the field case by Remark~\ref{rem:commutative algebra about nilpotent endomorphism}.
\end{proof}

\begin{proof}[Proof of Proposition~\ref{prop:bounding exponent of nilpotence}]
We work on the general case ($Z$ is only assumed to be proper). Let $Z=\bigcup_{i=1}^{m}U_{i}$ be a finite affine open cover. By a refined version of Chow's Lemma as in \cite[Cor.~2.6]{Conrad-Nagata}, there exists a projective $k$-scheme $P_{i}$ with a surjective morphism $q_{i}\colon P_{i}\rar Z$ such that $q_{i}$ restricted to $q_{i}^{-1}(U_{i})$ is an isomorphism onto $U_{i}$. Let $(q_{i})_{S}\colon (P_{i})_{S}\rar Z_{S}$ be the base-change of $q_{i}$, and let $\cE_{i}\coloneqq(q_{i})_{S}^{*}(\cE)$. As $(P_{i})_{S}\rar S$ is projective and $S$ is Noetherian, the flattening stratification of $S$ for $\cE_{i}$ exists (e.g. \cite[Thm.~2.1.5]{HL}), which is in particular a finite stratification $S=\bigcup_{j\in J_{i}}S_{j}$ by disjoint locally closed reduced subschemes $S_{j}\subset S$ such that, for each $j\in J_{i}$,  $(\cE_{i})_{(P_{i})_{S_{j}}}$ is $S_{j}$-flat. By the projective case (Corollary~\ref{cor:bounding exponent of nilpotence, projective}), there exists $N_{i,j}\in\bN$ such that, for any $s\in S_{j}$ and any $\cO_{(P_{i})_{s}}$-endomorphism $\alpha\colon(\cE_{i})_{(P_{i})_{s}}\rar(\cE_{i})_{(P_{i})_{s}}$, if $\alpha$ is nilpotent, then $\alpha^{N_{i,j}}=0$. Let $N_{i}\coloneqq\max_{j\in J_i}N_{i,j}$. Then, for any $s\in S$ and any $\cO_{(P_{i})_{s}}$-endomorphism $\alpha\colon(\cE_{i})_{(P_{i})_{s}}\rar(\cE_{i})_{(P_{i})_{s}}$, if $\alpha$ is nilpotent, then $\alpha^{N_{i}}=0$.

We claim that $N\coloneqq\max_{1\leq i\leq m}N_i$ gives the desired universal bound on the exponent of nilpotence. As in the last part of the proof of the projective case (Corollary~\ref{cor:bounding exponent of nilpotence, projective}), we only need to show that, if $s\in S$ and $\alpha\colon\cE_{Z_{s}}\rar\cE_{Z_{s}}$ is a nilpotent $\cO_{Z_{s}}$-linear endomorphism, then $\alpha^{N}=0$. For $1\le i\le m$, after pulling back by $(q_{i})_{s}\colon(P_{i})_{s}\rar Z_{s}$, we obtain a nilpotent $\cO_{(P_{i})_{s}}$-linear endomorphism $(q_{i})_{s}^{*}(\alpha)\colon(\cE_{i})_{(P_{i})_{s}}\rar(\cE_{i})_{(P_{i})_{s}}$. Hence $\left((q_{i})_{s}^{*}(\alpha)\right)^{N}=(q_{i})_{s}^{*}(\alpha^{N})=0$. As $(q_{i})_{s}$ yields an isomorphism $(q_{i})_{s}\colon (q_{i}^{-1}(U_{i}))_{s}\xrightarrow{\sim} (U_{i})_{s}$, this implies that $(q_{i})_{s}^{*}(\alpha^{N})\rvert_{(q_{i}^{-1}(U_{i}))_{s}}=0$, or $(\alpha^{N})\rvert_{(U_{i})_{s}}=0$. As $Z=\bigcup_{i=1}^{m}U_{i}$, we conclude $\alpha^{N}=0$.
\end{proof}

\subsection{Analytification of \texorpdfstring{$\calM_\dR$}{MdR} over a non-archimedean field}\label{sec:analytification of de Rham moduli}

In this subsection, we show that the analytification of the algebraic stack $\cM_{\dR}$ is the moduli stack of analytic integrable connections. 

Fix a complete non-archimedean field $K$ of mixed characteristic $(0,p)$. Let $Z_K$ be a smooth proper scheme over $S\coloneqq \Spec K$ with a relative normal crossings divisor $D_K$ and write $M_{Z_K}$ for the log structure on $Z_K$ given by $D_K$. Let $\calM_{\dR}\coloneqq \calM_{\dR}((Z_K,D_K)/\Spec K)\rightarrow (\Sch/\Spec K)_\fppf$ denote the moduli stack of integrable connections for $(Z_K,D_K)\rightarrow \Spec K$. 

A in \S\ref{sec:analytification}, for $T=\Spa(A,A^+)\in \Ob(\Afd_K)$ with $A$ complete, set 
\[
Z_A\coloneqq Z_K\times_{\Spec K}\Spec A\quad\text{and}\quad Z_T^\an\coloneqq Z_K^\an\times_{\Spa(K,\calO_K)}T=Z_A\times_{\Spec A}T.
\]
We use similar notation for other $K$-schemes such as $D_K$ and consider the log structure $M_{Z_T^\an}$ on $Z_T^\an$ given by $D_T^\an$ (see \cite{DLLZ-logadic} for the generalities of log adic spaces). Write $\omega^1_{Z_T^\an/D_T^\an}$ for the sheaf of log differentials. One can define the notion of (integrable) logarithmic connections in an obvious way.
Note that $A$ is Noetherian and every $\calO_{Z_A}$-module of finite presentation is coherent.

\begin{defn}\label{def:moduli stack of analytic integrable connections}
Define the \emph{moduli stack of analytic integrable connections},
\[
\calM_\andR\rightarrow \Afd_K,
\]
to be the category fibred in groupoids whose fiber over $T\in \Ob(\Afd_K)$ is the groupoid of coherent sheaves with integrable connection on $Z_{T}^\an/T$; namely, pairs $(\calE,\nabla)$ where $\calE$ is a $T$-flat coherent $
\calO_{Z_{T}^\an}$-module, and $\nabla\colon \calE\rightarrow \calE\otimes_{\calO_{Z_T^\an}}\omega_{Z_T^\an/T}^1$ is an  integrable $T$-connection. We simply call such $(\calE,\nabla)$ an \emph{integrable $T$-connection}.
\end{defn}

The main goal of this subsection is to prove the following.

\begin{thm}\label{thm:analytification of MdR}
The category fibered in groupoids $\cM_{\andR}$ is an adic stack over $\Spa(K,\cO_{K})$, and it is equivalent to the analytification $\cM_{\dR}^{\an}$ of $\calM_\dR$.
\end{thm}

Let $P^1$ denote the first log infinitesimal neighborhood of $(Z_K,M_{Z_K})$ in 
\[
(Z_K(1),M_{Z_K(1)})\coloneqq (Z_K,M_{Z_K})\times_{\Spec K}(Z_K,M_{Z_K}). 
\]
Concretely, $P^1\coloneqq \operatorname{\underline{Spec}}\calO_{Z_K(1)^{\mathrm{log}}}/\calI^2$ where $(Z_K,M_{Z_K})\hookrightarrow (Z_K(1)^{\mathrm{log}},M_{Z_K(1)^{\mathrm{log}}})$ is the exactification of $(Z_K,M_{Z_K})\hookrightarrow (Z_K(1),M_{Z_K(1)})$
 and $\calI\subset \calO_{Z_K(1)^{\mathrm{log}}}$ is the defining ideal of $Z_K\hookrightarrow Z_K(1)^{\mathrm{log}}$.
It is easy to see that $P_T^{1,\an}$ is the first log infinitesimal neighborhood of $(Z_T^\an,M_{Z_T^\an})$ in the self-product $(Z_T^\an(1),M_{Z_T^\an(1)})$ of $(Z_T^\an,M_{Z_T^\an})$ over $T$ (in the category of log adic spaces). Note that $P^1$ is proper over $S$, and that $\Delta\colon Z_T^\an\hookrightarrow P_T^{1,\an}$ is a closed immersion by a square-zero ideal, and the underlying topological spaces are identified. In particular, we often regard $\calO_{Z_T^\an}$ as a sheaf on $P_T^{1,\an}$.
Write $p_1,p_2$ for the two projections $P_T^{1,\an}\hookrightarrow Z_T^\an(1)\rightrightarrows  Z_T^\an$.

\begin{lem}\label{lem:connection and P1-linear isom}
For $\calE\in \Coh(\calO_{Z_T^\an})$, equipping $\calE$ with a connection $\nabla\colon \calE\rightarrow \calE\otimes\omega^{1}_{Z_T^\an/T}$ is the same as giving an $\calO_{P_T^{1,\an}}$-linear isomorphism $\epsilon\colon p_2^\ast\calE \xrightarrow{\cong}p_1^\ast \calE$ satisfying $\Delta^\ast\epsilon=\id_\calE$ via $\nabla(e)=\epsilon(1\otimes e)-e\otimes 1$ for a section $e$ of $\calE$. Moreover, morphisms also correspond in an obvious way.    
\end{lem}

\begin{proof}
This is standard as $\omega^{1}_{Z_T^\an/T}=\Ker(\calO_{P_T^{1,\an}}\rightarrow \calO_{Z_T^\an})$ (cf.~\cite[p.~94]{DLLZ-logadic}).
\end{proof}

\begin{lem}\label{lem:analytic MdR is a stack}
The functor $\calM_\andR\rightarrow (\Afd_K)_\Et$ is a stack in groupoids.
\end{lem}

\begin{proof}
First we show that for $x=(\calE,\nabla), x'=(\calE',\nabla')\in\Ob( (\calM_{\andR})_T)$, the presheaf $\Isom(x,x')$ on $\Afd_K/T$ is a sheaf. To see this, observe that $\Isom(x,x')$ is a subfunctor of the Hom-sheaf sending $(T'\rightarrow T)$ to $\Hom_{Z_{T'}^\an}(\calE|_{Z_{T'}^\an},\calE'|_{Z_{T'}^\an})$, consisting of isomorphisms compatible with $\nabla|_{Z_{T'}^\an}$ and $\nabla'|_{Z_{T'}^\an}$. Since these conditions are preserved by pullback and can be checked \'etale locally on $T'$, we see that $\Isom(x,x')$ is a sheaf.
It remains to show that the descent data are effective. Using Lemma~\ref{lem:connection and P1-linear isom}, one can easily deduce this from Propositions~\ref{prop:coherent sheaves on etale site of adic space}(ii) and \ref{prop:descent of coherent sheaves}, Remark~\ref{rem:flatnesss of coherent sheaf and morphism for adic spaces}, and the fact that the integrability of a connection can be checked \'etale locally. 
\end{proof}

\begin{proof}[Proof of Theorem~\ref{thm:analytification of MdR}]
Take $T=\Spa(A,A^+)\in \Ob(\Afd_K)$ with $A$ complete. We know from Theorem~\ref{thm:rigid relative GAGA} that the analytification functors
\[
\Coh(\calO_{Z_A})\rightarrow\Coh(\calO_{Z_T^{\an}})\quad\text{and}\quad
\Coh(\calO_{P^1_A})\rightarrow\Coh(\calO_{P_T^{1,\an}})
\]
are equivalences since $Z$ and $P^1$ are proper over $K$. Moreover, we know from Proposition~\ref{prop:analytification and flatness} that $A$-flatness corresponds to $T$-flatness under these equivalences. By Lemma~\ref{lem:connection and P1-linear isom}, one can associate to every $(\Spec A)$-connection $(\calE,\nabla)$ on $Z_A$ a $T$-connection $(\calE^\an,\nabla^\an)$ on $Z_T^\an$. Moreover, $(\calE,\nabla)$ is integrable if and only if so is $(\calE^\an,\nabla^\an)$: to see this, directly compare $\nabla^2$ and $\nabla^{\an,2}$; alternatively, use the fact that a connection in characteristic zero is integrable if and only if it is upgraded to a stratification (cf.~\cite[Prop.~1.2.7]{Shiho-II}). Combining this result with the preceding discussion, we obtain a natural equivalence $\alpha_T\colon (\calM_{\dR})_{\Spec A}\xrightarrow{\cong} (\calM_{\andR})_{T}$ that is functorial in $T$. By definition, $\calM_\dR^\an$ is the stackification of $T\mapsto (\calM_{\dR})_{\Spec A}$. By Lemma~\ref{lem:analytic MdR is a stack} (and \cite[042W]{stacks-project}), this is already a stack. Hence $\alpha_T$ yields the desired equivalence $\calM_\dR^\an\xrightarrow{\cong}\calM_{\andR}$.
\end{proof}

\begin{rem}
Keep the setup at the beginning of this subsection. The above arguments also show that the analytification of $\calCoh(Z_K/\Spec K)\rightarrow (\Sch/\Spec K)_\fppf$ is equivalent to the category fibered in groupoids over $\Afd_K$ whose fiber over $T$ is the groupoid of $T$-flat $\calO_{Z_T^\an}$-modules of finite presentation.
\end{rem}

\begin{defn}
Let $T=\Spa(A,A^+)\in \Ob(\Afd_K)$ with $A$ complete and keep the notation as in the proof of Theorem~\ref{thm:analytification of MdR}. For an integrable $(\Spec A)$-connection $(\calE,\nabla)$ on $Z_A$, we call the associated integrable $T$-connection on $Z_T^\an$ the \emph{analytification} of $(\calE,\nabla)$ and write $(\calE^\an,\nabla^\an)$.

Define the $i$-th de Rham cohomology $H^i_\dR((Z_A,M_{Z_A}),(\calE,\nabla))$ to be the $i$-th hypercohomology of the log de Rham complex
\[
H^i_\dR((Z_A,M_{Z_A}),(\calE,\nabla))\coloneqq \mathbb{H}^i(Z_A,(\calE\otimes_{\calO_{Z_A}}\omega^\bullet_{Z_A/A},\nabla)).
\]
Similarly, we define 
\[
H^i_\dR((Z_T^\an,M_{Z_T^\an}),(\calE^\an,\nabla^\an))\coloneqq \mathbb{H}^i(Z_T^\an,(\calE^\an\otimes_{\calO_{Z_T^\an}}\omega^\bullet_{Z_T^\an/T},\nabla^\an)).
\]
There is a natural isomorphism of $A$-modules
\[
H^i_\dR((Z_A,M_{Z_A}),(\calE,\nabla))\xrightarrow{\cong}H^i_\dR((Z_T^\an,M_{Z_T^\an}),(\calE^\an,\nabla^\an))
\]
induced by the spectral sequences
\begin{align*}
&E_1^{a,b}=H^b(Z_A,\calE\otimes_{\calO_{Z_A}}\omega^a_{Z_A/A})&&\implies H^{a+b}_\dR((Z_A,M_{Z_A}),(\calE,\nabla)),\\
&E_1^{a,b}=H^b(Z_T^\an,\calE^\an\otimes_{\calO_{Z_T^\an}}\omega^a_{Z_T^\an/T})&&\implies H^{a+b}_\dR((Z_T^\an,M_{Z_T^\an}),(\calE^\an,\nabla^\an))
\end{align*}
and the isomorphisms $H^b(Z_A,\calE\otimes_{\calO_{Z_A}}\omega^a_{Z_A/A})\cong H^b(Z_T^\an,\calE^\an\otimes_{\calO_{Z_T^\an}}\omega^a_{Z_T^\an/T})$.

\end{defn}

\subsection{Infinitesimal structures of the algebraic and analytic moduli}\label{sec:Infinitesimal}

In this subsection, keep the assumption $S=\Spec K$. We will identify the tangent spaces of the algebraic moduli space $\cM_{\dR}$ and the analytic moduli space $\cM_{\andR}$ with appropriate de Rham cohomology groups. 
See \cite[06I1]{stacks-project} and Definition~\ref{Tangent-Spaces-Defn} for the notion of tangent spaces and infinitesimal automorphisms of algebraic stacks and adic stacks, respectively.

\begin{lem}\label{lem:algebraic tangent space}
Let $L/K$ be a field extension, and let $(\cE,\nabla)\in\Ob((\cM_{\dR})_{\Spec L})$ be an integrable $L$-connection.
Write $\underline{\End}(\cE,\nabla)$ for the coherent sheaf $\underline{\End}(\cE)$ on $Z_{L}$ with an integrable connection $\nabla_{\underline{\End}(\cE,\nabla)}\colon \underline{\End}(\cE)\rightarrow\underline{\End}(\cE)\otimes \omega^{1}_{Z_{L}/L}$ defined by $\nabla_{\underline{\End}(\cE,\nabla)}(f)\coloneqq \nabla\circ f - (f\otimes 1)\circ \nabla$.
Then, there are natural isomorphisms of $L$-vector spaces,
\begin{align*}
T_{\cM_{\dR},(\cE,\nabla)}
& \cong H_{\dR}^{1}((Z_{L},M_{Z_{L}}),\underline{\End}(\cE,\nabla))    \qquad\text{and}\\
\Inf_{\cM_{\dR},(\cE,\nabla)}
&\cong H_{\dR}^{0}((Z_{L},M_{Z_{L}}),\underline{\End}(\cE,\nabla)).
\end{align*}
\end{lem}

\begin{proof}
By definition, $T_{\cM_{\dR},(\cE,\nabla)}$ is the set of isomorphism classes of integrable $L[\epsilon]/(\epsilon^{2})$-connections $(\widetilde{\cE},\widetilde{\nabla}\colon\widetilde{\cE}\rightarrow\widetilde{\cE}\otimes\omega_{(Z_{L[\epsilon]/(\epsilon^{2})})/(L[\epsilon]/(\epsilon^{2}))}^{1})$ whose mod $\epsilon$ reduction is isomorphic to $(\cE,\nabla)$. 
Let $q\colon Z_{L[\epsilon]/(\epsilon^{2})}\rightarrow Z_{L}$ denote the projection map. Then we obtain an integrable $L$-connection $(\cE',\nabla')\coloneqq(q_{\ast}\widetilde{\cE},q_{\ast}\widetilde{\nabla})$, where $\nabla'$ is regarded as a morphism of abelian sheaves
\[
\nabla'\colon\cE'\rightarrow q_{\ast}(\widetilde{\cE}\otimes\omega^{1}_{Z_{L[\epsilon]/(\epsilon^{2})}/(L[\epsilon]/(\epsilon^{2}))})=q_{\ast}(\widetilde{\cE}\otimes q^{\ast}\omega^{1}_{Z_{L}/L})=\cE'\otimes\omega^{1}_{Z_{L}/L}.
\] 

The multiplication-by-$\epsilon$ map defines an endomorphism $\times\epsilon\colon(\cE',\nabla')\rightarrow(\cE',\nabla')$ as an $L$-connection, and the cokernel is isomorphic to $(\cE,\nabla)$. 
Furthermore, the structure map $\Spec L[\epsilon]/(\epsilon^{2})\rightarrow \Spec L$ is a splitting of the first-order thickening $\Spec L\rightarrow\Spec L[\epsilon]/(\epsilon^{2})$ and it gives a splitting $\iota\colon(\cE,\nabla)\rightarrow(\cE',\nabla')$ to the projection $(\cE',\nabla')\rightarrow(\cE,\nabla)$. 
Hence $(\cE',\nabla')$ sits in a short exact sequence of $L$-connections
\[
0\rightarrow(\cE,\nabla)\rightarrow(\cE',\nabla')\rightarrow(\cE,\nabla)\rightarrow0,
\]
where the subobject is the image of $\iota(\cE,\nabla)$ via the map $\times\epsilon$. Conversely, given such a short exact sequence of integrable $L$-connections, one may recover the multiplication-by-$\epsilon$ map on $(\cE',\nabla')$ as the map that vanishes on the subconnection $(\cE,\nabla)$ and sends the quotient connection $(\cE,\nabla)$ to the subconnection $(\cE,\nabla)$. From this, we can recover the integrable $L[\epsilon]/(\epsilon^{2})$-connection $(\widetilde{\cE},\widetilde{\nabla})$. 
Furthermore, as $\nabla'$ is integrable, so is $\widetilde{\nabla}$.
Therefore, $T_{\cM_{\dR},(\cE,\nabla)}$ is the set of isomorphism classes of extensions of $(\cE,\nabla)$ by itself in the category of coherent sheaves with integrable connections on $Z_{L}$ over $L$. As the de Rham cohomology with coefficients is the right derived functor of the solutions functor on the category of coherent sheaves with integrable connections (e.g.~\cite[(4.4)]{Katz_Nilpotent}), the set of isomorphism classes of self-extensions is naturally identified with $H_{\dR}^{1}((Z_{L},M_{Z_L}),\underline{\End}(\cE,\nabla))$. One may check easily that the isomorphism identifies the two $L$-vector space structures.

For the second isomorphism, identify the trivial deformation of $(\cE,\nabla)$ with $(\cE'_{0},\nabla'_{0})\coloneqq(\cE,\nabla)^{\oplus 2}$ equipped with the $L[\epsilon]/(\epsilon^2)$-structure given by $\times\epsilon=\begin{pmatrix}0&1\\0&0\end{pmatrix}$. An endomorphism of $(\cE'_{0},\nabla'_{0})$ is represented by a matrix $\begin{pmatrix}a&b\\c&d\end{pmatrix}$ for $a,b,c,d\in H^{0}_{\dR}((Z_{L},M_{Z_{L}}),\underline{\End}(\cE,\nabla))$; it respects the $L[\epsilon]/(\epsilon^2)$-structure if and only if $c=0$ and $a=d$, in which case its mod-$\epsilon$ reduction is given by $a$. Thus we obtain
\begin{align*}
\Inf_{\cM_{\dR},(\cE,\nabla)}
&=\left\lbrace\begin{pmatrix}1&b\\0&1\end{pmatrix}~:~b\in H^{0}_{\dR}((Z_{L},M_{Z_{L}}),\underline{\End}(\cE,\nabla))\right\rbrace\\
&=H^{0}_{\dR}((Z_{L},M_{Z_{L}}),\underline{\End}(\cE,\nabla)).
\end{align*}
One may check that the identification is an $L$-linear isomorphism.
\end{proof}

\begin{prop}\label{prop:analytic Inf and T}
Let $(L,L^{+})$ be a complete analytic affinoid field over $(K,\cO_{K})$ and set $T=\Spa(L,L^+)$. Let $(\cE,\nabla)\in\Ob((\cM_{\dR})_{\Spec L})$ be an integrable $L$-connection, and let $(\calE^\an,\nabla^\an)\in\Ob((\cM_{\andR})_T)$ denote the corresponding integrable $T$-connection via Theorem~\ref{thm:analytification of MdR}. Then, there are natural commutative diagrams of $L$-vector space isomorphisms
\[
\xymatrix{
T_{\cM_{\dR},(\cE,\nabla)}\ar[r]^-\cong\ar[d]_-\cong &  H_{\dR}^{1}((Z_{L},M_{Z_{L}}),\underline{\End}(\cE,\nabla)) \ar[d]^-\cong\\
T_{\cM_{\andR},(\calE^\an,\nabla^\an)}\ar[r]^-\cong & H_{\dR}^{1}((Z_T^\an,M_{Z_T^\an}),\underline{\End}(\calE^\an,\nabla^\an))
}
\]
and
\[
\xymatrix{
\Inf_{\cM_{\dR},(\cE,\nabla)}\ar[r]^-\cong\ar[d]_-\cong &  H_{\dR}^{0}((Z_{L},M_{Z_{L}}),\underline{\End}(\cE,\nabla)) \ar[d]^-\cong\\
\Inf_{\cM_{\andR},(\calE^\an,\nabla^\an)}\ar[r]^-\cong & H_{\dR}^{0}((Z_T^\an,M_{Z_T^\an}),\underline{\End}(\calE^\an,\nabla^\an)).
}
\]
\end{prop}

\begin{proof}
By the proof of Theorem~\ref{thm:analytification of MdR}, the analytification gives an equivalence of the tangent groupoids $\fT_{\cM_{\dR},(\cE,\nabla)}\rightarrow\fT_{\cM_{\andR},(\cE^\an,\nabla^\an)}$. Therefore, the desired statements follow from Lemma~\ref{lem:algebraic tangent space}  and the analogous statement for $(\calE^\an,\nabla^\an)$ (which can be proved exactly in the same manner).
\end{proof}

\subsection{The generic fiber of the \texorpdfstring{$p$}{p}-adic completion}\label{sec:generic fiber of p-adic completion of de Rham moduli}
Keep the notation in \S\ref{sec:analytification of de Rham moduli}. We further assume that there exists a smooth proper scheme $Z\rightarrow \Spec \calO_K$ with relative normal crossings divisor $D$ such that the previous $(Z_K,D_K)$ is the generic fiber of $(Z,D)$. Let $\widehat{\calM}_\dR\rightarrow (\Sch/\calO_K)_\fppf$ be the formal algebraic stack obtained as the $p$-adic completion of $\calM_\dR((Z,D)/\Spec\calO_K)$ (see Definition~\ref{def:I-adic completion of algebraic stack}). Definition~\ref{defn:comparison map from generic fiber to analytification} gives a $1$-morphism of adic stacks over $\Afd_K$
\[
\widehat{\calM}_{\dR,\eta}\rightarrow \calM_{\dR}^\an\cong \calM_{\andR}.
\]
We will give an explicit description of this $1$-morphism. Take $T=\Spa(A,A^+)\in \Ob(\Afd_K)$. We assume that $A$ is complete and $(A,A^+)$ is of topologically finite type over a complete analytic affinoid field $(L,L^+)$. 

We start with a preliminary discussion. 
Fix an open $L^+$-subalgebra $A_{0}\subset A^+$ of topological finite type. Note $A_0[p^{-1}]=A$.

Consider the smooth proper morphism $Z_{A_0}\rightarrow \Spec A_0$ and its $p$-adic completion $\widehat{Z}_{A_0}\rightarrow \Spf A_0$. Since $A_{i,0}$ is $p$-torsion free of topologically finite type over the $p$-adically complete valuation ring $L^+$, the formal scheme $\widehat{Z}_{A_{i,0}}$ satisfies \cite[Situation I.10.1.1]{Fujiwara-Kato} by \cite[Cor.~0.9.2.8, 0.9.2.9]{Fujiwara-Kato}. Under this assumption, Fujiwara and Kato proved the following GFGA existence theorem. 

\begin{thm}[GFGA existence theorem]\label{thm:GFGA existence}
Let $i\colon \widehat{Z}_{A_{0}}\rightarrow Z_{A_{0}}$ denote the morphism of locally ringed spaces. Then $i^\ast\colon \Mod(\calO_{Z_{A_{0}}})\rightarrow \Mod(\calO_{\widehat{Z}_{A_{0}}})$ induces an equivalence of categories
$\Coh(\calO_{Z_{A_{0}}})\xrightarrow{\cong}\Coh(\calO_{\widehat{Z}_{A_{0}}})$, and the functor $i^\ast$ on $\Coh(\calO_{Z_{A_{0}}})$ agrees with the $p$-adic completion functor $\calF\mapsto \widehat{\calF}\coloneqq \varprojlim_{n}\calF/p^n\calF$.
\end{thm}

\begin{proof}
This is \cite[Thm.~I.10.1.2, \S{}I.9.1.(a)]{Fujiwara-Kato}.
\end{proof}

\begin{cor}\label{cor:GFGA existence for connections}
The $p$-adic completion functor induces an equivalence from the category of coherent sheaves with integrable $(\Spec A_0)$-connections on $Z_{A_0}$ to the category of coherent sheaves with integrable $(\Spf A_0)$-connections on $\widehat{Z}_{A_0}$.
\end{cor}

\begin{proof}
Apply Theorem~\ref{thm:GFGA existence} to $Z_{A_0}$ and  the first log infinitesimal neighborhood of $(Z_{A_0},M_{Z_{A_0}})$ in the self-fiber product of $(Z_{A_0},M_{Z_{A_0}})$ over $\Spec A_0$; we deduce from a stratification argument that the category of coherent sheaves with $(\Spec A_0)$-connections on $Z_{A_0}$ is equivalent to the category of coherent sheaves with integrable $(\Spf A_0)$-connections on $\widehat{Z}_{A_0}$ via the $p$-adic completion. It is easy to check that the integrability of a connection corresponds via this equivalence (e.g.~by taking a local coordinate).
\end{proof}

Let us have an extra discussion on the generic fiber, which will be used in \S\ref{sec:adic moduli stack of isocrystals}: assume further that $A^+$ is the integral closure of $A_0$ in $A$\footnote{If $L^+=L^\circ$, then every open $L^+$-subalgebra of $A^+$ of topological finite type satisfies this property by \cite[Lem.~4.4]{Huber-gen}.}. Note that such $A_0$'s are cofinal among the open $L^+$-subalgebras of $A^+$ of topological finite type (see Remark~\ref{rem:description of generic fiber of formal scheme over L}(i)).
Since $\Spf A_0$ and $\widehat{Z}_{A_0}$ satisfy the condition (S)(b) in \cite[\S1.9]{Huber-etale}, the functor $d$ in \cite[Prop.~1.9.1]{Huber-etale} is defined on these formal schemes, and it agrees with the generic fiber functor in \S\ref{sec:generic fiber} if $A_0$ is Noetherian; hence we obtain a morphism of analytic adic spaces $d(\widehat{Z}_{A_0})\rightarrow d(\Spf A_0)$. The additional assumption on $A_0$ implies $d(\Spf A_0)=\Spa(A,A^+)=T$. Moreover, on one hand, these analytic adic spaces also sit in the following commutative diagram of topologically ringed spaces
\[
\xymatrix{
(d(\widehat{Z}_{A_0}),\calO_{d(\widehat{Z}_{A_0})}) \ar[r]\ar[d] & (T,\calO_T) \ar[d]\\
(\widehat{Z}_{A_0},\calO_{\widehat{Z}_{A_0}}) \ar[r] & (\Spf A_0,\calO_{\Spf A_0}).
}
\]
On the other hand, we have a morphism $d(\widehat{Z}_{A_0}) \rightarrow Z_{A_0}\times_{\Spec A_0}T$ of adic spaces over $T$, which is an isomorphism by \cite[Prop.~1.9.6]{Huber-etale} and the fact that $Z_{A_0}\rightarrow \Spec A_0$ is proper. It is easy to see that $Z_{A_0}\times_{\Spec A_0}T=Z_A\times_{\Spec A}T=Z_T^\an$. Combining these arguments, we obtain a morphism of topologically ringed spaces 
\[
\lambda \colon (Z_T^\an,\calO_{Z_T^\an}) \rightarrow (\widehat{Z}_{A_0},\calO_{\widehat{Z}_{A_0}}).
\]

\begin{prop}\label{prop:analytification and completion of coheren sheaves and connections}
Assume that $A^+$ is the integral closure of $A_0$ in $A$.
\begin{enumerate}
 \item There is a natural identification $\lambda_\ast\calO_{Z_T^\an}=\calO_{\widehat{Z}_{A_0}}[p^{-1}]$.
 \item For $\calF\in \Coh(\calO_{Z_{A_0}})$, let $\widehat{\calF}\in \Coh(\calO_{\widehat{Z}_{A_0}})$ denote its $p$-adic completion. Write $\calF_T^\an\in \Coh(\calO_{Z_T^\an})$ for the analytification of $\calF|_{Z_A}\in\Coh(\calO_{Z_A})$. 
Then $\lambda_\ast \calF_T^\an\cong \widehat{\calF}[p^{-1}]$. A similar equality holds for a coherent sheaf with $A_0$-connection on $Z_{A_0}$ in place of $\calF$.
\end{enumerate}
\end{prop}

\begin{proof}
Part (i) follows from \cite[Prop.~1.9.1, Pf.]{Huber-etale}. For (ii), the first assertion follows from (i), \cite[1.9.4]{Huber-etale}, and the construction of these functors, and the second is deduced similarly via a stratification argument.
\end{proof}

\smallskip
\noindent
\textbf{The description of $\widehat{\calM}_{\dR,\eta}\rightarrow \calM_{\andR}$}.
Let $x\in \Ob((\widehat{\calM}_{\dR,\eta})_T)$. By Remark~\ref{rem:generic fiber of formal algebraic stacks}, there exist an \'etale covering $\{T_i=\Spa(A_i,A_i^+)\rightarrow T\}$ with $A_i$ complete and for each $i$, an open $L^+$-subalgebra $A_{i,0}$ of $A_i^+$ of topologically finite type such that $x|_{T_i}$ is in the essential image of 
\[
(\widehat{\calM}_\dR)_{\Spf A_{i,0}} \rightarrow \colim_{A_{i,0}'\subset A_i^+}(\widehat{\calM}_\dR)_{\Spf A_{i,0}'}\rightarrow (\widehat{\calM}_{\dR,\eta})_{T_i}.
\]
Take $x_i\in \Ob((\widehat{\calM}_\dR)_{\Spf A_{i,0}})$ that is mapped to $x|_{T_i}$. Let $\pi\in\calO_K$ be a uniformizer. By definition, $x_i$ corresponds to a compatible system $(\calE_{i,n},\nabla_{i,n})_n$ of integrable connections: 
\[
(\calE_{i,n},\nabla_{i,n})\in \Ob((\calM_\dR((Z,D)/\Spec\calO_K))_{\Spec A_{i,0}/(\pi^n)})
\]
such that $(\calE_{i,n+1},\nabla_{i,n+1})|_{Z_{A_{i,0}/(\pi^n)}}=(\calE_{i,n},\nabla_{i,n})$. Then $(\widehat{\calE}_i,\widehat{\nabla}_i)\coloneqq \varprojlim_n(\calE_{i,n},\nabla_{i,n})$ is a $p$-torsion free a.q.c.~sheaf of finite type with integrable $A_{i,0}$-connection on the $p$-adic formal scheme $\widehat{Z}_{A_{i,0}}$; $\widehat{\calE}_i$ is a coherent $\calO_{\widehat{Z}_{A_{i,0}}}$-module as every $p$-torsion free a.q.c.~sheaf of finite type on $\widehat{Z}_{A_{i,0}}$ is coherent by \cite[Prop.~0.4.1.8, I.3.5.10, I.7.2.3]{Fujiwara-Kato}. By Theorem~\ref{cor:GFGA existence for connections}, there exists a coherent sheaf with integrable $A_{i,0}$-connection $(\calE_i,\nabla_i)$ on $Z_{A_{i,0}}$ such that the $p$-adic completion of $(\calE_i,\nabla_i)$ is $(\widehat{\calE}_i,\widehat{\nabla}_i)$. Moreover, $\calE$
 is flat over $A_{i,0}$ by construction and \cite[Prop.~I.2.1.5, Cor.~I.4.8.2]{Fujiwara-Kato}.

Let $x_{T_i}^\an\in \Ob((\calM_\andR)_{T_i})$ denote the image of $x|_{T_i}$ under the $1$-morphism $\widehat{\calM}_{\dR,\eta}\rightarrow \calM_{\andR}$. We deduce from Remark~\ref{rem:description of map from generic fiber to analytification} that $x_{T_i}^\an$ is given by the analytification of $(\calE_i,\nabla_i)|_{Z_{A_i}}$.

Factor $\widehat{\calM}_{\dR,\eta}\rightarrow \calM_{\andR}$ as $\widehat{\calM}_{\dR,\eta}\rightarrow \calW\subset \calM_{\andR}$ with $\calW$ being an open substack of $\calM_{\andR}$ and the first $1$-morphism an epimorphism as in Proposition~\ref{prop:image of comparison map from generic fiber to analytification}.

\begin{prop}\label{prop:analytic connection coming from generic fiber of completion}
Consider the above factorization $\widehat{\calM}_{\dR,\eta}\rightarrow \calW\subset \calM_{\andR}$. Let $y=(\calE,\nabla)\in \Ob((\calM_{\andR})_T)$ be a $T$-flat coherent sheaf with integrable $T$-connection on $Z_T^\an$. Then $y$
lands in $\Ob(\calW_{T})$ if and only if there exist
\begin{itemize}
 \item an \'etale covering $\{T_i=\Spa(A_i,A_i^+)\rightarrow T\}$ with $A_i$ complete,
 \item an open $L^+$-subalgebra $A_{i,0}$ of $A_i^+$ of topologically finite type, and
 \item an $A_{i,0}$-flat coherent sheaf with integrable $A_{i,0}$-connection $(\calE_i,\nabla_i)$ on $Z_{A_{i,0}}$
\end{itemize}
such that $(\calE,\nabla)|_{Z_{T_i}^\an}$ is isomorphic to the analytification of $(\calE_i,\nabla_i)|_{Z_{A_i}}$ for each $i$.
\end{prop}

\begin{proof}
This easily follows from the discussion so far and Definition~\ref{def:1-morphism being mono or epi}.
\end{proof}

\section{Crystals and isocrystals on crystalline site}\label{sec:crystals and isocrystals on crystalline site}
We recall crystals and isocrystals on the (log) crystalline site. The standard references are \cite{Berthelot-book, BO, Berthelot-Breen-Messing, Kato-log, hyodo-kato, Ogus-F-crystals, Olsson-crystalline, Beilinson-crystalline-period-map}. Note that we will only need the logarithmic crystalline site in which the base PD log ($p$-adic formal) scheme or algebraic space has the trivial log structure and the canonical PD-structure on $(p)$.

Let $k$ be a perfect field of characteristic $p$. We set $W\coloneqq W(k)$, $W_n\coloneqq W_n(k)$, and $K \coloneqq W[p^{-1}]$. For a $p$-adic formal scheme or algebraic space $T$ over $W$, set $T_n\coloneqq T\times_WW_n$ for each $n\geq 1$; if $T$ has a log structure $M_T$, then $M_{T_n}$ denotes the pullback log structure on $T_n$.

\subsection{Crystalline site and crystals}\label{sec:crystalline site and crystals in p nilpotent case}

Fix an algebraic space $S$ over $W$ in which $p$ is locally nilpotent. 
Equip $S$ with the trivial log structure and the canonical PD-structure on the ideal $(p)$ (on the small \'etale site) and use the same symbol to denote the resulting PD algebraic space. 
Let $(Z,M_Z)$ be a fine log algebraic space over $S$. Let 
\[
((Z, M_Z)/S)_\cris
\]
denote the \emph{small \'etale logarithmic crystalline site} of $(Z,M_Z)$ relative to $S$: an object is a pair $((U,M_U)\hookrightarrow (T,M_T),\gamma)$ where $U$ is an \'etale \emph{scheme} over $Z$ (with $M_U\coloneqq M_Z|_{U_\et}$), $(U,M)\hookrightarrow (T,M_T)$ is an exact closed immersion of fine log schemes\footnote{One could allow log algebraic spaces, but for simplicity, we only consider log schemes.} over $S$, and $\gamma$ is a PD-structure on $J_T\coloneqq\Ker(\calO_T\rightarrow \calO_U)$ compatible with the PD-structure on $(p)$. We often write $(U,T)$ for an object. Note that $J_T$ is locally nilpotent and thus $U\hookrightarrow T$ is a nil-immersion.
Morphisms are the obvious ones, and we consider \'etale topology as in \cite[(5.2)]{Kato-log}: a collection $\{(U_i,T_i)\rightarrow (U,T)\}$ of morphisms is a \emph{covering} if $T_i\rightarrow T$ is \'etale with $M_{T_i}=(M_T)|_{T_{i,\et}}$ and $U_i=U\times_TT_i$ for each $i$, and the maps $T_i\rightarrow U$ are jointly surjective.

Let $((Z,M_Z)/S)_\Ncris\subset ((Z,M_Z)/S)_\cris$ denote the full subcategory consisting of $(U,T)$ whose PD-ideal is PD-nilpotent, namely, $\calJ_T^{[m]}=0$ for some $m$. It becomes a site (called \emph{PD-nilpotent logarithmic crystalline site}) with respect to the induced \'etale topology. Note that the canonical PD-structure on  $pW_n$ is PD-nilpotent if and only if $p\neq 2$ by \cite[Cor.~I.3.2.4]{Berthelot-book}. Finally, we remark that, when $M_Z$ is trivial, we recover the small \'etale (PD-nilpotent) crystalline site\footnote{Strictly speaking, we consider PD-structures on quasi-coherent ideals on the \emph{small \'etale site}.}, for which we simply write $(Z/S)_{\mathrm{(N)cris}}$.

Giving a sheaf $\calF$ on $((Z, M_Z)/S)_\cris$ is equivalent to giving a sheaf $\calF_{(U,T)}\in \Sh(T_\et)=\Sh(U_\et)$ for every $(U,T)\in ((Z,M_Z)/S)_\cris$  together with a transition morphism $g_\calF^{-1}\colon g^{-1}_\et\calF_{(U,T)}\rightarrow \calF_{(U',T')}$ for every morphism $g\colon (U',T')\rightarrow (U,T)$ satisfying the standard compatibilities for the identity and composition. 
In particular, the association 
\[
((U,M)\hookrightarrow (T,M_T),\gamma)\mapsto \calO_T(T)
\]
defines a sheaf $\calO_{Z/S}$ of rings on $((Z, M_Z)/S)_\cris$. A \emph{crystal (of $\calO_{Z/S}$-modules)} is an sheaf $\calE$ of $\calO_{Z/S}$-modules such that the induced map
\[
g_\calF^\ast\colon g_\et^\ast\calF_{(U,T)}\coloneqq \calO_{T'}\otimes_{g^{-1}_\et\calO_T}g^{-1}_\et\calF_{(U,T)}\rightarrow \calF_{(U',T')}
\]
is an isomorphism for every morphism $g\colon (U,T')\rightarrow (U,T)$ in $((Z, M_Z)/S)_\cris$. We say that a crystal $\calF$ is \emph{quasi-coherent (of finite type)} if each $\calF_{(U,T)}$ is a quasi-coherent sheaf of $\calO_T$-modules (of finite type). Similarly, we define the property of a crystal of $\calO_{Z/S}$-modules being \emph{of finite presentation} or \emph{locally free of finite rank}.

The log crystalline topoi enjoy the following functoriality: for another algebraic space $S'$ over $W$ and a fine log algebraic space $(Z',M_{Z'})$ over $S'$ together with the commutative diagram
\begin{equation}\label{eq:crystalline functoriality diagmra}
 \xymatrix{
(Z',M_{Z'})\ar[r]^f\ar[d] & (Z,M_Z)\ar[d]\\
S'\ar[r] & S,
}
\end{equation}
we have a morphism of topoi
\[
 f_\cris=(f_\cris^{-1},f_{\cris,\ast})\colon \Sh(((Z',M_{Z'})/S')_\cris)\rightarrow \Sh(((Z,M_Z)/S)_\cris).
\]
See \cite[(5.9)]{Kato-log} for the detail.\footnote{The discussion therein also works when $S$ and $S'$ are log algebraic spaces}
It comes with a morphism $f_\cris^{-1}\calO_{Z/S}\rightarrow\calO_{Z'/S'}$, which induces the pullback functor
\[
f_\cris^\ast\colon \Sh(((Z,M_Z)/S)_\cris,\calO_{Z/S}) \rightarrow \Sh(((Z',M_{Z'})/S')_\cris,\calO_{Z'/S'})
\]
sending $\calF$ to $f_\cris^\ast\calF\coloneqq \calO_{Z'/S'}\otimes_{f_\cris^{-1}\calO_{Z/S}}f_\cris^{-1}\calF$. Note that $f_\cris^\ast$ sends crystals of $\calO_{Z/S}$-modules to crystals of $\calO_{Z'/S'}$-modules and preserves the property of being quasi-coherent of finite type, of finite presentation, or locally free of finite rank (cf.~\cite[Cor.~IV.1.2.4]{Berthelot-book}).

The crystals on $(Z,M_Z)$ only depend on the mod $p$ fiber $(Z_1,M_{Z_1})$:

\begin{prop}\label{prop:crystals only depend on mod p fiber}
The pullback functor induces an equivalences of categories from the category of crystals of $\calO_{Z/S}$-modules to that of crystals of $\calO_{Z_1/S}$-modules.
\end{prop}

\begin{proof}
The proof of \cite[Thm.~IV.1.4.1]{Berthelot-book} also works in the logarithmic setup thanks to the existence of the log PD-envelope proved in \cite[Def.~5.4]{Kato-log}.
\end{proof}

Let us turn to the description of crystals in terms of quasi-nilpotent integrable connections. See \cite[Def.~4.10]{BO}, \cite[Thm.~6.2(iii)]{Kato-log}, or \cite[p.~19]{Ogus-F-crystals} for the definition of quasi-nilpotence of an integrable connection.

\begin{prop}\label{prop:crystal and integral connections}
Assume that there exists a log algebraic space $(Y,M_Y)$ smooth\footnote{Smoothness as in \cite[(3.3)]{Kato-log}, which is also referred to as log smoothness.} over $S$ together with an $S$-closed immersion $(Z,M_Z)\hookrightarrow(Y,M_Y)$. Let $(\calD,M_\calD)$ denote the log PD-envelope of $(Z,M_Z)$ in $(Y,M_Y)$. Then the evaluation on $(Z,\calD)$ gives rise to an equivalence of categories from the category of crystals of $\calO_{Z/S}$-modules to the category of $\calO_\calD$-modules with quasi-nilpotent integrable $S$-connections.

Moreover, it restricts to an equivalence between the full subcategories of objects with each of the following properties: quasi-coherent (of finite type); of finite presentation; locally free of finite rank.
\end{prop}

Let us give a few remarks before the proof: since $S$ has the trivial log structure, $(Y,M_Y)\rightarrow S$ is also integral. In particular, the underlying morphism $Y\rightarrow S$ is flat by \cite[Cor.~4.5]{Kato-log}; since the canonical PD-structure on $p\calO_S$ extends to $Y$, the pair $(Z,\calD)$ defines an object of $((Z,M_Z)/S)_\cris$ by \cite[(5.5.2)]{Kato-log}; if $(Y,M_Y)=(Z,M_Z)$ itself is smooth over $S$, then the log PD-envelope $(\calD,M_\calD)$ is $(Z,M_Z)$.

\begin{proof}
This follows from \cite[Thm.~6.2]{Kato-log} in the case where $S$ and $Y$ are schemes, and one can deduce the general case from the scheme case by descent: the descent of crystals is straightforward; for the descent of connections, one uses the stratification interpretation as in \cite[Lem.~2.2.14]{Olsson-crystalline}) in the case of trivial log structure, and the same argument works for the genera log structure thanks to \cite[(6.7)]{Kato-log}. 
\end{proof}

\begin{rem}\label{rem:compatibility of functoriality for crystals and connections}
Consider the functoriality setup as in \eqref{eq:crystalline functoriality diagmra} and suppose further that the assumption of Proposition~\ref{prop:crystal and integral connections} holds for both $(Z,M_Z)$ and $(Z',M_{Z'})$ together with the commutative diagram
\[
\xymatrix{
(Z',M_{Z'})\ar[r]\ar[d]_-f & (\calD',M_{\calD'}) \ar[r]\ar[d]_-{f_\calD}  & (Y',M_{Y'}) \ar[r]\ar[d] & S'\ar[d] \\
(Z,M_Z)\ar[r] & (\calD,M_\calD) \ar[r]  & (Y,M_Y) \ar[r] & S.
}
\]
Let $\calF$ be a crystal of $\calO_{Z/S}$-modules and let $(\calF_\calD,\nabla)$ denote the corresponding $\calO_\calD$-module with quasi-nilpotent integrable $S$-connection. By unwinding the proof of Proposition~\ref{prop:crystal and integral connections}, one can see that the crystal $f^\ast_\cris\calF$ of $\calO_{Z'/S'}$-module corresponds to the $\calO_{\calD'}$-module with quasi-nilpotent integrable $S'$-connection $f_\calD^\ast(\calF_\calD,\nabla)$ (cf.~\cite[Lem.~IV.1.6.1(ii)]{Berthelot-book}).
\end{rem}

\begin{lem}\label{lem:S-flatness of quasi-coherent crystal}
Assume that, for $i=1,2$, $(Z,M_Z)$ admits an $S$-closed immersion to a log algebraic space $(Y_{i},M_{Y_{i}})$ smooth over $S$ with log PD-envelope $(\calD_i,M_{\calD_i})$ as in Proposition~\ref{prop:crystal and integral connections}.
 Let $\calF$ be a quasi-coherent crystal of $\calO_{Z/S}$-modules and write $(\calF_i,\nabla_i)$ for the corresponding quasi-coherent $\calO_{\calD_{i}}$-module with quasi-nilpotent integrable $S$-connection. Then $\calF_1$ is flat over $\calO_S$ if and only if $\calF_2$ is flat over $\calO_S$.
\end{lem}

\begin{proof}
Let $(Y,M_Y)\coloneqq (Y_{1},M_{Y_{1}})\times_S(Y_{2},M_{Y_{2}})$ and write $(\calD,M_\calD)$ for the log PD-envelope of $(Z,M_Z)$ in $(Y,M_Y)$. Then one can see that the projection $\pi_i\colon (\calD,M_\calD)\rightarrow (\calD_i,M_{\calD_i})$ is faithfully flat since 
$\calO_\calD$ is locally a PD-polynomial algebra over $\calO_{\calD_i}$
(cf.~\cite[Lem.~2.22, Pf.]{hyodo-kato}). It follows from  Proposition~\ref{prop:crystal and integral connections} that $\pi_1^\ast \calF_1\cong\pi_2^\ast \calF_2$ on $\calD$, which implies the assertion.
\end{proof}

\begin{defn}\label{def:S-flat crystals}
Assume that $(Z,M_Z)$ admits an $S$-closed immersion to a log algebraic space $(Y,M_{Y})$ smooth over $S$ with log PD-envelope $(\calD,M_{\calD})$. We say that a quasi-coherent crystal of $\calO_{Z/S}$-modules is \emph{$S$-flat} if the underlying quasi-coherent $\calO_{\calD}$-module of the corresponding $S$-connection is flat over $\calO_S$. By Lemma~\ref{lem:S-flatness of quasi-coherent crystal}, this condition is independent of the choice of such $(Y,M_Y)$. In fact, we can and do define $S$-flatness when such a $(Y,M_Y)$ exists \'etale locally, thanks to fpqc descent.
\end{defn}

Let us also recall the relation between the quasi-nilpotence of an integrable connection and the $p$-curvature of its reduction modulo $p$ introduced in Definition~\ref{def:p-curvature}.

\begin{lem}\label{lem:quasi-nilpotence and nilpotence of p-curvature}
Assume that $(Z,M_Z)$ is smooth over $S$. Let $(\calF,\nabla\colon \calF\rightarrow \calF\otimes_{\calO_Z}\omega^1_{(Z,M_Z)/S})$ be an $\calO_Z$-module with integrable $S$-connection, and write $(\calF_1,\nabla_1)\coloneqq (\calF,\nabla)|_{Z_1}$ for its mod $p$ reduction. Then the following conditions are equivalent:
\begin{enumerate}
 \item $(\calF,\nabla)$ is quasi-nilpotent;
 \item $(\calF_1,\nabla_1)$ is quasi-nilpotent;
 \item $(\calF_1,\nabla_1)$ is nilpotent in the sense of Definition~\ref{def:p-curvature}.
\end{enumerate}
\end{lem}

\begin{proof}
This easily follows from the definitions and \cite[Cor.~5.5]{Katz_Nilpotent} when $Z$ is a scheme with trivial log structure. The general case can be proved similarly by using \cite[Rem.~1.2.2]{Ogus-F-crystals} and (a logarithmic version of) \cite[Lem.~2.3.25]{Olsson-crystalline}.
\end{proof}

\subsection{The case of a \texorpdfstring{$p$}{p}-adic base}\label{sec:The case of a $p$-adic base}

We extend the above discussion to the $p$-adic base: let $S$ be a $p$-adic formal algebraic space over $W$ and equip it with the trivial log structure and the canonical PD-structure on $(p)$. Let $(Z,M_Z)$ be a fine log algebraic space over $S_{n_0}$ for some $n_0$. For $n'\geq n\geq n_0$, we have an obvious fully faithful embedding $((Z,M_Z)/S_{n})_\cris\hookrightarrow ((Z,M_Z)/S_{n'})_\cris$. Let 
\[
((Z,M_Z)/S)_\cris
\]
denote the colimit of the sites $((Z,M_Z)/S_{n})_\cris$ ($n\geq n_0$) 
and call it the \emph{small \'etale logarithmic crystalline site}: each object of $((Z,M_Z)/S)_\cris$ comes from an object of $((Z,M_Z)/S_{n})_\cris$ for some $n\geq n_0$ (cf.~\cite[Def.~7.17]{BO} and \cite[\S2.7]{Olsson-crystalline}). When $M_Z$ is trivial, we simply write $(Z/S)_\cris$.

We define the structure sheaf $\calO_Z/S$ and the notion and properties of crystals of $\calO_{Z/S}$-modules in an obvious way (an exception is the definition of $S$-flatness, which is given below). 
It follows from Proposition~\ref{prop:crystals only depend on mod p fiber} that the category of crystals of $\calO_{Z/S}$-modules is equivalent to the category of crystals of $\calO_{Z_1/S}$-modules.

\begin{prop}\label{prop:crystals and connections in p-adic case}
Assume that there exists a log $p$-adic formal algebraic space $(Z_S,M_{Z_S})$ that is adic and smooth\footnote{By descent and Remark~\ref{rem:adically P vs P}, this is equivalent to the condition that $(Z_{S,n},M_{Z_{S,n}})\rightarrow S_n$ is smooth for every $n\geq n_0$.} over $S$ and satisfies $(Z,M_Z)=(Z_S,M_{Z_S})\times_SS_{n_0}$. Then there is a natural equivalence from the category of quasi-coherent crystals of $\calO_{Z/S}$-modules to the category of a.q.c.~$\calO_{Z_S}$-modules with topologically quasi-nilpotent integrable $S$-connections. Moreover, it restricts to an equivalence between the full subcategories of objects of finite type (resp.~of finite presentation, resp.~locally free of finite rank).
\end{prop}

An integrable connection on an a.q.c.~$\calO_Z$-module is said to be \emph{topologically quasi-nilpotent} if its base change to $Z_n$ is quasi-nilpotent for each $n$. This is equivalent to the condition that its mod $p$ reduction is nilpotent by Lemma~\ref{lem:quasi-nilpotence and nilpotence of p-curvature}.

\begin{proof}
Recall our notation $(Z_{S,n},M_{Z_{S,n}})\coloneqq (Z_S,M_{Z_S})\times_WW_n=(Z_S,M_{Z_S})\times_SS_n$. It follows from Proposition~\ref{prop:crystal and integral connections} that a quasi-coherent crystal of $\calO_{Z/S}$-module corresponds to a projective system $(\calF_n,\nabla_n)_{n\geq n_0}$ where $(\calF_n,\nabla_n)$ is a quasi-coherent sheaf with quasi-nilpotent $S_n$-connection on $Z_{S,n}$ such that $(\calF_n,\nabla_n)\otimes_{W_n}W_m=(\calF_m,\nabla_m)$ for $n\geq m\geq n_0$ under the identifications $(Z_S)_\et\cong (Z_{S,n})_\et\cong (Z_{S,m})_\et$. Then $(\calF,\nabla)\coloneqq \varprojlim_n (\calF_n,\nabla_n)$ is an a.q.c.~$\calO_{Z_S}$-module with topologically quasi-nilpotent integrable $S$-connection and satisfies $(\calF,\nabla)\otimes_WW_n=(\calF_n,\nabla_n)$: this follows from \cite[Prop.~I.3.4.1]{Fujiwara-Kato} and Proposition~\ref{prop:aqc sheaves on Zariski site and etale site} in the formal scheme case, and the general case is deduced from the formal scheme case by descent.
One can verify that this construction gives the desired equivalence and preserves the property of being of finite type and of finite presentation.
\end{proof}

\begin{defn}
Assume that $(Z,M_Z)\rightarrow S_{n_0}$ lifts, \'etale locally, to a log $p$-adic formal algebraic space $(Z_S,M_{Z_S})\rightarrow S$ that is adic and smooth as in Proposition~\ref{prop:crystals and connections in p-adic case}. Let $\calF$ be a quasi-coherent crystal of $\calO_{Z/S}$-modules and write $(\calF_{Z_S},\nabla)$ for the corresponding a.q.c.~$\calO_{Z_S}$-module with topologically quasi-nilpotent integrable $S$-connection.
We say that $\calF$ is \emph{$S$-flat} if $\calF_{Z_S}$ is adically flat over $\calO_S$ in the sense that $\calF_{Z_S}\otimes_WW_n$ is a flat $\calO_{S_n}$-module for every $n$. It follows from Lemma~\ref{lem:S-flatness of quasi-coherent crystal} and the fpqc descent that this condition is independent of the choice of $(Z_S,M_{Z_S})$.
\end{defn}

Let $S'$ be a $p$-adic formal algebraic space over $W$ with a morphism $S'\rightarrow S$. Suppose that we are given a fine log algebraic space $(Z',M_{Z'})$ over $S'_{n'_0}$ (for some $n'_0\leq n_0$) together with the commutative diagram
\[
 \xymatrix{
(Z',M_{Z'})\ar[r]^f\ar[d] & (Z,M_Z)\ar[d]\\
S'\ar[r] & S.
}
\]
The functoriality in \S\ref{sec:crystalline site and crystals in p nilpotent case} yields a morphism of topoi
\[
 f_\cris=(f_\cris^{-1},f_{\cris,\ast})\colon \Sh(((Z',M_{Z'})/S')_\cris)\rightarrow \Sh(((Z,M_Z)/S)_\cris).
\]
It comes with a morphism $f_\cris^{-1}\calO_{Z/S}\rightarrow\calO_{Z'/S'}$, which induces the pullback functor
\[
f_\cris^\ast\colon \Sh(((Z,M_Z)/S)_\cris,\calO_{Z/S}) \rightarrow \Sh(((Z',M_{Z'})/S')_\cris,\calO_{Z'/S'})
\]
sending $\calF$ to $f_\cris^\ast\calF\coloneqq \calO_{Z'/S'}\otimes_{f_\cris^{-1}\calO_{Z/S}}f_\cris^{-1}\calF$. Note that $f_\cris^\ast$ sends crystals of $\calO_{Z/S}$-modules to crystals of $\calO_{Z'/S'}$-modules and preserves the property of being quasi-coherent of finite type, of finite presentation, or locally free of finite rank.

\begin{lem}\label{lem:pulback stability of S-flatness of crystals}
Suppose that $(Z,M_Z)\rightarrow S_{n_0}$ lifts to a log $p$-adic formal algebraic space $(Z_S,M_{Z_S})\rightarrow S$ that is adic and smooth. Assume further one of the following:
\begin{enumerate}
\item the induced morphism $(Z',M_{Z'})\rightarrow (Z_S,M_{Z_S})\times_S S'_{n_0'}$ is an isomorphism;
\item $S=S'$ and \'etale locally, $(Z',M_{Z'})\rightarrow S_{n'_0}$ lifts to an adic and smooth log $p$-adic formal algebraic space $(Z'_S,M_{Z'_S})\rightarrow S$ and $f$ lifts to $f_S\colon (Z'_S,M_{Z'_S})\rightarrow (Z_S,M_S)$ such that $f_{S,n}\colon Z'_{S,n}\rightarrow Z_{S,n}$ is flat for every $n$. 
\end{enumerate}
Then for every $S$-flat quasi-coherent crystal $\calF$ of $\calO_{Z/S}$-modules, the quasi-coherent crystal $f^\ast_\cris\calF$ of $\calO_{Z'/S'}$-modules is $S'$-flat.
\end{lem}

\begin{proof}
Note that the notion of $S'$-flatness is defined in either case. The assertion follows from the definition of $S$-flatness, Remark~\ref{rem:compatibility of functoriality for crystals and connections}, and \cite[0584, 00HI]{stacks-project}. 
\end{proof}

\subsection{Isocrystals}

Let $S$ be a $p$-adic formal algebraic space over $W$ that is \emph{flat} over $W$.
Fix a morphism $(Z_S,M_{Z_S})\rightarrow S$ of log $p$-adic formal algebraic spaces that is adic and smooth, and set $(Z,M_Z)\coloneqq (Z_S,M_{Z_S})\times_Wk= (Z_S,M_{Z_S})\times_SS_1$.

\begin{defn}
The category $\Isoc((Z,M_Z)/S)$ of \emph{isocrystals} on $((Z,M_Z)/S)_\cris$ is defined to be the isogeny category of the category of $S$-flat crystals of $\calO_{Z/S}$-modules of finite presentation. Let $\Isoc^\lf((Z,M_Z)/S)$ denote the full subcategory given by the essential image of locally free crystals of $\calO_{Z/S}$-modules of finite rank. When we fix $S$ and $(Z,M_Z)$, these categories for different choices of $(Z_S,M_{Z_S})$ are equivalent to each other by Proposition~\ref{prop:crystals only depend on mod p fiber} and Lemma~\ref{lem:S-flatness of quasi-coherent crystal}.
We simply write $\Isoc^{(\lf)}(Z/S)$ if $M_Z$ is trivial.
\end{defn}

\begin{rem}
We put the $S$-flatness condition in the definition of isocrystals to consider a reasonable moduli problem; see Theorem~\ref{thm:moduli interpretation of Mcris}. 
\end{rem}

\begin{rem}
Assume $S=\Spf \calO_L$ for a finite extension $L$ over $K$. Then one may drop the $S$-flatness condition. In fact, every crystal $\calF$ of $\calO_{Z/S}$-module of finite presentation admits an $S$-flat quotient which gives rise to the same object as $\calF$ after passing to the isogeny category; to see this, use Proposition~\ref{prop:crystals and connections in p-adic case} and consider the $p$-torsion free quotient of the corresponding a.q.c.~$\calO_Z$-module with topologically quasi-nilpotent integrable $S$-connection.
Moreover, one can canonically associate to each isocrystal $\calF$ an $\calO_{Z_S}[p^{-1}]$-module with integrable connection $(\calF_S,\nabla)$; it is known that $\calF_S$ is a locally free $\calO_{Z_S}[p^{-1}]$-module of finite rank (cf.~\cite[Prop.~2.4]{Drinfeld-stacky}). We will refer to the \emph{rank} of $\calF$ as the rank of $\calF_S$.
\end{rem}

\begin{prop}\label{prop:Esnault-Shiho}
Assume $S=\Spf \calO_L$ for a finite extension $L$ over $K$. Every rank one isocrystal on $((Z,M_Z)/S)_\cris$ belongs to $\Isoc^\lf((Z,M_Z)/S)$. If $\dim Z=1$, then $\Isoc^\lf((Z,M_Z)/S)=\Isoc((Z,M_Z)/S)$.
\end{prop}

\begin{proof}
When $L=K$, these two statements are proved in \cite[Prop.~2.11]{Esnault-Shiho} and \cite[1.3]{Esnault-Shiho-119b}\footnote{Note that [1, Prop.~2.10] therein should read [1, Prop.~2.11].}, respectively. These proofs work verbatim in the general case.
\end{proof}

\begin{prop}\label{prop:pullback stability of isocrystals}
Let $S'$ be a $p$-adic formal algebraic space flat over $W$ and $(Z'_{S'},M_{Z'_{S'}})\rightarrow S'$ an adic and smooth morphism of log $p$-adic formal algebraic spaces fitting in the following commutative diagram
\[
 \xymatrix{
(Z'_{S'},M')\ar[r]^-{f_S}\ar[d] & (Z_S,M_S)\ar[d]\\
S'\ar[r] & S.
}
\]
Let $f$ denote the induced morphism $f_S\times_Wk\colon (Z',M_{Z'})\rightarrow (Z,M_Z)$.
\begin{enumerate}
 \item The pullback $f_\cris^\ast$ induces a functor $\Isoc^\lf((Z,M_Z)/S)\rightarrow\Isoc^\lf((Z',M_{Z'})/S')$.
 \item If $Z'\rightarrow Z\times_{S}S'$ is flat, then $f_\cris^\ast$ induces a functor $\Isoc((Z,M_Z)/S)\rightarrow\Isoc((Z',M_{Z'})/S')$.
\end{enumerate}
\end{prop}

\begin{proof}
Note that $f_\cris^\ast$ preserves the property of two quasi-coherent crystals being isogenous. Since it also preserves local freeness of finite rank, part (i) follows. For (ii), we claim that $Z'_{S',n}\coloneqq Z'_{S'}\times_WW_n\rightarrow Z_S\times_SS'_n$ is flat for every $n\geq 1$: to see this, first observe that $(Z_S,M_{Z_S})\times_SS'\rightarrow S'$ is a smooth morphism of log $p$-adic formal algebraic spaces with $S'$ having trivial log structure. Hence $Z_S\times_SS'$ is flat over $W$ (cf.~a remark after Proposition~\ref{prop:crystal and integral connections}), and the claim is deduced from our assumption (the case $n=1$) by \cite[Prop.~7.13]{Ogus-crystalline-prism}.  By writing $f_S$ as the composite $(Z'_{S'},M_{Z'_{S'}})\rightarrow (Z_S,M_{Z_S})\times_SS'\rightarrow (Z_S,M_{Z_S})$, we deduce from the claim and Lemma~\ref{lem:pulback stability of S-flatness of crystals} that $f_\cris^\ast$ also preserves $S$-flatness. Hence (ii) follows.
\end{proof}

\subsection{Frobenius and \texorpdfstring{$F$}{F}-isocrystals}\label{sec:Frobenius pullback of crystals}

We consider the above general framework in the following case:
Let $(Z_W,M_{Z_W})$ be a fine log scheme over $W$ such that the underlying morphism $Z_W\rightarrow \Spec W$ is smooth. Set $(Z,M_Z)\coloneqq (Z_W,M_{Z_W})\times_Wk$ and write $(Z_W^\wedge,M_{Z_W^\wedge})$ for the $p$-adic completion of $(Z_W,M_{Z_W})$.
For a $p$-adic formal algebraic space $S$ over $W$ (with trivial log structure), set $(Z_S,M_{Z_S})\coloneqq (Z_W^\wedge,M_{Z_W}^\wedge)\times_{\Spf W}S$. Note $Z_{S_1}=Z_W\times_WS_1=Z\times_kS_1$.

Assume, moreover, that $k=\F_q$ is a finite field with $q=p^a$. Then the $a$th iterate of the absolute Frobenius $\Fr_Z\colon (Z,M_Z)\rightarrow (Z,M_Z)$ is a $k$-morphism. Note that $\Fr_Z\colon Z\rightarrow Z$ is flat by \cite[0EC0]{stacks-project} since $Z\rightarrow \Spec k$ is smooth.
We have the commutative diagram
\begin{equation}\label{eq:diagram for Frobenius pullback on crystalline site}
\xymatrix@C+2pc{
(Z,M_Z)\times_{\F_q}S_1 \ar[d]\ar[r]^{\Fr_Z^a\times \id_{S_1}}
& (Z,M_Z)\times_{\F_q}S_1\ar[d]\\
S\ar@{=}[r] &S, 
}    
\end{equation}
and this induces a morphism of topoi
\begin{equation}\label{eq:Frobenius on crystallie site}
\Fr^a_{Z/S,\cris}\colon 
\Sh(((Z_{S_1},M_{Z_{S_1}})/S)_\cris)\rightarrow \Sh(((Z_{S_1},M_{Z_{S_1}})/S)_\cris),    
\end{equation}
which is functorial in $S$. For a crystal $\calF$ of $\calO_{Z_{S_1}/S}$-modules on $((Z_{S_1},M_{Z_{S_1}})/S)_\cris$, the pullback $\Fr_{Z/S,\cris}^{a,\ast} \calF$ is a crystal of $\calO_{Z_{S_1}/S}$-modules on $((Z_{S_1},M_{Z_{S_1}})/S)_\cris$. If $\calF$ is locally finite free, then so is $\Fr_{Z/S,\cris}^{a,\ast} \calF$. In the case when $S$ is flat over $W$, Proposition~\ref{prop:pullback stability of isocrystals} also gives
\[
\Fr_{Z/S,\cris}^{a,\ast}\colon \Isoc^\lf((Z_{S_1},M_{Z_{S_1}})/S)\rightarrow\Isoc^\lf((Z_{S_1},M_{Z_{S_1}})/S).
\]
\[
\Fr_{Z/S,\cris}^{a,\ast}\colon \Isoc((Z_{S_1},M_{Z_{S_1}})/S)\rightarrow\Isoc((Z_{S_1},M_{Z_{S_1}})/S).
\]

\begin{lem}\label{lem:Frobenius pullback functoriality}
For a morphism $f\colon S'\rightarrow S$ of $p$-adic formal algebraic spaces over $W$, we have $f_\cris^\ast\circ \Fr_{Z/S,\cris}^{a,\ast}\cong \Fr_{Z/S',\cris}^{a,\ast}\circ f_\cris^\ast$, where $f_\cris^\ast$ denotes the pullback functor $\Sh(((Z_{S_1},M_{Z_{S_1}})/S)_\cris)\rightarrow \Sh(((Z_{S'_1},M_{Z_{S'_1}})/S')_\cris)$.
\end{lem}

\begin{proof}
This follows from $(\id_Z\times f_1)\circ(\Fr_Z^a\times \id_{S'_1})=(\Fr_Z^a\times \id_{S_1})\circ (\id_Z\times f_1)$.
\end{proof}

\begin{lem}\label{lem:S-linearity of Frobenius pullback}
The Frobenius pullback functor $\Fr_{Z/S,\cris}^{a,\ast}$ is $\Gamma(S,\calO_S)$-linear on the category of crystals of $\calO_{Z_{S_1}/S}$-modules as well as on
$\Isoc((Z_{S_1},M_{Z_{S_1}})/S)$.
\end{lem}

\begin{proof}
This is clear from the definition.
\end{proof}

Let $\overline{\Q}_p$ denote an algebraic closure of $K$.

\begin{defn}\label{defn:L-F-isocrystals}
An \emph{$F$-isocrystal} on (the crystalline site of) $(Z,M_Z)/W$ is an isocrystal $\calF$ on $((Z,M_Z)/W)_\cris$ together with an isomorphism $\varphi_\calF\colon \Fr^{a,\ast}_{Z/W,\cris}\calF\xrightarrow{\cong}\calF$; a morphism $(\calF,\varphi_\calF)\rightarrow (\calF',\varphi_{\calF'})$ of $F$-isocrystals is a morphism $\calF\rightarrow\calF'$ of isocrystals that is compatible with $\varphi_\calF$ and $\varphi_{\calF'}$.
More generally, for a finite extension $L$ of $K=W[p^{-1}]$, an \emph{$L$-$F$-isocrystal} on $(Z,M_Z)/W$ is an $F$-isocrystal $(\calF,\varphi_\calF)$ together with a $K$-algebra map $\iota_{\calF,L}\colon L\rightarrow \End((\calF,\varphi_\calF))$. 
Write $\FIsoc(Z,M_Z)_L$ for the category of $L$-$F$-isocrystals. For a finite extension $L'$ of $L$, it is straightforward to define the scalar extension functor $\FIsoc(Z,M_Z)_L\rightarrow \FIsoc(Z,M_Z)_{L'}$ sending $(\calF,\varphi_\calF,\iota_{\calF,L})$ to $(\calF\otimes_LL',\varphi_\calF\otimes_LL',\iota_{\calF,L}\otimes_LL')$\footnote{The underlying isocrystal is isomorphic to $\calF^{\oplus [L':L]}$. See also \cite[1.4.1]{Abe-Langlands}}. A \emph{$\overline{\Q}_p$-$F$-isocrystal} refers to an object of the $2$-colimit $\varinjlim_{L}\FIsoc(Z,M_Z)_L$.

An \emph{$F$-crystal} on  $(Z,M_Z)/W$ is a quasi-coherent crystal $\calF$ of $\calO_{Z/W}$-modules of finite presentation together with an isomorphism $\varphi_\calF\colon \Fr^{a,\ast}_{Z/W,\cris}\calF\xrightarrow{\cong}\calF$. We similarly define the notions of morphisms of $F$-crystals, $L$-$F$-crystals, and $\overline{\Z}_p$-$F$-crystals.  \end{defn}

\begin{rem}\label{rem:two ways to define L-F-isocrystals}
One can also define an $L$-$F$-isocrystal to be an isocrystal  $\calF'$ on $((Z,M_Z)/\calO_L)_\cris$ together with an isomorphism $\varphi_{\calF'}\colon \Fr^{a,\ast}_{Z/\calO_L,\cris}\calF'\xrightarrow{\cong}\calF'$. The equivalence of these two definitions follows from Lemma~\ref{lem:S-linearity of Frobenius pullback} and the fact that every $K$-homomorphism $L\rightarrow \End(\calF)$ for an isocrystal $\calF$ on $((Z,M_Z)/W)_\cris$ arises from a $W$-homomorphism $\calO_L\rightarrow \End(\calF_0)$ for some crystal $\calF_0$ underlying $\calF$. The latter is standard and is also proved in the second half of the proof of Corollary~\ref{cor:moduli interpretation of classical points of Misoc}.
\end{rem}

\begin{defn}
Let $(\calF,\varphi_{\calF},\iota_{\calF,L})$ be an $L$-$F$-isocrystal. Define the $i$-th crystalline cohomology $H^i_\cris((Z,M_Z),\calF)$ by
\[
H^i_\cris((Z,M_Z),\calF)\coloneqq H^i(((Z,M_Z)/W)_\cris,\calF_0)[p^{-1}],
\]
where $\calF_0$ is a crystal underlying $\calF$ such that $\iota_{L,\calF}$ restricts to a $W$-algebra map $\calO_L\rightarrow \End(\calF_0)$. It is clear that this is independent of the choice of $\calF_0$. Moreover, $H^i_\cris((Z,M_Z),\calF)$ is an $L$-vector space by $\iota_{L,F}$ and is equipped with an $L$-linear automorphism $\varphi$ given by
\[
H^i_\cris((Z,M_Z),\calF)\xrightarrow{\Fr_{Z/W,\cris}^{a,\ast}}
 H^i_\cris((Z,M_Z),\Fr_{Z/W,\cris}^{a,\ast}\calF) 
\xrightarrow{H^i_\cris(\varphi_\calF)} H^i_\cris((Z,M_Z),\calF).
\]
For a finite extension $L'$ of $L$, we have an $L'$-linear isomorphism 
$H^i_\cris((Z,M_Z),\calF)\otimes_LL'\xrightarrow{\cong}H^i_\cris((Z,M_Z),\calF\otimes_LL')$ that is compatible with $L'$-linear Frobenii. 

Hence to a $\overline{\Q}_p$-$F$-isocrystal $\calG$, we associate an $\overline{\Q}_p$-vector space $H^i_\cris((Z,M_Z),\calG)$ together with a $\overline{\Q}_p$-linear Frobenius $\varphi$ by taking an underlying $L$-$F$-isocrystal $\calF$ and setting $H^i_\cris((Z,M_Z),\calG)\coloneqq H^i_\cris((Z,M_Z),\calF)\otimes_L\overline{\Q}_p$.

If $M_Z$ is trivial, we often write $H^i_\cris(Z,\calF)$ for simplicity.
\end{defn}

For the rest of the subsection, we assume that $Z_W\rightarrow \Spec W$ is smooth proper and $M_{Z_W}$ is associated to a relative normal crossings divisor $D_W\rightarrow \Spec W$.

\begin{rem}\label{rem:F-isocrystals and overconvergent/convergent ones}
When $M_Z$ is trivial (i.e., $D=\emptyset$), our category of $F$-isocrystals agrees with the category of overconvergent $F$-isocrystals or the category of convergent $F$-isocrystals in the literature by the Dwork trick (see \cite[2.3.7, Thm.~2.4.2]{Berthelot-rigid}); our category of $\overline{\Q}_p$-$F$-isocrystals agrees with the category $\Isoc^\dagger(Z/(\overline{\Q}_p)_F)$ in \cite[p.~929, 1.4.10]{Abe-Langlands} (see also \cite[1.2.14]{Abe-Caro-weights}).

In general, our category of $F$-isocrystals is fully faithfully embedded into the category of overconvergent $F$-isocrystals on $Z\smallsetminus D$; this follows from the log version of \cite[Thm.~2.4.2]{Berthelot-rigid} and \cite[Thm.~6.4.5]{Kedlaya-semistablereduction1} (cf.~\cite[\S7]{Kedlaya-notesonisocrystals}).
\end{rem}

As stated in the next theorem, the crystalline cohomology of an $F$-isocrystal coincides with the rigid cohomology of the associated overconvergent $F$-isocrystal on $Z\smallsetminus D$.
For an $F$-isocrystal $\calF$ on $(Z,M_Z)/W$, take a crystal $\calF_0$ underlying $\calF$. Then $\calF_0$ corresponds to a coherent $\calO_{Z_W^\wedge}$-module $(\widehat{\calF}_0,\widehat{\nabla}_0)$ with topologically quasi-nilpotent integrable $W$-connection on $Z_W^\wedge$ by Proposition~\ref{prop:crystals and connections in p-adic case}. As in Proposition~\ref{prop:analytification and completion of coheren sheaves and connections} and the preceding paragraph, $(\widehat{\calF}_0,\widehat{\nabla}_0)[p^{-1}]$ defines a coherent sheaf $(\calF_K^\an,\nabla_K^\an)$ with integrable $K$-connection on the analytification $(Z_K^\an,M_{Z_K^\an})$ of the generic fiber $(Z_K,M_{Z_K})$ of $(Z_W,M_{Z_W})$, which is independent of the choice $\calF_0$.

\begin{thm}\label{thm:rigid-crystalline comparison for F-isocrystals}
Let $\calF$ be an $F$-isocrystal on $(Z,M_Z)/W$. Write $\calF_\rig$ for the associated overconvergent $F$-isocrystal on $Z\smallsetminus D$ and $(\calF_K^\an,\nabla_K^\an)$ for the associated coherent sheaf with integrable $K$-connection on $(Z_K^\an,M_{Z_K^\an})$. 
There there are natural $K$-linear isomorphisms
\[
H^i_\rig(Z\smallsetminus D,\calF_\rig)\cong H^i_\cris((Z,M_Z),\calF)\cong H^i_\dR((Z_K^\an,M_{Z_K^\an}), (\calF_K^\an,\nabla_K^\an)).
\]
Moreover, the first isomorphism is compatible with Frobenii. Similar statements hold for $L$-$F$-isocrystals and $\overline{\Q}_p$-$F$-isocrystals.
\end{thm}

\begin{proof}
This can be deduced from \cite[Cor.~2.4.13, Thm.~3.1.1]{Shiho-I}; the former corollary concerns the comparison between rigid cohomology and log-convergent cohomology, and the latter theorem concerns the comparison between log-convergent cohomology and log-crystalline cohomology (which is our crystalline cohomology). In both of the comparison results, the cohomology groups are shown to be naturally isomorphic to the analytic cohomology defined in \cite[Def.~2.2.12]{Shiho-II}, which coincides with $H^i_\dR((Z_K^\an,M_{Z_K^\an}), (\calF_K^\an,\nabla_K^\an))$ since $Z_W$ is proper over $W$. In fact, the comparison results with the analytic cohomology hold without the properness assumption and satisfy \'etale hyperdescent. When an \'etale scheme $Z'$ over $Z$ lifts to $\Spec W$ together with a Frobenius lift, one can check that these isomorphisms for $Z'$ are compatible with Frobenii, from which one can deduce that the isomorphism $H^i_\rig(Z\smallsetminus D,\calF_\rig)\cong H^i_\cris((Z,M_Z),\calF)$ is compatible with Frobenii.
\end{proof}

\begin{rem}\label{rem:two de Rham cohomologies for L-F-isocrystals}
With the notation as in Remark~\ref{rem:two ways to define L-F-isocrystals} and Theorem~\ref{thm:rigid-crystalline comparison for F-isocrystals}, one can obtain, for an $L$-$F$-isocrystal $\calF$, a natural $L$-linear isomorphism
\[
H^i_\dR((Z_K^\an,M_{Z_K^\an}), (\calF_K^\an,\nabla_K^\an))\cong H^i_\dR((Z_L^\an,M_{Z_L^\an}), (\calF'^\an_L,\nabla'^\an_L)).
\]
\end{rem}

The theory of weights exists for overconvergent $F$-isocrystals; after earlier works \cite{Crew-finiteness-theorem, Kedlaya-p-adic-Weil-II,Caro-Devissages}, the fundamental properties of the theory of weights in the framework of overholonomic $F$-complexes of arithmetic $\calD$-modules are proved in \cite{Abe-Caro-weights}. Here we only need the notion of overconvergent $F$-isocrystals being pure of weight $i$ for $i\in \Z$ (i.e., $\iota$-pure of weight $i$ for every isomorphism $\iota\colon \overline{\Q}_p\cong\C$), whose definition is explained, for example, in \cite[Def.~10.2]{Kedlaya-notesonisocrystals}. We say that a $\overline{\Q}_p$-$F$-isocrystal on $(Z,M_Z)/W$ is \emph{pure of weight $i$} if it is so when regarded as an overconvergent $\overline{\Q}_p$-$F$-isocrystal on $Z\smallsetminus D$.

\begin{thm}\label{thm:weight theory for F-isocrystals}
Let $\calF$ be an irreducible $\overline{\Q}_p$-$F$-isocrystal on $(Z,M_Z)/W$.
\begin{enumerate}
 \item There exists a $\overline{\Q}_p$-$F$-isocrystal $\calL$ of rank $1$ on $\Spec k/W$ such that the determinant $\det (\calF\otimes f_\cris^\ast\calL)$ is of finite order.
 \item If $\det\calF$ is of finite order, then $\calF$ is pure of weight $0$.
 \item Assume that $\calF$ is pure of weight $0$. If $D=\emptyset$, then $H^i_\cris(Z,\calF)$ is pure of weight $i$. In general, $H^i_\cris((Z,M_Z),\calF)$ is $\iota$-mixed of weight at least $i$ for every $\iota\colon \overline{\Q}_p\cong \C$.
\end{enumerate}
\end{thm}

\begin{proof}
Part (i) is \cite[Lem.~6.1(ii)]{Abe-petits-camarades}. Part (ii) is obtained independently in \cite[4.6. Cor.]{Abe-Esnault} and \cite[Thm.~4.11]{Kedlaya-companionI}; note that the crucial input is \cite{Abe-Langlands}. Part (iii) is \cite[Thm.~5.4.1]{Kedlaya-p-adic-Weil-II} and Theorem~\ref{thm:rigid-crystalline comparison for F-isocrystals}; note that a vast generalization is \cite[Thm.~4.1.3]{Abe-Caro-weights}.
\end{proof}

\subsection{Crystals and isocrystals on the nilpotent crystalline site}
Let us briefly discuss crystals and isocrystals on the nilpotent crystalline site $((Z,M_Z)/S)_\Ncris$. Since this is not the main topic of our work, the details are left to the reader.
Assume $p\neq 2$.

With the setup and notation as in \S\ref{sec:crystalline site and crystals in p nilpotent case}, one defines the notion of \emph{crystals of $\calO_{Z/S}$-modules (of finite presentation)} on $((Z,M_Z)/S)_\Ncris$ exactly in the same way. The nilpotent crystalline topoi enjoy the functoriality compatible with that of crystalline topoi, and the pullback functor induces an equivalence from the category of crystals on $((Z,M_Z)/S)_\Ncris$ to that of crystals on $((Z_1,M_{Z_1})/S)_\Ncris$; the proof of Proposition~\ref{prop:crystals only depend on mod p fiber} works verbatim.

\begin{prop}\label{prop:crystal and integral connections on nilpotent case}
Assume that there exists a log algebraic space $(Y,M_Y)$ smooth over $S$ together with an $S$-closed immersion $(Z,M_Z)\hookrightarrow(Y,M_Y)$. Let $(\calD,M_\calD)$ denote the log PD-envelope of $(Z,M_Z)$ in $(Y,M_Y)$. Then the evaluation on $(Z,\calD)$ gives rise to an equivalence from the category of crystals on $((Z,M_Z)/S)_\Ncris$ to the category of $\calO_\calD$-modules with integrable $S$-connections.
\end{prop}

\begin{proof}
When $S$ and $Y$ are schemes and $M_Y$ is trivial, this is essentially \cite[Rem.~IV.1.6.6]{Berthelot-book} (when $M_Z$ is trivial) and \cite[Thm.~1.1.8]{Ogus-F-crystals}. The general case is reduced to this case as in Proposition~\ref{prop:crystal and integral connections}.
\end{proof}

Under the assumption of this proposition, one defines the notion of \emph{$S$-flat} crystals of finite presentation on $((Z,M_Z)/S)_\Ncris$ as in Definition~\ref{def:S-flat crystals}.
These definitions and results are naturally generalized to the $p$-adic base case as in \S\ref{sec:The case of a $p$-adic base}.

Returning to the special setup of \S\ref{sec:Frobenius pullback of crystals}, let $Z_W$ be a smooth proper scheme over $W$ equipped with a normal crossings divisor $D_W$ relative to $W$, and let $M_{Z_W}$ denote the associated log structure. We define $(Z,M_Z)\coloneqq (Z_W,M_{Z_W})\times_Wk$. 

For a $p$-adic formal algebraic space $S$ over $W$, the diagram \eqref{eq:diagram for Frobenius pullback on crystalline site} defines the Frobenius pullback endofunctor $\Fr_{Z/S,\Ncris}^a$ on the category of $S$-flat crystals of finite presentation on $((Z_{S_1},M_{Z_{S_1}})/S)_\Ncris$. If $S$ is flat over $W$, we define the category of \emph{isocrystals} on $((Z_{S_1},M_{Z_{S_1}})/S)_\Ncris$ to be the isogeny category of the category of $S$-flat crystals of finite presentation on $((Z_{S_1},M_{Z_{S_1}})/S)_\Ncris$, which also admits the endofunctor $\Fr_{Z/S,\Ncris}^a$. An \emph{$F$-isocrystal} on $((Z,M_Z)/W)_\Ncris$ is an isocrystal $\calG$ on $((Z,M_Z)/W)_\Ncris$ together with an isomorphism $\varphi_\calG\colon \Fr^{a,\ast}_{Z/W,\Ncris}\calG\xrightarrow{\cong}\calG$. We also define \emph{$L$-$F$-isocrystals} (for $L/K$ finite) and \emph{$\overline{\Q}_p$-$F$-isocrystals} on $((Z,M_Z)/W)_\Ncris$ as in Definition~\ref{defn:L-F-isocrystals}.

\begin{prop}\label{prop:comparing F-isocrystals on crystalline site and nilpotent crystalline site}
The restriction functor from the category of $F$-isocrystals on $(Z,M_Z)/W$ to that of $F$-isocrystals on $((Z,M_Z)/W)_\Ncris$ is an equivalence. The same holds for $L$-$F$-isocrystals and $\overline{\Q}_p$-$F$-isocrystals.
\end{prop}

\begin{proof}
We briefly sketch the first assertion and leave the details to the reader. These two categories of $F$-isocrystals admit the internal $\Hom$, and the restriction functor is compatible with taking the internal $\Hom$. Moreover, for an $F$-isocrystal $(\calF,\varphi_\calF)$ on $(Z,M_Z)/W$ with an underlying crystal $\calF_0$, the induced map
\[
H^i(((Z,M_Z)/W)_\cris,\calF_0)[p^{-1}]\rightarrow H^i(((Z,M_Z)/W)_\Ncris,\calF_0)[p^{-1}]
\]
is an isomorphism (even without inverting $p$): both sides can be computed by the same de Rham complex. Since this isomorphism is compatible with the Frobenius action coming from $\varphi_\calF$, taking the Frobenius invariant part of $H^0$ of the internal $\Hom$ yields the full faithfulness. For the essential surjectivity, take an $F$-isocrystal $(\calG,\varphi_\calG)$ on $((Z,M_Z)/W)_\Ncris$ with an underlying crystal $\calG_0$. Twisting $\varphi_\calG$, we may assume that $\varphi_\calG(\Fr^{a,\ast}_{Z/W,\Ncris}\calG_0)\subset \calG_0$. 
Note that $\calG_0$ corresponds to a coherent $\calO_{Z_W^\wedge}$-module $(\widehat{\calG}_0,\widehat{\nabla}_0)$ with integrable $W$-connection on $Z_W^\wedge$. We need to show that $\widehat{\nabla}_0$ is topologically quasi-nilpotent, which can be checked \'etale locally on $Z_W$, so we may assume that $Z_W$ is affine and admits an \'etale map to $\mathbb{A}^n_W$ with coordinates $(z_1,\ldots,z_n)$ such that $D_W$ is defined by $z_1\cdots z_d=0$. We can then lift $\Fr_Z^a$ to a $W$-formal scheme endomorphism $F\colon Z_W^\wedge\rightarrow Z_W^\wedge$ such that $F^\ast z_i=z_i^{q}$. As in Remark~\ref{rem:compatibility of functoriality for crystals and connections}, we find that $F^\ast(\widehat{\calG}_0,\widehat{\nabla}_0)$ gives a crystal on $((Z,M_Z)/W)_\Ncris$ underlying $\Fr^{a,\ast}_{Z/W,\Ncris}\calG$. It is straightforward to see that $F^\ast\widehat{\nabla}_0$ has vanishing $p$-curvature after modulo $p$ and thus is topologically quasi-nilpotent. Since $\varphi_\calG$ is an isomorphism, one may regard $F^\ast(\widehat{\calG}_0,\widehat{\nabla}_0)$ as a subconnection of  $(\widehat{\calG}_0,\widehat{\nabla}_0)$ that contains $p^r\widehat{\calG}_0$ for some $r>0$, from which we deduce that $\widehat{\nabla}_0$ is topologically quasi-nilpotent.
\end{proof}

\section{Moduli stacks of crystals and isocrystals}\label{section:moduli of crystals and isocrystals}

In this section, we construct the moduli stack of crystals $\cM_{\cris}$ and the moduli stack of isocrystals $\calM_\isoc$; $\calM_\cris$ is a formal algebraic stack locally formally of finite type, and $\calM_\isoc$ is an adic stack.

Fix a perfect field $k$ of characteristic $p$. We set $W\coloneqq W(k)$, $W_n\coloneqq W_n(k)$ and $K \coloneqq W[p^{-1}]$. Let $Z_W$ be a smooth proper scheme over $W$ equipped with a normal crossings divisor $D_W$ relative to $W$ and let $M_{Z_W}$ denote the associated log structure. For any algebraic space $S$ over $W$, let $(Z_S,D_S)$ (resp.~$M_{Z_S}$) denote the base change of $(Z_W,D_W)$ along $S\rightarrow \Spec W$ (resp.~the pullback log structure on $Z_S$ from $M_{Z_W}$). When $S=\Spec k$, we simply write $(Z,D)$ and $M_Z$.

\subsection{Formal moduli stack of crystals}\label{sec:formal moduli of crystals}

Let
\[
\calM_{\dR,W}\coloneqq \calM_\dR((Z_W,D_W)/W)\rightarrow (\Sch/W)_\fppf
\]
denote the moduli stack of integrable connections as in Definition~\ref{def:moduli stack of integrable connections}.
We write $\calM_{\dR,k}\rightarrow (\Sch/k)_\fppf\subset (\Sch/W)_\fppf$ for the mod $p$ fiber $\calM_{\dR,W}\times_{\Spec W}\Spec k$ and $\calM_{\dR,W}^\wedge\rightarrow (\Sch/W)_\fppf$ for the $p$-adic completion of $\calM_{\dR,W}$.
Let $\calM_{\dR,k}^\nilp$ denote the substack of $\calM_{\dR,k}$ of nilpotent connections; this is a formal algebraic stack over $k$ and its reduced substack $(\calM_{\dR,k}^\nilp)_\red$ is a closed substack of $\calM_{\dR,W}$ by Theorem~\ref{thm:reduced closed substack of nilpotent connections} and Corollary~\ref{cor:nilpotent locus is formal completion}.

\begin{defn}[Formal moduli stack of crystals]\label{defn:Mcris}
Define a formal algebraic stack 
\[
\calM_\cris\coloneqq \calM_\cris(Z,M_Z)\rightarrow (\Sch/W)_\fppf
\]
to be the completion of $\calM_{\dR,W}$ along the closed subset $\lvert (\calM_{\dR,k}^\nilp)_\red\rvert \subset \lvert \calM_{\dR,W}\rvert$. By construction, $\calM_\cris$ is a locally Noetherian formal algebraic stack, locally formally of finite type over $\Spf W$; it is also obtained as the completion of $\calM_{\dR,W}^\wedge$ along $\lvert (\calM_{\dR,k}^\nilp)_\red\rvert$.
\end{defn}

The formal algebraic stack $\calM_\cris$ has the following moduli interpretation.

\begin{thm}\label{thm:moduli interpretation of Mcris}
 Let $S$ be a $p$-adic formal algebraic space over $W$, and set $S_1\coloneqq S\times_Wk$. Then $\Mor_{(\Sch/W)_\fppf}(S,\calM_\cris)$ is functorially identified with the groupoid of $S$-flat crystals of $\calO_{Z_{S_1}/S}$-modules of finite presentation on the crystalline site $((Z_{S_1},M_{Z_{S_1}})/S)_\cris$. 
\end{thm}

\begin{proof}
For each $n\geq 1$, set $S_n=S\times_WW_n$.
On one hand, we know from Remark~\ref{rem:enhanced moduli interpretation of completion of algebraic stack} that giving a $1$-morphism $\calS_S\rightarrow \calM_\cris$ is the same as giving a $1$-morphism $\calS_S\rightarrow \calM_{\dR,W}$ such that the composite $\calS_{S_\red}\rightarrow \calS_S\rightarrow \calM_{\dR,W}$ factors through $(\calM_{\dR,k}^{\nilp})_\red$; the latter is equivalent to the condition that the composite $\calS_{S_1}\rightarrow \calS_S\rightarrow \calM_{\dR,W}$ factors through $\calM_{\dR,k}^{\nilp}$ by Lemma~\ref{lem:nilpotence can be checked after restricting to reduction}. 
On the other hand, since $S$ is a $p$-adic formal algebraic space, every $1$-morphism $\calS_S\rightarrow \calM_{\dR,W}$ automatically factors through  $\calM_{\dR,W}^\wedge$. It follows from Remarks~\ref{rem:formal scheme valued point of formal algebraic stack} and \ref{rem:morphism from algebraic space to MdR} that giving a $1$-morphism $\calS_S\rightarrow \calM_{\dR,W}^\wedge$ is equivalent to giving a system $(\calE_n,\nabla_n)_n$ where $(\calE_n,\nabla_n)$ is an $S_n$-flat $\calO_{Z_{S_n}}$-module of finite presentation with integrable $S_n$-connection together with an isomorphism $(\calE_{n+1},\nabla_{n+1})|_{Z_{S_n}}\cong (\calE_n,\nabla_n)$. 
Now the assertions follow from these observations, Proposition~\ref{prop:crystals and connections in p-adic case} (and the proof therein); moreover, these identifications are obviously functorial in $S$.
\end{proof}

\begin{thm}\label{thm:Verschiebung on Mcris}
Assume that $k$ is the finite field $\F_{p^a}$. Then there exists a unique $1$-morphism
\[
V\colon \calM_\cris\rightarrow \calM_\cris
\]
over $(\Sch/W)_\fppf$ such that for every $p$-adic formal algebraic space $S$ over $W$, the induced functor $V\colon \Mor(S,\calM_\cris)\rightarrow \Mor(S,\calM_\cris)$ is identified with the Frobenius pullback $\Fr_{Z_{S_1}/S,\cris}^{a,\ast}$ of crystals.
We call it the \emph{Verschiebung map} on $\calM_\cris$.
\end{thm}

\begin{proof}
We know from Remark~\ref{rem:locally Noetherian formal algebraic stack as colimit} that there is an inductive system of substacks $\calM_n$ of $\calM_{\cris}$ such that
$\calM_{1}=(\calM_{\dR,k}^\nilp)_\red$, $\calM_n$ is an algebraic stack over $W_n$ with $\calM_{n-1}\hookrightarrow\calM_n$ a thickening, and $\varinjlim_n\calM_n\rightarrow \calM_\cris$ is an equivalence.
First fix $n$ and choose a presentation $[U/R]\xrightarrow{\cong}\calM_n$ by a smooth groupoid $(U,R,s,t,c)$ in algebraic spaces over $W_n$. By Theorem~\ref{thm:moduli interpretation of Mcris}, the structure $1$-morphism $\pi\colon \calS_U\rightarrow \calM_\cris$ defines a $U$-flat crystal $\calE$ of $\calO_{Z_{U_1}/U}$-modules of finite presentation on $((Z_{U_1},M_{Z_{U_1}})/U)_\cris$ together with an isomorphism $s_\cris^\ast \calE\cong t_\cris^\ast\calE$ on $((Z_{R_1},M_{Z_{R_1}})/R)_\cris$.
By Lemma~\ref{lem:Frobenius pullback functoriality}, the crystal $\Fr_{Z_{U_1}/U,\cris}^{a,\ast}\calE$ on $((Z_{U_1},M_{Z_{U_1}})/U)_\cris$ is equipped with an isomorphism $s_\cris^\ast \Fr_{Z_{U_1}/U,\cris}^{a,\ast}\calE\cong t_\cris^\ast\Fr_{Z_{U_1}/U,\cris}^{a,\ast}\calE$ on $((Z_{R_1},M_{Z_{R_1}})/R)_\cris$, which yields a $1$-morphism $\pi'\colon \calS_U\rightarrow \calM_\cris$ such that $\pi'\circ s=\pi'\circ t$. By Lemma~\ref{lem:2-coequalizer property of quotient stack}(ii), we obtain a $1$-morphism $V\colon \calM_n\rightarrow \calM_\cris$ such that $V\circ \pi=\pi'$. It is now straightforward to check that $V$ is independent of the choice of presentation and glues to a $1$-morphism $V\colon \calM_\cris\rightarrow \calM_\cris$ by varying $n$ and that for every $p$-adic formal space $S$ over $W$, the induced functor $V\colon \Mor(S,\calM_\cris)\rightarrow \Mor(S,\calM_\cris)$ is identified with the Frobenius pullback $\Fr_{Z_{S_1}/S,\cris}^{a,\ast}$ of crystals.
\end{proof}

The $p$-adic formal stack $\calM_{\dR,W}^\wedge$ also has a moduli interpretation in terms of crystals on the \emph{nilpotent} crystalline site.

\begin{thm}
Assume $p>2$. 
\begin{enumerate}
 \item Let $S$ be a $p$-adic formal algebraic space over $W$, and set $S_1\coloneqq S\times_Wk$. Then $\Mor_{(\Sch/W)_\fppf}(S,\calM_{\dR,W}^\wedge)$ is functorially identified with the groupoid of $S$-flat crystals of $\calO_{Z_{S_1}/S}$-modules of finite presentation on the nilpotent crystalline site $((Z_{S_1},M_{Z_{S_1}})/S)_\Ncris$. 
 \item If $k$ is the finite field $\F_{p^a}$, then there exists a unique $1$-morphism
\[
V\colon \calM_{\dR,W}^\wedge\rightarrow \calM_{\dR,W}^\wedge
\]
over $(\Sch/W)_\fppf$ such that for every $p$-adic formal algebraic space $S$ over $W$, the functor $V\colon \Mor(S,\calM_{\dR,W}^\wedge)\rightarrow \Mor(S,\calM_{\dR,W}^\wedge)$ is identified with the Frobenius pullback $\Fr_{Z_{S_1}/S,\Ncris}^{a,\ast}$ of crystals on $((Z_{S_1},M_{Z_{S_1}})/S)_\Ncris$.
We call it the \emph{Verschiebung map} on $\calM_{\dR,W}^\wedge$.
\end{enumerate}
\end{thm}

\begin{proof}
The proof is similar to those of Theorems~\ref{thm:moduli interpretation of Mcris} and \ref{thm:Verschiebung on Mcris} using Proposition~\ref{prop:crystal and integral connections on nilpotent case}. In fact, for (ii), since $\calM_{\dR,W}^\wedge$ is a $p$-adic formal stack, one may take a presentation of $\calM_{\dR,W}^\wedge$ by a smooth groupoid in $p$-adic formal algebraic spaces and apply (i).
\end{proof}

\subsection{Adic moduli stack of isocrystals}\label{sec:adic moduli stack of isocrystals}
Keep the notation in \S\ref{sec:formal moduli of crystals}.
We have $1$-morphisms
\[
\calM_\cris\rightarrow \calM_{\dR,W}^\wedge\rightarrow \calM_{\dR,W}
\]
of stacks in groupoids over $(\Sch/\calO_K)_\fppf$. These yield $1$-morphisms of adic stacks over $K$
\[
\calM_{\cris,\eta}\rightarrow (\calM_{\dR,W}^\wedge)_\eta\rightarrow \calM_{\dR,K}^\an\xrightarrow{\cong}\calM_{\andR},
\]
where the first $1$-morphism is an open immersion by Proposition~\ref{prop:generic fiber of formal completion along subset in the case of formal algebraic stacks}, the second is given in Definition~\ref{defn:comparison map from generic fiber to analytification}, and the last equivalence is Theorem~\ref{thm:analytification of MdR}.

\begin{defn}[Adic moduli stack of isocrystals]\label{defn:Misoc}
Define the open substack 
\[
\calM_{\isoc}\subset \calM_{\andR} 
\]
by applying Proposition~\ref{prop:image of comparison map from generic fiber to analytification} to the open immersion $\calM_{\cris,\eta}\hookrightarrow (\calM_{\dR,W}^\wedge)_\eta$. 
\end{defn}

\begin{prop}\label{prop:moduli interpretation of Misoc}
Let $T=\Spa(A,A^+)\in\Ob(\Afd_K)$ and assume that $A$ is complete and $(A,A^+)$ is of topologically finite type over a complete analytic affinoid field $(L,L^+)$.  Let $(\calE,\nabla)\in\Ob((\calM_{\andR})_T)$ be an integrable $T$-connection on $Z_T^\an$. Then $(\calE,\nabla)$ is in $\calM_{\isoc,T}$ if and only if the following holds:
there exist 
\begin{itemize}
 \item an \'etale covering $\{T_i=\Spa(A_i,A_i^+)\rightarrow T\}$ with $A_i$ complete,
 \item an open $L^+$-subalgebra $A_{i,0}\subset A_i^+$ of topological finite type,
 \item an $A_{i,0}$-flat coherent sheaf with topologically quasi-nilpotent integrable $A_{i,0}$-connection $(\calE_{i,0},\nabla_{i,0})$ on $Z_{A_{i,0}}$
\end{itemize}
such that $(\calE,\nabla)|_{Z_{T_i}}$ is isomorphic to the analytification of $(\calE_i,\nabla_i)|_{Z_{A_i}}$ as integrable $T_i$-connections.
\end{prop}

\begin{proof}
It follows from Propositions~\ref{prop:analytic connection coming from generic fiber of completion} and \ref{prop:crystals and connections in p-adic case}, and the proof of Theorem~\ref{thm:moduli interpretation of Mcris}.
\end{proof}

Before explaining further properties of $\calM_{\isoc}$, let us start with a preliminary discussion, following \S\ref{sec:generic fiber of p-adic completion of de Rham moduli}: take any $T=\Spa(A,A^+)\in\Ob(\Afd_K)$ with $A$ complete and let $(L,L^+)$ be a complete analytic affinoid field  over which $(A,A^+)$ is of topologically finite type. Let $A_{0}$ be an open $L^+$-subalgebra of $A^+$ of topological finite type such that $A^+$ is the integral closure of $A_0$ in $A$. With the notation therein, we have a morphism of topologically ringed spaces $\lambda \colon (Z_T^\an,\calO_{Z_T^\an}) \rightarrow (\widehat{Z}_{A_0},\calO_{\widehat{Z}_{A_0}})$.

A topologically quasi-nilpotent integrable $A_0$-connection $(\calE_0,\nabla_0)$ on the scheme $Z_{A_0}$ gives rise to a topologically quasi-nilpotent integrable $A_0$-connection $(\widehat{\calE_0},\widehat{\nabla_0})$ on $\widehat{Z}_{A_0}$ via $p$-adic completion and an integrable $T$-connection $(\calE,\nabla)$ on $Z_T^\an$ via the analytification of the generic fiber. We know from Proposition~\ref{prop:analytification and completion of coheren sheaves and connections}(ii) that they are related by $\lambda_\ast (\calE,\nabla)\cong (\widehat{\calE_0}[p^{-1}],\widehat{\nabla_0}[p^{-1}])$.

Now we can give a simpler description of $K$-valued points of $\calM_{\isoc}$.

\begin{prop}\label{prop:moduli interpretation of classical points of Misoc}
Let $T=\Spa(L,L^\circ)$ for a complete non-archimedean field $L$ over $K$ and let $(\calE,\nabla)\in\Ob((\calM_{\andR})_T)$ be an integrable $T$-connection on $Z_T^\an$. Then $(\calE,\nabla)$ is in $\left(\calM_{\isoc}\right)_{T}$ if and only if there exists a topologically quasi-nilpotent integrable connection $(\calE_0,\nabla_0)$ on the scheme $Z_{L^\circ}$ such that $(\calE,\nabla)$ is isomorphic to the analytification of $(\calE_0,\nabla_0)|_{Z_{L}}$ as integrable $T$-connections.
\end{prop}

\begin{proof}
Every \'etale covering of $T$ can be refined to a covering of the form $\{T'\coloneqq \Spa(L',(L')^\circ)\rightarrow T\}$ for a finite Galois extension $L'/L$. Hence we know from Proposition~\ref{prop:moduli interpretation of Misoc} that there exist a finite Galois extension $L'/L$ and a topologically quasi-nilpotent $(L')^\circ$-flat integrable connection $(\calE_0',\nabla_0')$ on the scheme $Z_{(L')^\circ}$ such that $(\calE,\nabla)|_{Z_{T'}}$ is isomorphic to the analytification of $(\calE_0',\nabla_0')|_{Z_{L'}}$ as integrable $T'$-connections. 

Let $\lambda'$ denote $(Z_{T'}^\an,\calO_{Z_{T'}^\an}) \rightarrow (\widehat{Z}_{(L')^\circ},\calO_{\widehat{Z}_{(L')^\circ}})$ and identify $\lambda'_\ast ((\calE,\nabla)|_{Z_{T'}})$ with $(\widehat{\calE_0'}[p^{-1}],\widehat{\nabla_0'}[p^{-1}])$. The Galois group $G\coloneqq\Gal(L'/L)$ acts on $\widehat{\calE_0'}[p^{-1}]$ and we define $\widehat{\calE_0''}$ to be the $\calO_{\widehat{Z}_{(L')^\circ}}$-submodule $\sum_{g\in G}g\widehat{\calE_0'}$ of $\widehat{\calE_0'}[p^{-1}]$. Since the connection $\widehat{\nabla_0'}[p^{-1}]$ is $G$-equivariant, it yields an integrable connection $\widehat{\nabla_0''}$ on $\widehat{\calE_0''}$ by restriction. We are going to equip the pair $(\widehat{\calE_0''},\widehat{\nabla_0''})$ with a descent datum with respect to the finite flat morphism $\widehat{Z}_{(L')^\circ}\rightarrow \widehat{Z}_{L^\circ}$. For this, note that $(\widehat{\calE_0'}[p^{-1}],\widehat{\nabla_0'}[p^{-1}])$ comes with a natural descent datum, which is described by the action of $G$ on $((\calE,\nabla)|_{Z_{T'\times_TT'}})\cong\prod_G ((\calE,\nabla)|_{Z_{T'}})$. By construction, this descent datum yields a descent datum on $(\widehat{\calE_0''},\widehat{\nabla_0''})$. By the descent theory (\cite[Thm.~2.1]{Bosch-Gortz} or \cite[Prop.~I.6.1.3]{Fujiwara-Kato}) and the GFGA existence theorem (Corollary~\ref{cor:GFGA existence for connections}), we obtain an $L^0$-flat coherent sheaf with integrable $(\Spec L^\circ)$-connection $(\calE_0,\nabla_0)$ on $Z_{L^\circ}$; use the interpretation of a connection in terms of a stratification for the descent of $\widehat{\nabla_0''}$. One can directly check that $(\widehat{\calE_0''},\widehat{\nabla_0''})$ is  topologically quasi-nilpotent (defined as in \cite[Thm.~6.2(iii)]{Kato-log}).  It follows that $(\calE_0,\nabla_0)$ is topologically quasi-nilpotent and the analytification of $(\calE_0,\nabla_0)|_{Z_L}$ is isomorphic to $(\calE,\nabla)$.
\end{proof}

\begin{cor}\label{cor:moduli interpretation of classical points of Misoc}
Let $T=\Spa(L,\calO_L)$ for a finite extension $L$ over $K$. Then the groupoid $(\calM_\isoc)_T$ is naturally identified with the groupoid $\Isoc(Z,M_Z)_L$ of pairs $(\calF,\iota_{\calF,L})$ where $\calF$ is an isocrystal on $((Z,M_Z)/W)_\cris$ and $\iota_{\calF,L}\colon L\rightarrow \End(\calF)$ is a $K$-algebra map to the endomorphism ring of $\calF$ as an isocrystal.
\end{cor}

\begin{proof}
Observe that the pushforward along $Z_T^\an\rightarrow (Z_K)^\an$ identifies the groupoid $(\calM_{\andR})_T$ with the groupoid of pairs consisting of an integrable $K$-connection $(\calE,\nabla)$ on $(Z_K)^\an$ and a $K$-algebra map $L\rightarrow \End((\calE,\nabla))$. This together with Proposition~\ref{prop:moduli interpretation of classical points of Misoc} defines a well-defined fully faithful functor $(\calM_\isoc)_T\rightarrow \Isoc(Z,M_Z)_L$.
To see that it is essentially surjective, take an object $(\calF,\iota_{\calF,L})$ of $\Isoc(Z,M_Z)_L$; it is represented by a topologically quasi-nilpotent $W$-flat integrable connection $(\widehat{\calE_0},\widehat{\nabla_0})$ on $\widehat{Z}_W$ together with a $K$-algebra map $\iota\colon L\rightarrow\End_K((\widehat{\calE_0}[p^{-1}],\widehat{\nabla_0}[p^{-1}]))$.
Take a $W$-basis $(a_i)$ of $\calO_L$ and define 
$\widehat{\calE_0'}$ to be the $\calO_{\widehat{Z}_{W}}$-submodule $\sum_{i}\iota(a_i)(\widehat{\calE_0})$ of $\widehat{\calE_0}[p^{-1}]$. Since the connection $\widehat{\nabla_0}[p^{-1}]$ is $L$-equivariant, it yields an integrable connection $\widehat{\nabla_0'}$ on $\widehat{\calE_0'}$ by restriction, and we can check that it is topologically quasi-nilpotent. Via the $\calO_L$-action on $(\widehat{\calE_0'},\widehat{\nabla_0'})$, we regard it as a topologically quasi-nilpotent $\calO_L$-flat integrable connection on $\widehat{Z}_{\calO_L}$, which defines a object of $(\calM_\andR)_T$ which maps to $(\calF,\iota_{\calF,L})$. 
\end{proof}

Next we show that the Verschiebung $V$ on $\calM_{\cris}$ defines a $1$-endomorphism on  $\calM_{\isoc}$ over $(\Afd_K)_\et$. 

\begin{thm}\label{thm:Verschiebung on Misoc}
The $1$-morphism $V_\eta\colon \calM_{\cris,\eta}\rightarrow \calM_{\cris,\eta}$ descends to a $1$-morphism
\[
V\colon \calM_{\isoc}\rightarrow \calM_{\isoc}
\]
over $(\Afd_K)_\et$. We call $V$ the \emph{Verschiebung map} on $\calM_\isoc$.
\end{thm}

\begin{proof}
Keep the notation as above. Proposition~\ref{prop:crystals and connections in p-adic case} and \S\ref{sec:Frobenius pullback of crystals} define a new topologically quasi-nilpotent integrable $A_0$-connection $\Fr^{a,\ast}_{Z/\Spf A_0,\cris}(\widehat{\calE_0},\widehat{\nabla_0})$ on $\widehat{Z}_{A_0}$, which is algebraized to an integrable $A_0$-connection $(\Fr_\cris^{a,\ast}\calE,\Fr_\cris^{a,\ast}\nabla)$ on $Z_{A_0}$ and then gives an integrable $T$-connection $(\Fr_\cris^{a,\ast}\calE_T^\an,\Fr_\cris^{a,\ast}\nabla_T^\an)$ on $Z_T^\an$. It is straightforward to see that $(\Fr_\cris^{a,\ast}\calE_T^\an,\Fr_\cris^{a,\ast}\nabla_T^\an)$ only depends on $(\widehat{\calE_0}[p^{-1}],\widehat{\nabla_0}[p^{-1}])$ or $(\calE_T^\an,\nabla_T^\an)$, and the formation is functorial in $A_0$ and $T$.
Now the theorem easily follows from Theorem~\ref{thm:Verschiebung on Mcris}, Proposition~\ref{prop:moduli interpretation of Misoc}, and the preceding discussion.
\end{proof}

We define a few variants for spaces that parametrize \emph{$F$-isocrystals}. The first is the locus in $\cM_{\isoc}$ fixed by $V$.

\begin{defn}\label{def:Fix(V)}
The adic stack $\Fix V\coloneqq \Fix(V/\calM_\isoc)\rightarrow (\Afd_K)_\et$ is defined as the $2$-fiber product 
\[
\xymatrix{
\Fix(V/\calM_{\isoc})\ar[r]\ar[d]&
{\cM_{\isoc}}\ar[d]^-{(\id,V)}\\
{\cM_{\isoc}}\ar[r]^-{\Delta}&
{\cM_{\isoc}}\times_{\Spa K}\cM_{\isoc}.
}
\]
\end{defn}

\begin{prop}\label{prop:Fix V and F-isocrystals}
Let $T=\Spa(L,\calO_L)$ for a finite extension $L$ of $K$. Then the groupoid $(\Fix V)_T$ is naturally identified with the category $\FIsoc(Z,M_Z)_L$ of $L$-$F$-isocrystals (whose morphisms are isomorphisms between them).
\end{prop}

\begin{proof}
This follows from Lemma~\ref{lem:S-linearity of Frobenius pullback}, Corollary~\ref{cor:moduli interpretation of classical points of Misoc}, and Theorem~\ref{thm:Verschiebung on Misoc}.
\end{proof}

In \S\ref{section:Misoc irr}, we show that the open substack $(\Fix 
V)_\irr\subset \Fix V$ consisting of (absolutely) irreducible objects is a $\Gm^\an$-gerbe and describe its structure (when $D=\emptyset$).

Finally, let us briefly discuss the nilpotent crystalline variant.

\begin{defn} Assume $p>2$.
Define the open substack 
\[
\calM_{\mathrm{Nisoc}}\subset \calM_{\andR} 
\]
by applying Proposition~\ref{prop:image of comparison map from generic fiber to analytification} to the open immersion $(\calM_{\dR,W}^\wedge)_\eta\hookrightarrow (\calM_{\dR,W}^\wedge)_\eta$. 
Note that $\cM_{\isoc}$ is an open substack of $\calM_{\mathrm{Nisoc}}$.
\end{defn}

\begin{prop}
Assume $p>2$.
\begin{enumerate}
 \item The $1$-morphism $V_\eta\colon (\calM_{\dR,W}^\wedge)_\eta\rightarrow (\calM_{\dR,W}^\wedge)_\eta$ descends to a $1$-morphism
\[
V\colon \calM_{\mathrm{Nisoc}}\rightarrow \calM_{\mathrm{Nisoc}}
\]
over $(\Afd_K)_\et$, and its restriction to $\cM_{\isoc}$ coincides with the Verschiebung given in Theorem~\ref{thm:Verschiebung on Misoc}.
 \item 
 Define the adic stack $\Fix(V/\calM_\mathrm{Nisoc})\rightarrow (\Afd_K)_\et$  as the $2$-fiber product 
\[
\xymatrix{
\Fix(V/\calM_{\mathrm{Nisoc}})\ar[r]\ar[d]&
\calM_{\mathrm{Nisoc}}\ar[d]^-{(\id,V)}\\
\calM_{\mathrm{Nisoc}}\ar[r]^-{\Delta}&
\calM_{\mathrm{Nisoc}}\times_{\Spa K} \calM_{\mathrm{Nisoc}}.
}
\]
Then, the induced $1$-morphism $\Fix(V/\calM_\isoc)\rightarrow\Fix(V/\calM_{\mathrm{Nisoc}})$ is an equivalence over $T=\Spa(L,\calO_L)$ for every finite extension $L$ of $K$.
\end{enumerate}
\end{prop}

\begin{proof}
Part (i) is proved as in Theorem~\ref{thm:Verschiebung on Misoc}: the current situation is in fact simpler, as the nilpotent crystalline case does not concern the topological quasi-nilpotence of an integrable connection. For (ii), one can show that an analogue of Proposition~\ref{prop:Fix V and F-isocrystals} holds for $\Fix(V/\calM_{\mathrm{Nisoc}})$ and $L$-$F$-isocrystals on $((Z,M_Z)/W)_\Ncris$. Hence the assertion follows from Proposition~\ref{prop:comparing F-isocrystals on crystalline site and nilpotent crystalline site}.
\end{proof}

\section{Geometry of \texorpdfstring{$\calM_\isoc$}{Misoc} and \texorpdfstring{$\Fix V$}{Fix V} over the irreducible locus}\label{sec:geometry of irreducible locus}

In this section, we introduce an open adic substack $\calM_{\isoc,\irr}\subset \calM_\isoc$ of (geometrically) irreducible isocrystals, which we call the \emph{irreducible locus} of $\cM_{\isoc}$, and study the geometry of $\Fix V\times_{\calM_\isoc}\calM_{\isoc,\irr}$. The advantage of the irreducible locus lies in the fact that $\calM_{\isoc,\irr}$ is a $\Gm^\an$-gerbe in the sense of Definition~\ref{def:Gm an-gerbe}. For this, we first discuss the stable locus $\calM_{\dR,\st}$ of the algebraic stack $\calM_{\dR}$ of integrable connections.
The arguments are essentially identical for the non-log case and the log case, except the fact that a stable connection may not be irreducible in the log case. To simplify the notation, we discuss the non-log case in \S\ref{section:MdRst}, \ref{section:Misoc irr} and briefly mention the log case in \S\ref{section:stable locus in the log case}.

\subsection{Gerbes and the stable moduli stack \texorpdfstring{$\cM_{\dR,\st}$}{MdRst}}\label{section:MdRst}
The goal of this subsection is to define an open substack $\cM_{\dR,\st}\subset\cM_{\dR}$, the \emph{moduli space of stable integrable connections}, and show that it is a \emph{$\Gm$-gerbe}.

To achieve this, for this section, we need to assume that $Z\rar S$ is a smooth \emph{projective}\footnote{Projective in the sense of EGA.} scheme with geometrically connected fibers and fix an $S$-ample invertible sheaf $L$ on $Z$. For simplicity, we also only work with the non-log setup.

\begin{defn}\label{defn:stable}
Let $u\colon\Spec \kappa\rar S$ be a geometric point. An integrable $\kappa$-connection $\nabla\colon\cE_{\kappa}\rar\cE_{\kappa}\otimes_{\cO_{Z_{\kappa}}}\Omega^{1}_{Z_{\kappa}/\kappa}$ on $Z_\kappa$ is called \emph{stable}\footnote{This stability is also called the \emph{Gieseker stability} or \emph{$p$-stability}.} if the following two conditions are satisfied.
\begin{enumerate}
\item The coherent sheaf $\cE_{\kappa}$ is torsion-free.
\item For any nonzero proper $\kappa$-subconnection $(\cE'_{\kappa},\nabla')\subset (\cE_{\kappa},\nabla)$, we have 
\[
p_{\cE'_{\kappa}}(N)<p_{\cE_{\kappa}}(N)\quad\text{for any}\quad N\gg0, 
\]
where $p(\cE_{\kappa}),p(\cE'_{\kappa})$ denote the normalized Hilbert polynomials with respect to $L_{\kappa}$ as in Definition~\ref{defn:Hilbert polynomial}.
\end{enumerate}
More generally, for an $S$-scheme $T\rar S$, an integrable $T$-connection $(\calE,\nabla)$ on $Z_T$ is \emph{stable} if the pullback along every geometric point of $T$ is stable.
\end{defn}

The following is well-known.

\begin{prop}\label{prop:stable=irreducible}
Let $u\colon\Spec \kappa\rar S$ be a geometric point of characteristic zero. Then every stable integrable $\kappa$-connection is irreducible as a $\kappa$-connection and its underlying coherent sheaf is locally free. Conversely, any irreducible integrable $\kappa$-connection is stable as a $\kappa$-connection.
\end{prop}

\begin{proof}
The statement about local freeness is Remark~\ref{rem:connection is locally free over char 0}. Furthermore, for an irreducible integrable connection, the condition (ii) of Definition~\ref{defn:stable} is trivially satisfied, as there is no nonzero proper subconnection. This shows the reverse direction.

For the forward direction, as $Z_{\kappa}$ is Noetherian, we may reduce to the case where $\kappa=\bC$. Then, as mentioned in \cite[p. 24]{Simpson2}, any integrable connection $(\cE_{\bC},\nabla_{\bC})$ on $Z_{\bC}$ satisfies $p_{\cE_{\bC}}(t)=p_{\cO_{Z_{\bC}}}(t)$. More precisely, the curvature of the complex analytic connection $(\cE_{\bC}^{\an},\nabla_{\bC}^{\an})$ vanishes, which implies that $c_{i}(\cE^{\an}_{\bC})_{\bC}\in H^{2i}(Z_{\bC}^{\an},\bC)$  vanishes for $i\ge1$ (e.g. \cite[p.~306, Thm.]{Milnor-Stasheff}), where $c_{i}(\cE^{\an}_{\bC})_{\bC}$ is the image of the $i$-th topological Chern class $c_{i}(\cE^{\an}_{\bC})\in H^{2i}(Z_{\bC}^{\an},\bZ)$ in the singular cohomology with complex coefficients. This implies that $c_{i}(\cE_{\bC}^{\an})$ is torsion, or its image $c_{i}(\cE_{\bC}^{\an})_{\bQ}\in H^{2i}(Z_{\bC}^{\an},\bQ)$ vanishes. Therefore, the Chern character $\mathrm{ch}(\cE^{\an}_{\bC})=\rank\cE$. By the product formula for Chern characters (e.g. \cite[Thm.~4.4.3]{Hirzebruch}), $\mathrm{ch}(\cE(n)_{\bC}^{\an})=\rank(\cE)\mathrm{ch}(\cO_{Z_{\bC}}(n)^{\an})$ for any $n\in\bZ$. By the Hirzebruch--Riemann--Roch theorem (e.g. \cite[Thm.~21.1.1]{Hirzebruch}) and GAGA, this implies that $P_{\cE_{\bC}}(n)=\rank(\cE) P_{\cO_{Z_{\bC}}}(n)$ for any $n\in\bZ$, or $p_{\cE_{\bC}}(t)=p_{\cO_{Z_{\bC}}}(t)$.

Now suppose that $(\cE_{\bC},\nabla_{\bC})$ is a stable integrable connection on $Z_{\bC}$. If there is any nonzero proper subconnection $(\cE'_{\bC},\nabla'_{\bC})$, then by the above paragraph, $p_{\cE'_{\bC}}(N)=p_{\cE_{\bC}}(N)$ for any $N\in\bZ$, violating the condition (ii) of Definition~\ref{defn:stable}. Therefore, this implies that $(\cE_{\bC},\nabla_{\bC})$ is irreducible as a connection.
\end{proof}

\begin{rem}
The proof also shows that every integrable connection is \emph{semistable}, i.e. the conditions of Definition~\ref{defn:stable} are satisfied except that, for condition (ii), you also allow its subconnections to have the same normalized Hilbert polynomial.
\end{rem}

\begin{defn}
The \emph{moduli stack of stable integrable connections}
\[
\cM_{\dR,\st}(Z/S)\rar(\Sch/S)_{\fppf}
\]
is defined to be the strictly full subcategory of $\cM_{\dR}(Z/S)$ whose objects over $T\rar S$ are the stable integrable $T$-connections  $(\cE,\nabla)$ on $Z_{T}$.

For $r\in\bN$, we may define $\cM_{\dR,\st}^{r}(Z/S)\rar(\Sch/S)_{\fppf}$ as the strictly full subcategory of objects of $\cM_{\dR,\st}(Z/S)$ whose underlying families of coherent sheaves have constant rank $r$ with respect to $L$ (see Remark~\ref{rem:rank}).
\end{defn}

\begin{prop}\label{prop:MdRst}
The natural functors
\[
\cM_{\dR,\st}^{r}(Z/S)\rightarrow \cM_{\dR,\st}(Z/S)\rightarrow\cM_{\dR}(Z/S)
\]
are open immersions. In particular, $\cM_{\dR,\st}(Z/S)$ and $\cM_{\dR,\st}^{r}(Z/S)$ are algebraic stacks with affine diagonal and is locally of finite presentation over $S$. 
Moreover, we have a decomposition $\cM_{\dR,\st}(Z/S)=\coprod_{r\in\bN}\cM_{\dR,\st}^{r}(Z/S)$.
\end{prop}

\begin{proof}
We first claim the following: given a $T$-connection $\nabla\colon\cE\rar\cE\otimes_{\cO_{Z_{T}}}\Omega_{Z_{T}/T}^{1}$ for an $S$-scheme $T$, there exists an open subscheme $U\subset T$ such that a geometric point $x\colon\Spec k\rar T$ is contained in $U$ if and only if $(\cE_{x},\nabla_{x})$ is stable. This statement is shown in \cite[Lem.~3.7]{Simpson1} and deduced from the properness of the Quot schemes and the finiteness of the set of possible Hilbert polynomials of the quotient. The claim implies that $\cM_{\dR,\st}(Z/S)\rightarrow\cM_{\dR}(Z/S)$ is an open immersion. Since the  Hilbert polynomial (hence the multiplicity) of a flat family of pure sheaves is locally constant (e.g., \cite[Prop.~2.1.2]{HL}), we also see that $\cM_{\dR,\st}^{r}\rar\cM_{\dR,\st}$ is an open substack and the desired decomposition holds for $\cM_{\dR,\st}$.
\end{proof}

In what follows, $\Gm$ denotes the multiplicative group over $S$ and set $B\Gm\coloneqq [S/\Gm]\rightarrow S$.

\begin{defn}\label{def:gerbe}
An algebraic stack $\calX\rightarrow (\Sch/S)_\fppf$ is said to be a \emph{$\Gm$-gerbe} if there exist an $S$-algebraic space $X$ and a $1$-morphism $f\colon \cX\rar \calS_X$ over $S$ satisfying the following conditions:
\begin{enumerate}
 \item $f$ is an epimorphism;
 \item for every $T\in \Sch/S$ and $x,x'\in \Ob(\calX_T)$ with $f(x)=f(x')$ in $X(T)$, there exists an fppf covering $\{T_i\rightarrow T\}$ such that the pullbacks $x|_{T_i},x'|_{T_i}\in \Ob(\calX_{T_i})$ are isomorphic for every $i$;
 \item for every $T\in \Sch/S$ and $x\in \Ob(\calX_T)$, there is an isomorphism of sheaves $\iota_x\colon (\Gm)_T\rightarrow \Isom_{\calX_T}(x,x)$;
  \item for every $T\in \Sch/S$ and an isomorphism $\alpha\colon x\xrightarrow{\cong}x'$ in $\Ob(\calX_T)$, we have $\iota_{x'}=\mathrm{Inn}_\alpha\circ \iota_x$ where $\mathrm{Inn}_\alpha\colon \Isom_{\calX_T}(x,x)\rightarrow\Isom_{\calX_T}(x',x')$ is the isomorphism sending $\beta$ to $\alpha\circ\beta\circ \alpha^{-1}$.
\end{enumerate}
We call $X$ the \emph{coarse moduli space} of the $\Gm$-gerbe $\cX$ (see Remark~\ref{rem:gerbe} below). By abuse of notation, we simply write $\calX\rightarrow X$ for $f$ and say that $\calX\rightarrow X$ realizes $\calX$ as a $\Gm$-gerbe.
\end{defn}

\begin{rem}\label{rem:gerbe}
There is a redundancy in Definition~\ref{def:gerbe}. Namely, condition (iii) implies (iv) (after possibly modifying $\iota_x$) by applying \cite[0CJY]{stacks-project} to $\calC\coloneqq (\Sch/S)_\fppf/X$; (iii) also implies (i) and (ii) by \cite[06QJ]{stacks-project}.
In general, if a $1$-morphism $f\colon \calX\rightarrow \calS_X$ to an algebraic space $X$ satisfies (i) and (ii), then $X$ is canonically identified with the fppf (or \'etale) sheafification of the presheaf
\[
T\mapsto \Ob(\cX_{T})/\text{isomorphisms}
\]
by \cite[06QD]{stacks-project}. 
\end{rem}

\begin{prop}\label{prop:basic properties of Gm-gerbes}
Let $\calX$ be a $\Gm$-gerbe over $S$ and let $f\colon \calX\rightarrow X$ be a $1$-morphism to an algebraic space realizing $\calX$ as a $\Gm$-gerbe. Then the following properties hold.
\begin{enumerate}
 \item If $f$ admits a section $x\colon X\rightarrow\calX$, $x$ induces an equivalence $\iota \colon X\times_SB\Gm\xrightarrow{\cong} \calX$ such that $f\circ \iota=\pr_X$.
 \item The $1$-morphism $f$ is smooth.
 \item There exists a surjective \'etale morphism $U\rightarrow X$ from a scheme $U$ such that $U\times_X\calX\cong U\times_SB\Gm$.
\end{enumerate}
Conversely, if $f\colon \calX\rightarrow \calS_X$ is a $1$-morphism over $S$ from an algebraic stack $\calX$ to an algebraic space $X$ satisfying part (iii), then $f$ realizes $\calX$ as a $\Gm$-gerbe. 
\end{prop}

\begin{proof}
(i) Observe that there is a natural equivalence $X\times_S\Gm\cong X\times_{\Delta\circ x,\calX\times_S\calX,\Delta\circ x} X$; if $X$ is a scheme, this follows from condition (iii) of Definition~\ref{def:gerbe}, and the general case is deduced from a scheme case by taking a presentation. Now the assertion is proved as in \cite[06QG(1)]{stacks-project}; one first shows that there is a natural identification between the groupoid in algebraic spaces induced by $(X,R\coloneqq X\times_{x,\calX,x}X)$ and the one induced by the trivial action of $\Gm$ on $X$. This identification induces a fully faithful $1$-morphism $[X/\Gm]=X\times_SB\Gm\rightarrow \calX$, which is essentially surjective by condition (ii) of Definition~\ref{def:gerbe}. 

(ii) Since $f$ is an epimorphism, there exists an fppf cover $U\rightarrow X$ from a scheme $U$ such that $U\rightarrow X$ factors through $f$. It is straightforward to see that $U\times_X\calX\rightarrow U$ realizes $f_U\colon U\times_X\calX$ as a $\Gm$-gerbe. By (i), $f_U$ is identified with $U\times B\Gm\rightarrow U$. Since $B\Gm\rightarrow S$ is smooth, we conclude that $f$ is smooth.

(iii) By (ii), there exists a surjective \'etale morphism $U\rightarrow X$ from a scheme $U$ such that $U\rightarrow X$ factors through $f$. Now (iii) follows from (i).

Finally, it is straightforward to see the last assertion by directly verifying the conditions of Definition~\ref{def:gerbe}.
\end{proof}

\begin{prop}\label{prop:MdRst is Gmgerbe}
The stable moduli stack $\cM_{\dR,\st}$ is a $\Gm$-gerbe.
\end{prop}

We will denote the coarse moduli space of $\cM_{\dR,\st}$ as $M_{\dR,\st}$.

\begin{proof}
Let $T$ be an $S$-scheme and let $x\in\Ob((\cM_{\dR,\st})_{T})$. Write $\nabla\colon\cE\rar\cE\otimes_{\cO_{Z_{T}}}\Omega_{Z_{T}/T}^{1}$ for a stable integrable $T$-connection corresponding to $x$.  As multiplying on $\cE$ by an element of $\cO_{T}^{\times}$ gives an isomorphism of $(\calE,\nabla)$, we have a natural morphism $\iota_{x}\colon\Gm\rar\Isom_{(\cM_{\dR,\st})_{T}}(x,x)$. It is enough to show that $\iota_x$ is an isomorphism; condition (iv) in Definition~\ref{def:gerbe} then follows easily, and by Remark~\ref{rem:gerbe}, we conclude that $\calM_{\dR,\st}$ is a $\Gm$-gerbe. 

By Proposition~\ref{prop:MdRst}, $\Isom_{(\cM_{\dR,\st})_{T}}(x,x)$ is an algebraic space locally of finite presentation over $T$. Therefore, by \cite[Cor.~17.9.5]{EGAIV4} and faithfully flat descent, we may assume that $T=\Spec k$ for an algebraically closed field $k$. 

 Suppose that $\varphi\colon\cE\rar\cE$ is a homomorphism that respects the connection. If $\varphi$ is not zero, then $\Ker\varphi\subset\cE$ is a proper subconnection. If $\Ker\varphi\neq 0$, we have $p_{\Ker\varphi}(N)< p_{\cE}(N)$ for $N\gg0$, or $p_{\Image\varphi}(N)=p_{\cE/\Ker\varphi}(N)>p_{\cE}(N)$ for $N\gg0$, which contradicts the stability of $(\cE,\nabla)$. Therefore, $\varphi$ is either $0$ or injective. As $\cE$ is finitely generated and $k$ is algebraically closed, there exists $\lambda\in k$ such that $\ker(\varphi-\lambda\id_{\cE})\ne0$. This  implies that $\varphi=\lambda\id_{\cE}$. Therefore, $\iota_{x}$ is bijective on the $k$-points. Note that, as $\Isom_{(\cM_{\dR,\st})_{k}}(x,x)$ is a $k$-group algebraic space locally of finite presentation, by \cite[Lem.~4.2]{ArtinAlgebraizationI}, $\Isom_{(\cM_{\dR,\st})_{k}}(x,x)$ is a $k$-group scheme. Therefore,  $\iota_{x}$ is universally injective, as the diagonal $\Delta_{\iota_{x}}$ is surjective (on closed points).

For an Artin local $k$-algebra $A$ with $\dim_{k}A=n$, consider the integrable $A$-connection $(\cE_{A},\nabla_{A})$ on $Z_{A}$. Then, an $\cO_{Z_{A}}$-module homomorphism $\varphi_{A}\colon \cE_{A}\rar\cE_{A}$ respecting the connection $\nabla_{A}$ is equivalent to an $\cO_{Z}$-module homomorphism $\varphi\colon\cE^{\oplus n}\rar\cE^{\oplus n}$, respecting the connection $(\cE^{\oplus n},\nabla^{\oplus n})$, and also respecting the $A$-module structure on $\cE^{\oplus n}$ after identifying it with $\cE_{A}$, using a certain fixed $k$-basis of $A$. Thus, $\varphi$ can be expressed as an $n\times n$ matrix with entries in $\End_{\cO_{Z}}(\cE)$. As each matrix entry is an endomorphism that also respects the connection $\nabla$, by the above paragraph, it follows that each entry is of the form $\lambda\id_{\cE}$ for some $\lambda\in k$. As this is an $A$-module homomorphism, this matrix must be the matrix corresponding to multiplication by a scalar $a\in A$, and for this to be an isomorphism, $a$ must be invertible in $A$. Thus, it follows that $\iota_{x}(A)\colon A^{\times}\rar\Isom_{(\cM_{\dR,\st})_{k}}(x,x)(A)$ is an isomorphism. This shows that $\iota_{x}$ induces isomorphisms between the formal completions at the identity of $\bG_{m,k}$ and $\Isom_{(\cM_{\dR,\st})_{k}}(x,x)$. Therefore, $\iota_{x}$ is (formally) \'etale at the identity, so it is an \'etale morphism. As $\iota_{x}$ is already shown to be universally injective, it is an open immersion. As $\iota_{x}$ is bijective on closed points, it follows that $\iota_{x}$ is an isomorphism.
\end{proof}

\subsection{Irreducible locus of \texorpdfstring{$\calM_\isoc$}{Misoc}}\label{section:Misoc irr}

We keep the notation in \S\ref{section:moduli of crystals and isocrystals} and further assume that $Z_W\rightarrow \Spec W$ is projective and $D_W=\emptyset$ as in the previous subsection.
Let $\calM_\andR\rightarrow (\Afd_K)_\Et$ be the adic stack of analytic integrable connections on $Z_K^\an$ in Definition~\ref{def:moduli stack of analytic integrable connections}. Recall from \S\ref{sec:analytification of de Rham moduli} that there is an equivalence $\alpha\colon \calM_\dR^\an\xrightarrow{\cong}\calM_{\andR}$
from the analytification of $\calM_\dR\coloneqq \calM_\dR(Z_K/K)$ in such a way that for every $T=\Spa(A,A^+)\in\Ob(\Afd_K)$ with $A$ complete, $\alpha_T$ is induced from the analytification $(\calM_\dR)_{\Spec A}\xrightarrow{\cong}(\calM_{\andR})_T$.

\begin{defn}\label{defn:Misocirr}
Define a category fibered in groupoids $\calM_{\andR,\irr}\rightarrow \Afd_K$ as the strictly full subcategory $\calM_{\andR,\irr}\subset \calM_{\andR}$ whose objects over $T\rightarrow S$ are integrable $T$-connections $(\calE,\nabla)$ on $Z_T^\an$ such that for every algebraically closed affinoid field $\overline{t}=\Spa(L,L^+)\in\Ob(\Afd_K)$ and every map $\overline{t}\rightarrow T$, the pullback $(\calE,\nabla)|_{Z_{\overline{t}}^{\an}}$ is irreducible as a $\overline{t}$-connection.
It is straightforward to see that $\calM_{\andR,\irr}$ is a substack of $\calM_{\andR}$. Moreover, it follows from Proposition~\ref{prop:stable=irreducible} that $\alpha\colon \calM_\dR^\an\xrightarrow{\cong}\calM_{\andR}$ restricts to an equivalence $\calM_{\dR,\st}^\an\xrightarrow{\cong}\calM_{\andR,\irr}$. In particular, $\calM_{\andR,\irr}$ is an open substack of $\calM_\andR$.

More generally, for any open substack $\calV$ of $\calM_\andR$, we define an open substack $\calV_\irr$ of $\calV$ in a similar way.
\end{defn}

For simplicity, write $\Spa K$ for $\Spa(K,\calO_K)$ and $\times$ for $\times_{\Spa K}$. Let $\Gm^\an$ denote the analytification of $(\Gm)_{\Spec K}$ and set $B\Gm^\an\coloneqq [\Spa K/\Gm^\an]$. We are going to show that $\calM_{\andR,\irr}$ is a $\Gm^\an$-gerbe. For this, we first define $\Gm^\an$-gerbes as in Definition~\ref{def:gerbe} and discuss the basic properties.

\begin{defn}\label{def:Gm an-gerbe}
An adic stack $\calX\rightarrow (\Afd_K)_\Et$ is said to be a \emph{$\Gm^\an$-gerbe} if there exist an \'etale algebraic space $X$ over $K$ and a $1$-morphism $f\colon \cX\rar \calS_X$ over $K$ satisfying the following conditions:
\begin{enumerate}
 \item $f$ is an epimorphism;
 \item for every $T\in \Ob(\Afd_K)$ and $x,x'\in \Ob(\calX_T)$ with $f(x)=f(x')$ in $X(T)$, there exists a covering $\{T_i\rightarrow T\}$ in $(\Afd_K)_\Et$ such that the pullbacks $x|_{T_i},x'|_{T_i}\in \Ob(\calX_{T_i})$ are isomorphic for every $i$;
 \item for every $T\in \Ob(\Afd_K)$ and $x\in \Ob(\calX_T)$, there is an isomorphism of sheaves $\iota_x\colon (\Gm^\an)_T\rightarrow \Isom_{\calX_T}(x,x)$;
  \item for every $T\in \Ob(\Afd_K)$ and an isomorphism $\alpha\colon x\xrightarrow{\cong}x'$ in $\Ob(\calX_T)$, we have $\iota_{x'}=\mathrm{Inn}_\alpha\circ \iota_x$ where $\mathrm{Inn}_\alpha\colon \Isom_{\calX_T}(x,x)\rightarrow\Isom_{\calX_T}(x',x')$ is the isomorphism sending $\beta$ to $\alpha\circ\beta\circ \alpha^{-1}$.
\end{enumerate}
We call $X$ the \emph{coarse moduli space} of the $\Gm^\an$-gerbe $\cX$. By abuse of notation, we simply write $\calX\rightarrow X$ for $f$ and say that $\calX\rightarrow X$ realizes $\calX$ as a $\Gm^\an$-gerbe.
We say that $\calX$ is \emph{split} if $\calX$ is equivalent to the $\Gm^\an$-gerbe of the form $X\times B\Gm^\an$ for an \'etale space $X$.
\end{defn}

\begin{prop}\label{prop:basic properties of Gm an-gerbes}
Let $\calX\rightarrow (\Afd_K)_\Et$ be an adic stack and let $f\colon \calX\rightarrow \calS_X$ be a $1$-morphism to an \'etale space $X$ over $K$.
\begin{enumerate}
 \item If $f$ realizes $\calX$ as a $\Gm^\an$-gerbe, then $X$ is canonically identified with the \'etale sheafification of the presheaf
\[
T\mapsto \Ob(\cX_{T})/\text{isomorphisms}.
\]
 \item The $1$-morphism $f$ realizes $\calX$ as a $\Gm^\an$-gerbe if and only if there exists a surjective \'etale morphism $U\rightarrow X$ from an \'etale space $U$ over $K$ such that $U\times_X\calX\cong U\times B\Gm^\an$ is a split $\Gm^\an$-gerbe. 
\end{enumerate}
\end{prop}

Note that one can even take $U$ in (ii) to be an adic space of finite type over $K$.

\begin{proof}
Part (i) is proved as in \cite[06QD]{stacks-project}, and one can easily modify the proof of Proposition~\ref{prop:basic properties of Gm-gerbes} to show (ii). 
\end{proof}

We need the following variant of (ii) in a later discussion.

\begin{prop}\label{prop:bootstrap for Gm an-gerbe}
Let $\calX\rightarrow (\Afd_K)_\Et$ be an adic stack satisfying conditions (iii) and (iv) of Definition~\ref{def:Gm an-gerbe} and let $f\colon \calX\rightarrow \calS_Y$ be a $1$-morphism to an \'etale space $Y$. Assume that there exists an \'etale surjection $Y'\rightarrow Y$ of \'etale spaces such that 
\begin{enumerate}
  \item $\calX'\coloneqq \calX\times_{\calS_Y}\calS_{Y'}$ is a $\Gm^\an$-gerbe, and
  \item the $1$-morphism $\pi\colon \calX'\rightarrow\calX$ is compatible with $\Gm^\an$-structures on isomorphisms (i.e., for every $T\in \Ob(\Afd_K)$ and $x'\in \Ob(\calX'_T)$, the map $\Isom_{\calX'_T}(x',x')\rightarrow \Isom_{\calX_T}(\pi(x'),\pi(x'))$ is compatible with $\iota_{x'}$ and $\iota_{\pi(x')}$).
\end{enumerate}
 Then $\calX$ is a $\Gm^\an$-gerbe.
\end{prop}

\begin{proof}
Let $X$ denote the sheaf on $(\Afd_K)_\Et$ obtained as the sheafification of the presheaf $T\mapsto \Ob(\cX_{T})/\text{isomorphisms}$ and define $X'$ similarly. We have the following commutative diagram
\[
\xymatrix{
\calX' \ar[r]\ar[d]_-\pi& \calS_{X'}\ar[r]\ar[rd] & \calS_{X\times_YY'}\ar[r]\ar[d] & \calS_{Y'}\ar[d]\\
\calX\ar[rr]&& \calS_X \ar[r]& \calS_Y.
}
\]
Observe that the induced $1$-morphism $\calX'\rightarrow \calX\times_{\calS_X}\calS_{X'}$ is an equivalence; one can verify the conditions in Lemma~\ref{lem:criteria for equivalence of stacks} thanks to the assumptions on $\calX$, $\calX'$, and $\pi$. In particular, $\calX\times_{\calS_X}\calS_{X'}\rightarrow \calX\times_{\calS_X}\calS_{X\times_YY'}$ is an equivalence. It follows that $X'\rightarrow X\times_YY'$ is an isomorphism. Hence $X'\rightarrow X$ is representable by \'etale spaces and \'etale and surjective. It follows from Proposition~\ref{prop:bootstrap for etale spaces (sheaf admitting etale surjective map from etale sapce)} that $X$ is an \'etale space.  Finally, we conclude that $\calX$ is a $\Gm^\an$-gerbe by Proposition~\ref{prop:basic properties of Gm an-gerbes}.
\end{proof}

\begin{prop}\label{prop:open sub of gerbe is gerbe}
Let $\calX\rightarrow (\Afd_K)_\Et$ be a $\Gm^\an$-gerbe. Then every open substack $\calV$ of $\calX$ is also a $\Gm^\an$-gerbe.
\end{prop}

\begin{proof}
Choose a $1$-morphism $f\colon \calX\rightarrow X$ to an \'etale space $X$ that realizes $\calX$ as a $\Gm^\an$-gerbe. By Proposition~\ref{prop:basic properties of Gm an-gerbes}(ii), there exists a surjective \'etale morphism $U\rightarrow X$ from an adic space $U$ of finite type over $K$ such that $U\times_X\calX\cong U\times_SB\Gm^\an$.

The pullback of the open immersion $U\times_X\calV\rightarrow U\times_X\calX$ along $U\rightarrow U\times_SB\Gm^\an\cong U\times_X\calX$ defines an open subspace $W\subset U$. By constructions, two open immersions 
$U\times_X\calV\rightarrow U\times_X\calX\cong U\times_SB\Gm^\an$ and $W\times_SB\Gm^\an\rightarrow U\times_SB\Gm^\an$ define the same open substack of $U\times_SB\Gm^\an$. Hence there is an equivalence $U\times_X\calV\xrightarrow{\cong}W\times_SB\Gm^\an$, which also implies that the open immersion $W\times_X\calV\rightarrow U\times_X\calV$ is an equivalence.

Define $V$ to be the image of $W\subset U\rightarrow X$. It is straightforward to see that $V$ is an open subspace of $X$ and the $1$-morphism $\calV\rightarrow \calX\xrightarrow{f}X$ factors as $\calV\xrightarrow{f_\calV}V\rightarrow X$ with $f_\calV$ being an epimorphism. Note that the open immersion $W\hookrightarrow U\times_XV$ is an equality; for this, first observe that the induced $1$-morphism $W\times_V\calV\rightarrow (U\times_XV)\times_V\calV$ is identified with the open immersion $W\times_X\calV\rightarrow U\times_X\calV$. Since the latter is an equivalence and $f_\calV$ is an epimorphism, we conclude $W= U\times_XV$. In particular, the surjective \'etale morphism $W\rightarrow V$ induces an equivalence $W\times_V\calV\cong W\times_SB\Gm^\an$. Hence $\calV$ is a $\Gm^\an$-gerbe by Proposition~\ref{prop:basic properties of Gm an-gerbes}(ii).
\end{proof}

\begin{thm}\label{thm:Misocirr is Gm-gerbe}
The adic stack $\calM_{\andR,\irr}\rightarrow(\Afd_K)_\Et$ is a $\Gm^\an$-gerbe and its coarse moduli space is identified with the analytification $M_{\dR,\st}^\an$ of the coarse moduli space $M_{\dR,\st}$ of $\calM_{\dR,\st}$. In particular, $\calM_{\isoc,\irr}\rightarrow(\Afd_K)_\Et$ is a $\Gm^\an$-gerbe.
\end{thm}

\begin{proof}
It is straightforward to see that the analytification of a $\Gm$-gerbe over $(\Sch/K)_\fppf$ is a $\Gm^\an$-gerbe, and the coarse moduli of the $\Gm^\an$-gerbe is given by the analytification of the coarse moduli of the $\Gm$-gerbe. Hence the first assertion follows from Proposition~\ref{prop:MdRst is Gmgerbe}. 
The second assertion follows from Proposition~\ref{prop:open sub of gerbe is gerbe}.
\end{proof}

Recall from Definition~\ref{def:Fix(V)} that $\Fix V$ is defined as the $2$-fiber product 
\[
\xymatrix{
\Fix V\ar[r]\ar[d]_q&
{\cM_{\isoc}}\ar[d]^-{(\id,V)}\\
{\cM_{\isoc}}\ar[r]^-{\Delta}&
{\cM_{\isoc}}\times_{\Spa K}\cM_{\isoc}.
}
\]

\begin{defn}
Define $(\Fix V)_{\irr}$ to be the open substack of $\Fix V$ that is equivalent to the open immersion $\Fix V\times_{q,\calM_\isoc}\calM_{\isoc,\irr}\rightarrow \Fix V$.
\end{defn}

\begin{prop}\label{prop:Fix(V)st is Gm-gerbe}
The adic stack $(\Fix V)_\irr\rightarrow (\Afd_K)_\Et$ is a $\Gm^\an$-gerbe and the projection $(\Fix V)_\irr\rightarrow \calM_{\isoc,\irr}$  is a $1$-morphism of $\Gm^\an$-gerbes.\footnote{Namely, it is a $1$-morphism over $\Afd_K$ that is compatible with the $\Gm^\an$-structures on isomorphisms in condition (iii) of Definition~\ref{def:Gm an-gerbe}.}
\end{prop}

\begin{proof}
Let $V^{-1}(\calM_{\isoc,\irr})$ denote the open substack of $\calM_{\isoc}$ that is equivalent to the open immersion $\calM_{\isoc,\irr}\times_{\calM_{\isoc},V}\calM_{\isoc}\rightarrow \calM_{\isoc}$. It is immediate to see that $V^{-1}(\calM_{\isoc,\irr})$ is an open substack of $\calM_{\isoc,\irr}$, the $1$-morphism $(\Fix V)_\irr\rightarrow \calM_{\isoc,\irr}$ factors through $V^{-1}(\calM_{\isoc,\irr})$, and $(\Fix V)_\irr$ is equivalent to the $2$-fiber product $\calY$ sitting in the diagram
\[
\xymatrix{
\calY\ar[r]\ar[d]
& V^{-1}(\cM_{\isoc,\irr})\ar[d]^-{(\id,V)}
\\
V^{-1}(\cM_{\isoc,\irr})\ar[r]^-{(\id,\mathrm{incl})} 
& V^{-1}(\cM_{\isoc,\irr})\times_{\Spa K}\cM_{\isoc,\irr}.
}
\]
Since $V^{-1}(\cM_{\isoc,\irr})$ is a $\Gm^\an$-gerbe by Proposition~\ref{prop:open sub of gerbe is gerbe} and $V\colon V^{-1}(\cM_{\isoc,\irr})\rightarrow \cM_{\isoc,\irr}$ is a $1$-morphism of $\Gm^\an$-gerbes, the assertion follows from Lemma~\ref{lem:intersection of gerbes is gerbe} below.
\end{proof}

\begin{lem}\label{lem:intersection of gerbes is gerbe}
Let $\calX,\calY$ be $\Gm^\an$-gerbes and let $f_1,f_2\colon \calX\rightarrow\calY$ be $1$-morphisms of $\Gm^\an$-gerbes. Write $\Gamma_{f_1},\Gamma_{f_2}\colon \calX\rightarrow \calX\times_{\Spa K}\calY$ be the corresponding graph $1$-morphisms and define the adic stack $\mathrm{Eq}(f_1,f_2)$ by the following 2-fiber product
\[
\xymatrix{
\mathrm{Eq}(f_1,f_2)\ar[r]\ar[d] & \calX \ar[d]_-{\Gamma_{f_2}}\\
\calX\ar[r]^-{\Gamma_{f_1}}& \calX\times_{\Spa K}\calY.
}
\]
Then $\mathrm{Eq}(f_1,f_2)$ is also a $\Gm^\an$-gerbe.
\end{lem}

\begin{proof}
We will use the notation $\mathrm{Eq}(-,-)$ to denote the $2$-fiber product of the graph $1$-morphisms.
Let $\calX\rightarrow X$ and $\calY\rightarrow Y$ be $1$-morphisms to algebraic spaces realizing $\calX$ and $\calY$ as $\Gm^\an$-gerbes. Since $f_i$ is a $1$-morphism of $\Gm^\an$-gerbes, the composite $\calX\xrightarrow{f_i}\calY\rightarrow Y$ factors as $\calX\rightarrow X\xrightarrow{f'_i}Y$ for some morphism $f'_i$ of \'etale spaces. We first claim that the induced $1$-morphism $\calX\rightarrow X\times_{f'_i,Y}\calY$ is an equivalence; one can easily verify the conditions in Lemma~\ref{lem:criteria for equivalence of stacks} using the fact that $\calX$ and $\calY$ are $\Gm^\an$-gerbes and  $f_i$ is a $1$-morphism of $\Gm^\an$-gerbes.

We claim that we may assume that $\calX$ is split. In fact, take an \'etale surjection $\alpha'\colon U\rightarrow X$ such that $U\times_X\calX$ is split as in Proposition~\ref{prop:basic properties of Gm an-gerbes}(ii). Write $\alpha$ for the $1$-morphism $U\times_X\calX\rightarrow \calX$. It is straightforward to see the following diagram
\[
\xymatrix{
\mathrm{Eq}(f_1\circ \alpha,f_2\circ \alpha)\ar[r]\ar[d]& U\ar[d]\\
\mathrm{Eq}(f_1,f_2)\ar[r]  & X
}
\]
is Cartesian. It is easy to see that $\mathrm{Eq}(f_1,f_2)$ satisfies conditions (iii) and (iv) of Definition~\ref{def:Gm an-gerbe} since $f$ and $g$ are $1$-morphisms of $\Gm^\an$-gerbes; moreover, the above left vertical $1$-morphism is compatible with the $\Gm^\an$-sturctures on isomorphisms. Hence the claim follows from Proposition~\ref{prop:bootstrap for Gm an-gerbe}.

Now take an \'etale surjection $\beta'\colon V\rightarrow Y$ such that $V\times_Y\calY$ is split. We are going to show that $\mathrm{Eq}(f_1,f_2)\times_YV$ is a split $\Gm^\an$-gerbe, which will complete the proof by Proposition~\ref{prop:bootstrap for Gm an-gerbe}.

Define an \'etale space $X_i$ and an adic stack $\calX_i$ and the map $g_i$ as the ($2$-)fiber products and the map appearing in the following commutative diagram 
\[
\xymatrix{
\calX_i\ar[r]\ar[d] & X_i\ar@/^1.0pc/[rr]^-{g_i}\ar[r]_-{\Gamma_{g_i}}\ar[d] & X\times V\ar[r]\ar[d]& V\ar[d]\\
\calX \ar[r] & X \ar[r]^-{\Gamma_{f_i'}} \ar@/_1.0pc/[rr]_{f_i'}& X\times Y \ar[r]& Y
}
\]
with every square being Cartesian. Let $Z$ denote the \'etale space $X_1\times_{\Gamma_{g_1},X\times V,\Gamma_{g_2}}X_2$. By construction, $Z$ is isomorphic to $\mathrm{Eq}(f_1',f_2')\times_YV$. 

By the splitting assumption, we identify $V\times_Y\calY\cong V\times\Gm^\an$, $\calX\cong X\times B\Gm^\an$, and $\calX_i\cong X_i\times\Gm^\an$. Via these equivalences, the base change of $f_i\colon \calX\rightarrow \calY$ along $V\rightarrow Y$ is identified with $g_i\times \id\colon X_i\times B\Gm^\an \rar V\times B\Gm^{\an}$. Hence $\mathrm{Eq}(f_1,f_2)\times_YV$ sits in the following Cartesian diagram
\[
\xymatrix{
\mathrm{Eq}(f_1,f_2)\times_YV \ar[r]\ar[d] & X_2\times B\Gm^\an\ar[d]_-{\Gamma_{g_2\times\id}}\\
X_1\times B\Gm^\an \ar[r]^-{\Gamma_{g_1\times\id}}& (X\times B\Gm^\an)\times (V\times B\Gm^\an).
}
\]
It is easy to see that the above $2$-fiber product is equivalent to the product of $Z=X_1\times_{\Gamma_{g_1},X\times V,\Gamma_{g_2}}X_2$ and $B\Gm^\an\times_{\Delta,B\Gm^\an\times B\Gm^\an,\Delta}B\Gm^\an$. The latter adic stack is equivalent to the inertia stack of $B\Gm^\an$ by \cite[034H]{stacks-project}, which is $\Gm^\an\times B\Gm^\an$ as $\Gm^\an$ is commutative. We conclude that $\mathrm{Eq}(f_1,f_2)\times_YV$ is equivalent to $Z\times\Gm^\an\times B\Gm^\an$, which is obviously a split $\Gm^\an$-gerbe.
\end{proof}

\begin{defn}\label{def:MFisoc}
Let $M_{\isoc,\irr}$ denote the coarse moduli space of the $\Gm^\an$-gerbe $\calM_{\isoc,\irr}$ and define an open \'etale subspace $V^{-1}(M_{\isoc,\irr})\subset M_{\isoc,\irr}$ similarly.
Define the \'etale space $M_{\isoc,\irr}^{V=\id}$ by the following fiber product
\[
\xymatrix{
M_{\isoc,\irr}^{V=\id}\ar[r]\ar[d]
& V^{-1}(M_{\isoc,\irr})\ar[d]^-{(\id,V)}
\\
V^{-1}(M_{\isoc,\irr})\ar[r]^-{(\id,\mathrm{incl})} 
& V^{-1}(M_{\isoc,\irr})\times_{\Spa K}M_{\isoc,\irr}.
}
\]Corollary~\ref{cor:moduli interpretation of classical points of Misoc} (together with the proof of Proposition~\ref{prop:moduli interpretation of classical points of Misoc}) implies that for every finite extension $L$ over $K$, the set $M_{\isoc,\irr}^{V=\id}(\Spa L)$ is naturally identified with the set of isomorphism classes of \emph{absolutely irreducible $L$-isocrystals on $Z$ underlying $L$-$F$-isocrystals}; more precisely, $M_{\isoc,\irr}^{V=\id}(\Spa L)$ is the set of isomorphism classes of pairs $(\calF,\iota_{\calF,L})$ where $\calF$ is an isocrystal on $(Z/W)_\cris$ and $\iota_{\calF,L}\colon L\rightarrow \End(\calF)$ is a $K$-algebra map such that
\begin{enumerate}
  \item $(\calF,\iota_{\calF,L})$ underlies an $L$-$F$-isocrystal; and
  \item $(\calF,\iota_{\calF,L})$ is absolutely irreducible (i.e., the scalar extension $(\calF\otimes_{L}L', \iota_{\calF,L}\otimes_LL')$ is irreducible in an obvious sense for every finite extension $L'$ of $L$).
\end{enumerate}
We have a natural $1$-morphism $(\Fix V)_\irr\rightarrow M_{\isoc,\irr}^{V=\id}$, and the proof of Lemma~\ref{lem:intersection of gerbes is gerbe} shows that \'etale locally on $(\Fix V)_\irr$ and $M_{\isoc,\irr}^{V=\id}$, the $1$-morphism is identified with $M_{\isoc,\irr}^{V=\id}\times \Gm^\an\times B\Gm^\an\rightarrow M_{\isoc,\irr}^{V=\id}$.
\end{defn}

Using the weight theory of $F$-isocrystals, we can now describe the tangent spaces of $(\Fix V)_\irr$ and $M_{\isoc,\irr}^{V=\id}$.

\begin{prop}\label{prop:tangent space of MFisoc}
Let $L$ be a finite extension of $K$. Take $x\in \Ob((\Fix V)_\irr)_{\Spa L})$ and continue to write $x$ for the image in  $M_{\isoc,\irr}^{V=\id}(\Spa L)$. Then we have
\[
\dim_L T_{(\Fix V)_\irr,x}=1\quad\text{and}\quad T_{M_{\isoc,\irr}^{V=\id},x}=0.
\]
\end{prop}

\begin{proof}
By abuse of notation, we also write $x$ for the image in $\Ob((\calM_{\isoc,\irr})_{\Spa L})$. Obviously, the map $T_{V^{-1}(\calM_{\isoc,\irr}),x}\rightarrow T_{\calM_{\isoc,\irr},x}$ is an isomorphism and the same holds for $\Inf$. Lemma~\ref{lem:T exact sequence} gives a long exact sequence
\[
\xymatrix{
0 \ar[r] & \Inf_{(\Fix V)_\irr,x}\ar[r] & (\Inf_{\calM_{\isoc,\irr},x})^{\oplus 2}\ar[r]^-\alpha& (\Inf_{\calM_{\isoc,\irr},x})^{\oplus 2}\ar[dll]\\
  & T_{(\Fix V)_\irr,x}\ar[r] & (T_{\calM_{\isoc,\irr},x})^{\oplus 2}\ar[r]^-\beta& (T_{\calM_{\isoc,\irr},x})^{\oplus 2},
}
\]
where both $\alpha$ and $\beta$ are given by $\begin{psmallmatrix}\id&\id\\ \id&dV_x\end{psmallmatrix}$.

By Proposition~\ref{prop:Fix V and F-isocrystals}, $x$ corresponds to an $L$-$F$-isocrystal $\calE$ such that the associated $\overline{\Q}_p$-$F$-isocrystal $\calE\otimes_L\overline{\Q}_p$ is irreducible. 
The image of $x$ in $\Ob((\calM_{\andR})_{\Spa L})$ defines an integrable $L$-connection $\nabla_L^\an\colon\cE_L^\an\rar\cE_L^\an\otimes_{\cO_{Z_{L}^\an}}\Omega^{1}_{Z_{L}^\an/L}$. From Proposition~\ref{prop:analytic Inf and T}, Theorem~\ref{thm:rigid-crystalline comparison for F-isocrystals}, Remark~\ref{rem:two de Rham cohomologies for L-F-isocrystals} and the constructions, we deduce natural $L$-vector space isomorphisms
\begin{align*}
&\Inf_{\calM_{\isoc,\irr},x}=\Inf_{\calM_{\andR},x}\cong H^{0}_{\dR}(Z_{L}^\an,\underline{\End}(\cE_L^\an,\nabla_L^\an))\cong H^0_\cris(Z/W,\underline{\End}(\calE))\;\text{and}\\
&T_{\calM_{\isoc,\irr},x}=T_{\calM_{\andR},x}\cong H^{1}_{\dR}(Z_{L}^\an,\underline{\End}(\cE_L^\an,\nabla_L^\an))\cong H^1_\cris(Z/W,\underline{\End}(\calE))    
\end{align*}
such that $dV_x$ on the leftmost side corresponds to the Frobenius on the rightmost side. It is straightforward to see that $dV_x$ on $\Inf_{\calM_{\isoc,\irr},x}$ is the identity. We know from Theorem~\ref{thm:weight theory for F-isocrystals} that the endomorphism $\underline{\End}(\calF\otimes_L\overline{\Q}_p)$ is pure of weight $0$ and thus the eigenvalues of $dV_x$ on $T_{\calM_{\isoc,\irr},x}$ are pure of weight $1$. We conclude from these that $\dim_L T_{(\Fix V)_\irr,x}=1$.

Similarly, Lemma~\ref{lem:T exact sequence} gives an exact sequence
\[
(\Inf_{M_{\isoc,\irr},x})^{\oplus 2}
\rar 
T_{M_{\isoc,\irr}^{V=\id},x}
\rar
(T_{M_{\isoc,\irr},x})^{\oplus 2}
\xrar{\left(\begin{smallmatrix}\id&\id\\ \id&dV_x\end{smallmatrix}\right)}
(T_{M_{\isoc,\irr},x})^{\oplus 2}.
\]
Since $M_{\isoc,\irr}$ is the coarse moduli space of $\calM_{\isoc,\irr}$ (and thus an \'etale space), we have $\Inf_{M_{\isoc,\irr},x}=0$ and $T_{\calM_{\isoc,\irr},x}\cong T_{M_{\isoc,\irr},x}$. Hence  $T_{M_{\isoc,\irr}^{V=\id},x}=0$.
\end{proof}

\begin{thm}\label{thm:irreducible F-isocrystals discrete}
The \'etale space $M_{\isoc,\irr}^{V=\id}$ is of the form $\coprod_{i\in I}\Spa(L_{i},\cO_{L_{i}})$, where each $L_{i}$ is a finite extension of $K$.
\end{thm}

\begin{proof}
We need to show that $M_{\isoc,\irr}^{V=\id}$ is an adic space \'etale over $K$. Take a presentation $U/R\xrightarrow{\cong}M_{\isoc,\irr}^{V=\id}$ by an \'etale equivalence relation $s,t\colon R\rightrightarrows U$ in $\Rig_K$. We first show that the adic space $U$ is \'etale over $K$. Since $\Omega_{U/K}$ is a coherent sheaf on $U$, it suffices to show that the stalk $(\Omega_{U/K})_u$ is zero for every $u\in U$ such that $k(u)$ is finite over $K$, which follows from Lemma~\ref{lem:tangent space of etale spaces} and Proposition~\ref{prop:tangent space of MFisoc}. Hence $R$ is also \'etale over $K$. It follows that $(t,s)\colon R\rightarrow U\times_{\Spa K}U$ is a closed immersion. By
Theorem~\ref{thm:effective-etale-equivalence-relation}, we conclude that $M_{\isoc,\irr}^{V=\id}$ is an adic space \'etale over $K$.
\end{proof}

\begin{rem}
Fix $r\geq 1$. Let $M_{\isoc,\irr,\rank r}^{V=\id}\subset M_{\isoc,\irr}^{V=\id}$ denote the open subspace of rank $r$ $F$-isocrystals. Then $M_{\isoc,\irr,\rank r}^{V=\id}$ is quasi-compact, namely, written as a \emph{finite} disjoint union of $\Spa(L_i,L_i^+)$'s. This follows from \cite[Cor.~4.3]{Abe-Esnault} or \cite[Cor.~4.40]{Kedlaya-companionI}; note that the cited result claims finiteness up to the twist by rank one $F$-isocrystals on $\Spec k$, but two objects in $(\Fix V)_\irr$ define the same point in $M_{\isoc,\irr}^{V=\id}$ if and only if they differ by such a twist.
\end{rem}

\begin{rem}
The geometry of Verschiebung on $\cM_{\cris}$ is expected to be more complicated. For example, \cite{Laszlo} and \cite{Yang} show that there is a non-trivial positive-dimensional flat family of sheaves on a genus $2$ curve over a finite field of characteristic $2$ fixed by a power of Verschiebung.
\end{rem}

\subsection{The stable locus in the log case}\label{section:stable locus in the log case}

Let us briefly discuss the log case. 
Keep the notation in \S\ref{section:moduli of crystals and isocrystals} and further assume that $Z_W\rightarrow \Spec W$ is projective (but possibly $D_W\neq \emptyset$).
We set $\calM_\dR\coloneqq \calM_\dR((Z_K.D_K)/K)$ and $\calM_\andR\coloneqq \calM_\andR((Z_K^\an,D_K^\an))$.

We start with the stable locus of the algebraic stack $\calM_\dR$ by following \cite{Nitsure-Log}.

\begin{defn}[{\cite[Def.~2.2, 2.6]{Nitsure-Log}}]
Let $u\colon\Spec \kappa\rightarrow \Spec K$ be a geometric point. An integrable $\kappa$-connection $\nabla\colon\cE_{\kappa}\rar\cE_{\kappa}\otimes_{\cO_{Z_{\kappa}}}\omega^{1}_{Z_{\kappa}/\kappa}$ on $Z_\kappa$ is called \emph{stable} if $\cE_{\kappa}$ is torsion-free, and for any nonzero proper $\kappa$-subconnection $(\cE'_{\kappa},\nabla')\subset (\cE_{\kappa},\nabla)$, we have $p_{\cE'_{\kappa}}(N)<p_{\cE_{\kappa}}(N)$ for any $N\gg0$.
More generally, for an $K$-scheme $T$, an integrable $T$-connection $(\calE,\nabla)$ on $Z_T$ is \emph{stable} if the pullback along every geometric point of $T$ is stable.

Define the \emph{moduli stack of stable integrable connections}
\[
\calM_{\dR,\st}\coloneqq \cM_{\dR,\st}((Z_K,D_K)/K)\rar(\Sch/K)_{\fppf}
\]
to be the strictly full subcategory of $\cM_{\dR}$ whose objects over $T$ are the stable integrable $T$-connections.
\end{defn}

\begin{lem}
Let $\kappa$ is an algebraically closed field over $K$. Every irreducible integrable $\kappa$-connection on $Z_\kappa$ is stable.
\end{lem}

\begin{proof}
It is easy to see that the torsion-subsheaf of the underlying coherent sheaf $\calE_\kappa$ is stable under the connection, so it must be zero. Hence $\calE_\kappa$ is torsion-free, and the stability follows from \cite[Prop.~2.3]{Nitsure-Log} by reducing it to the case $\kappa=\C$.
\end{proof}

\begin{prop}\label{prop:Mdr st in log case}
The natural $1$-morphism $\cM_{\dR,\st}\rightarrow\cM_{\dR}$ is an open immersion. Moreover, $\cM_{\dR,\st}$ is a $\Gm$-gerbe.
\end{prop}

\begin{proof}
The proofs of Propositions~\ref{prop:MdRst} and \ref{prop:MdRst is Gmgerbe} work, thanks to \cite[Prop.~2.5, 2.15]{Nitsure-Log}.
\end{proof}

\begin{defn}
Define $\calM_{\andR,\st}\rightarrow \Afd_K$ to be the strictly full subcategory $\calM_{\andR,\st}\subset \calM_{\andR}$ whose objects over $T\in\Ob(\Afd_{K})$ are integrable $T$-connections $(\calE,\nabla)$ on $Z_T^\an$ such that for every algebraically closed affinoid field $\overline{t}=\Spa(L,L^+)\in\Ob(\Afd_K)$ and every map $\overline{t}\rightarrow T$, the algebraization of the pullback $(\calE,\nabla)|_{Z_{\overline{t}}^{\an}}$ is stable as a $\overline{t}$-connection.
It follows from Proposition~\ref{prop:Mdr st in log case} that $\alpha\colon \calM_\dR^\an\xrightarrow{\cong}\calM_{\andR}$ restricts to an equivalence $\calM_{\dR,\st}^\an\xrightarrow{\cong}\calM_{\andR,\st}$. In particular, $\calM_{\andR,\st}$ is an open substack of $\calM_\andR$.
More generally, for any open substack $\calV$ of $\calM_\andR$, we define an open substack $\calV_\st$ of $\calV$ in a similar way.
\end{defn}

\begin{defn}
Define $(\Fix V)_{\st}$ to be the open substack of $\Fix V$ that is equivalent to the open immersion $\Fix V\times_{q,\calM_\isoc}\calM_{\isoc,\st}\rightarrow \Fix V$.
It follows from Lemma~\ref{lem:intersection of gerbes is gerbe} that  $(\Fix V)_\st$ is a $\Gm^\an$-gerbe and the projection $(\Fix V)_\st\rightarrow \calM_{\isoc,\st}$ is a $1$-morphism of $\Gm^\an$-gerbes.

Let $M_{\isoc,\st}$ denote the coarse moduli space of the $\Gm^\an$-gerbe $\calM_{\isoc,\st}$ and define an open \'etale subspace $V^{-1}(M_{\isoc,\st})\subset M_{\isoc,\st}$ similarly.
Define the \'etale space $M_{\isoc,\st}^{V=\id}$ by the following fiber product
\[
\xymatrix{
M_{\isoc,\st}^{V=\id}\ar[r]\ar[d]
& V^{-1}(M_{\isoc,\st})\ar[d]^-{(\id,V)}
\\
V^{-1}(M_{\isoc,\st})\ar[r]^-{(\id,\mathrm{incl})} 
& V^{-1}(M_{\isoc,\st})\times_{\Spa K}M_{\isoc,\st}.
}
\]
For every finite extension $L$ over $K$, the set $M_{\isoc,\st}^{V=\id}(\Spa L)$ is naturally identified with the set of isomorphism classes of pairs $(\calF,\iota_{\calF,L})$ where $\calF$ is an isocrystal on $((Z,M_Z)/W)_\cris$ and $\iota_{\calF,L}\colon L\rightarrow \End(\calF)$ is a $K$-algebra map such that
\begin{enumerate}
  \item $(\calF,\iota_{\calF,L})$ underlies an $L$-$F$-isocrystal; and
  \item the $L$-connection on $Z_L^\an$ corresponding to $(\calF,\iota_{\calF,L})$ is stable.
\end{enumerate}
We have a natural $1$-morphism $(\Fix V)_\st\rightarrow M_{\isoc,\st}^{V=\id}$.
\end{defn}

\begin{prop}\label{prop:tangent space of MFisoc in log case}
Let $L$ be a finite extension of $K$. Take $x\in \Ob((\Fix V)_\st)_{\Spa L})$ and continue to write $x$ for the image in  $\Ob((M_{\isoc,\st}^{V=\id})_{\Spa L})$. If the $\overline{\Q}_p$-$F$-isocrystal attached to $x$ is irreducible, then
\[
\dim_L T_{(\Fix V)_\st,x}=1\quad\text{and}\quad T_{M_{\isoc,\st}^{V=\id},x}=0.
\]
\end{prop}

\begin{proof}
The proof of Proposition~\ref{prop:tangent space of MFisoc} works, thanks to the second part of Theorem~\ref{thm:weight theory for F-isocrystals}(ii).
\end{proof}

\begin{thm}
The \'etale space $M_{\isoc,\st}^{V=\id}$ admits an open and closed subspace $M_{\isoc,\irr}^{V=\id}$ satisfying the following properties:
\begin{enumerate}
 \item for every finite extension $L$ over $K$, a point $x\in M_{\isoc,\st}^{V=\id}(\Spa L)$ lies in $M_{\isoc,\irr}^{V=\id}$ if and only if the $\overline{\Q}_p$-$F$-isocrystal attached to $x$ is irreducible;
 \item $M_{\isoc,\irr}^{V=\id}$ is of the form $\coprod_{i\in I}\Spa(L_{i},\cO_{L_{i}})$, where each $L_{i}$ is a finite extension of $K$.
\end{enumerate}
\end{thm}

\begin{proof}
Take a presentation $U/R\xrightarrow{\cong}M_{\isoc,\st}^{V=\id}$ by an \'etale equivalence relation $s,t\colon R\rightrightarrows U$ in $\Rig_K$. Take any $u\in U$ such that $k(u)$ is finite over $K$ and the $\overline{\Q}_p$-$F$-isocrystal attached to $u$ is irreducible. 
It follows from Lemma~\ref{lem:tangent space of etale spaces} and Proposition~\ref{prop:tangent space of MFisoc in log case} that $(\Omega_{U/K})_u=0$, which implies that $u$ has an open neighborhood $V_u\subset U$ of the form $\Spa(k(u),\calO_{k(u)})$. Set $U_u\coloneqq s(t^{-1}(V_u))\subset U$. Then $U_u$ is an open neighborhood of $u$ that is \'etale over $K$ and it satisfies $s^{-1}(U_u)=t^{-1}(U_u)$ (cf.~\cite[044G]{stacks-project}). Define $M_{\isoc,\irr}^{V=\id}$ to be the union of open subspaces $U_u/s^{-1}(U_u)$ of $M_{\isoc,\st}^{V=\id}$ where $u$ runs over the points that yield irreducible $\overline{\Q}_p$-$F$-isocrystals. By construction, $M_{\isoc,\irr}^{V=\id}$ satisfies conditions (i) and (ii). So it remains to show that $M_{\isoc,\irr}^{V=\id}$ is also a closed subspace of $M_{\isoc,\st}^{V=\id}$. 
To see this, observe that $M_{\isoc,\st}^{V=\id}$ is written as the disjoint union $\coprod M_{\isoc,\st,\rank r}^{V=\id}$ where $M_{\isoc,\st,\rank r}^{V=\id}$ is an open and closed subspace of points whose underlying coherent sheaves are geometric fiberwise generically locally free of rank $r$.
It follows from \cite[Cor.~4.3]{Abe-Esnault} or \cite[Cor.~4.40]{Kedlaya-companionI} that $M_{\isoc,\st,\rank r}^{V=\id}\cap M_{\isoc,\irr}^{V=\id}$ is finite \'etale over $K$ and thus it is also closed in $M_{\isoc,\st,\rank r}^{V=\id}$.
\end{proof}

\section{Rank one case}\label{sec:rank 1}

In this section, we describe the case of rank $1$ isocrystals in more detail. 
Let $k=\bF_{q}$ be a finite field, and let $W=W(k)$. Suppose that $Z_{W}\rar W$ is either an abelian scheme or a relative curve (i.e., a smooth proper morphism of relative dimension $1$ with geometrically connected fibers). Set $Z\coloneqq Z_{W}\times_{\Spec W}\Spec k$.

Recall from Definition~\ref{def:MFisoc} that there exists an \'etale space $M_{\isoc,\irr}^{V=\id}$ such that, for every finite extension $L$ over $K$, the set $M_{\isoc,\irr}^{V=\id}(\Spa L)$ is naturally identified with the set of isomorphism classes of absolutely irreducible $L$-isocrystals on $Z$ underlying $L$-$F$-isocrystals.

\begin{thm}\label{thm:counting rank 1 F-isocrystals}
\hfill
\begin{enumerate}
 \item[(i)] The open subspace $M_{\isoc,1}^{V=\id}\subset M_{\isoc,\irr}^{V=\id}$ consisting of rank one objects is finite \'etale over $\Spa K$ of degree equal to $\lvert\Pic^{0}(Z/k)(k)\rvert$.
 \item[(ii)] The number of rank one $\overline{\Q}_p$-$F$-isocrystals on $Z$ up to twists by rank one $\overline{\bQ}_{p}$-$F$-isocrystals on $\Spec k$ is equal to $\lvert\Pic^{0}(Z/k)(k)\rvert$.
\end{enumerate}
\end{thm}

We note that (i) and (ii) are equivalent. Indeed, we know from Theorem~\ref{thm:irreducible F-isocrystals discrete} that $M_{\isoc,1}^{V=\id}$ is \'etale over $K$. Moreover, the direct limit
\[
\varinjlim_{L/K\;\text{finite}}M_{\isoc,1}^{V=\id}(\Spa L) 
\]
represents the set of isomorphism classes of rank one $\overline{\Q}_p$-isocrystals on $Z$ underlying $\overline{\Q}_p$-$F$-isocrystals; it is easy to identify this set with the set of rank one $\overline{\Q}_p$-$F$-isocrystals on $Z$ up to twists by rank one $\overline{\bQ}_{p}$-$F$-isocrystals on $\Spec k$.

In fact, the following strengthening of Theorem~\ref{thm:counting rank 1 F-isocrystals}(i) holds, which was mentioned in Introduction.

\begin{thm}
\label{thm:counting rank 1 F-isocrystals Fqm}
For every $m\geq 1$, the open subspace $M_{\isoc,1}^{V^m=\id}\subset M_{\isoc,\irr}^{V^{m}=\id}$ (defined in a similar way) is finite \'etale over $K$ of degree equal to $\lvert\Pic^{0}(Z/k)(\bF_{q^{m}})\rvert$.
\end{thm}

This section studies the explicit geometry of the moduli spaces of rank one crystals and isocrystals. As a consequence, we will show Theorem~\ref{thm:counting rank 1 F-isocrystals}(i), hence also (ii) and Theorem~\ref{thm:counting rank 1 F-isocrystals Fqm}, in \S\ref{sec:Counting rank 1 F-isocrystals using geometry}.

\begin{rem}
Of course, the geometric class field theory implies Theorem~\ref{thm:counting rank 1 F-isocrystals}(ii): upon the choice of a geometric point $\overline{z}$ of $Z$, the image of $\pi_{1}^{\et}(Z_{\overline{k}},\overline{z})$ in the abelianization $\pi_{1}^{\et}(Z,\overline{z})^{\mathrm{ab}}$ is naturally identified with $\Pic^{0}(Z/k)(k)$ (see \cite[Thm.~1.3.1]{Deligne-WeilII}). 
Hence, combined with Theorem~\ref{thm:irreducible F-isocrystals discrete}, we get Theorem~\ref{thm:counting rank 1 F-isocrystals}(i). The main purpose of this section is to illustrate that the geometry involving rank one objects can be made very explicit.
\end{rem}

\begin{rem}
In his IHES lectures \cite{Deligne-IHES}, Deligne sketched a proof of a result analogous to Theorem~\ref{thm:counting rank 1 F-isocrystals} in the case where $p>2$ and $Z_{W}/W$ is a relative curve. This used the nilpotent crystalline site interpretation of the universal vector extensions.
\end{rem}

In \S\ref{sec:universal vector extension}, we review the notion of the \emph{universal vector extension} of an abelian scheme, and relate it with $\cM_{\dR,\Bun_{1}}$ of both an abelian scheme and a relative curve. Subsection~\ref{sec:Moduli stacks of rank 1 crystals and isocrystals} study $\calM_{\cris,1}$ and $\calM_{\isoc,1}$ as well as their coarse moduli spaces $M_{\cris,1}$ and $M_{\isoc,1}$ in terms of the universal vector extension. Finally, Theorem~\ref{thm:geometric description of M1cris,BunV=id} in \S\ref{sec:Counting rank 1 F-isocrystals using geometry} will show that $M_{\cris,1}^{V=\id}$ is finite flat over $\Spf W$ of degree $\left|\Pic^{0}(Z/k)(k)\right|$, which easily gives Theorem~\ref{thm:counting rank 1 F-isocrystals}(i).

\subsection{Universal vector extensions}\label{sec:universal vector extension}

We review the universal vector extension. For this subsection only, let $S$ be a Noetherian scheme.

\smallskip
\noindent\textbf{Universal vector extensions of abelian schemes.} 
Let $f\colon \cA\rar S$ be a projective abelian scheme, and let $e\colon S\rar\cA$ be the identity section. Note that the isomorphism $\calO_S\cong f_\ast\calO_\calA$ induces an isomorphism $e^\ast \Omega_{\calA/S}^1\cong f_\ast\Omega_{\calA/S}^1$ of locally free $\calO_S$-modules which is compatible with arbitrary base change $T\rightarrow S$.

Consider the Picard functor $\Pic(\calA/S)$ for $\calA/S$: 
for an $S$-scheme $T$, 
\[
\Pic(\calA/S)(T)\coloneqq\Pic(\calA_T)/f_T^\ast\Pic(T),
\]
which will often be identified with the $e$-rigidified Picard functor evaluated at $T$, namely, the set of isomorphism classes of pairs $(\calL,u)$ where $\calL$ is an invertible sheaf on $\calA_T$ and $u\colon \calO_T\xrightarrow{\cong}e_T^\ast\calL$ is an $\calO_T$-linear isomorphism; $\Pic(\calA/S)$ is representable by a scheme. Let $\cA^{\ast}\coloneqq \Pic^{0}(\calA/S)$ denote the (relative) identity component. Then $\Pic^{0}(\calA/S)$ is representable by a projective abelian scheme over $S$ and $\Pic^{0}(\calA/S)(T)$ consists of isomorphism classes of pairs $(\calL,u)$  such that for every geometric point $\overline{t}$ of $T$, the pullback $\calL_{\calA_{\overline{t}}}$ is algebraically equivalent to $\calO_{\calA_{\overline{t}}}$.
We call $\calA^\ast$ the \emph{dual abelian scheme}; there is a canonical isomorphism $\cA\xrightarrow{\sim}(\cA^{\ast})^{\ast}$ given by the Poincar\'e bundle (e.g., \cite[Cor.~6.8]{GIT}, \cite[Thm.~8.2.5]{BLR}, \cite[Lem.~9.2.9, Prop.~9.5.10]{FGA}).

\begin{defn}[{\cite[p.~2, Prob.~B]{MM}}]
The \emph{universal vector extension} $\univext{\cA^{\ast}/S}$ of $\cA^{\ast}$ over $S$ is a sheaf on $(\Sch/S)_{\fppf}$ which is universal among extensions of $\cA^{\ast}$ by quasi-coherent $\cO_{S}$-modules. 
\end{defn}

\begin{thm}\label{thm:moduli interpretation of universal vector extension}
The universal vector extension $\univext{\cA^{\ast}/S}$ exists and is representable by a smooth $S$-group scheme, sitting in a short exact sequence of fppf sheaves
\[
0\rar e^{\ast}\Omega^{1}_{\cA/S}\rar\univext{\cA^{\ast}/S}\rar \cA^{\ast}\rar0.
\]
Moreover, it represents the functor assigning to an $S$-scheme $T$ the set of isomorphism classes of triples $(\calL,u,\nabla)$ where $(\calL,u)\in \Pic^{0}(\calA/S)(T)$ and $\nabla\colon\cL\rar\cL\otimes_{\cO_{\cA_{T}}}\Omega^{1}_{\cA_{T}/T}$ is an integrable connection, 
\end{thm}

\begin{proof}
The first part is \cite[1.10, 2.6]{MM}, and the second is essentially the combination of \cite[Prop.~2.6.7, 3.2.3, 4.2.1, (4.6.3)]{MM}. 
\end{proof}

\begin{rem}\label{rem:identification of universal vector extension and MdR10}
In fact, a closer look at \cite{MM} also shows the following compatibility. Let $M$ denote the functor in Theorem~\ref{thm:moduli interpretation of universal vector extension}. The kernel $\Ker(M\rightarrow \Pic^0(\calA/S)=\calA^\ast)$ (as an fppf abelian sheaf) is canonically identified with $f_\ast\Omega_{\calA/S}^1$, and the identification $E(\calA^\ast/S)\cong M$ is compatible with the projection to $\calA^\ast$ under the natural identification $e^\ast \Omega_{\calA/S}^1\cong f_\ast\Omega_{\calA/S}^1$; see, in particular, \cite[p.~3, pp.~15-16, p.~40]{MM}.
\end{rem}

Recall from Definition~\ref{def:moduli stack of vector bundles with integrable connections} and Proposition~\ref{prop:MdRBun} the $1$-morphism $\cM_{\dR,\Bun_{1}}(\cA/S)\rightarrow \Bun_1(\calA/S)$ of algebraic stacks over $S$. 
As in Proposition~\ref{prop:MdRst is Gmgerbe}, we see that the multiplication by an invertible scalar on the underlying invertible sheaf makes $\cM_{\dR,\Bun_{1}}(\cA/S)$ and $\Bun_{1}(\cA/S)$ $\Gm$-gerbes. Moreover, it is clear that the above $1$-morphism is compatible with $\Gm$-structures. 

\begin{lem}\label{lem:MdRBun1 split Gm gerbe, abelian scheme case}
The coarse moduli space of $\Bun_1(\calA/S)$ is $\Pic(\calA/S)$, and the $1$-morphism $\cM_{\dR,\Bun_{1}}(\cA/S) \rightarrow\Bun_1(\calA/S)$ is a $1$-morphism of split $\Gm$-gerbes. 
\end{lem}

The lemma claims the equivalences
\[
\cM_{\dR,\Bun_{1}}(\cA/S)\cong M_{\dR,\Bun_{1}}(\cA/S)\times_{S}B\Gm\quad\text{and}\quad \Bun_1(\calA/S)\cong \Pic_{\calA/S}\times_{S}B\Gm,
\]
where $M_{\dR,\Bun_{1}}(\cA/S)$ is the coarse moduli space of $\cM_{\dR,\Bun_{1}}(\cA/S)$. 

\begin{proof}
By \cite[Prop. 4.2.5]{Wang}\footnote{\emph{Op.~cit.} assumes that $S$ is of characteristic $0$, but the proof works in general.}, the $1$-morphism $\Pic(\calA/S)\times_{S}B\Gm\rightarrow \Bun_{1}(\cA/S)$ sending the pair $((\calL,u),\calN)\in \Pic(\calA/S)(T)\times \Ob(B\Gm)_T$ to $\calL\otimes_{\calO_{\calA_T}} f_T^\ast\calN$ over an $S$-scheme $T$ is an equivalence of $\Gm$-gerbes; namely, $\Bun_{1}(\cA/S)$ is a split $\Gm$-gerbe, whose coarse moduli space is $\Pic(\cA/S)$. It follows that $\cM_{\dR,\Bun_{1}}(\cA/S)$ is also a split $\Gm$-gerbe.
\end{proof}

Define the algebraic stacks $\cM_{\dR,\Bun_{1}}(\cA/S)^{0}$ and $\Bun_{1}(\cA/S)^{0}$ over $S$ by the open substacks sitting in the following Cartesian diagram
\[
\xymatrix{
\cM_{\dR,\Bun_{1}}(\cA/S)^{0}\ar[r]\ar@{^{(}->}[d] & \Bun_{1}(\cA/S)^{0}\ar[r]\ar@{^{(}->}[d] & \Pic^0(\calA/S)\ar@{^{(}->}[d]\\
\cM_{\dR,\Bun_{1}}(\cA/S)\ar[r]& \Bun_{1}(\cA/S)\ar[r] & \Pic(\calA/S).
}
\]

\begin{lem}\label{lem:coarse moduli is universal vector extension, abelian scheme case}
The coarse moduli space $M_{\dR,\Bun_{1}}(\cA/S)^{0}$ of $\cM_{\dR,\Bun_{1}}(\cA/S)^{0}$ is identified with the universal vector extension $\univext{\cA^{\ast}/S}$.
\end{lem}

\begin{proof}
This follows from Theorem~\ref{thm:moduli interpretation of universal vector extension} and Lemma~\ref{lem:MdRBun1 split Gm gerbe, abelian scheme case}.
\end{proof}

\smallskip
\noindent\textbf{Universal vector extensions of Jacobians.} 
Let $f\colon \cC\rar S$ be a relative curve together with a section $e\colon S\rar \cC$. Define the Picard functor $\Pic(\calC/S)$ as in the abelian scheme case; $\Pic(\calS/S)$ is representable by an $S$-scheme, and the (relative) identity component $\Jac({\cC/S})\coloneqq\Pic^{0}({\cC/S})$, the \emph{Jacobian} of $\cC$ over $S$, is a projective abelian scheme over $S$ equipped with the canonical isomorphism $\varphi_\Theta\colon \Jac({\cC/S})\xrightarrow{\cong}\Jac({\cC/S})^\ast$ given by the theta divisor. It also comes with the \emph{Abel--Jacobi embedding} $\iota\colon \cC\rar \Jac({\cC/S})$, sending $c\in \calC(T)$ to $\calO_{\calC_T}(c-e_T)\in \Jac({\cC/S})(T)$ for any $S$-scheme $T$. Then $\iota^\ast\colon \Jac(\calC/S)^\ast=\Pic^0(\Jac(\calC/S)/S)\rightarrow \Jac(\calC/S)=\Pic^0(\calC/S)$ is the inverse to $-\varphi_\Theta$. For these statements, see \cite[Thm.~9.3.1, Prop.~9.4.4, Pf.]{BLR} and \cite[Lem.~6.9, Rem.~6.10]{MilneJV}.

As in the abelian scheme case, consider the $1$-morphism $\cM_{\dR,\Bun_{1}}(\cC/S)\rightarrow \Bun_1(\calC/S)$ of $\Gm$-gerbes over $S$.
Let $M_{\dR,\Bun_{1}}(\cA/S)$ denote the coarse moduli space of $\cM_{\dR,\Bun_{1}}(\cA/S)$. 

\begin{lem}\label{lem:MdRBun1 split Gm gerbe, curve case}
The coarse moduli space of $\Bun_1(\calC/S)$ is $\Pic(\calC/S)$, and the $1$-morphism $\cM_{\dR,\Bun_{1}}(\cC/S) \rightarrow\Bun_1(\calC/S)$ is a $1$-morphism of split $\Gm$-gerbes. 
\end{lem}

\begin{proof}
The proof of Lemma~\ref{lem:MdRBun1 split Gm gerbe, abelian scheme case} works verbatim in this case.
\end{proof}

Define the algebraic stacks $\cM_{\dR,\Bun_{1}}(\cC/S)^{0}$ and $\Bun_{1}(\cC/S)^{0}$ over $S$ by the restriction to $\Pic^0(\calC/S)\subset \Pic(\calS)$ as in the abelian scheme case. The pullback along $\iota\colon \calC\rightarrow \Jac(\calC/S)$ gives the commutative diagram
\[
\xymatrix{
\calM_{\dR,\Bun_1}(\Jac(\calC/S)/S)^0\ar[r]\ar[d]_-{\iota^\ast} & \Bun_1(\Jac(\calC/S)/S)^0 \ar[r]\ar[d]_-{\iota^\ast} & \Jac(\calC/S)^\ast \ar[d]_-{\iota^\ast}^-\cong\\
\calM_{\dR,\Bun_1}(\calC/S)^0\ar[r] & \Bun_1(\calC/S)^0\ar[r] & \Jac(\calC/S).
}
\]
Since the right vertical map is an isomorphism, the middle one is an equivalence by Lemma~\ref{lem:MdRBun1 split Gm gerbe, curve case}.

\begin{lem}\label{lem:coarse moduli is universal vector extension, relative curve case}
The pullback $\iota^{\ast}\colon \cM_{\dR,\Bun_{1}}(\Jac({\cC/S})/S)^{0}\xrightarrow{\sim}\cM_{\dR,\Bun_{1}}(\cC/S)^{0}$ is also an equivalence, and the coarse moduli space $M_{\dR,\Bun_{1}}(\cC/S)^{0}$ of $\cM_{\dR,\Bun_{1}}(\cC/S)^{0}$ is naturally identified with the universal vector extension $\univext{\Jac(\cC/S)/S}$.
\end{lem}

Note that $\iota^\ast \univext{\Jac(\cC/S)^\ast/S}\xrightarrow{\cong}\univext{\Jac(\cC/S)/S}$ is an isomorphism, so one may also identify $M_{\dR,\Bun_{1}}(\cC/S)^{0}$ with the former.

\begin{proof}
For simplicity, write $J\coloneqq \Jac({\cC/S})$ in this proof. 
For the first assertion, it is enough to show that the induced morphism $\iota^\ast\colon M_{\dR,\Bun_1}(J/S)^0\rightarrow M_{\dR,\Bun_1}(\calC/S)^0$ is an isomorphism of fppf sheaves on $S$.
Let $f^J\colon J\rightarrow S$ denote the structure map and write $e^J\colon S\rightarrow J$ for the identity section.
Consider the diagram of fppf sheaves
\[
\xymatrix{
0\ar[r] & f^J_\ast \Omega_{J/S}^1 \ar[r]\ar[d]_-\alpha & M_{\dR,\Bun_1}(J/S)^0 \ar[r]^-\beta\ar[d]_-{\iota^\ast} & J^\ast \ar[r]\ar[d]_-{\iota^\ast}^-\cong  & 0\\
0\ar[r] & f_\ast \Omega_{\calC/S}^1 \ar[r] & M_{\dR,\Bun_1}(\calC/S)^0 \ar[r]^-\gamma & J \ar[r]  & 0,
}
\]
where the horizontal maps are the obvious ones making the rows left exact, and $\alpha$ is induced from the maps $\Omega_{J/S}^1\rightarrow \iota_\ast\iota^\ast \Omega_{J/S}^1\rightarrow \iota_\ast\Omega_{\calC/S}^1$. It is straightforward to see that the left square is commutative. Since we know that $\beta$ is surjective, so is $\gamma$. To show that the middle vertical map is an isomorphism, it is enough to show that $\alpha$ is an isomorphism. Since $\alpha$ is a morphism of finite locally free $\calO_S$-modules, it suffices to check that $\alpha$ is an isomorphism after restricting to points of $S$, which is covered in \cite[Prop.~5.3(a)]{MilneJV}.

The second assertion follows by identifying the above diagram of fppf sheaves with the diagram
\[
\xymatrix{
0 \ar[r]& e^{J,\ast}\Omega_{J/S}^1 \ar[r]\ar[d] & E(J^\ast/S)\ar[r]\ar[d]_-{\iota^\ast}^-\cong & J^\ast \ar[r]\ar[d]_-{\iota^\ast}^-\cong & 0\\
0 \ar[r]& e^{J^\ast,\ast}\Omega_{J^\ast/S}^1\ar[r] & E(J/S) \ar[r] & J\ar[r] &0.
}
\]
This is tedious but straightforward as in Remark~\ref{rem:identification of universal vector extension and MdR10}, so we leave it to the reader.

\end{proof}

\subsection{Moduli stacks of rank one crystals and isocrystals}\label{sec:Moduli stacks of rank 1 crystals and isocrystals}

In this subsection, we define the moduli stacks of rank $1$ crystals and isocrystals and describe them using the universal vector extension.

\smallskip
\noindent
\textbf{Setup.}
Let $k=\bF_{q}$ be a finite field with $q=p^{a}$, $W=W(k)$, and $K=\Frac W$. Let $S=\Spec W$ and $f_W\colon Z_{W}\rightarrow S$ be either a relative curve or an abelian scheme.    

For a $W$-scheme $T$, we let $f_T\colon Z_{T}\coloneqq Z_{W}\times_{W}T\rightarrow T$. For simplicity, we use $Z\coloneqq Z_{k}$ and $f\coloneqq f_k$. 
Let $\Bun_1(Z_T/T)^0$ denote the open substack of $\Bun_1(Z_T/T)$ consisting of invertible sheaves whose geometric fibers are algebraically equivalent to zero. In other words, $\Bun_1(Z_T/T)^0\cong \Bun_1(Z_T/T)\times_{\Pic(Z_T/T)}\Pic^0(Z_T/T)$, which is a $\Gm$-gerbe with $\Pic^0(Z_T/T)$ being the coarse moduli space.\footnote{The definitions of $\Pic(Z_T/T)$ and $\Pic^0(Z_T/T)$ in \S\ref{sec:universal vector extension} are valid under the assumption that $f$ has a section. In general, the former is defined as the fppf sheafification of $T'\mapsto \Pic(Z_{T'})$. See \cite[Chap.~8,9]{BLR}.} 
Similarly, we define the open substack $\cM_{\dR,\Bun_{1}}(Z_{T}/T)^0$ of $\cM_{\dR,\Bun_{1}}(Z_{T}/T)$. The forgetful $1$-morphism $\cM_{\dR,\Bun_{1}}(Z_{T}/T)^0\rightarrow \Bun_1(Z_T/T)^0$ is a $1$-morphism of $\Gm$-gerbes. Let $M_{\dR,\Bun_{1}}(Z_{T}/T)^{0}$ denote the coarse moduli space of $\cM_{\dR,\Bun_{1}}(Z_{T}/T)^0$.

In \S\ref{section:moduli of crystals and isocrystals}, the formal moduli stack of crystals $\cM_{\cris}$ was defined as the completion of $\cM_{\dR}(Z_{W}/W)$ along $\lvert \calM_{\dR}(Z/k)^\nilp_\red\rvert$. In the same section, we also constructed the adic moduli stack of isocrystals $\cM_{\isoc}$ over $(Z/W)_\cris$ from $\cM_{\cris}$.

\begin{defn}
We define the formal algebraic stack $\cM_{\cris,1}$ as the formal algebraic open substack of $\cM_{\cris}$ equivalent to the $2$-fiber product
\[
\xymatrix{
\cM_{\cris,1}\ar@{^(->}[r]\ar[d]&
\cM_{\cris}\ar[d]\\
\cM_{\dR,\Bun_{1}}(Z_{W}/W)^{0}\ar@{^(->}[r]&
\cM_{\dR}(Z_{W}/W).
}
\]
Equivalently, $\cM_{\cris,1}$ is the completion of $\cM_{\dR,\Bun_{1}}(Z_{W}/W)^{0}$ along the closed subset $\lvert \calM_{\dR}(Z/k)^\nilp_\red\rvert \cap \lvert \cM_{\dR,\Bun_{1}}(Z_{W}/W)^{0}\rvert$.

We define the open substack $\calM_{\isoc,1}\subset \calM_{\andR} $
by applying Proposition~\ref{prop:image of comparison map from generic fiber to analytification} to the open immersion $(\calM_{\cris,1})_\eta\hookrightarrow (\calM_{\dR,W}^\wedge)_\eta$. 
We call the adic stack $\calM_{\isoc,1}$ the \emph{moduli stack of rank $1$ isocrystals} over $(Z/W)_\cris$.
\end{defn}

\begin{rem}\label{rem:moduli of rank 1 isocrystals}
If $T=\Spa(L,\calO_L)$ for a finite extension $L$ over $K$, then $\Ob((\calM_{\isoc,1})_T)$ is the groupoid consisting of rank $1$ locally free isocrystals on $(Z_{(\calO_L/p)}/\calO_L)_\cris$; this follows from Proposition~\ref{prop:Esnault-Shiho}. In general, we do not know whether every rank $1$ locally free isocrystals admits an underlying crystal that is locally free. To remedy this, one could consider an open substack $(\calM_{\isoc,1})'$ of $\cM_{\isoc}$ defined by the $2$-fiber product
\[
\xymatrix{
(\cM_{\isoc,1})'\ar@{^(->}[rr]\ar[d]&&
\cM_{\isoc}\ar[d]\\
(\cM_{\dR,\Bun_{1}}(Z_{K}/K)^{0})^{\an}\ar@{^(->}[r]&
\cM_{\dR}(Z_{K}/K)^{\an}\ar[r]^-\cong& \calM_\andR(Z_K^\an/K).
}
\]
Note that
the underlying invertible sheaf of any integrable connection over a field of characteristic $0$ is algebraically equivalent to zero, hence one obtains the same stack $(\cM_{\isoc,1})'$ even if one replaces $(\cM_{\dR,\Bun_{1}}(Z_{K}/K)^{0})^{\an}$ in the above diagram with $(\cM_{\dR,\Bun_{1}}(Z_{K}/K))^{\an}$.
Obviously, $\calM_{\isoc,1}\subset (\cM_{\isoc,1})'$. By the remark at the beginning, these two moduli stacks agree over $T=\Spa(L,\calO_L)$ for any finite extension $L$ over $K$, and thus yield the same results for $F$-isocrystals.
\end{rem}

\begin{prop}\label{prop:V on algebraically trivial locus}
The Verschiebung endomorphisms $V$ of $\cM_{\cris}$ and $\cM_{\isoc}$ restrict to $V\colon \cM_{\cris,1}\rar\cM_{\cris,1}$ and $V\colon \cM_{\isoc,1}\rar\cM_{\isoc,1}$, respectively.
\end{prop}

\begin{proof}
We know from Remark~\ref{rem:compatibility of functoriality for crystals and connections} that the $1$-morphism $\calM_\cris\times_{\Spf W}\Spec k\rightarrow \calM_\dR(Z/k)$, regarded as the inclusion of strictly full subcategories of $\calM_\dR(Z_W/W)$, is compatible with the restriction of $V$ on the source and the Frobenius pullback $\Fr_{Z}^{a,\ast}$ on the target.
Since $\Fr_{Z}^{a,\ast}$ preserves $\cM_{\dR}(Z/k)^{\nilp}$ and $\cM_{\dR,\Bun_{1}}(Z/k)^{0}$, we deduce that the Verschiebung on $\calM_\cris$ restricts to $V\colon \cM_{\cris,1}\rar\cM_{\cris,1}$. The assertion for $\calM_{\isoc,1}$ follows from this and Theorem~\ref{thm:Verschiebung on Misoc}.
\end{proof}

Let $\bG_{\mathrm{m},W}^{\wedge}$ denote the $p$-adic completion of $\bG_{\mathrm{m},W}$.

\begin{lem}
The formal algebraic stack $\calM_{\cris,1}$ is a $\bG_{\mathrm{m},W}^{\wedge}$-gerbe,\footnote{In the sense of Definition~\ref{def:gerbe} with ``algebraic'' and $\Gm$ replaced with ``formal algebraic'' and $\bG_{\mathrm{m},W}^{\wedge}$.}
and $\calM_{\isoc,1}$ is a $\Gm^{\an}$-gerbe.
\end{lem}

\begin{proof}
The first assertion follows easily from the fact that $\cM_{\dR,\Bun_{1}}(Z_{W}/W)^{0}$ is a $\bG_{\mathrm{m},W}$-gerbe and $\calM_{\dR}(Z/k)^\nilp_\red$ lives over $k$. The second follows from Proposition~\ref{prop:open sub of gerbe is gerbe} and Theorem~\ref{thm:Misocirr is Gm-gerbe}.
\end{proof}

\begin{defn}\label{defn:coarse moduli of rank 1 isocrystals}
Let $M_{\cris,1}$ denote the coarse moduli space of $\calM_{\cris,1}$ and let $M_{\isoc,1}$ denote the coarse moduli space of $\calM_{\isoc,1}$.
By definition, $M_{\cris,1}$ is a formal algebraic space over $\Spf W$, and $M_{\isoc,1}$ is an \'etale space over $K$. It is clear that the Verschiebung endomorphisms $V\colon \cM_{\cris,1}\rar\cM_{\cris,1}$ and $V\colon\cM_{\isoc,1}\rar\cM_{\isoc,1}$ descend to $M_{\cris,1}$ and $M_{\isoc,1}$, respectively. 
\end{defn}

We will have a more concrete description of $\cM_{\cris,1}$ and $\cM_{\isoc,1}$. We start by describing the reduced nilpotent locus 
\[
(\cM_{\dR,\Bun_{1}}(Z/k)^{0})^{\nilp}_{\red}\coloneqq\cM_{\dR}(Z/k)^{\nilp}_{\red}\cap\cM_{\dR,\Bun_{1}}(Z/k)^{0}.
\]

For an $\F_p$-scheme $T$, let $\Fr_T$ denote the absolute Frobenius on $T$. Consider the vector bundle $\mathbb{V}((f_\ast \Fr_{Z}^{\ast}\Omega^{1}_{Z/k})^\vee)\rightarrow\Spec k$ associated with $H^{0}(Z,\Fr^{\ast}_{Z}\Omega^{1}_{Z/k})$. 

\begin{defn}\label{defn:Hitchin map}
Define the \emph{Hitchin map} 
\[
\cH\colon \cM_{\dR,\Bun_{1}}(Z/k)^{0}\rightarrow \mathbb{V}((f_\ast \Fr_{Z}^{\ast}\Omega^{1}_{Z/k})^\vee)
\]
by assigning to $T\in \Ob(\Sch/k)$ and $(\cL,\nabla)\in \Ob((\cM_{\dR,\Bun_{1}}(Z/k)^{0})_T)$ the $p$-curvature $\cH(\cL,\nabla)\coloneqq\psi(\cL,\nabla)\in \Hom(\Der(Z_{T}/T),\Fr_{Z_{T},*}\underline{\End}(\calL))$; here, for the structure map $\pi\colon Z_{T}\rightarrow Z$, we have used identifications 
\begin{align*}
&\hspace{5mm} \mathbb{V}((f_\ast \Fr_{Z}^{\ast}\Omega^{1}_{Z/k})^\vee)(T)=\Hom(\cO_{Z_{T}},\pi^{\ast}\Fr_{Z}^{\ast}\Omega^{1}_{Z/k})
\\
&=\Hom(\cO_{Z_{T}},\Fr_{Z_T}^{\ast}\Omega^{1}_{Z_T/T}\otimes\underline{\End}(\cL))=\Hom(\underline{\Hom}(\Fr_{Z_T}^{\ast}\Omega^{1}_{Z_T/T},\cO_{Z_{T}}),\underline{\End}(\cL))
\\
&=\Hom(\Fr_{Z_T}^\ast\Der(Z_{T}/T),\underline{\End}(\cL))=\Hom(\Der(Z_{T}/T),\Fr_{Z_{T},*}\underline{\End}(\calL)),
\end{align*}
where $\Hom$, $\otimes$, $\underline{\End}$, and $\underline{\Hom}$ are all taken as $\calO_{Z_T}$-modules.

 This induces a map from the coarse moduli space 
\[
H\colon M_{\dR,\Bun_{1}}(Z/k)^{0}\rar \mathbb{V}((f_\ast \Fr_{Z}^{\ast}\Omega^{1}_{Z/k})^\vee), 
\]
which will also be referred to as the Hitchin map.

Let $0\colon\Spec k\rar \mathbb{V}((f_\ast \Fr_{Z}^{\ast}\Omega^{1}_{Z/k})^\vee)$ be the zero section, and let $\cH^{-1}(0)$ and $H^{-1}(0)$ be the pullback of the zero section along $\cH$ and $H$, respectively. 
\end{defn}

\begin{lem}
For any $k$-scheme $T$ and $(\calL,\nabla),(\calL',\nabla')\in \Ob((\cM_{\dR,\Bun_{1}}(Z/k)^{0})_T)$, we have
\[
\psi((\calL,\nabla)\otimes (\calL',\nabla'))=\psi(\calL,\nabla)+\psi(\calL',\nabla').
\]
In particular, $H$ is a morphism of commutative group schemes over $k$.
\end{lem}

\begin{proof}
This follows from a simple computation, which is left to the reader.
\end{proof}

\begin{lem}\label{lem:reduced nilpotent locus. rank 1}
As algebraic closed substacks of $\cM_{\dR,\Bun_{1}}(Z/k)^{0}$, we have
\[
(\cM_{\dR,\Bun_{1}}(Z/k)^{0})^{\nilp}_{\red}=\cH^{-1}(0).
\]
\end{lem}

\begin{proof}
For any $k$-scheme $T$, $(\cH^{-1}(0))_{T}\subset (\cM_{\dR,\Bun_{1}}(Z/k)^{0})_{T}$ is the strictly full subcategory consisting of the objects $(\cL,\nabla)$ with vanishing $p$-curvature. This implies that $\cH^{-1}(0)$ is a substack of $(\cM_{\dR,\Bun_{1}}(Z/k)^{0})_{\red}^{\nilp}$. Conversely, if $T$ is a \emph{reduced} $k$-scheme, then $((\cM_{\dR,\Bun_{1}}(Z/k)^{0})^{\nilp}_{\red})_{T}$ is the strictly full subcategory of $(\cM_{\dR,\Bun_{1}}(Z/k)^{0})_{T}$ consisting of the objects $(\cL,\nabla)$ with nilpotent $p$-curvature. By the locally free case of Proposition~\ref{prop:bounding exponent of nilpotence}, any nilpotent endomorphism of $\cL$ is $0$, which implies that $(\cM_{\dR,\Bun_{1}}(Z/k)^{0})_{\red}^{\nilp}$ is a substack of $\cH^{-1}(0)$. Therefore, both coincide as closed substacks of $\cM_{\dR,\Bun_{1}}(Z/k)^{0}$. 
\end{proof}

The Cartier descent describes the reduced nilpotent locus. Write $\Fr_{Z_T/T}\colon Z_T\rar Z_T^{(p)}\coloneqq Z_T\times_{T,\Fr_T}T$ for the relative Frobenius map with respect to a $k$-scheme $T$. 
Proposition~\ref{prop:Cartier descent for general MdR} defines the $1$-morphism $C_{Z/k}\colon\calCoh(Z^{(p)}/k)\rightarrow \calM_{\dR,k}$, sending $\calF\in \Ob(\calCoh(Z^{(p)}/k)_T)$ to $(\Fr_{Z_T/T}^{\ast}\calF,\nabla_{\can})\in \Ob((\calM_{\dR,k})_T)$.

\begin{lem}\label{lem:Cartier descent}
The $1$-morphism $C_{Z/k}$ restricts to an equivalence of algebraic stacks 
\[
C_{Z/k}\colon \Bun_1(Z^{(p)}/k)^0\xrightarrow{\cong} \calH^{-1}(0).
\]
\end{lem}

\begin{proof}
This follows from Proposition~\ref{prop:Cartier descent for general MdR} and Lemma~\ref{lem:reduced nilpotent locus. rank 1}.
\end{proof}

In the next proposition, we assume that there is a $W$-rational point $e\colon W\rar Z_{W}$, where we take $e$ to be the identity section in the abelian scheme case. For a $W$-scheme $T$, we let $e_{T}\colon T\rar Z_{T}$ be the pullback of $e$ to $T$. In this case, we know from Lemmas~\ref{lem:coarse moduli is universal vector extension, abelian scheme case} and \ref{lem:coarse moduli is universal vector extension, relative curve case} that $\cM_{\dR,\Bun_{1}}(Z_{T}/T)^{0}$ is a split $\Gm$-gerbe, and its coarse moduli space 
$M_{\dR,\Bun_{1}}(Z_{T}/T)^{0}$ is identified with the universal vector extension of the Jacobian $\Pic^{0}(Z_{T}/T)$:
\[
M_{\dR,\Bun_{1}}(Z_{T}/T)^{0}\cong\univext{\Pic^{0}(Z_{T}/T)/T}.
\]
The Cartier descent gives an isomorphism of commutative group schemes over $k$
\[
C_{Z/k}\colon\Pic^{0}({Z^{(p)}/k})\xrightarrow{\sim} H^{-1}(0).
\]
We also have an equivalence $\cH^{-1}(0)\cong H^{-1}(0)\times_{\Spec k}B\Gm$.
Regarding $H^{-1}(0)$ as a closed subscheme of $\univext{\Pic^{0}({Z_{W}/W})/W}$, we define  $\univext{\Pic^{0}({Z_{W}/W})/W}^{\wedge}_{H^{-1}(0)}$ to be the formal completion along $H^{-1}(0)$. Since $k=\F_{p^a}$, we have the Frobenius pullback $\Fr_{Z}^{a,\ast}\colon\Pic^{0}({Z/k})\rar\Pic^{0}({Z/k})$, which induces an endomorphism $\Fr_{Z}^{a,\ast}\colon \univext{\Pic^{0}({Z/k})/k}^{\wedge}_{H^{-1}(0)}\rar \univext{\Pic^{0}({Z/k})/k}^{\wedge}_{H^{-1}(0)}$. Under the Cartier descent isomorphism $C_{Z/k}$, the restriction of $\Fr_{Z}^{a,\ast}$ to $H^{-1}(0)$ corresponds to the Frobenius pullback $\Fr_{Z^{(p)}}^{a,\ast}\colon\Pic^{0}({Z^{(p)}/k})\rar\Pic^{0}({Z^{(p)}/k})$; this follows from the equality $\Fr_{Z_T^{(p)}}^a\circ \Fr_{Z_T/T}=\Fr_{Z_T/T}\circ \Fr_{Z_T}^a\colon Z_T\rar Z_T^{(p)}$. Lastly, let $B\bG_{\mathrm{m},W}^{\wedge}\coloneqq[\Spf W/\bG_{\mathrm{m},W}^{\wedge}]$.

\begin{prop}\label{prop:description of rank 1 Mcris and Misoc using section}
Assume that $Z_W\rightarrow \Spec W$ has a section $e$. Using the notation as above, there is an equivalence of formal algebraic stacks over $W$
\[
\cM_{\cris,1}\cong \univext{\Pic^{0}({Z_{W}/W})/W}^{\wedge}_{H^{-1}(0)}\times_{\Spf W}B\bG_{\mathrm{m},W}^{\wedge}
\]
such that the mod $p$ reduction of $V\colon \cM_{\cris,1}\rar\cM_{\cris,1}$ corresponds to 
\[
\Fr_{Z}^{a,\ast}\times \id\colon \univext{\Pic^{0}({Z/k})/k}^{\wedge}_{H^{-1}(0)}\times_{\Spec k}B\Gm\rar\univext{\Pic^{0}({Z/k})/k}^{\wedge}_{H^{-1}(0)}\times_{\Spec k}B\Gm.
\]
Similarly, there is an equivalence of adic stacks over $K$
\[
\cM_{\isoc,1}\cong (\univext{\Pic^{0}({Z_{W}/W})/W}^{\wedge}_{H^{-1}(0)})_{\eta}\times_{\Spa K}B\Gm^{\an}.
\]
In particular, we have identifications
\[
M_{\cris,1}\cong \univext{\Pic^{0}({Z_{W}/W})/W}^{\wedge}_{H^{-1}(0)}\;\text{and}\;
M_{\isoc,1}\cong (\univext{\Pic^{0}({Z_{W}/W})/W}^{\wedge}_{H^{-1}(0)})_{\eta}.
\]
\end{prop}

\begin{proof}
The first assertion follows from the discussion before the proposition. Since $(B\bG_{\mathrm{m},W}^{\wedge})_{\eta}\rar B\Gm^{\an}$ is an epimorphism, the second assertion follows from the first and Proposition~\ref{prop:image of comparison map from generic fiber to analytification}.
\end{proof}

\begin{cor}\label{cor:Misoc1 is the generic fiber of Mcris1}
The $1$-morphism $(\calM_{\cris,1})_\eta\rightarrow \calM_{\isoc,1}$ is an equivalence of adic stacks over $K$. Moreover, $V$ on $\cM_{\isoc,1}$ agrees with the generic fiber $V_\eta$ applied to $V$ on $\cM_{\cris,1}$. Similarly, the morphism $(M_{\cris,1})_\eta\rightarrow M_{\isoc,1}$ is an isomorphism of \'etale spaces over $K$ compatible with the Verschiebung endomorphisms.
\end{cor}

\begin{proof}
When $Z_W\rightarrow\Spec W$ admits a section, this follows from Proposition~\ref{prop:description of rank 1 Mcris and Misoc using section}. In the general case, there exists a finite extension $k'$ of $k$ such that $Z_{W(k')}\rightarrow\Spec W(k')$ admits a section since  $Z_W\rightarrow\Spec W$ is smooth. One can easily see that it is enough to check the assertions after the base change to a finite extension $K'$ of $K$ and one may take $K'=W(k')[p^{-1}]$; note that the construction involves passing to completion and generic fiber and applying Proposition~\ref{prop:image of comparison map from generic fiber to analytification}, which are all compatible with base change along $W\rightarrow W(k')$.
\end{proof}

\subsection{Counting rank one \texorpdfstring{$F$}{F}-isocrystals using geometry}\label{sec:Counting rank 1 F-isocrystals using geometry}

In this subsection, we give a geometric proof of Theorem~\ref{thm:counting rank 1 F-isocrystals}. 

\begin{defn}
We define the formal scheme $M_{\cris,1}^{V=\id}$ by the fiber product
\[
\xymatrix{
M_{\cris,1}^{V=\id}\ar[r]\ar[d]
&
M_{\cris,1}\ar[d]^-{(\id,V)}
\\
M_{\cris,1}\ar[r]^-\Delta
&
M_{\cris,1}\times_{\Spf W}M_{\cris,1}.
}
\]
Observe that $M_{\cris,1}$ is a commutative formal group scheme, and $V$ is a group homomorphism. Hence we have $M_{\cris,1}^{V=\id}=\Ker(V-\id\colon M_{\cris,1}\rar M_{\cris,1})$, which makes $M_{\cris,1}^{V=\id}$ a commutative formal group scheme. We set $(M_{\cris,1}^{V=\id})_{k}\coloneqq M_{\cris,1}^{V=\id}\times_{\Spf W}\Spec k$.
\end{defn}

The following is a geometric description of $M_{\cris,1}^{V=\id}$.
Let 
\[
\Pic^{0}({Z^{(p)}/k})[\Fr_{Z^{(p)}}^{a,\ast}-\id]\coloneqq\Ker(\Fr_{Z^{(p)}}^{a,\ast}-\id\colon \Pic^{0}({Z^{(p)}/k})\rightarrow\Pic^{0}({Z^{(p)}/k})).
\]

\begin{thm}\label{thm:geometric description of M1cris,BunV=id}
The Cartier descent induces an isomorphism
\[
C_{Z/k}\colon \Pic^{0}({Z^{(p)}/k})[\Fr_{Z^{(p)}}^{a,\ast}-\id]\xrightarrow{\cong}
(M_{\cris,1}^{V=\id})_{k}.
\]
Moreover, the commutative formal group scheme $M_{\cris,1}^{V=\id}$ is finite and flat of degree $\lvert\Pic^{0}({Z/k})(k)\rvert$ over $\Spf W$.
\end{thm}

\begin{proof}
For simplicity, write $F_{T}$ for $\Fr_{T}^{a}$. 
First we assume that $f_W\colon Z_W\rightarrow \Spec W$ has a section.
We know from Proposition~\ref{prop:description of rank 1 Mcris and Misoc using section} that
\[
(M_{\cris,1}^{V=\id})_k\cong\univext{\Pic^{0}({Z/k})/k}^{\wedge}_{H^{-1}(0)}[F_Z^\ast-\id].
\]
Moreover, \cite[Cor. I.1.2.3]{Fujiwara-Kato} shows that the latter is isomorphic to the formal completion of the kernel of $F_Z^\ast-\id$:
\[
\univext{\Pic^{0}({Z/k})/k}^{\wedge}_{H^{-1}(0)}[F_{Z}^{\ast}-\id]\cong (\univext{\Pic^{0}({Z/k})/k}[F_{Z}^{\ast}-\id])^{\wedge}_{H^{-1}(0)[F_{Z}^{\ast}-\id]}.
\]

Let us study the $k$-group scheme $\univext{\Pic^{0}({Z/k})/k}[F_{Z}^{\ast}-\id]$.
Consider the commutative diagram
\[
\xymatrix{
0\ar[r]& f_\ast\Omega_{Z/k}^{1}\ar[r]\ar[d]^-{F_Z^\ast-\id} & \univext{\Pic^{0}({Z/k})/k}\ar[r]\ar[d]^-{F_Z^\ast-\id}& \Pic^{0}(Z/k)\ar[r]\ar[d]^-{F_Z^\ast-\id} &0\\
0\ar[r]& f_\ast\Omega_{Z/k}^{1}\ar[r] & \univext{\Pic^{0}({Z/k})/k}\ar[r]& \Pic^{0}(Z/k)\ar[r] &0
}
\]
of sheaves on $(\Sch/k)_\fppf$ with exact rows. Since $\Fr_{Z_T}^\ast\omega=0$ for any $k$-scheme $T$ and $\omega\in H^0(T,\Omega_{Z_T/T}^1)=\Gamma(T,f_\ast\Omega_{Z/k})$, we see $F_{Z}^{\ast}-\id=-\id$ on $f_\ast\Omega_{Z/k}^{1}$, and therefore the induced map $\univext{\Pic^{0}({Z/k})/k}[F_{Z}^{\ast}-\id]\rar \Pic^{0}({Z/k})[F_{Z}^{\ast}-\id]$ is an isomorphism.
Moreover, Proposition~\ref{prop:description of rank 1 Mcris and Misoc using section} yields an isomorphism $C_{Z/k}\colon \Pic^{0}({Z^{(p)}/k})[F_{Z^{(p)}}^{\ast}-\id]\xrightarrow{\cong} H^{-1}(0)[F_{Z}^{\ast}-\id]$.
We thus obtain the commutative diagram
\[
\xymatrix{
H^{-1}(0)[F_{Z}^{\ast}-\id]\ar@{^{(}->}[r]^-\alpha & \univext{\Pic^{0}({Z/k})/k}[F_{Z}^{\ast}-\id]\ar[d]^-\cong \\
\Pic^{0}({Z^{(p)}/k})[F_{Z^{(p)}}^{\ast}-\id]\ar[u]^-{C_{Z/k}}_-\cong\ar[r]^-{\Fr_{Z/k}^\ast} & \Pic^{0}({Z/k})[F_{Z}^{\ast}-\id],
}
\]
where $\alpha$ is the canonical closed immersion.
We claim that $\Fr_{Z/k}^\ast$ is an isomorphism and thus so is $\alpha$; 
to see this, observe that $F_{Z^{(p)}}=\Fr_{Z^{(p)}}^a$ factors as $Z^{(p)}\xrightarrow{g} Z\xrightarrow{\Fr_{Z/k}} Z^{(p)}$ for some $g$. Then $F_Z=g\circ \Fr_{Z/k}$, and $g^\ast$ gives the inverse to $\Fr_{Z/k}^\ast$ since $F_{Z^{(p)}}^\ast=\id$ on the source and $F_Z^\ast=\id$ on the target.

We obtain from the claim the identifications
\[
(M_{\cris,1}^{V=\id})_k\cong(\univext{\Pic^{0}({Z/k})/k}[F_{Z}^{\ast}-\id])^{\wedge}_{H^{-1}(0)[F_{Z}^{\ast}-\id]}\underset{\cong}{\xleftarrow{C_{Z/k}}} \Pic^{0}({Z^{(p)}/k})[F_{Z^{(p)}}^{\ast}-\id].
\]
Note that $\Pic^{0}({Z^{(p)}/k})[F_{Z^{(p)}}^{\ast}-\id]\cong\Pic^{0}({Z/k})[F_{Z}^{\ast}-\id]$ is a finite $k$-group scheme of degree $\left|\Pic^{0}({Z/k})(k)\right|$. Indeed, $F_{\Pic^{0}({Z/k})}=\Fr_{\Pic^{0}({Z/k})}^{a}\colon\Pic^{0}({Z/k})\rar\Pic^{0}({Z/k})$ is the dual isogeny of $F_Z^\ast$, and $\Pic^{0}({Z/k})[F_{Z}^{\ast}-\id]$ is Cartier dual to $\Pic^{0}({Z/k})[F_{\Pic^{0}({Z/k})}-\id]$. As $\Pic^{0}({Z/k})[F_{\Pic^{0}({Z/k})}-\id]$ is the constant $k$-group scheme associated with the group $\Pic^{0}({Z/k})(k)$, we see that $\Pic^{0}({Z/k})[F_{Z}^{\ast}-\id]$ is a finite $k$-group scheme of degree $\left|\Pic^{0}({Z/k})(k)\right|$.

Since the formal scheme $(M_{\cris,1}^{V=\id})_{k}$ is indeed a finite $k$-scheme, we also see from \cite[Prop.~I.4.2.3]{Fujiwara-Kato} that $M_{\cris,1}^{V=\id}$ is a $p$-adic formal scheme that is finite over $\Spf W$. So it remains to show that $M_{\cris,1}^{V=\id}$ is flat over $\Spf W$, which will follow if we show that $V-\id\colon M_{\cris,1}\rar M_{\cris,1}$ is adically flat. We know from Proposition~\ref{prop:description of rank 1 Mcris and Misoc using section} that $V-\id\colon M_{\cris,1}\rar M_{\cris,1}$ is adic. 

Take any affine open formal subschemes $\Spf A\subset M_{\cris,1}$ and $\Spf B\subset (V-\id)^{-1}(\Spf A)$. We name the restriction of $V-\id$ to $\Spf B$ as $h\colon \Spf B\rar \Spf A$ and also refer to the corresponding maps of rings $A\rar B$ or affine schemes $\Spec B\rar \Spec A$ as $h$ for simplicity. Recall that $M_{\cris,1}= \univext{\Pic^{0}({Z_{W}/W})/W}^{\wedge}_{H^{-1}(0)}$ and $H^{-1}(0)\cong \Pic^0(Z^{(p)}/k)$ is smooth over $k$. Hence $A$ and $B$ are Noetherian and regular. 

Let $I\subset A$ and $J\subset B$ be the ideal of $A$ and $B$ given by the restriction of the ideal sheaf defining $H^{-1}(0)\hookrightarrow M_{\cris,1}$, respectively; $I$ is a prime ideal and an ideal of definition of the adic ring $A$, and the same holds for $J$. We also have $\operatorname{ht}(I)=\dim \univext{\Pic^{0}({Z_{W}/W})/W}^{\wedge}-\dim H^{-1}(0)=\operatorname{ht}(J)$. Obviously, $h(I)\subset J$, and the map $h\colon\Spec B/J\rar\Spec A/I$ is identified with the restriction of $F_{Z}^{\ast}-\id$ to an affine open subset of $H^{-1}(0)$. Since we saw that $H^{-1}(0)[F_{Z}^{\ast}-\id]$ is finite over $k$, we conclude that $F_{Z}^{\ast}-\id\colon H^{-1}(0)\rar H^{-1}(0)$ is finite and flat. Therefore, $h\colon \Spec B/J\rar\Spec A/I$ is quasi-finite and flat.

Since $h\colon\Spf B\rar\Spf A$ is adic, it suffices to show that $h\colon\Spec B\rar \Spec A$ is flat at all maximal ideals $J\subset \fq\subset B$. Fix such $\fq$ and let $\fp\coloneqq h^{-1}(\fq)\subset A$, which necessarily contains $I$. Since $h\colon \Spec B/J\rar\Spec A/I$ is quasi-finite and flat, we deduce that $\fp$ is also a maximal ideal. On one hand, since $A_\fp$ is regular local and $I$ is prime, we have $\dim A_\fp=\dim A_\fp/IA_\fp+\operatorname{ht}(I)$. On the other hand, $A/I$ is an integral domain of finite type over $k$ with $\fp/I$ a maximal ideal. Hence we have $\dim A_\fp/IA_\fp=\dim A/I$ by \cite[00OS]{stacks-project}, which gives $\dim A_\fp=\dim A/I+\operatorname{ht}(I)$.
Similarly, we have $\dim B_\fq=\dim B/J+\operatorname{ht}(J)$. As $\dim A/I=\dim B/J$ and $\operatorname{ht}(I)=\operatorname{ht}(J)$, we conclude $\dim A_\fp=\dim B_\fq$. We finally deduce from the miracle flatness criterion (e.g. \cite[00R4]{stacks-project}) that $\Spec B\rar \Spec A$ is flat at $\fq$.
This completes the proof when $f_W$ has a section.

We now prove the general case, where we do not assume that $f_{W}\colon Z_{W}\rar \Spec W$ has a section. Let $k'=\F_{p^{ab}}$ be a finite extension of $k$ such that $f_W$ admits a section over $W'\coloneqq W(k')$. Then, we have an identification $M_{\cris,1}\times_WW'\cong \univext{\Pic^{0}({Z_{W'}/W'})/W'}^{\wedge}_{(H^{-1}(0)\times_kk')}$.
Observe that the canonical map $(Z_{k'})^{(p)}\rightarrow (Z^{(p)})\times_kk'$ is a $k'$-scheme isomorphism. Therefore, there is the commutative diagram of $k'$-schemes
\[
\xymatrix{
\Pic^0((Z_{k'})^{(p)}/k') \ar[d]^-\cong\ar[r]^-{\Fr_{Z_{k'}/k'}} &\Pic^0(Z_{k'}/k') \ar[d]^-\cong\\
\Pic^0(Z^{(p)}/k)\times_kk' \ar[r]^-{\Fr_{Z/k}\times \id} & \Pic^0(Z/k)\times_kk'.
}
\]
In particular, the $k'$-endomorphism $\Fr_Z^a\times\id$ on $\Pic^0(Z/k)\times_kk'$ defines an endomorphism $(\Fr_Z^a\times\id)^\ast$ on $\univext{\Pic^{0}({Z_{k'}/k'})/k'}$ and yields an isomorphism
\[
(M_{\cris,1}^{V=\id})_{k}\times_kk'\cong\univext{\Pic^{0}({Z_{k'}/k'})/k'}^{\wedge}_{(H^{-1}(0)\times_kk')}[(\Fr_Z^a\times\id)^\ast-\id].
\]
Arguing as above, we see that 
\[
(H^{-1}(0)\times_kk')[(\Fr_Z^a\times\id)^\ast-\id]\hookrightarrow \univext{\Pic^{0}({Z_{k'}/k'})/k'}[(\Fr_Z^a\times\id)^\ast-\id]
\]
is an isomorphism. Therefore, $C_{Z/k}\colon \Pic^{0}({Z^{(p)}/k})[\Fr_{Z^{(p)}}^{a,\ast}-\id]\rightarrow
(M_{\cris,1}^{V=\id})_{k}$ is an isomorphism, since it is so after the base change to $k'$.
For the second assertion, it is enough to show that $M_{\cris,1}^{V=\id}\times_WW'$ is finite flat over $W'$ of degree $\lvert\Pic^{0}({Z/k})(k)\rvert=\lvert (H^{-1}(0)\times_kk')[(\Fr_Z^a\times\id)^\ast-\id]\rvert$, which is proved using the first assertion as above.
\end{proof}

\begin{defn}
Define the adic space $M_{\isoc,1}^{V=\id}$ by the following fiber product
\[
\xymatrix{
M_{\isoc,1}^{V=\id}\ar[r]\ar[d]
&
M_{\isoc,1}\ar[d]^-{(\id,V)}
\\
M_{\isoc,1}\ar[r]^-{\Delta}
&
M_{\isoc,1}\times_{\Spa K}M_{\isoc,1}.
}
\]
Since $M_{\isoc,1}$ is a commutative group adic space, and $V$ is a group homomorphism, we have $M_{\isoc,1}^{V=\id}=\Ker(V-\id\colon M_{\isoc,1}\rar M_{\isoc,1})$, which makes $M_{\isoc,1}^{V=\id}$ a commutative group adic space over $K$. 
\end{defn}

\begin{proof}[Proof of Theorems~\ref{thm:counting rank 1 F-isocrystals}(i) and \ref{thm:counting rank 1 F-isocrystals Fqm}]
On one hand, it follows from Theorem~\ref{thm:irreducible F-isocrystals discrete} and Remark~\ref{rem:moduli of rank 1 isocrystals} that $M_{\isoc,1}^{V=\id}$ defined above is the open subspace of $M_{\isoc,\irr}^{V=\id}$ consisting of rank one objects in Theorem~\ref{thm:counting rank 1 F-isocrystals}(i). On the other hand, we know from Proposition~\ref{prop: properties of generic fiber of formal algebraic spaces}(i) and Corollary~\ref{cor:Misoc1 is the generic fiber of Mcris1} that $M_{\isoc,1}^{V=\id}$ is the generic fiber of $M_{\cris,1}^{V=\id}$. Hence Theorem~\ref{thm:counting rank 1 F-isocrystals}(i) follows from Theorem~\ref{thm:geometric description of M1cris,BunV=id}.
We also note that the proof of Theorem~\ref{thm:geometric description of M1cris,BunV=id} applies verbatim to the analogous result for $M_{\cris,1}^{V^m=\id}$ for any $m\ge1$. Namely, the Cartier descent induces an isomorphism
\[
\Pic^{0}(Z^{(p)}/k)[\Fr_{Z^{(p)}}^{am,\ast}-\id]\xrightarrow{\cong}(M_{\cris,1}^{V^m=\id})_{k}
\]
for any $m\ge1$, and $M_{\cris,1}^{V=\id}$ is finite flat of degree $\lvert\Pic^{0}(Z/k)(\bF_{q^{m}})\rvert$ (recall that $k=\bF_{q}$, $q={p^{a}}$). This implies Theorem~\ref{thm:counting rank 1 F-isocrystals Fqm}.
\end{proof}

\appendix

\section{Generalities on sites and stacks}\label{sec:appendix}

We recall several definitions and results about sites and stacks from \cite{stacks-project}. A full subcategory $\calC'$ of a category $\calC$ is said to be \emph{strictly full} if every object of $\calC$ that is isomorphic to an object of $\calC'$ is also in $\calC'$.

\smallskip
\noindent
\textbf{Special cocontinuous functors}.

\begin{defn}[{\cite[03CG]{stacks-project}}]
A functor $u\colon \calC\rightarrow\calD$ between sites is said to be \emph{special cocontinuous} if it satisfies the following conditions:
\begin{enumerate}
    \item $u$ is continuous and cocontinuous;
    \item given morphisms $a,b\colon U\rightarrow U'$ in $\calC$ with $u(a)=u(b)$, there exists a covering $\{f_i\colon U_i\rightarrow U\}$ in $\calC$ such that $a\circ f_i=b\circ f_i$ for every $i$;
    \item given $U,U'\in\Ob(\calC)$ and a morphism $d\colon u(U)\rightarrow u(U')$ in $\calD$, there exist a covering $\{f_i\colon U_i\rightarrow U\}$ and morphisms $c_i\colon U_i\rightarrow U'$ in $\calC$ such that $u(c_i)=d\circ u(f_i)$ for every $i$;
    \item for every $V\in\Ob(\calD)$, there exists a covering in $\calD$ of the form $\{u(U_i)\rightarrow V\}$.
\end{enumerate}
In this case, the cocontinuous functor $u$, via \cite[00XO]{stacks-project}, yields the morphism of topoi $\Sh(\calC)\rightarrow \Sh(\calD)$, which is an equivalence by \cite[03A0]{stacks-project}.
\end{defn}

\smallskip
\noindent
\textbf{Categories fibered in setoids and groupoids}.
A $1$-morphism of categories fibered in groupoids over a category $\calC$ is referred to as a $1$-morphism over $\calC$ for simplicity.

\begin{defn}[{\cite[04TM]{stacks-project}}]\label{def:category fibered in setoids associated to presheaf}
Let $\calC$ be a category and let $F$ be a presheaf on $\calC$. Define a category fibered in groupoids $\calS_{F}\rightarrow \calC$ as follows:
\begin{itemize}
\item an object in $\calS_{F}$ is a pair $(U,x)$ where $U\in\Ob(\calC)$ and $x\in F(U)$;
\item a morphism $(U,x)\rightarrow (V,y)$ is a morphism $f\colon U\rightarrow V$ in $\calC$ with $f^\ast y=x$;
\item the functor $\calS_{F}\rightarrow\calC$ sends $(U,x)\in\Ob(\calS_\calC)$ to $U\in\Ob(\calC)$ and $(f\colon (U,x)\rightarrow (V,y))\in \Mor_{\calS_F}((U,x),(V,y))$ to $(f\colon U\rightarrow V)\in \Mor_{\calC}(U,V)$.
\end{itemize}
\end{defn}

\begin{rem}\label{rem:category fibered in setoids associated to presheaf}
Keep the notation as in Definition~\ref{def:category fibered in setoids associated to presheaf}. Then $\calS_F$ is indeed a category fibered in setoids. Moreover, every category fibered in setoids over $\calC$ is equivalent over $\calC$ to $\calS_F$ for some presheaf $F$ (see \cite[02Y2, 0045]{stacks-project}).
If $G$ is another presheaf, we have a natural identification
\[
\Mor_{\PSh(\calC)}(F,G)=\Mor(\calS_F,\calS_G)/\text{$2$-isomorphism},
\]
where the latter is the categories of $1$-morphisms $\calS_F\rightarrow\calS_G$ of categories over $\calC$ up to $2$-isomorphisms (see \cite[04SC]{stacks-project}).
If $F=h_{U}$ for $U\in\Ob(\calC)$, then $\calS_{F}=\calC/U$.
If $\calC$ is a site, then $\calS_F$ is a stack if and only if $F$ is a sheaf (see \cite[0430]{stacks-project}).
\end{rem}

\begin{rem}\label{rem:equivalent to a split category fibered in groupoids}
Even more generally, given a functor $\calF\colon (\calC)^\mathrm{op}\rightarrow \Grpd$, one can similarly define a category fibered in groupoids $\calS_\calF\rightarrow \calC$ (see \cite[0049]{stacks-project}). A category fibered in groupoids obtained in this way is called \emph{split}. Every category fibered in groupoids over $\calC$ is equivalent to a split one by \cite[02XY]{stacks-project}.
\end{rem}

\smallskip
\noindent
\textbf{Stackification of categories fibered in groupoids}.

\begin{lem}[{\cite[02ZP, 0436]{stacks-project}, \cite[Rem.~4.6.6]{Olsson-book}}]\label{lem:stackification}
Let $\calC$ be a site and let $\calX\rightarrow \calC$ be a category fibered in groupoids over $\calC$. Then there exist a stack in groupoids $\calX'\rightarrow \calC$ and a $1$-morphism $f\colon \calX\rightarrow \calX'$  over $\calC$ satisfying the following two properties:
\begin{enumerate}
 \item for any $U\in\Ob(\calC)$ and any $x_1,x_2\in \Ob(\calX_U)$, the morphism
 \[
 \Mor_\calX(x_1,x_2)\rightarrow \Mor_{\calX'}(f(x_1),f(x_2))
 \]
 induced by $f$ makes the target the sheafification of the source;
 \item for any $U\in\Ob(\calC)$ and $x'\in \Ob(\calX'_U)$, there exists a covering $\{U_i\rightarrow U\}$ such that $x'|_{U_i}$ is in the essential image of $f\colon \calX_{U_i}\rightarrow\calX'_{U_i}$ for every $i$.
\end{enumerate}
Moreover, it satisfies the following universal property: for every stack in groupoids $\calY\rightarrow\calC$ and every $1$-morphism $g\colon \calX\rightarrow\calY$ over $\calC$, there exists a $1$-morphism $g'\colon \calX'\rightarrow\calY$ over $C$ such that the diagram
\[
\xymatrix{
\calX\ar[rr]^-g\ar[rd]_-f&&\calY\\
&\calX'\ar[ru]_-{g'}&
}
\]
is $2$-commutative. Such an $\calX'$ is unique up to a $1$-isomorphism that is unique up to unique $2$-isomorphism and called the \emph{stackification} of $\calX$. 
\end{lem}

\begin{lem}\label{lem:sheafification and stackification}
 Let $F$ be a presheaf on a site $\calC$ and write $F^a$ for its sheafification. Then the induced $1$-morphism $\calS_F\rightarrow \calS_{F^a}$ over $\calC$ identifies $\calS_{F^a}$ with the stackification of $\calS_F$.
\end{lem}

\begin{proof}
Using Remark~\ref{rem:category fibered in setoids associated to presheaf}, one can easily check that $\calS_F\rightarrow \calS_{F^a}$ satisfies the properties in Lemma~\ref{lem:stackification}. 
\end{proof}

The following is also used to prove the equivalence of stacks whose definition involves stackification.

\begin{lem}[{\cite[04WQ, 046N]{stacks-project}}]\label{lem:criteria for equivalence of stacks}
A $1$-morphism $f\colon \calX\rightarrow\calY$ of stacks over a site $\calC$ is an equivalence if and only if the following two conditions hold:
\begin{enumerate}
 \item for any $U\in\Ob(\calX)$ and $x_1,x_2\in\Ob(\calX_U)$, the morphism
 \[
f\colon \Mor_\calX(x_1,x_2)\rightarrow \Mor_\calY(f(x_1),f(x_2))
 \]
 is an isomorphism of sheaves on $\calC/U$;
 \item for any $U\in\Ob(\calC)$ and any $y\in \Ob(\calY_U)$, there exists a covering $\{U_i\rightarrow U\}$ such that $y|_{U_i}$ is in the essential image of $f\colon \calX_{U_i}\rightarrow\calY_{U_i}$ for every $i$.
\end{enumerate}
More generally, the full faithfulness of $f$ is equivalent to (i).
\end{lem}

\smallskip
\noindent
\textbf{Equivalence relations and quotient stacks}.
Let $\calC$ be a site.

\begin{defn}\label{def:equivalence relation for site}
An \emph{equivalence relation} on a sheaf $U$ on $\calC$ consists of a sheaf $R$ on $\calC$ and a pair of morphisms $s,t\colon R\rightrightarrows  U$ such that the map $(t,s)\colon R\rightarrow U\times U$ induces an injection $R(T)\rightarrow U(T)\times U(T)$ and defines an equivalence relation on $U(T)$ for each $T\in \Ob(\calC)$. 
In this case, write $U/R$ for the sheaf on $\calC$ associated to the presheaf $T\mapsto U(T)/R(T)$. When $\calC$ is subcanonical, an equivalence relation on $U\in\Ob(\calX)$ refers to an equivalent relation on the sheaf $h_U$ represented by $U$.
\end{defn}

\begin{example}\label{eg:coequalizer diagram of sheaves}
If $U\rightarrow X$ is an epimorphism of sheaves on $\calC$, then $R\coloneqq U\times_XU$ yields an equivalence relation on $U$ and the induced morphism $U/R\rightarrow X$ is an isomorphism by \cite[086K]{stacks-project}.
\end{example}

\begin{defn}\label{def:groupoid for site}
A \emph{groupoid in sheaves} on $\calC$ is a quintuple $(U,R,s,t,c)$, where $U,R$ are sheaves on $\calC$, and $s,t\colon R\rightarrow U$ and $c\colon R\times_{s,U,t}R\rightarrow R$ are morphisms of sheaves
such that, for every $T\in\Ob(\calC)$, the quintuple $(U(T),R(T),s,t,c)$ defines a groupoid category $[U(T)/R(T)]$ (cf. \cite[0230]{stacks-project}). Note that $s,t$ are necessarily epimorphisms of sheaves.

 For a groupoid in sheaves $(U,R,s,t,c)$, define a category fibered in groupoids $[U/_{p}R]\rightarrow \calC$ as follows: set
\[
\Ob([U/_{p}R])\coloneqq\{(T,x) \mid T\in\Ob(\calC),~x\in U(T)\},
\]
and for $(T,x),(T',x')\in\Ob([U/_{p}R])$, 
\begin{align*}
\Mor_{[U/_{p}R]}&((T,x),(T',x'))\\
&\coloneqq\{(f,\phi)\mid f\in \Mor_\calC(T,T'), \phi\in\Mor_{[U(T)/R(T)]}(x,f^\ast x')\},
\end{align*}
where $f^\ast x'\in U(T)$ is the pullback of $x'\in U(T')$. The structural functor $[U/_{p}R]\rightarrow\calC$ is given by $(T,x)\mapsto T$. Note that the fiber category over $T$ is $[U(T)/R(T)]$.

Let $[U/R]\rightarrow \calC$ denote the stackification of $[U/_pR]$.
\end{defn}

\begin{lem}\label{lem:2-coequalizer property of quotient stack}
Let $(U,R,s,t,c)$ be a groupoid in sheaves on $\calC$.
\begin{enumerate}
 \item There exist a canonical $1$-morphism $\pi\colon \calS_U\rightarrow [U/R]$ over $\calC$ and a canonical $2$-morphism $\alpha\colon \pi\circ s\rightarrow \pi\circ t$ making the diagram
\[
\xymatrix{
\calS_R \ar[r]^-s\ar[d]_-t & \calS_U \ar[d]_-\pi \\
\calS_U \ar[r]^-\pi & [U/R]
}
\]
$2$-commutative. If $\calC$ is subcanonical, then the diagram is a $2$-fiber product of stacks in groupoids over $\calC$.
 \item For every stack $\calX\rightarrow \calC$ in groupoids with a $1$-morphism $f\colon \calS_U\rightarrow \calX$ and a $2$-morphism $\beta\colon f\circ s\rightarrow f\circ t$ such that $\beta\star \id_c=(\beta\star \id_{\pr_0})\circ(\beta\star \id_{\pr_1})$, there exists a $1$-morphism $[U/R]\rightarrow \calX$ over $\calC$ which makes the diagram 
 \[
\xymatrix{
\calS_R \ar[r]^-s\ar[d]_-t & \calS_U \ar[d]_-\pi \ar[rdd]^-f&\\
\calS_U \ar[r]^-\pi\ar[rrd]_-f & [U/R]\ar[rd]&\\
&& \calX
}
\]
$2$-commutative.
\end{enumerate}
\end{lem}

\begin{proof}
(i) For the first part, the proof of \cite[044R, 044S]{stacks-project} works verbatim. The second part is proved as in \cite[044V, 04M9]{stacks-project}, which only uses the subcanonicality of the site in question.

(ii) The proof of \cite[044T, 044U]{stacks-project} works verbatim.
\end{proof}

\begin{defn}[{cf.~\cite[04MB]{stacks-project}}]
Let $(U,R,s,t,c)$ be a groupoid in sheaves on $\calC$. A \emph{$[U/R]$-descent datum} relative to a covering $\calT=\{T_i\rightarrow T\}$ in $\calC$ consists of a system $(u_i,r_{ij})$ where $u_i\in U(T_i)$ and $r_{ij}\in R(T_i\times_TT_j)$ for each $i$ and $j$ such that\footnote{Unlike \textit{loc.~cit.}, we use the convention $\pr_1$, $\pr_2$, and $\pr_{ij}$ ($i,j\in\{1,2,3\}$) for the projections.}
\begin{itemize}
 \item $s(r_{ij})=\pr_1^\ast u_i$ and $t(r_{ij})=\pr_2^\ast u_j$ in $U(T_i\times_TT_j)$;
  \item $c(\pr_{23}^\ast r_{jk}, \pr_{12}^\ast r_{ij})=\pr_{13}^\ast r_{ij}$ in $R(T_i\times_TT_j\times_TT_k)$.
\end{itemize}
We have obvious notions of morphisms between $[U/T]$-descent data relative to $\calT$ and the pullback of $[U/T]$-descent data relative to a refinement $\calT'\rightarrow\calT$ of coverings.
\end{defn}

The quotient stack $[U/R]$ is described by $[U/R]$-descent data as follows. We refer the reader to \cite[044X]{stacks-project} for the precise formulation and proof.

\begin{lem}\label{lem:explicit description of quotient stack}
Assume that $\calC$ is subcanonical. Let $(U,R,s,t,c)$ be a groupoid in sheaves on $\calC$ and write $\pi\colon \calS_U\rightarrow [U/R]$ for the natural $1$-morphism. For every $T\in\Ob(\calC)$ and every $x\in \Ob([U/R]_T)$, there exists a covering $\{f_i\colon T_i\rightarrow T\}$ such that $f_i^\ast x\cong \pi(u_i)$ for some $u_i\in U(T_i)$. Moreover, the composition of isomorphisms
\[
\pi(\pr_1^\ast u_i)\cong \pr_1^\ast f_i^\ast x\cong \pr_2^\ast f_j^\ast x \cong \pi(\pr_2^\ast u_j)
\]
is of the form $\pi(r_{ij})$ for some $r_{ij}\in R(T_i\times_TT_j)$.
In this way, $x$ gives rise to a $[U/R]$-descent datum $(u_i,r_{ij})$.
Any $[U/R]$-descent datum arises in this way, and under this correspondence, morphisms in $[U/R]_T$ correspond to morphisms between $[U/R]$-descent data. These correspondences of objects and morphisms are compatible with refinements of coverings.
\end{lem}

\newcommand{\etalchar}[1]{$^{#1}$}
\providecommand{\bysame}{\leavevmode\hbox to3em{\hrulefill}\thinspace}
\providecommand{\MR}{\relax\ifhmode\unskip\space\fi MR }
\providecommand{\MRhref}[2]{%
  \href{http://www.ams.org/mathscinet-getitem?mr=#1}{#2}
}
\providecommand{\href}[2]{#2}

\end{document}